\documentclass[12pt,leqno,oneside]{article} 

\usepackage{amsfonts}
\usepackage{latexsym}
\usepackage{amsthm}  
\usepackage{amsmath}
\usepackage{amssymb}
\usepackage[dvips]{graphicx,color}
\usepackage{makeidx}

\usepackage[all]{xy}

\usepackage{bbm}
\usepackage{enumerate}
\usepackage{geometry}
\usepackage{rotating}

\numberwithin{equation}{section}

{\theoremstyle{plain}\newtheorem{Def}{Definition}[section]}
\newtheorem{Lemma}[Def]{Lemma}
\newtheorem{Th}[Def]{Theorem}
\newtheorem{Cor}[Def]{Corollary}
\newtheorem{Prop}[Def]{Proposition}
{\theoremstyle{remark} \newtheorem{Rem}[Def]{Remark}}
{\theoremstyle{remark} \newtheorem{Exa}[Def]{Example}}

\newtheorem{Assumption}[Def]{Assumption}
\newtheorem{Convention}[Def]{Product Convention}

\newcommand{\ZZ}{\mathbb{Z}}
\newcommand{\QQ}{\mathbb{Q}}

\newcommand{\RR}{\mathbb{R}}

\newcommand{\CC}{\mathbb{C}}

\newcommand{\PP}{\mathbb{P}}

\newcommand{\s}{\hspace{3pt}}

\newcommand{\laa}{\langle}
\newcommand{\raa}{\rangle}
\newcommand{\bs}{\backslash}

\newcommand{\fin}{\hspace{4 pt} $\square$}
\newcommand{\ov}{\overline}

\newcommand{\ca}{\mathcal}
\newcommand{\mk}{\mathfrak}
\newcommand{\wt}{\widetilde}
\newcommand{\wh}{\widehat}

\newcommand{\M}[4]
{\left(
\begin{array}{cc}
#1 & #2 \\
#3 & #4
\end{array} 
\right)}
\newcommand{\V}[2]
{\left(
\begin{array}{c}
#1\\
#2
\end{array} \right)}
\newcommand{\R}[2]{\left(
\begin{array}{cc}
#1 & #2\\
\end{array} 
\right)}

\newcommand{\biglslant}[2]{{\raisebox{-.20em}{$#1$}\backslash\raisebox{+.20em}{$#2$}}}
\newcommand{\bigrslant}[2]{{\raisebox{.2em}{$#1$}/\raisebox{-.2em}{$#2$}}}

\newcommand*{\lcdot}{\raisebox{-0.25ex}{\scalebox{1.5}{$\cdot$}}}

\DeclareMathOperator{\sign}{sign}
\DeclareMathOperator{\Norm}{{\bf N}}
\DeclareMathOperator{\Tr}{Tr}

\newcommand{\bz}{\scalebox{0.6}[0.9]{$\mathbbm{O}$}}
\newcommand{\bzt}{\scalebox{0.4}[0.6]{$\mathbbm{O}$}}

\DeclareMathOperator{\bu}{\mathbbm{1}}

\DeclareMathOperator{\ii}{i}

\DeclareMathOperator{\Hom}{Hom}
\DeclareMathOperator{\cov}{cov}
\DeclareMathOperator{\End}{End}
\DeclareMathOperator{\res}{res}
\DeclareMathOperator{\Res}{Res}
\DeclareMathOperator{\Imm}{Im}
\DeclareMathOperator{\Ree}{Re}

\DeclareMathOperator{\Maps}{Maps}
\DeclareMathOperator{\Stab}{Stab}
\DeclareMathOperator{\Latt}{Latt}
\DeclareMathOperator{\ord}{ord}
\DeclareMathOperator{\signb}{{sign}}

\DeclareMathOperator{\sg}{ {sg}}
\DeclareMathOperator{\Isom}{ {Isom}}
\DeclareMathOperator{\Hol}{ {Hol}}

\DeclareMathOperator{\lcm}{ {lcm}}
\DeclareMathOperator{\im}{ {im}}
\DeclareMathOperator{\Sym}{ {Sym}}
\DeclareMathOperator{\SPD}{ {SPD}}
\DeclareMathOperator{\Pic}{Pic}
\DeclareMathOperator{\fish}{{part}}

\DeclareMathOperator{\modd}{mod}
\DeclareMathOperator{\Gal}{Gal}
\DeclareMathOperator{\cond}{cond}

\newcommand{\spart}[2]{\mbox{$[#1;#2]$-$\fish$-}}

\newcommand{\spartsing}[2]{\mbox{$[#1;#2]$-$\fish$}}

\begin{document}
\title{$GL_2$-real analytic Eisenstein series twisted by parameter matrices and multiplicative integral quasi-characters}
\author{Hugo Chapdelaine}
\date{July, 2016}
\maketitle
\begin{abstract}
Let $K$ be a totally real field of dimension $g$ over $\QQ$ and let $\ca{O}_K$ be its ring of integers. Consider the hermitian symmetric domain 
$\mk{h}^g$ consisting of the cartesian product of $g$ copies of 
the Poincar\'e upper half-plane. The group $SL_2(K)$ acts naturally on $\mk{h}^{g}$ by M\"obius transformations. 
In this work, we make a detailed study of certain families of Eisenstein series $\{G(z,s)\}_{s\in\CC}$ 
where $z\in\mk{h}^{g}$ and $s\in\CC$. The function $G(z,s)$ is real analytic in the variable $z\in\mk{h}^{g}$ and holomorphic in the variable $s\in\CC$. Moreover, 
it is modular in the variable $z$ with respect to a discrete subgroup of $SL_2(K)$ which is commensurable to
$SL_2(\ca{O}_K)$. The construction of $G(z,s)$ consists in taking a sum over the direct sum of 
two lattices $\mk{m}\oplus \mk{n}\subseteq K^2$, where the general term of the
defining summation
is \lq\lq twisted\rq\rq\; simultaneously by a parameter matrix $U\in M_2(K)$ and by an integral quasi-character
of $(K\otimes_{\QQ}\CC)^{\times}$.
The first main result of this work gives an analytic characterization of $G(z,s)$ in terms of the $\{\spart{c}{s}G\}_{c\in\PP^1(K)}$, where
 $\spart{c}{s}G$ may be viewed as the non-square-integrable part of $[z\mapsto G(z,s)]$
in a neighborhood of the cusp $c$. The second main result provides an explicit description of the Fourier series expansion of $[z\mapsto G(z,s)]$
which leads to a proof of the  meromorphic continuation of $[s\mapsto G(z,s)]$ to all of $\CC$.
The third main result gives two proofs of a functional equation which relates $G(z,s)$ to $G^*(z,1-s)$, where
$G^*(z,s)$ is the \lq\lq dual Eisenstein series\rq\rq\s associated to $G(z,s)$. Finally, the fourth main result of this monograph, gives a new 
proof of the meromorphic continuation and of the functional equation of a class of partial zeta functions that had been studied previously by the
author. 
\end{abstract}
\tableofcontents
\section{Introduction}
Let $K=\QQ$ be the field of rational numbers. The symmetric space associated to the algebraic group $SL_2/\QQ$
is  $SL_2(\RR)/SO(2)\simeq\mk{h}^+$, where $\mk{h}^+=\mk{h}:=\{x+\ii y\in\CC:y>0\}$ corresponds to the Poincar\'e upper half-plane. Let also 
$\mk{h}^-:=\{x+\ii y\in\CC:y<0\}$ denote the lower half-plane, and $\mk{h}^{\pm}:=\mk{h}\cup\mk{h}^{-}$ be the disjoint union
of the upper and lower half-planes.
The spaces $\mk{h}$ and $\mk{h}^-$ are isomorphic as real analytic manifolds through the complex conjugation.
Note that $\mk{h}^{\pm}\simeq GL_2(\RR)/SO(2)\RR^{\times}$, where $SO(2)\RR^{\times}=\{\gamma\in GL_2(\RR):\det(\gamma)>0,\gamma \ii=\ii\}$.
Therefore, one may interpret $SO(2)\RR^{\times}$ as the \lq\lq positive determinant stabilizer\rq\rq\; of $\ii$ with respect to 
$GL_2(\RR)$. By analogy with the group $SL_2/\QQ$, we call 
$\mk{h}^{\pm}$ the {\it $\pm$-symmetric space associated to $GL_2/\QQ$}. Note though, that it is not a symmetric space of $GL_2/\QQ$ 
in the usual sense, since $SO(2)\RR^{\times}$ is not a finite index subgroup of $O(2)$ (a maximal compact subgroup of $GL_2(\RR)$). The (disconnected) space
$\mk{h}^{\pm}$ is equipped with the Poincar\'e metric $ds^2=\frac{dx^2+dy^2}{y^2}$, so in particular, it is a Riemannian space.
The group $PSL_2(\RR)=SL_2(\RR)/\{\pm I_2\}$ acts faithfully,  properly discontinuously and isometrically
on $\mk{h}$ and $\mk{h}^{-}$ by M\"obius transformations. 
For $z\in\mk{h}^{\pm}$ and $s\in\CC$ with $\Ree(s)>1$, 
the {\it classical real analytic Eisenstein series} associated to the discrete group $PSL_2(\ZZ)\leq PSL_2(\RR)$ is defined by 
\begin{align}\label{couv0}
E(z,s):=\sum_{\Gamma_{\infty}\bs PSL_2(\ZZ)}|\Imm(\gamma z)|^s=\frac{1}{2}\sum_{\substack{(m,n)\in\ZZ^2\\ \gcd(m,n)=1}} \frac{|y|^s}{|mz+n|^{2s}},
\end{align}
where $\Gamma_{\infty}:=\Stab_{PSL_2(\ZZ)}(\infty)=\left\{\pm\M{1}{n}{0}{1}:n\in\ZZ\right\}$ corresponds to the isotropy
group of the cusp $\infty=\frac{1}{0}$ with respect to the discrete group $PSL_2(\ZZ)$. Note that the function 
$E(z,s)$ is real analytic in $z\in\mk{h}^{\pm}$, even though its defining series involves the absolute value of $y$, since
the $y$ coordinate never crosses the real axis when $z$ varies inside $\mk{h}^{\pm}$. If one thinks of an Eisenstein 
series as a function in $z$, the it is probably more accurate to say that $E(z,s)$ is a {\it family} of real analytic Eisenstein series
of weight $0$, where $s\in\Pi_1:=\{s\in\CC:\Ree(s)>1\}$. Sometimes in order to emphasize this point of view we may use
the suggestive notation $\{E(z,s)\}_{s\in\Pi_1}$.

Consider now the following {\it modified real analytic Eisenstein series}
\begin{align}\label{couv}
G(z,s):=\sum_{(0,0)\neq (m,n)\in\ZZ^2}\frac{|y|^s}{|mz+n|^{2s}},
\end{align}
where $z\in\mk{h}^{\pm}$ and $s\in\Pi_1$. Using the fact that $\ZZ$ is a unique factorization domain, it is straightforward to see that 
\begin{align*}
G(z,s)=2\zeta(2s)\cdot E(z,s),
\end{align*}
where $\zeta(s)$ corresponds to the Riemann zeta function. Note that, in comparison to $E(z,s)$, the defining
summation of the modified Eisenstein series $G(z,s)$ 
goes over the {\it full punctured} lattice $\ZZ^2\bs(0,0)$ rather than just the set of ordered pairs of coprime integers 
(i.e., the set of right cosets of $\biglslant{\Gamma_{\infty}}{ PSL_2(\ZZ)}$).
The Eisenstein series constructed in this manuscript, when the totally real field $K$ is equal to $\QQ$, may be 
viewed as a natural generalization of $G(z,s)$. Let us explain 
in what sense they generalize \eqref{couv}. Let $U:=\M{u_1}{v_1}{u_2}{v_2}\in M_2(\QQ)$ be
a {\it parameter matrix}, $p,w\in\ZZ$ be {\it integral weights}, and
$\mk{m},\mk{n}\subseteq\QQ$ be a pair of two lattices (discrete $\ZZ$-modules of rank one).
To each quadruple $\ca{Q}=((\mk{m},\mk{n}),U,p,w)$, the $(GL_2/\QQ)$-Eisenstein series considered in this work can be defined explicitly  as 
\begin{align}\label{couv1}
&G_{\ca{Q}}(z,s)=G_{(\mk{m},\mk{n})}^w(U,p\;;z,s):=\\[2mm] \notag
&\sum_{(-v_1,-v_2)\neq(m,n)\in(\mk{m},\mk{n})} 
\frac{\omega_p((m+v_1)z+(n+v_2))}{\left((m+v_1)z+(n+v_2)\right)^{w}}\cdot\frac{e^{2\pi\ii (u_1(m+v_1)+u_2(m+v_2))}}{|(m+v_1)z+(n+v_2)|^{2s}}
\cdot |y|^s,
\end{align}
where  $z\in\mk{h}^{\pm}$ and $s\in\Pi_1$. Note that the parameter $w$ is associated to the {\it integral quasi-character} of $\CC^{\times}$
of weight $w$: $\Norm^w:\CC^{\times}\rightarrow\CC^{\times}$, $z\mapsto z^w$; while  
the integer $p$ is associated to the {\it integral unitary character} of $\CC^{\times}$: $\omega_p:\CC^{\times} \rightarrow S^1$,
$z\mapsto \left(\frac{z}{|z|}\right)^p$.
Also, the map $(m,n)\mapsto e^{2\pi\ii u_1 m+u_2 n}$, for $(m,n)\in(\mk{m},\mk{n})$, 
may be viewed as a finite order character of the abelian group $\mk{m}\oplus\mk{n}$.
Therefore, the general term of the summation \eqref{couv1} is twisted {\it simultaneously} by a finite order character of
$\mk{m}\oplus\mk{n}$ and by the \lq\lq integral quasi-character\rq\rq\; of $(\QQ\otimes_\QQ\CC)^{\times}\simeq\CC^{\times}$ given
by $z\mapsto \omega_p(z)\Norm(z)^{w}$. Moreover, the quantity $mz+n$  appears in the summation \eqref{couv1} after having been
shifted, additively, by the quantity $v_1z+v_2$.

When $\mk{m}=\mk{n}=\ZZ$, $U=\M{u_1}{v_1}{u_2}{v_2}\in M_2(\QQ)$, $p=0$ and $w=0$, 
the Eisenstein series in \eqref{couv1} appears in equation (2) on p. 622 of \cite{Eps} (see also \cite{Eps2} and \cite{Sie1}).
More precisely, for $z=x+\ii y\in\mk{h}^{\pm}$ fixed, it corresponds to the Epstein zeta function of degree $2$ associated to the positive definite quadratic form
$(m,n)\mapsto Q_z(m,n):=|y|^{-1}|mz+n|^2$ with characteristics $\Big|\begin{matrix} v_1 & v_2\\ u_1 & u_2\end{matrix}\Big|$. 
Let us explain in more details the origin of this equivalence. Recall that 
an Epstein zeta function of degree $r$ and characteristics $\Big|\begin{matrix} v \\ u\end{matrix}\Big|$\; ($u,v\in\RR^r$) associated
to a symmetric positive definite $r$ by $r$ matrix $Q$,  is defined as
\begin{align}\label{tupp}
Z_r\Big|\begin{matrix} v\\ u\end{matrix}\Big|(Q,s):=\sum_{m\in\ZZ^r\bs\{(0,\ldots,0)\}} \frac{e^{2\pi\ii m u^t}}{Q[m+v]^{s}}.
\end{align}
Here $s$ is a complex number with $\Ree(s)>\frac{r}{2}$, the element $m\in\ZZ^r$ is viewed as a row vector and
$Q[m]:=m Q m^t$, where $m^t$ corresponds to the transpose of $m$. Note that the function $[Q\mapsto Z_g\Big|\begin{matrix} v\\ u\end{matrix}\Big|(Q,s)]$
may be viewed as a function on the usual symmetric space associated to $GL_r(\RR)$ which can be naturally identified with the space
\begin{align*}
\SPD_r(\RR):=\{Q\in M_r(\RR):\mbox{$Q$ is symmetric and positive definite}\}.
\end{align*}
Recall that for $\gamma\in GL_r(\RR)$ and $Q\in\SPD_r(\RR)$ (resp. $AO(r)\in GL_r(\RR)/O(r)$) the
left $GL_r(\RR)$-action is given explicitly by $Q\mapsto\gamma Q\gamma^t$ (resp. $AO(r)\mapsto\gamma AO(r)$). 
It is straightforward to check that the map $\phi:\SPD_r(\RR)\rightarrow GL_r(\RR)/O(r)$,
$ AA^t\mapsto A O(r)$ is an isomorphism of left $GL_r(\RR)$-homogeneous spaces. 

In the special case where the characteristics $u,v\in\QQ^r$ are rationals, it follows directly from the definition
of the Epstein zeta function that
$[Q\mapsto Z_r\Big|\begin{matrix} v\\ u\end{matrix}\Big|(Q,s)]$ is modular of weight $0$ with respect to a suitable
discrete subgroup $\Gamma_{v,u}\leq GL_r(\ZZ)$ which depends only on the class of the pair $(v,u)$ inside $\QQ^r/\ZZ^r\times\QQ^r/\ZZ^r$.
In particular, the restricted function $[Q\mapsto Z_r\Big|\begin{matrix} v\\ u\end{matrix}\Big|(Q,s)]$ is an example of a 
{\it real analytic Eisenstein series} (in the sense of \cite{Borel1})
on the symmetric space $GL_r(\RR)/O(r)$.  When $r=2$, we have the sequence of maps
\begin{align*}
\mk{h}\stackrel{\psi}{\rightarrow} SL_2(\RR)/SO(2)\leq GL_2(\RR)/O(2)\stackrel{\phi^{-1}}{\rightarrow} \SPD_2(\RR),
\end{align*}
where $\psi$ and $\phi^{-1}$ are isomorphisms of left $SL_2(\RR)$-homogeneous spaces. The map
$\psi$ is given explicitly by $z=x+\ii y\mapsto \M{|y|^{1/2}}{x|y|^{-1/2}}{0}{|y|^{-1/2}} SO(2)$.
Note that $\phi^{-1}\circ\psi(z)$ is nothing else than the symmetric matrix associated to the quadratic form $Q_z$.
In particular, if we restrict $[Q\mapsto Z_r\Big|\begin{matrix} v\\ u\end{matrix}\Big|(Q,s)]$ to $\mk{h}$ via 
$\phi^{-1}\circ\psi$, we readily see that it coincides with \eqref{couv1}.

The main goal of this manuscript is a detailed study of \eqref{couv1} and its generalizations when one replaces the $\QQ$-algebraic group
$GL_2/\QQ$ by the $\QQ$-algebraic group $\Res_{\QQ}^{K}(GL_2)$, where $K$ is an arbitrary totally real field. The generalization
of \eqref{couv1}, when $K$ is an arbitrary totally real field, is given in Definition \ref{eis_ser} of Section \ref{def_eis}. 
It is formally equivalent to \eqref{couv1}, except
that one is required to \lq\lq quotient the summation\rq\rq\; by a certain finite index subgroup of $\ca{O}_K^{\times}$, in order 
to insure the convergence of the series. Note that, when $g=[K:\QQ]>1$, the Eisenstein series in Definition \ref{eis_ser}, 
is no more equivalent to an Epstein zeta function. 
Instead, it should be viewed as 
a generalization of the real analytic Eisenstein series studied by Asai in \cite{Asai}. 
As the Eisenstein series in \eqref{couv1} for $GL_2/\QQ$, our Eisenstein series for $\Res_{\QQ}^{K}(GL_2)$ will depend 
on some additional data: a pair lattices $\mk{m},\mk{n}\subseteq K$, 
a parameter matrix $U\in M_2(K)$, an integral unitrary weight $p\in\ZZ^g$, and a \lq\lq parallel weight 
shift\rq\rq\; $w\in\ZZ$. We view these Eisenstein series as functions 
$(z,s)\mapsto G_{\ca{Q}}(z,s)=G_{(\mk{m},\mk{n})}^w(U,p\,;z,s)$, where $(z,s)\in(\mk{h}^{\pm})^g\times\CC$. We call
$(z,s)$ the {\it arguments} of $G_{\ca{Q}}(z,s)$, and the data $\ca{Q}=((\mk{m},\mk{n}),U,p,w)$, the {\it parameters} of 
$G_{\ca{Q}}(z,s)$. Usually, we  think of the parameters
as being fixed. When the parameters $U,p$ and $w$ are trivial and $\mk{m}=\mk{n}$, we recover the Eisenstein series studied by Asai in \cite{Asai}.
Note though, that Asai did not restrict himself to totally real fields $K$, as we do in the present work;
mainly in order to simplify the presentation. There is a technical aspect which we improve on \cite{Asai}. Contrary to 
\cite{Asai}, we are not assuming that $K$ has class number one. Let us also explain a difference of points of view 
between the present work and \cite{Asai}. Let $(\mk{h}^{\pm})^g\simeq (GL_2(\RR)/SO(2)\RR^{\times})^g$ 
be the $\pm$-symmetric space associated to $\ca{G}:=\Res_{\QQ}^{K}(GL_2)$.
In \cite{Asai}, the author decides to work with the connected component $\mk{h}^g$ which corresponds
to the symmetric space associated to $\ca{G}_1:=\Res_{\QQ}^{K}(SL_2)$. However, in this manuscript, 
we decided to work systematically with $(\mk{h}^{\pm})^g$.
We would like to emphasize that the natural domain of definition of $G_{\ca{Q}}(z,s)$, when viewed as a function in $z$, is really
the disconnected space $(\mk{h}^{\pm})^g$ rather than the smaller connected space $\mk{h}^g$. It is true though, that the space
$(\mk{h}^{\pm})^g$ is a disjoint union of $2^g$ components, where each component is {\it real analytically} isomorphic to $\mk{h}^g$.
Based on this observation, it is possible to associate a collection of $2^g$ functions on $\mk{h}^g$ to the function 
$[z\mapsto G_{\ca{Q}}(z,s), z\in(\mk{h}^{\pm})^g]$. Note that 
the corresponding  collection of functions is entirely determined by the restriction of
$G_{\ca{Q}}(z,s)$ to anyone of the $2^g$ connected components of $(\mk{h}^{\pm})^g$.
The main advantage of working with $(\mk{h}^{\pm})^g$, rather than with $\mk{h}^g$, is that it allows us to interpret the function 
$z\mapsto G_{\ca{Q}}(z,s)$ as being modular 
with respect to an arithmetic subgroup of $GL_2(K)$, rather than an arithmetic subgroup of the smaller group 
$GL_2^+(K)$ (invertible matrices with a totally positive determinant). Moreover, when $z$ is allowed to be in 
$(\mk{h}^{\pm})^g$, the expression $G_{\ca{Q}}(z,s)$
satisfies additional symmetries (e.g. symmetries induced from sign changes of the entries of $U$) 
which are not visible when $z$ is restricted to be in $\mk{h}^g$.

This monograph is divided in 9 sections. Each of the section is divided in subsections.
For the convenience of the reader, we also added, at the end of the paper, an appendix which contains 
some background material and various calculations that did not fit well in the main body of the text.
Let us finish this introduction by giving a brief description of each section.

In Section \ref{nota}, 
we introduce the notation and the necessary background that will be in effect for 
the rest of the work. In Section \ref{def_eis}, we define the real analytic Eisenstein series.
In Section \ref{sec_beg}, we present
the four main results ({\bf A}, {\bf B}, {\bf C} and {\bf D}) of this work.

In Section \ref{frip}, we introduce special hypergeometric functions and prove various formulas for them which are essential 
for the explicit computation of the non-zero
Fourier series coefficients of $[x\mapsto G_{\ca{Q}}(x+\ii y,s)]$. In Section \ref{mero_cont}, 
we discuss in details the analytic continuation of these special hypergeometric functions when their parameters
vary in certain regions of the $3$-dimensional complex space. In Section \ref{ref_for}, we prove reflection formulas that are 
needed for the second proof of the functional equation of $G_{\ca{Q}}(z,s)$.

In Section \ref{lat_nu2}, 
we record, for the convenience of the reader, some general results about lattices in number fields that we could not easily find
in the literature. Note that this section is largely independent of the rest of the work.  
In Section \ref{Four}, we record classical results about the Fourier series of functions on the standard $g$-dimensional torus which are needed
for the sequel. In Section \ref{mon_con}, we introduce the monomial in $z$: $P(\alpha,\beta;z)$. This monomial, is the building block
of our Eisenstein series. We also introduce the {\it Product Convention \ref{conv}} which
gives a precise meaning to the expression $P(\alpha,\beta;z)$ when some of the coordinates of $z$ are {\it negative real numbers}. 
In Section \ref{dos}, we provide an explicit formula for the Fourier series expansion of the periodic function 
$[x\mapsto R_{\mathcal{L}}(\alpha,\beta;x+\ii y)]$, which can be viewed as the average of the monomial $P(\alpha,\beta;z)$
over the lattice $\ca{L}$.

In Section \ref{hubo}, we introduce various subgroups
of $GL_1(K)$ and $GL_2(K)$ that naturally appear in the definition and in the study
of $G_{\ca{Q}}(z,s)$. We also present and prove their main properties. In Section \ref{real_an}, we choose to rewrite the Eisenstein series 
$G_{\ca{Q}}(z,s)=G_{(\mk{m},\mk{n})}^w(U,p\,;z,s)$ in terms of a bi-weight $[\alpha(s),\beta(s)]\in\CC^g\times\CC^g$ which depends
on $p\in\ZZ^g$ and $s\in\CC$. The auxiliary notation $G_{(\mk{m},\mk{n})}^{\alpha(s),\beta(s)}(U\,;z):=G_{(\mk{m},\mk{n})}^w(U,p\,;z,s)$ 
is conducive to some subsequent calculations.
In Section \ref{real_an2}, we introduce the notion of modular forms of bi-weight $[\alpha,\beta]$. In Section \ref{symm_a}, we prove
various symmetries (in the arguments and in the parameters) of the expression $G_{(\mk{m},\mk{n})}^w(U,p\,;z,s)$.
In Section \ref{Fourier}, we give the definition of the Fourier series expansion
of a modular form of bi-weight $[\alpha,\beta]$ at a cusp $c\in\PP^1(K)$. In Section \ref{groo}, we give some standard upper bound
for the function $(z,s)\mapsto G_{(\mk{m},\mk{n})}^0(U,p\,;z,s)$ ($w=0$) when $\Ree(s)>1$ and $z$ tends to
a cusp. 

In Section \ref{def_mas}, we define Maa\ss\; graded operators and partial-graded Laplacians. 
In Section \ref{eigen1}, we show that $[z\mapsto G_{(\mk{m},\mk{n})}^0(U,p\,;z,s)]$ is an eigenvector, with eigenvalue
$s(1-s)$, with respect to each of the partial-graded Laplacian. In Section \ref{Hil}, 
we introduce the Hilbert space $L^2(\mk{h};\Gamma;p)$ endowed with its Petersson inner product.

In Section \ref{qual}, we give a qualitative description of the Fourier series coefficients of certain
families of real analytic modular forms of integral unitary weight $p\in\ZZ^g$. In Section \ref{exis},
we prove that such non-trivial families admit at least one non-zero Fourier coefficient. In Section 
\ref{unicity}, we show that certain cuspidal real analytic families of modular forms of integral unitary weight
$p$ don't exist. In Section \ref{car}, we provide an {\it analytic characterization} of certain real analytic families
of modular forms, which include the family $\{G_{(\mk{m},\mk{n})}^0(U,p;z,s)\}_{s\in\CC}$ as a special case.
In Section \ref{party0}, we introduce the useful notation
$\spart{c}{s}G$ and $\spart{c}{1-s}G$, whenever $c\in\PP^1(K)$ is a cusp and $\{G(z,s)\}_{s\in\CC}$
is a real analytic family of modular forms.
Finally, in Section \ref{writing_app_7}, we explain, under the assumption of a technical condition,
how to write explicitly the real analytic Eisenstein series $G_{(\mk{m},\mk{n})}^0(U,p\, ;z,s)$ as a sum of classical real analytic Poincar\'e-Eisenstein series.

In Section \ref{nico0}, we compute the Fourier series expansion of $[z\mapsto G_{(\mk{m},\mk{n})}^w(U,p\, ;z,s)]$
at the cusp $\infty$ and write it as a sum of three main terms: $\left( T_1+T_2+T_3\right)|\Norm(y)|^s$. 
In Section \ref{lun}, we prove the meromorphic continuation of $[s\mapsto G_{(\mk{m},\mk{n})}^w(U,p\, ;z,s)]$.
In Section \ref{sec_h}, we compute the standard Fourier series
expansion of the holomorphic Eisenstein series $[z\mapsto G_{(\mk{m},\mk{n})}^w(U,\bz\, ;z,s)]$ when $w\geq 3$.

In Section \ref{pars}, we define the 
{\it completed Eisenstein series}. In Section \ref{fun_eqn}, we give two proofs of a functional
equation which is satisfied by the completed Eisenstein series. In Section \ref{der}, we rewrite this functional equation
in terms of the uncompleted Eisenstein series $G_{\ca{Q}}(z,s)$. 
In Section \ref{Colmez_trick}, using an idea of Colmez, we explain how the calculation done in the second proof of 
Theorem \ref{nice_fn_eq} leads to a new proof of the meromorphic continuation and the functional equation of the partial 
zeta functions that were studied in \cite{Ch4}. In Section \ref{weighted}, we prove a functional equation for certain weighted sums of the (uncompleted) Eisenstein series
$G_{(\mk{m}_i,\mk{n}_i)}^w(U_i,p\;z,s)$'s, where the parameters $(p,w)$ are being fixed, and the pairs
$((\mk{m}_i,\mk{n}_i),U_i)$'s are allowed to vary. Moreover, when the parameter matrices $U_i$'s have a particular shape and the pair of
lattices are {\it diagonal pairs}, we also explain how the previous weighted sums lead to sums of {\it ray class invariant of $K$}
which also satisfy a functional equation.

In Appendix \ref{app_1a}, we present the main properties of partial zeta functions which are 
twisted simultaneously by a finite additive character of a lattice, and a (multiplicative) 
signature character of $K^{\times}$.  We also
rewrite the functional equation proved in \cite{Ch4} in a different way which is needed for the Proof of
Theorem \ref{nice_fn_eq}. In Appendix \ref{app_3}, we study a particular one-dimensional space of a linear system of ODEs of order $2$ in $g$-variables.
In Appendix \ref{app_4}, we discuss some recurrence relations satisfied by the coefficients of the 
Taylor series expansion of $[s\mapsto G_{(\mk{m},\mk{n})}^0(U,p\, ;z,s)]$ at $s=1$. In Appendix \ref{app_5},
we introduce some background material on the Riemannian symmetric space $\mk{h}^g$ and on the set of
cusps associated to a subgroup $\Gamma\leq (SL_2(\RR))^g$ which is commensurable to $SL_2(\ca{O}_K)$.
We also introduce the notion of distance between a point $z\in\mk{h}^g$ and a cusp $c\in\PP^1(K)$. 
Finally, in Appendix \ref{app_6}, we provide a proof of Proposition \ref{poids} by using the 
{\it point-pair invariant kernel method} due to Selberg.

\section{Notation, background and statement of the main results}\label{not_pre}

\subsection{Notation and background}\label{nota}
We would like first to introduce some notation that will be in effect for the whole monograph. Let $K$
be a totally real field. If $A$ is a $\QQ$-algebra, we let $K_A:=K\otimes_\QQ A$. 
In particular, $K_A$ has a natural structure of a left 
$A$-module given by $a(x\otimes b)=x\otimes ab$ where $x\in K$ and $a,b\in A$. We may view $\CC$ as a $\QQ$-algebra
and, therefore, consider the $\CC$-algebra $K_\CC=K\otimes_\QQ \CC$. Let $\{\sigma_{i}\}_{i=1}^g$ 
be the set of real embeddings of $K$ into $\RR$. If $x\in K$, then $x^{(i)}$ is taken to mean $\sigma_i(x)$.
We have a natural $\CC$-algebra isomorphism  
$\theta:K_{\CC}\rightarrow \CC^{g}$ given on simple tensors by $x\otimes w\mapsto (x^{(i)}w)_{i=1}^g$.
Note that the isomorphism class of $K_\CC$, as a $\CC$-algebra, depends only on the degree $[K:\QQ]$
and not on the field $K$ itself. We endow $K_\CC$ with the topology induced from $\CC^g$ via $\theta$. 
We will always view $K$ as a $\QQ$-subalgebra of $K_\CC$ via $x\mapsto x\otimes 1$.
Notice that $\theta|_{K_\RR}:K_\RR\rightarrow\RR^g$ gives, in a similar way, an $\RR$-algebra isomorphism such that $K$ is a dense subset of $K_\RR$. 
If $w=u+\ii v\in\CC$, then $\ov{w}=u-\ii v$ denotes the complex conjugate. By convention, we use the symbol $\ii:=\sqrt{-1}$, whereas 
the symbol $i$ will be used usually as an index. For $z=(z_i)_{i=1}^g\in\CC^g$, we let 
\begin{enumerate}[(a)]
 \item  $\Norm(z):=\prod_{i=1}^g z_i\in\CC$ (the norm of $z$),
 \item $\Tr(z):=\sum_{i=1}^g z_i\in\CC$ (the trace of $z$),
 \item  $\Ree(z):=(\Ree(z_i))_{i=1}^g\in\RR^g$ (the real part vector of $z$),
 \item $\Imm(z):=(\Imm(z_i))_{i=1}^g\in \RR^g$ (the imaginary part vector of $z$),
 \item $|z|:=(|z_1|,\ldots,|z_g|)\in\RR^g$ (the absolute value vector of $z$).
\end{enumerate}
From now on, we identify the $\CC$-algebras $K_\CC$ with $\CC^g$ through the isomorphism $\theta$ without any further mention.
Therefore, if $z\in K_\CC$, then $\Tr(z)$, $\Norm(z)$, $\Ree(z)$ and $\Imm(z)$ are taken to mean $\Tr(\theta(z))$,
$\Norm(\theta(z))$, $\Ree(\theta(z))$, $\Imm(\theta(z))$. Moreover, if $z\in K_\CC$, then $z_i$ is taken to mean $\theta(z)_i$, the 
$i$-th coordinate of $\theta(z)$. We note that if $x\in K\subseteq K_\CC$, then $\Tr(x)$ and $\Norm(x)$ coincide with the usual definitions of
the trace and the norm of an algebraic number.
It is also convenient to introduce the following shorthand notation: for $\alpha\in\CC^g$ and $z\in\CC^g$, such that for all $i$, $z_i\neq 0$, we define
\begin{align*}
z^{\alpha}:=\prod_{j=1}^g z_j^{\alpha_j},
\end{align*}
where $z_j^{\alpha_j}:=e^{\alpha_j\log z_j}$. Here, $\log z_j$
is computed with respect to the {\it principal branch} of the 
logarithm, i.e., for $w\in\CC\bs\{0\}$, $\log w:=\log|w|+\ii\arg(w)$ where $-\pi <\arg(w)\leq \pi$. 
Note that the principal branch of the logarithm satisfies the rule
\begin{align}\label{chum0}
\ov{\log w}=\log \ov{w}+\chi_w,
\end{align}
where $\chi_w=0$, if $w\notin\RR_{<0}$ and $\chi_w=-2\pi\ii$, if $w\in\RR_{<0}$. It follows that
for $w_1,w_2\in\CC\bs\{0\}$ and $\alpha\in\CC$, that one has the rule
\begin{align}\label{chum}
(z_1z_2)^{\alpha}z_1^{-\alpha}z_2^{-\alpha}=e^{\alpha(\arg(w_1w_2)-\arg(w_1)-\arg(w_2))}.
\end{align}
If $h\in K\bs\{0\}$ and $\alpha\in\CC$, then $h^{\alpha}$ is taken to mean $\theta(h)^{\alpha}$. 
We let $\bu:=(1,1,\ldots,1)\in\CC^g$ and $\bz:=(0,0,\ldots,0)\in\CC^g$. We also let $e_j=(0,\ldots,1,\ldots,0)$, where we place
$1$ in the $j$-th coordinate and $0$ elsewhere. Let us illustrate the notation introduced before with two examples. 
Let $z=x+\ii y\in\CC^g$, where $x=\Ree(z),y=\Imm(z)\in\RR^g$, and let $s\in\CC$. Then the
notation above give us the two identities:
$|\Norm(\Imm(z))|^s=|y|^{\bu\cdot s}$ and $|\Norm(z)|^s=|z|^{\bu\cdot s}$. Both of these identities will be used freely
in the rest of the paper.

\subsubsection{Quasi-characters of $K_{\CC}^{\times}$}

We would like now to give some notation in order to handle {\it quasi-characters} of locally compact abelian groups. 
Let $G$ be a locally compact abelian group. A quasi-character $\chi$ of $G$ is defined as a {\it continuous} group homomorphism
$\chi:G\rightarrow\CC^{\times}$. A character $\chi$ of $G$ is defined as a quasi-character $\chi$ of $G$ such that $\im(\chi)\subseteq S^1$,
where $S^1$ corresponds to the unit circle. We define $\wh{G}:=\Hom_{\mbox{\tiny{cont}}}(G,\CC^{\times})$ to be
the group of quasi-characters of $G$. We endow $\wh{G}$ with the compact-open topology. 
We also let ${}^{u}\wh{G}:=\{\chi\in\wh{G}:\im(\chi)\subseteq S^1\}$
be group of characters of $G$ (i.e., unitary quasi-characters of $G$) which is often called the {\it Pontryagin dual of $G$}.
It is also convenient to define ${}^{p}\wh{G}:=\{\chi\in\wh{G}:\im(\chi)\subseteq \RR_{>0}\}$,
the {\it group of positive quasi-characters of $G$}. From the topological group isomorphism $\CC^{\times}\rightarrow S^1\times\RR_{>0}$ given
by $z\mapsto (\frac{z}{|z|},|z|)$, it follows that every quasi-character $\chi\in\wh{G}$ can be written uniquely as
$\chi=\omega\cdot\eta$, where $\omega\in {}^{u}\wh{G}$ and $\eta\in {}^{p}\wh{G}$. For example, if $G=\CC^{\times}$, then 
${}^{u}\wh{\CC^{\times}}\simeq \ZZ\times \RR$ and ${}^{p}\wh{\CC^{\times}}\simeq \RR$. In particular,
$\wh{\CC}\simeq\ZZ\times\RR\times\RR$. In general, using the Pontryagin-van Kampen structure theorem for locally compact abelian 
groups (see Theorem 25 of \cite{Morris}) it is not difficult to prove that ${}^{p}\wh{G}$ is again a locally compact abelian group. 
Since the Pontryagin dual ${}^{u}\wh{G}$ is well-known to be locally compact, it follows that 
$\wh{G}$ is locally compact. For the whole paper, we should be mainly interested in the set of quasi-characters 
on $K_{\CC}^{\times}$ (a locally compact abelian group)
which are {\it invariant} under the action of a given finite index subgroup $\mathcal{V}^+$ of $\ca{O}_K^{\times}(\infty)$, where
$\ca{O}_K^{\times}(\infty)$ denotes the group of totally positive units of $\ca{O}_K$. It is precisely this set of quasi-characters that
can be used to construct real analytic Eisenstein series. 

Since $K_{\CC}^{\times}$ is isomorphic to $(\CC^{\times})^g$, it
follows that $\wh{K_{\CC}^{\times}}$ is isomorphic, as a topological group, to $\ZZ^g\times \RR^g\times\RR^g$. Let us choose (arbitrarily) 
such a topological isomorphism $\varphi:\wh{K_{\CC}^{\times}}\simeq \ZZ^g\times \RR^g\times\RR^g$. 
\begin{Def}
We define $\mk{X}_{0}\leq \wh{K_{\CC}^{\times}}$ to 
be the subgroup of quasi-characters corresponding to $\ZZ^g\times\{0\}\times\{0\}$ under the isomorphism $\varphi$. 
\end{Def}
Note that, the topological group $\mk{X}_0$ is independent of $\varphi$ and that it depends only on the degree $[K:\QQ]=g$ and not
on $K$ itself. One has the following characterization of characters in $\mk{X}_0$: 
if $\chi\in\wh{K_\CC^{\times}}$, one may check that $\chi\in\mk{X}_0$ if and only if, for all 
$\lambda\in(\RR^{\times})^{g}\subseteq K_\CC^{\times}$ and all $z\in K_\CC^{\times}$, $\chi(\lambda z)=\chi(z)$.
It follows from the previous observation that all the quasi-characters in $\mk{X}_0$ are unitary 
and, therefore, lie in fact in ${}^{u}\wh{K_{\CC}^{\times}}$.

We choose to identify the integral lattice $\mk{S}:=\ZZ^g$ with the subgroup of characters $\mk{X}_{0}$ 
by the explicit map $p\mapsto \omega_p$, given by
\begin{align*}
\omega_{p}:K_{\CC}^{\times}&\rightarrow S^1\\
            (z_{i})_{i=1}^g&\mapsto \prod_{i=1}^g\left(\frac{z_i}{|z_i|}\right)^{p_i}.
\end{align*}
\begin{Def}
We call the set $\mk{S}$ the {\it integral weight lattice} associated to the $\CC$-algebra $K_\CC$, and the elements of $\mk{S}$ are called {\it integral weights}. 
A character $\omega_p\in\mk{X}_0$ is called an {\it integral unitary character of $K_{\CC}^{\times}$}. 
\end{Def}

\vspace{-0.5cm}
As before, we denote the {\it trivial weight} by $\bz:=(0,0,\ldots,0)\in\mk{S}$ and the {\it unit weight}
by $\bu:=(1,1,\ldots,1)\in\mk{S}$.
\begin{Rem}
Let $L/K$ be a totally imaginary quadratic extension over $K$, so $L$ is a CM field. Let $\Phi$ be a fixed CM type of $L$ and
let $L\hookrightarrow K_\CC$ be the embedding obtained from $\Phi$. One may check that $\omega_p|L^{\times}$, for $p\neq\bz$, corresponds to the infinite part of a 
non-trivial unitary character of type $A$ in the sense of Weil (see \cite{Weil3}). Moreover,  
$\omega_p^2=\omega_{2p}$ becomes the infinite part of a non-trivial unitary character of type $A_0$.
\end{Rem}

\subsubsection{The signature group and the space $K_{\CC}^{\pm}$}

Let $x\in\RR^{\times}$ and $\ov{p}\in\ZZ/2\ZZ$, we define
\begin{enumerate}
 \item $\sign(x)=1$ if $x>0$, and $\sign(x)=-1$ if $x<0$, 
 \item $\sg(x)=\ov{0}\in\ZZ/2\ZZ$ if $x>0$, and $\sg(x)=\ov{1}\in\ZZ/2\ZZ$ if $x<0$,
 \item $[\ov{p}]=0\in\ZZ$ if $\ov{p}=\ov{0}$, and $[\ov{p}]=1\in\ZZ$ if $\ov{p}=\ov{1}$.
\end{enumerate}
\begin{Def}
The \textit{signature group}
associated to the $\RR$-algebra $K_\RR$ is defined to be the set $\ov{\mk{S}}:=(\ZZ/2\ZZ)^{g}$. Elements of $\ov{\mk{S}}$ are called signatures, or signature elements. 
\end{Def}
We let $\ov{\bz}:=(\ov{0},\ov{0},\ldots,\ov{0})\in\ov{\mk{S}}$ 
be the zero signature and $\ov{\bu}=(\ov{1},\ov{1},\ldots,\ov{1})\in\ov{\mk{S}}$ be the unit signature. 
We define the function 
\begin{align*}
\sg:(\RR^{\times})^g\rightarrow\ov{\mk{S}}
\end{align*}
by $\sg(x)=\ov{p}\in\ov{\mk{S}}$ where $\ov{p}$ is such that $\sign(x_i)=(-1)^{[\ov{p}_i]}$ for all $i\in\{1,\ldots,g\}$. Note that, when $g=1$, the function
$\sg$ defined above agrees with our previous definition of the function $\sg$ on $\RR^{\times}$. An element 
$x\in (\RR^{\times})^{g}$, which satisfies
$\sg(x)=\ov{p}$, is said to have the {\it signature $\ov{p}$}. Since $K^{\times}
\subseteq(\RR^{\times})^g\subseteq K_\CC^{\times}$, we may restrict the function $\sg$ to $K^{\times}$. As usual, an element
$x\in (\RR^{\times})^g$ with signature $\ov{\bz}$ is said to be totally positive. Sometimes, we may also use
the more conventional notation $x\gg 0$ to denote that $x$ is totally positive. 
\begin{Def}
A group homomorphism $\omega:K_{\RR}^{\times}\rightarrow\{\pm 1\}$ is called a 
{\it sign character of $K_\RR^{\times}$}. 
\end{Def}
\vspace{-0.5cm}
Note that the sign characters of $K_\RR^{\times}$ are automatically continuous.
For each embedding $\sigma_i$ of $K$, the map $s_i:K^{\times}\rightarrow\{\pm 1\}$, where $s_i(x):=\signb(\sigma_i(x))$, 
is an example of a sign character.
A sign character $\omega$ is completely determined by its restriction to $K^{\times}$. 
In particular, we have $\omega|_{K^{\times}}=\prod_{i=1}^g s_i^{[\ov{p}_i]}$ for a unique
signature $\ov{p}=(\ov{p}_i)_{i=1}^g\in\ov{\mk{S}}$ which we call the signature of $\omega$. Moreover,
every sign character $\omega$ of $K^{\times}$ extends uniquely to a character of $K_{\RR}^{\times}$ and, therefore, 
gives rise to an element in $\wh{K_{\RR}^{\times}}$. Let us choose arbitrarily a continuous group
isomorphism $\psi:{}^{u}\wh{K_\RR^{\times}}\rightarrow \{\pm 1\}^g\times\RR^g$. 
\begin{Def}
We define $\ov{\mk{X}}_{0}$ as 
$\psi^{-1}(\{\pm 1\}^g\times\{0\})\leq\wh{K_{\RR}^{\times}}$.
\end{Def}
\vspace{-0.5cm}
One may check that $\ov{\mk{X}}_{0}$ is a subgroup of $\wh{K_{\RR}^{\times}}$ which is independent of the choice of $\psi$. 
We have a natural restriction map $res:\wh{K_\CC^{\times}}\rightarrow \wh{K_\RR^{\times}}$, 
$\omega_p\in\mk{X}_0\mapsto\res(\omega_p):=\omega_p|_{K_\RR^{\times}}=\omega_{\ov{p}}\in\ov{\mk{X}}_0$. Here the map $p\mapsto \ov{p}$ corresponds 
to the natural projection $\ZZ^g\rightarrow (\ZZ/2\ZZ)^{g}$.
We note that $\omega_{\ov{\bzt}}$ is the {\it trivial sign character} of $K_{\RR}^{\times}$,
and that $\omega_{\ov{\bu}}=\signb\circ\Norm$ is the {\it unit weight} sign character of $K_{\RR}^{\times}$. 

Let 
\begin{align*}
K_{\CC}^{\pm}:=\{z=(z_i)_{i=1}^g\in K_{\CC}:\mbox{for all $i$, $\Imm(z_i)\neq 0$}\}\subseteq K_{\CC}^{\times}.
\end{align*}
Note that $K_\CC^{\pm}$ may be identified with $(\mk{h}^{\pm})^g$.
For each $z\in K_{\CC}^{\pm}$, one may associate a signature element of $\ov{\mk{S}}$ to the element $z$ in the following way:
$z\mapsto\sg(\Imm(z))\in\ov{\mk{S}}$. Using this observation, one may decompose the space $K_{\CC}^{\pm}$ as follows
\begin{align*}
K_\CC^{\pm}=\coprod_{\ov{p}\in\ov{\mk{S}}}\mk{h}^{\ov{p}}, \s\s\mbox{where}\s\s\s \mk{h}^{\ov{p}}:=\{z\in K_\CC^{\pm}:\sg(\Imm(z))=\ov{p}\}.
\end{align*} 
In particular, the subset $K_\CC^{\pm}$ is an open set of $K_\CC^{\times}$ with $2^g$ components sharing the same boundary. 
More precisely, for each $\ov{p}\in\ov{\mk{S}}$, the boundary of $\mk{h}^{\ov{p}}$ corresponds to $K_\RR=\RR^g+\bz\ii\subseteq K_\CC$.
Finally, note that the group $\ov{\mk{S}}$ acts naturally on $K_{\CC}$ in the following way: for $\ov{p}\in\ov{\mk{S}}$ and $z=(z_j)_{j=1}^g\in K_\CC$, we
let, for $j\in\{1,\ldots,g\}$, $(z^{\ov{p}})_j:=c_{\infty}^{[\ov{p}_j]}(z_j)$, where $c_{\infty}$ denotes the complex conjugation in $\CC$.
Note that the subsets $K_\RR$ and $K_{\CC}^{\pm}$ are stable under the action of $\ov{\mk{S}}$.

\subsubsection{The sign character $\omega_{\ov{p}}$ as the limit value of $\omega_p$}

Let $\ov{p}\in\ov{\mk{S}}$ be a fixed signature element and let 
$q\in\mk{S}$ be such that $\ov{q}=\ov{p}$. We would like now to explain how the sign character $\omega_{\ov{p}}$ of 
$K_{\RR}^{\times}$ may be viewed as the {\it limit value} of the 
integral unitary character $\omega_q$ of $K_\CC^{\times}$, as $z\in K_{\CC}^{\pm}$ tends to a point in $(\RR^{\times})^g\subseteq K_{\CC}$. 
Note that the sets $K_\RR$ and $K_\CC^{\pm}$, when viewed as subsets of $K_\CC$, are disjoint, and $K_\RR$ corresponds to the boundary
of $K_{\CC}^{\pm}$. For $w\in\CC\bs\RR_{\leq 0}$ and $n\in\ZZ$, 
we have
\begin{align}\label{lim}
w^{n/2}\cdot (\ov{w})^{-n/2}=e^{n\ii\arg(w)}=\left(\frac{w}{|w|}\right)^n.
\end{align}
Note that the first equality in \eqref{lim} {\it does not necessarily hold true} if $w\in\RR_{<0}$. Assume, now 
that $z=u+\ii v\in\CC\bs\{\RR\cup\ii\RR\}$. On may check that
\begin{align}\label{lim0}
\lim\limits_{v\rightarrow 0}w^{n/2}(\ov{w})^{-n/2}=\left(\sign(u)\right)^{[\ov{n}]}.
\end{align}
Now let   $z\in K_{\CC}^{\pm}$ be such that $\Ree(z)=u=(u_i)_{i=1}^g\in\RR^g$ with $u_i\neq 0$ for all $i$. Let also 
$p\in \ZZ^g$. It follows from \eqref{lim0} that
\begin{align}\label{feve}
\lim_{\substack{\Ree(z)=u\\ \Imm(z)\rightarrow\bz}} z^{q/2}z^{-q/2}=\omega_q(u)=\omega_{\ov{p}}(u).
\end{align}
The limit formula above was an initial motivation for the present paper. In particular, note that 
limit formula \eqref{feve} holds true even if there exists a coordinate $u_i\in\RR_{<0}$. However, note that
$u^{q/2}u^{-q/2}$ may fail to compute $\omega_{\ov{p}}(u)$ if there exists a coordinate $u_i\in\RR_{<0}$. This last observation is at the origin of
the Product Convention \ref{conv} introduced in Section \ref{mon_con}

\subsubsection{Automorphic factors of the group $\ca{G}(\RR)$}

The Eisenstein series $G(z,s)$ defined in Section \ref{def_eis} will be modular functions in the variable $z$ with respect to certain discrete subgroups of
$GL_2(K)$ which are commensurable to $GL_2(\ca{O}_K)$. In order to formulate precisely their modularity behavior, we introduce
in this section some standard notations related to arithmetic subgroups and automorphic factors. All of our arithmetic subgroups
will be subgroups of the $\RR$-points of the $\QQ$-algebraic group
$\ca{G}:=(\Res_{\QQ}^K GL_2)$. We let $\ca{G}^{o}$ denote the connected component of the identity of $\ca{G}$,
$\ca{G}_1:=(\Res_{\QQ}^K SL_2)$ and $\ca{U}$ the maximal unipotent subgroup of $\ca{G}$ of upper triangular matrices.
We have the topological group isomorphisms $\ca{G}(\RR)\simeq GL_2(\RR)^g$, $\ca{G}^o(\RR)\simeq GL_2^+(\RR)^g$,
$\ca{G}_1(\RR)\simeq SL_2(\RR)^g$, $\ca{U}(\RR)\simeq\RR^g$, and the obvious inclusions 
$\ca{U}(\RR)\subseteq\ca{G}_1(\RR)\subseteq \ca{G}^o(\RR)\subseteq \ca{G}(\RR)$. The group
$\ca{G}(\RR)$ (resp. $\ca{G}^o(\RR)$) acts on $K_{\CC}^{\pm}$ (resp. on $\mk{h}^{\ov{p}}$ for any $\ov{p}\in\ov{\mk{S}}$) via
M\"obius transformations in the usual way: for $\gamma=(\gamma_j)_{j=1}^g\in \ca{G}(\RR)$ and $z=(z_j)_{j=1}^g\in K_\CC^{\pm}$ we let
$\gamma z=(\gamma_j z_j)_{j=1}^g$. We also have a determinant vector map:
\begin{align*}
\det:\ca{G}(\RR) &\rightarrow GL_1(\RR)^g\\
            (\gamma_i)_i &\mapsto (\det(\gamma_i))_i
\end{align*}
We fix once and for all an embedding $\iota:K\hookrightarrow\RR$. Therefore, we get an inclusion 
$GL_2(K)\subseteq \ca{G}(\RR)$. We also define
\begin{align}\label{amp}
GL_2^+(K):=\{\gamma\in GL_2(K),\sg(\det(\gamma))=\ov{\bz}\}. 
\end{align}
We note that the M\"{o}bius transformation induced by the matrix $\gamma\in GL_2(K)$ will permute the various components of $K_\CC^{\pm}$, 
according to the signature element $\sg(\det(\gamma))\in\ov{\mk{S}}$.
For $z\in K_\CC^{\pm}$ and 
\begin{align*}
\gamma=(\gamma_i)_{i=1}^g=\M{a_i}{b_i}{c_i}{d_i}_{i=1}^g\in \ca{G}(\RR),
\end{align*}
we  define  
the automorphic factor 
$$
j(\gamma,z):=(c_iz_i+d_i))_{i=1}^g\in K_\CC^{\pm}\leq K_\CC^{\times}.
$$
Unlike some authors, we decided not to put the factor $\det(\gamma_i)^{-1/2}$ in the $i$-th coordinate 
of $j(\gamma,z)$, mainly, in order to avoid choosing between one of the two square roots. The application 
$\gamma\mapsto [z\mapsto j(\gamma,z)]$ may be viewed as 
a $1$-cocycle in $Z^1(\ca{G}(\RR),\Maps_{\mbox{\tiny{cont}}}(K_\CC^{\pm},K_{\CC}^{\pm}))$, 
where the set $\Maps_{\mbox{\tiny{cont}}}(K_{\CC}^{\pm},K_{\CC}^{\pm})$ is the set of $K_{\CC}^{\pm}$-valued continuous maps on $K_{\CC}^{\pm}$, 
which is endowed with the following right action of 
$\ca{G}(\RR)$: for $f\in\Maps_{\mbox{\tiny{cont}}}(K_{\CC}^{\pm},K_{\CC}^{\pm})$
and $\gamma\in \ca{G}(\RR)$, $z\in K_{\CC}^{\pm}$, we let $f^{\gamma}(z):=f(\gamma z)$. In particular, we have the usual $1$-cocycle identity 
\begin{align}\label{abeil}
j(\gamma_1\gamma_2,z)=j(\gamma_1,\gamma_2 z)\cdot j(\gamma_2,z),
\end{align}
for all $z\in K_{\CC}^{\pm}$ and $\gamma_1,\gamma_2\in \ca{G}(\RR)$. 

\subsubsection{Modular forms of unitary weight $\{p;s\}$}\label{unit_weight}

Let $p\in\ZZ^g$ be fixed. We endow 
the $\CC$-vector space $\Maps_{\mbox{\tiny{cont}}}(K_{\CC}^{\pm},\CC)$ ($\CC$-valued continuous maps of $K_{\CC}^{\pm}$) 
with the following right $\ca{G}(\RR)$-action: for
$f\in\Maps_{\mbox{\tiny{cont}}}(K_{\CC}^{\pm},\CC)$, $\gamma\in \ca{G}(\RR)$, $z\in K_\CC^{\pm}$ and $s\in\CC$ we set:
\begin{align}\label{slash}
f\big|_{\{p;s\},\gamma}(z):=\omega_p(j(\gamma,z))^{-1}\cdot |\det(\gamma)|^{-\bu\cdot s}\cdot f(\gamma z).
\end{align}
We call $\big|_{\{p;s\}}$ the $\{p;s\}$-slash action.
\begin{Def}
Let $\Gamma\leq \ca{G}(\RR)$ be a subgroup and let $f\in\Maps_{\mbox{\tiny{cont}}}(K_{\CC}^{\pm},\CC)$. If for all $\gamma\in\Gamma$
one has that $f\big|_{\{p;s\},\gamma}(z)=f(z)$, we say that $f$ {\it has unitary weight $\{p;s\}$ relative to the group} $\Gamma$.
\end{Def}
\vspace{-0.5cm}
The word unitary in the above definition refers to the first coordinate of $\{p;s\}$ which corresponds to the integral unitary character $\omega_p$ of $K_{\CC}^{\times}$.
When the value $s$ is clear from the context, we may sometimes say for simplicity that $f(z)$ has unitrary weight $p$ (relative
to $\Gamma$)  rather than weight $\{p\,;s\}$ (relative to $\Gamma$). If the subgroup $\Gamma$ is chosen inside 
$\ca{G}_1(\RR)$, then the factor $ |\det(\gamma)|^{-\bu\cdot s}\equiv 1$, where here the symbol $\equiv$ means \lq\lq identically equal to\rq\rq.  
In this case, the presence of the $s$ coordinate in the notation $\{p;s\}$ is irrelevant and, for this reason, we simply write $f\big|_{\{p\},\gamma}(z)$. 
In Section \ref{real_an2}, the more general notion of a modular function of bi-weight $\{[\alpha,\beta];\mu\}$ (relative to $\Gamma$) will be introduced.

\subsection{Definition of a family of real analytic Eisenstein series}\label{def_eis}
In this section, we introduce the class of $GL_2$-real analytic Eisenstein series that will be studied in this monograph.
We first introduce some auxiliary data on which these will depend.
Let $\mk{m},\mk{n}$ be two lattices of $K$ of maximal rank. In particular, $\mk{m}$ and $\mk{n}$ can be viewed as discrete free $\ZZ$-modules of rank $g=[K:\QQ]$
inside $K_\RR$. Let $p=(p_i)_{i=1}^g\in\mk{S}=\ZZ^g$
be an integral unitary weight and let 
\begin{align*}
U:=\M{u_1}{v_1}{u_2}{v_2}\in M_2(K),
\end{align*}
be a matrix which we call the {\it parameter matrix}. Finally, let $w\in\ZZ$ be an integer which we call the {\it holomorphic parallel 
weight shift}. Taking into account these auxiliary data, this paper aims to study
in detail the Eisenstein series associated to every {\it standard quadruple} $\ca{Q}=((\mk{m},\mk{n}),U,p,w)$. In order to define these Eisenstein series in a
compact way, we need to introduce certain subsets that will be indexing the defining summation of $G(z,s)$. Let
$\mathcal{V}^+:=\mathcal{V}_U^+(\mk{m},\mk{n})$ be the finite index subgroup of $\ca{O}_K^{\times}(\infty)$ which appears 
in Definition \ref{fromage}. According to the definition of $\mathcal{V}^+$, for all $(m,n)\in\mk{m}\oplus\mk{n}$ and all 
$\epsilon\in\mathcal{V}_U^+(\mk{m},\mk{n})$, one has that
\begin{align*}
(\epsilon(m+v_1),\epsilon(n+v_2))\equiv(m+v_1,n+v_2)\pmod{\mk{m}\oplus\mk{n}}.
\end{align*}
In particular, $\mathcal{V}_U^+(\mk{m},\mk{n})$ leaves the set $(\mk{m}+v_1,\mk{n}+v_2)$ stable under the diagonal action of $\ca{V}^+$.
\begin{Def}\label{set_def}
We define $\mathcal{R}$ to be an {\it arbitrarily} chosen complete set of representatives of the set $(\mk{m}+v_1,\mk{n}+v_2)\bs\{(0,0)\}$ under the diagonal
action of $\mathcal{V}^+$. For each representative $(m+v_1,n+v_2)\in\mathcal{R}$ we choose, arbitrarily, 
a matrix  $\gamma_{m,n}:=\M{*}{*}{m+v_1}{n+v_2}\in GL_2(K)$,
and we let $\mathcal{T}:=\{\gamma_{m,n}\in GL_2(K):(m+v_1,n+v_2)\in\mathcal{R}\}$.
\end{Def}

For $a\in\RR$, we define the open right half-plane
\begin{align}\label{plane}
\Pi_a:=\{s\in\CC:\Ree(s)>a\}.
\end{align}
We are now ready to define the Eisenstein series. 
\begin{Def}\label{eis_ser}
Let $\ca{Q}=((\mk{m},\mk{n}),U,p,w)$ be a standard quadruple. 
For $z\in K_\CC^{\pm}$ and
$s\in\Pi_{1-\frac{w}{2}}$, we define
\begin{align}\label{mou1}
&G_{(\mk{m},\mk{n})}^w(U,p\, ;z,s):=
\sum_{\mathcal{R}}
\frac{\omega_{p}((m+v_1)z+(n+v_2))\cdot e^{{2\pi\ii } \Tr(u_1(m+v_1)+u_2(n+v_2))}}{\Norm((m+v_1)z+(n+v_2))^w\cdot|\Norm((m+v_1)z+(n+v_2))|^{2s}}
\cdot|y|^{\bu\cdot s}\\[5mm] \notag
&=\sum_{\mathcal{T}}
\omega_{p}(j(\gamma_{m,n},z))\cdot\Norm(j(\gamma_{m,n},z))^w\cdot e^{2\pi\ii\Tr(u_1(m+v_1)+u_2(n+v_2))} 
|\det(\gamma_{m,n})|^{-\bu\cdot s}\cdot|\Imm(\gamma_{m,n}z)|^{\bu\cdot s}.
\end{align}
\end{Def}
To obtain the second equality from the first, we have used the identity
\begin{align*}
\Imm(\gamma z)=\det(\gamma)\cdot\frac{\Imm(z)}{|j(\gamma,z)|^2},
\end{align*}
which is valid for all $z\in K_{\CC}^{\pm}$ and all $\gamma\in \ca{G}(\RR)$. 
\begin{Rem}
The absolute convergence of the summation on the right-hand side of $\eqref{mou1}$, when $s\in\Pi_{1-\frac{w}{2}}$,
can be proved by comparing $G_{(\mk{m},\mk{n})}^w(U,p\, ;z,s)$ with a finite sum of classical 
real analytic Poincar\'e-Eisenstein series of weight $0$. See Section \ref{groo} where this is explained in more details.
In Theorem \ref{gr_est} of Section \ref{groo}, we also 
give some precise growth estimates of $G_{(\mk{m},\mk{n})}^w(U,p\, ;z,s)$ when $s\in\Pi_{1-\frac{w}{2}}$ and $z$ tends to a cusp.
\end{Rem}
We claim that the definition of $G_{(\mk{m},\mk{n})}^w(U,p\, ;z,s)$
is {\it independent} of the chosen set of representatives $\mathcal{R}$ (or $\mathcal{T}$).
In order to verify this assertion, it is enough to check that the general term of the summation is invariant under the diagonal action of $\mathcal{V}^+=
\mathcal{V}_U^+(\mk{m},\mk{n})$. This is indeed the case. First, note that the image of $\mathcal{V}^+$ inside $(\RR^{\times})^g\subseteq K_\CC^{\times}$ 
gives rise to a lattice of rank $g-1$. By
definition, the integral unitary character $z\mapsto\omega_{p}(z)$ is invariant under $(\RR^{\times})^g$ and, therefore, a fortiori, 
invariant under $\mathcal{V}^+$. Moreover, the quasi-character $z\mapsto \Norm(z)^w|\Norm(z)|^{2s}$  (for a fixed value of $s$) is obviously 
invariant under $\mathcal{V}^+$. Finally, for any $\epsilon\in\mathcal{V}^+$ and $(m,n)\in(\mk{m},\mk{n})$, it follows
from the definition of $\ca{V}^+$ that the absolute trace of the element
\begin{align}\label{obstruc}
&\epsilon u_1(m+v_1)+\epsilon u_2(n+v_2)- u_1(m+v_1)+u_2(n+v_2)\\[2mm]
&=(\epsilon-1)u_1v_1+(\epsilon-1)u_2v_2+(\epsilon-1)u_1m+(\epsilon-1)u_2n,
\end{align}
actually lies in $\ZZ$. It, thus, follows that the right-hand side of the second equality of \eqref{mou1} is indeed independent 
of the set of representatives $\mathcal{R}$. 

\begin{Rem}
On the right hand side of the 
second equality of \eqref{mou1}, the map $(m,n)\mapsto e^{2\pi\ii\Tr(u_1m+u_2 n)}$ is a finite
order character of $\mk{m}\oplus\mk{n}$; the map $z\mapsto \omega_{p}(z)$ is an integral unitary
character of $K_\CC^{\times}$; the map $z\mapsto\Norm(z)^w$ is an integral quasi-character of $K_{\CC}^{\times}$
in the sense of Section \ref{galina}.
\end{Rem}
\begin{Rem}
A priori, the series $G_{(\mk{m},\mk{n})}^w(U,p\, ;z,s)$ makes sense only when $s\in\Pi_{1-\frac{w}{2}}$. However, 
in Section \ref{pars} (see Theorem \ref{lune}), it will be shown that it
admits a meromorphic continuation to all of $\CC$.
\end{Rem}
\begin{Rem}
The Eisenstein series $G_{(\mk{m},\mk{n})}^w(U,p\, ;z,s)$ satisfies many symmetries with respect to the standard quadruple 
$\ca{Q}=((\mk{m},\mk{n}),U,p,w)$. In particular, up to a root of unity, it is possible to view $G_{(\mk{m},\mk{n})}^w(U,p\, ;z,s)$
as a function of the parameter matrix $U$, when $U$ varies in a certain subset of the {\it torus matrix}
$\mathcal{T}_{\mk{m},\mk{n}}$. See Section \ref{symm_a} for more details on that. 
\end{Rem}

\begin{Rem}
In general, if $K$ is a number field with $r_1$ real embeddings and $2r_2$ complex embeddings,
one may associate the symmetric domain $\mk{h}^{r_1}\times \mk{h}_q^{r_2}$ (no more hermitian if $r_2\geq 1$) 
to the $\QQ$-algebraic group $\Res_\QQ^{K}(SL_2)$. Here $\mk{h}_q$ stands for the quaternion upper half-space; see \cite{Asai}
for more details on that. It is possible to construct Eisenstein series (similar to the one constructed in this paper) 
on such symmetric domains based on our explicit approach. However, the author did not work out all the details. In \cite{Ch-Kl}, 
using a different approach, the meromorphic continuation of this larger class of Eisenstein series is proved.
\end{Rem}

\subsubsection{$\ca{V}^{+}$-integral quasi-characters of $K_{\CC}^{\times}$}\label{galina}
In general, it is possible to twist the {\it basic function}
\begin{align}\label{bami}
|\det(\gamma_{m,n})|^{-\bu\cdot s}\cdot|\Imm(\gamma_{m,n}z)|^{\bu\cdot s}
\end{align} 
(for $\gamma_{m,n}\in\mathcal{T}$) which appears in the defining 
summation \eqref{mou1} of $G_{(\mk{m},\mk{n})}^w(U,p\, ;z,s)$, by a quasi-character $\chi\in\wh{K_{\CC}^{\times}}$ of
the form $\chi(z)=\omega_p(z)\cdot\Norm(z)^w\cdot\eta_m(z)$,
where $p\in\ZZ^g$, $w\in\ZZ$ and $m\in\ZZ^{g-1}$. Here, the character $\eta_m(z)$ may be thought of as being induced by the infinite part
of the restriction of a Hecke Gr\"o\ss encharacter to the set of principal fractional ideals of $K$. 
Following the presentation on p. 209 of \cite{Asai},
the character $\eta_m(z)$ may be given explicitly by the formula
\begin{align*}
\eta_m(z)=\prod_{j=1}^{g}|z_j|^{2\pi\ii\sum_{k=1}^{g-1} m_k e_{j}^{\laa k\raa}}, 
\end{align*}
where the $e_{j}^{\laa k\raa}$ are given by
\begin{align*}
\left(
\begin{matrix}
\frac{1}{g} & \frac{1}{g} & \ldots & \frac{1}{g}  \\[2mm]
 e_1^{\laa 1\raa} & e_{2}^{\laa 1\raa} & \ldots & e_{g}^{\laa 1\raa}\\[2mm]
 \vdots           & \vdots & \vdots             & \vdots \\[2mm]
 e_1^{\laa g-1\raa}  & e_{2}^{\laa g-1\raa} & \ldots    & e_{g}^{\laa g-1\raa}
\end{matrix}
\right)=
\left(
\begin{matrix}
1 & \log|\epsilon_1^{(1)}| & \ldots & \log|\epsilon_{g-1}^{(1)}|  \\[2mm] 
1 & \log|\epsilon_1^{(2)}| & \ldots & \log|\epsilon_{g-1}^{(2)}|\\[2mm]
 \vdots           & \vdots & \vdots             & \vdots \\[2mm]
1  & \log|\epsilon_1^{(g)}| & \ldots & \log|\epsilon_{g-1}^{(g)}|
\end{matrix}
\right)^{-1},
\end{align*}
and $\{\epsilon_1,\ldots,\epsilon_{g-1}\}$ is a $\ZZ$-basis of $\mathcal{V}^+\subseteq K^{\times}\subseteq K_\CC^{\times}$. Note that 
the quasi-character $\chi=\omega_p\cdot\Norm(z)^w\cdot\eta_m(z)$ may be viewed as an element $\wh{(K_{\CC}^{\times}/\mathcal{V}^+)}$. 
\begin{Def}
We define
\begin{align*}
\mk{X}_{\ca{V}^+}:=\{\omega_p\cdot\Norm^w\cdot\eta_m\in \wh{K_{\CC}^{\times}}:p\in\ZZ^g,w\in\ZZ,m\in\ZZ^{g-1}\},
\end{align*}
and call the elements of $\mk{X}_{\ca{V}^+}$, {\it $\mathcal{V}^+$-integral quasi-characters} of $K$.  
\end{Def}
Note that by definition, for all $\ca{V}^+$, we always have $\mk{X}_0\leq \mk{X}_{\ca{V}^+}$.
One can show that there exists a topological group isomorphism
\begin{align*}
\varphi:\ZZ^{2g-1}\times\CC \stackrel{\simeq}{\longrightarrow}    \wh{(K_\CC^{\times}/\mathcal{V}^+)},
\end{align*}
normalized so that $s\mapsto\varphi((\bz,s))=|\Norm(\_)|^s$, for all $s\in\CC$.
In particular, the group $\mk{X}_{\ca{V}^+}$ may be viewed as a {\it certain lattice} of rank $2g$ inside 
$\wh{(K_\CC^{\times}/\mathcal{V}^+)}\simeq\ZZ^{2g-1}\times\CC$.

Note that the map $s\mapsto\varphi((\bz,s))=|\Norm(\_)|^s$, for $s\in\CC$, 
gives rise to a one (complex) dimensional family of quasi-characters of $K_\CC^{\times}/\mathcal{V}^+$. Taking
this point of view, we see that the basic function which appears in \eqref{bami} can be rewritten as
\begin{align}\label{bami2}
|\det(\gamma_{m,n})|^{-\bu\cdot s}\cdot|\Imm(\gamma_{m,n}z)|^{\bu\cdot s}=|\Imm(z)|^{\bu\cdot s}\cdot|\Norm(j(\gamma_{m,n},z))|^{-2s}.
\end{align}
Looking at the identity \eqref{bami2}, we thus see that the the Eisenstein series in \eqref{mou1} is equal
to $|\Imm(z)|^{\bu\cdot s}$ times the average over the set $\{j(\gamma_{m,n},z)\in K_{\CC}^{\pm}\}_{\gamma_{m,n}\in\ca{T}}$  
of the following {\it continuous family} of quasi-characters of $K_\CC^{\times}/\mathcal{V}^+$:
\begin{align*}
\left\{\omega_p(\_)\cdot\Norm^{-w}(\_)\cdot|\Norm(\_)|^{-2s}\right\}_{s\in\Pi_{1-\frac{w}{2}}}.
\end{align*}
In the present work, we have only considered Eisenstein series twisted by the $\ca{V}^+$-integral 
quasi-characters of the form $\omega_p\cdot\Norm^w$. The reason behind this restriction is twofold. 
First, it simplifies some of the formulas which are already a bit involved. Second, our method
of the proof of the meromorphic continuation of $[s\mapsto G_{(\mk{m},\mk{n})}^w(U,p\, ;z,s)]$ requires a priori
the knowledge of the meromorphic continuation of the following zeta function (see Appendix \ref{app_1a}):
\begin{align}\label{toil}
Z_V(a,b;\omega_{\ov{p}}\eta_m\Norm^w;s):=[\ca{O}_K:V]^s\sum_{\substack{a+v\in\mathfrak{R}\\ a+v\neq 0}}
\omega_{\ov{p}}\eta_m\Norm^w(a+v)\cdot\frac{e^{{2\pi\ii } \Tr_{K/\QQ}(b(a+v))}}{|\Norm_{K/\QQ}(a+v)|^{s}},
\end{align}
to all of $\CC$. Recall that series \eqref{toil} converges absolutely only when $s\in\Pi_{1+w}$. 
In \cite{Ch4}, partial zeta functions (twisted by characters of $\mk{X}_{\ca{V}^+}$) were constructed {\it only} 
for the characters in $\mk{X}_0$ (i.e., the $\omega_{\ov{p}}$'s) rather than 
the more general class of quasi-characters in $\mk{X}_{\ca{V}^+}$. 
It is likely to assume that the technique used in \cite{Ch4} can also be used 
to handle such quasi-characters, although, the author did not check all the details. Note here that the presence of
the quasi-character $\Norm^w$ creates no additional difficulty, since it simply involves a shift in the arguments of $s$ as the 
following formula indicates:
\begin{align*}
Z_V(a,b;\omega_{\ov{p}}\eta_m\Norm^w;s)=Z_V(a,b;\omega_{\ov{p}+w\cdot\ov{\bu}}\eta_m;s-w).
\end{align*}
\begin{Rem}
In \cite{Ch-Kl}, a different method (based on ideas of Colin de Verdi\`ere)
is used to show the meromorphic continuation of Eisenstein series twisted by {\it any} quasi-character in $\mk{X}_{\ca{V}^+}$
and, in the same paper, it is explained how this meromorphic continuation also leads to the meromorphic continuation and functional equation 
of the corresponding partial zeta functions.
\end{Rem}
\begin{Rem}
Using an idea of Colmez, we explain in Section \ref{Colmez_trick} how the explicit determination of the {\it non-constant} Fourier series coefficients 
$(a_{\xi}(y,s))_{\xi\in\ca{L}^*\bs\{0\}}$ of $[z\mapsto G_{(\mk{m},\mk{n})}^w(U,p\, ;z,s)]$ is enough to prove, a posteriori,
the meromorphic continuation and the functional equation of $[s\mapsto Z_V(a,b;\omega_{\ov{p}};s)]$. We also expect this approach to
generalize to all quasi-characters in $\mk{X}_{\ca{V}^+}$.
\end{Rem}
\subsubsection{Transformation formula of Eisenstein series under M\"{o}bius transformations}

We would like now to describe the transformation formula under which the Eisenstein series goes 
under a M\"{o}bius transformation of $GL_2(K)$. In order to simplify the presentation we will assume that $w=0$. Under this assumption, 
not much is lost in virtue of the following identity:
\begin{align}\label{tram}
G_{(\mk{m},\mk{n})}^{w}(U,p;z;s)=\frac{G_{(\mk{m},\mk{n})}^{0}(U,p-w\cdot\bu;z,s+\frac{w}{2})}{|y|^{\frac{w}{2}\cdot\bu}}.
\end{align}
In order to treat in a satisfactory way the case when $w\neq 0$, it is preferable (but not essential) to introduce the more
general notion of modular forms of bi-weight $[\alpha,\beta]$, where $\alpha,\beta\in\CC^g$ are weights subjected to the condition 
$\alpha-\beta\in\ZZ^g$. Here the weight 
$\alpha$ corresponds to the holomorphic weight (i.e., for the variable $z$) and $\beta$ to the anti-holomorphic weight 
(i.e., for the variable $\ov{z}$). This 
is accomplished in Section \ref{real_an2}. 

It is proved in Proposition \ref{pleut},  that the function $[z\mapsto G_{(\mk{m},\mk{n})}^0(U,p\, ;z,s)]$
satisfies the following transformation formula:
\begin{align}\label{tr_for}
G_{(\mk{m},\mk{n})}^0\big|_{\{-p,s\},\gamma}(U,p\,;z,s)=f_{\gamma}\cdot G_{(\mk{m},\mk{n})\gamma}^{0}(U^\gamma,p,z,s),
\end{align}
for all $\gamma\in GL_2(K)$ and $z\in K_\CC^{\pm}$. The map $(\mk{m},\mk{n})\mapsto(\mk{m},\mk{n})\gamma$ corresponds to the natural 
right action of $GL_2(K)$ on row vectors. However, the map $U\mapsto U^{\gamma}$ corresponds to the right action which appears in Definition \ref{chan} (this
right action {\it does not} correspond to the right matrix multiplication of $U$ by $\gamma$). We call it the {\it upper right action}.
The quantity $f_{\gamma}$ is a positive rational number which depends only on $U$, $\gamma$ and $(\mk{m},\mk{n})$ (see Definition \ref{chien} and \eqref{indexi}).

From the transformation formula \eqref{tr_for}, it is reasonable to expect that the function 
\begin{align*}
[z\mapsto G_{(\mk{m},\mk{n})}^0(U,p\,;z,s)],
\end{align*}
is a modular with respect to a suitable subgroup of $GL_2(K)$. This is indeed the case.
One may show that there exists a discrete subgroup $\Gamma_U(\mk{m},\mk{n})\leq GL_2(K)\leq \ca{G}(\RR)$
(see Definition \ref{grp}), such that for all $\gamma\in \Gamma_U(\mk{m},\mk{n})$
\begin{align}\label{tr0}
G_{(\mk{m},\mk{n})}^0\big|_{\{-p\,;s\},\gamma}(U,p\,;z,s)=\zeta_\gamma\cdot G_{(\mk{m},\mk{n})}^0(U,p\, ;z,s),
\end{align}
where $\zeta_{\gamma}$ is some explicit root of unity. The proof of \eqref{tr0} is a consequence of Proposition \ref{mod_for}.
In order to get rid of this root of unity $\zeta_{\gamma}$, one needs to
restrict the slash action $\big|_{\{p\,;s\},\gamma}$ to a smaller subgroup of $\Gamma_U(\mk{m},\mk{n})$. For the rest of the introduction, let us assume that 
$U\in\frac{1}{N}\M{\mk{m}^*}{\mk{m}}{\mk{n}^*}{\mk{n}}$ for some $N\in\ZZ_{\geq 1}$. 
We let $\Gamma_U(\mk{m},\mk{n};N)\leq \Gamma_U(\mk{m},\mk{n})$ be the discrete subgroup appearing in Definition \ref{grp}.
Again, it follows from Proposition \ref{mod_for}, that for all $\gamma\in \Gamma_U(\mk{m},\mk{n};N)$, one has that
\begin{align}\label{tr2}
G_{(\mk{m},\mk{n})}^0\big|_{\{-p\,;s\},\gamma}(U,p\,;z,s)=G_{(\mk{m},\mk{n})}^0(U,p\, ;z,s).
\end{align}
Therefore, the functions $[z\mapsto G_{(\mk{m},\mk{n})}^0(U,p\, ;z,s)]$ may be viewed as a real analytic modular form 
of unitary weight $\{-p\,;s\}$ relative to the discrete group $\Gamma_U(\mk{m},\mk{n};N)$.
\begin{Rem}
Suppose that there exists $\epsilon\in \mathcal{V}_U(\mk{m},\mk{n})\leq \ca{O}_K^{\times}$ 
(see Definition \ref{fromage}; the units in $\mathcal{V}_U(\mk{m},\mk{n})$ are no more required to be totally positive)
such that $\omega_{\ov{p}}(\epsilon)\Norm(\epsilon)^k=-1$. Considering the matrix $\M{\epsilon}{0}{0}{1}$,
one may deduce from \eqref{tr2}, that $G_{(\mk{m},\mk{n})}^0(U,p\, ;z,s)$  is identically equal to $0$ 
(cf. Remark \ref{vanish}). It is therefore important to assume that no such unit exists in order to avoid working with the zero function. 
The author does not know if the non-existence
of such units automatically implies that $[(z,s)\mapsto G_{(\mk{m},\mk{n})}^0(U,p\, ;z,s)]\not\equiv 0$.
\end{Rem}

\subsubsection{Eisenstein series as $(I',I'')$-fold differential forms}
In this section, we would like to give a \lq\lq differential form interpretation\rq\rq\s of the Eisenstein series
$[z\mapsto G_{(\mk{m},\mk{n})}^0(U,p\, ;z,s)]$ when the weight $p$ is even. 
First, note that the natural image of $\Gamma:=\Gamma_U(\mk{m},\mk{n};N)$ inside $PGL_2(K)\leq PGL_2(\RR)^g$, which we denote by
$\ov{\Gamma}$, is a discrete subgroup which acts {\it faithfully} and {\it properly discontinuously} on $K_{\CC}^{\pm}$ (see Appendix \ref{app_5}).
In particular, the quotient space $K_\CC^{\pm}/\ov{\Gamma}$ is
an orbifold. Note that the orbifold $K_\CC^{\pm}/\ov{\Gamma}$ will be connected if and only if $\sg(\det(\Gamma))=\ov{\mk{S}}$. 
Now, let $\Gamma_1:=\Gamma\cap \ca{G}_1(\RR)$ and denote by $\ov{\Gamma}_1$ the image of $\Gamma_1$ inside $PSL_2(\RR)^g$.
The group $\ov{\Gamma}_1$ acts faithfully and properly discontinuously on the connected component $\mk{h}^{\ov{\bz}}=\mk{h}^g$
of $K_\CC^{\pm}$.
For $p=(p_j)_{j=1}^g\in\ZZ^g$, we let $V(p;\ov{\Gamma}_1)$ be the $\CC$-vector space of all 
real analytic meromorphic functions on $\mk{h}^g$ 
which satisfy the formula \eqref{slash} under all $\gamma\in\ov{\Gamma}_1$, i.e., modular forms of
unitary weight $p$ relative to $\ov{\Gamma}_1$. Let $I:=\{1,2,\ldots,g\}$ and consider a fixed partition $I=I'\bigsqcup I''$ of $I$ into two sets.
Consider the following \lq\lq real analytic $(I',I'')$-fold form on $\mk{h}^g$\rq\rq:
\begin{align*}
\eta(z):=f(z) \cdot\left(\bigotimes_{i\in I'}\left(\frac{dz_i}{y_i}\right)^{\otimes p_i}\right)
\cdot\left(\bigotimes_{i\in I''}\left(\frac{d\ov{z}_i}{y_i}\right)^{\otimes p_i}\right).
\end{align*}
Such expressions may be viewed as elements of $\Sym_\CC^*(A^1(\mk{h}^g),\CC)$, where $A^{1}(\mk{h}^g,\CC)$ denotes 
the smooth differential $1$-forms on $\mk{h}^g$ and $\Sym_\CC^*$ corresponds to the symmetric algebra functor on $\CC$-vector spaces. 
For each coordinate $z_i=x_i+\ii y_i$ of $z\in\mk{h}^g$ and $\gamma=\M{a}{b}{c}{d}\in \ca{G}_1(\RR)$, one has the transformation formula
\begin{align}\label{tapi}
\frac{d(\gamma z_i)}{\Imm(\gamma z_i)}=\frac{|cz_i+d|^2}{(cz_i+d)^2}\cdot\frac{dz_i}{\Imm(z_i)}.
\end{align}
From \eqref{tapi}, it follows that $\eta(z)$ descends to a differential form on the orbifold $Y_{\ov{\Gamma}_1}:=\mk{h}^g/\ov{\Gamma}_1$,
if and only if $f(z)$ lies in the vector space 
$V(\ov{\Gamma}_1,\wt{p})$, where $\wt{p}$ is such that $\wt{p}_j=2p_j$ if $j\in I'$,
and $\wt{p}_j=-2p_j$ if $j\in I''$. So, in general, it is only when the weight $p$ is {\it even} (i.e., $p\in(2\ZZ)^g$) 
that $[z\mapsto G_{(\mk{m},\mk{n})}^0(U,p\, ;z,s)]$  may be thought of as $(I',I'')$-fold differential forms on $Y_{\ov{\Gamma}_1}$.

\subsection{Main results}\label{sec_beg}
In this section we present the four main results of this manuscript. 

\noindent {\bf A.} Our first main result is Theorem \ref{clop} which provides an analytic characterization of certain 
real analytic families $\{G(z,s)\}_{s\in\Pi_1}$  of modular forms of a {\it fixed} integral unitary weight $p\in\ZZ^g$. 
Since the statement of the theorem
is itself a bit technical, let us just mention that functions $\{G(z,s)\}_{s\in\Pi_1}$ are required to satisfy 
essentially two properties:
\begin{enumerate}[(1)]
 \item  they must be eigenvectors with respect to the partial-graded Laplacians associated to $p$ (see Section \ref{def_mas}), 
 \item they must satisfy some growth conditions at the cusps.
\end{enumerate}

\noindent {\bf B.} The second main result of this work is the explicit computation of the Fourier series expansion of 
$[z\mapsto G_{(\mk{m},\mk{n})}^0(U,p\, ;z,s)]$
at the cusp $\infty:=\frac{1}{0}$ (Theorem \ref{key_thmm}) and a proof of the meromorphic continuation of $[s\mapsto G_{(\mk{m},\mk{n})}^0(U,p\, ;z,s)]$
(Theorem \ref{lune}). For the precise meaning of the Fourier series expansion at a cusp $c\in\PP^1(K)$ 
we refer the reader to Section \ref{Fourier}.
The \textit{non-constant} Fourier coefficients of $G_{(\mk{m},\mk{n})}^0(U,p\, ;z,s)$ can be written in 
terms of the {\it Tricomi's confluent hypergeometric function} (see equations \eqref{tornade} and \eqref{trico} for the precise definitions) that was considered in \cite{Sie5}. 
On the other hand, the \textit{constant term} of the Fourier series expansion of $G_{(\mk{m},\mk{n})}^0(U,p\, ;z,s)$,
at the cusp $\infty$, can be written in terms of a special class of zeta functions associated to $K$ that was studied in \cite{Ch4} and \cite{Ch8}. 
Let $\ov{p}:=p\pmod{2}\in\ov{\mk{S}}=(\ZZ/2\ZZ)^g$. For a lattice $\mk{n}\subseteq K$, $u,v\in K$ and $s\in\CC$ with $\Ree(s)>1$, we define the {\it partial zeta
function}
\begin{align*}
Z_{\mk{n}}(u,v;\omega_{\ov{p}};s):=[\ca{O}_K:\mk{n}]^s\sum_{\left\{0\neq(n+u)\in(\mk{n}+u)\right\}/\mathcal{V}_{u,v,\mk{n}}^+}
\frac{\omega_{\ov{p}}(n+u)e^{{2\pi\ii } \Tr(v(n+u))}}{|\Norm_{K/\QQ}(n+u)|^{s}},
\end{align*} 
and, its associated {\it normalized zeta function}, 
\begin{align*}
\wt{Z}_{\mk{n}}(u,v,\omega_{\ov{p}};s):=\frac{Z_{\mk{n}}(u,v;\omega_{\ov{p}};s)}{[\ca{O}_K:\mk{n}]^s},
\end{align*}
where $\mathcal{V}_{u,v,\mk{n}}^+$  is the finite index subgroup of $\ca{O}_K^{\times}(\infty)$ which appears in Definition \ref{unit_def}.
We would like to emphasize that function $[s\mapsto\wt{Z}_{\mk{n}}(u,v,\omega_{\ov{p}};s)]$ could be identically equal to zero. See
Remark \ref{vanish} for an example where such a vanishing happens.

The constant term of the Fourier series expansion of 
$G_{(\mk{m},\mk{n})}^0(U,p\, ;z,s)$ at the cusp $\infty$ is given explicitly by the formula
\begin{align}\label{fati}
& e_1\cdot \delta_{\mk{m}}(v_1)\cdot \wt{Z}_{\mk{n}}(v_2,u_2,\omega_{\ov{p}};2s)\cdot |y|^{\bu\cdot s}\\[2mm] \notag
& \hspace{1cm}+e_2\cdot \delta_{\mk{n}^*}(u_2)\cdot\theta_0\cdot 2^{g(1-2s)}\cdot\varphi_{p}(1-s)\cdot\psi_{p}(s)
\cdot \wt{Z}_{\mk{m}}(v_1,u_1,\omega_{\ov{p}};2s-1)\cdot |y|^{\bu\cdot (1-s)},
\end{align}
where 
\begin{enumerate}
 \item $\varphi_{p}(s)$ and  $\psi_{p}(s)$ are some explicit products involving the gamma function and the fourth root of unity 
 $\sqrt{-1}=\ii$, (see Appendix \ref{app_1e} for the definitions of $\varphi_{p}(s)$ and $\psi_{p}(s)$ ),
 \item $\theta_0=\theta_0((\mk{m},\mk{n}),U,p;z):=\cov(\mk{n})^{-1}e^{2\pi\ii\Tr(u_2v_2)}\cdot(2\pi)^g\cdot\ii^{\Tr(p)}\cdot (-1)^{\Tr(\sg(z)\cdot\ov{p})}$ is a constant,
 \item $\delta_{\mk{m}}(v_1)=1$ if $v_1\in\mk{m}$ and $0$ otherwise; and $\delta_{\mk{n}^*}(u_2)=1$ if
$u_2\in\mk{n}^*$ and $0$ otherwise. 
\item The factors $e_1:=e_1((\mk{m},\mk{n}),U)$ and $e_2:=e_2((\mk{m},\mk{n});U)$ are positive integers measuring the discrepancy 
between certain subgroups of units of $\ca{O}_K^{\times}(\infty)$ (see Definition \ref{chien}).
\end{enumerate}
The expression
\eqref{fati} is what we call {\it the constant term of the Fourier series expansion of $G_{(\mk{m},\mk{n})}^0(U,p\, ;z,s)$  
at the cusp $\infty$ (relative to the arithmetic group $\Gamma_U(\mk{m},\mk{n})$)}.
The reader should be careful here: the terminology \lq\lq constant term\rq\rq\; is a bit abusive,
since it really depends on the variables $s$ and $\Imm(z)$. It is also convenient to define
\begin{enumerate}[(i)]
 \item $\spart{\infty}{s}G_{(\mk{m},\mk{n})}^0(U,p\, ;z,s):=e_1\cdot \delta_{\mk{m}}(v_1)\cdot \wt{Z}_{\mk{n}}(v_2,u_2,\omega_{\ov{p}};2s)$,
 \item $\spart{\infty}{1-s}G_{(\mk{m},\mk{n})}^0(U,p\, ;z,s):=\\[2mm] \mbox{\hspace{4cm}}
 e_2\cdot \delta_{\mk{n}^*}(u_2)\cdot\theta_0\cdot 2^{g(1-2s)}\cdot\varphi_{p}(1-s)\cdot\psi_{p}(s)
 \cdot \wt{Z}_{\mk{m}}(v_1,u_1,\omega_{\ov{p}};2s-1)$.
\end{enumerate}
We would like to emphasize here that the function 
$\psi_{p}(s)$ and $\varphi_{p}(s)$ only depend on the weight parameter $p\in\ZZ^g$
and not on the data $U,\mk{m},\mk{n}$. Moreover, the expression 
$\psi_{p}(s)$ {\it depends only} on $\ov{p}\in\ov{\mk{S}}$ rather than $p$ itself.
Therefore, it is legitimate to write $\psi_{\ov{p}}(s)$.
However, the function $\varphi_{p}(s)$ really depends on $p$ itself and not just on $\ov{p}$. 
It follows from these observations, that the constant term of the Fourier
series expansion of $[z\mapsto G_{(\mk{m},\mk{n})}^0(U,p\, ;z,s)]$, up to the function $\varphi_{p}(s)$, 
\lq\lq only sees\rq\rq\s the weight $p$ modulo $2$.
\begin{Rem}
From \eqref{fati}, it follows from the definitions of $\delta_{\mk{m}}(v_1)$ and $\delta_{\mk{n}^*}(u_2)$ 
that the constant term of the Fourier series of $G_{(\mk{m},\mk{n})}^0(U,p\, ;z,s)$
at the cusp $\infty$ vanishes identically in the variables $s$ and $\Imm(z)$, if $v_1\notin\mk{m}$ and $u_2\notin\mk{n}^*$.
\end{Rem}
\begin{Rem}
Let $\ov{q}\in\ov{\mk{S}}$. If we replace $z$ by $z^{\ov{q}}$ in \eqref{fati}, 
then the second term of \eqref{fati} gets multiplied by  $(-1)^{\Tr(\ov{q}\cdot\ov{p})}$. 
It follows from this observation that, up to the sign, the dependence of $\theta_0((\mk{m},\mk{n}),U,p;z)$
on $z$ is only a dependence on $\sg(z)=\sg(y)$.
\end{Rem}

Let $u,v\in K^{\times}$. In Appendices \ref{app_1a} and \ref{app_1b}, it is explained how a previous result of the author, proved in \cite{Ch4}, implies that 
\begin{align}\label{lor}
\wt{Z}_{\mk{n}}(u,v,\omega_{\ov{p}};s)=\cov(\mk{n}^*)\cdot e^{2\pi\ii\Tr_{K/\QQ}(uv)}\cdot\lambda_{\ov{p}}(1-s)\cdot
\wt{Z}_{\mk{n}^*}(-v,u,\omega_{\ov{p}};1-s), 
\end{align}
for some explicit function $\lambda_{\ov{p}}(s)$, where $\ov{p}\in\ov{\mk{S}}$. Note that since 
\begin{align}\label{clement}
Z_{\mk{n}}(-u,-v;\omega_{\ov{p}};s)=(-1)^{\Tr(\ov{p})}\cdot Z_{\mk{n}}(u,v;\omega_{\ov{p}};s),
\end{align}
the functional equation \eqref{lor} implies the symmetry
\begin{align}\label{tyyr}
\lambda_{\ov{p}}(1-s)=(-1)^{\Tr(\ov{p})}\lambda_{\ov{p}}(s)^{-1}.
\end{align}
For an integral weight $p\in\ZZ^g$, it is convenient to define $\lambda_p(s):=\lambda_{\ov{p}}(s)$. It is explained in Appendix \ref{app_1e}, that the following 
relations hold true
\begin{enumerate}
 \item[(1.a)]  $\varphi_{p}(1-s)=(-1)^{\Tr(p)}\cdot\varphi_{p}(s)^{-1}$,
 \item[(2.a)] $\varphi_{-p}(s)=\varphi_p(s)$,
 \item[(3.a)] $\psi_p(s)\cdot(2\pi)^{2g(1-s)}=\lambda_p(2s-1)$.
 \item[(4.a)] $\varphi_p(1-s)\cdot\psi_p(s)=\prod_{j=1}^g\frac{\Gamma(2s-1)}{\Gamma(s-p_j/2)\Gamma(s+p_j/2)}$.
 \end{enumerate}

\noindent {\bf C.} The third main result of this work (see Theorem \ref{nice_fn_eq} below) is a proof of a functional equation for the completed Eisenstein series
$\wh{G}_{\mk{m},\mk{n}}^0(U,p\,;z,s)$ (see Section \ref{pars} for the precise definition). In order
to state the functional equation, we need to introduce the \lq\lq Cartan involution\rq\rq\; $*$ of the
$K$-algebra $M_2(K)$. For $U=\M{u_1}{v_1}{u_2}{v_2}$, let
\begin{align}\label{ulu}
\M{u_1}{v_1}{u_2}{v_2}^*:=\M{v_2}{-u_2}{-v_1}{u_1}.
\end{align}
The Cartan involution, satisfies some nice compatibilities properties with respect to the upper
right action, see Proposition \ref{compa}. For a matrix $U=\M{a}{b}{c}{d}\in M_2(K)$, we also need to define 
\begin{align*}
\ell_U:=u_1v_1+u_2v_2.
\end{align*}
Direct computations show that 
\begin{enumerate}
 \item[(1.b)] $\ell_{U^*}=-\ell_U$ and,
 \item[(2.b)] for all $\gamma\in SL_2(K)$, $\ell_{U^{\gamma}}=\ell_U$.
\end{enumerate}

We may now state the functional equation.

{\bf Theorem \ref{nice_fn_eq}}
{\it
Let $\ca{Q}=((\mk{m},\mk{n}),U,p,0)$ be a standard quadruple. Let
\begin{align*}
\wh{G}_{(\mk{m},\mk{n})}^0(U,p\,;z,s):=C(\alpha(s),\beta(s);z)\cdot G_{(\mk{m},\mk{n})}^0(U,p\,;z,s),
\end{align*}
be the completed Eisenstein series, where $C(\alpha(s),\beta(s);z)$ is the \lq\lq Euler factor at infinity\rq\rq\;\;
as defined in Section \ref{pars}. Then
\begin{align}\label{tiry}
\wh{G}_{(\mk{m},\mk{n})}^0(U,p\,;z,s)= (-1)^{\Tr(p)}\cdot e^{2\pi\ii\Tr(\ell_U)}\cdot\frac{\cov(\mk{n}^*)}{\cov(\mk{m})}\cdot\wh{G}_{(\mk{n}^*,\mk{m}^*)}^0(U^*,p\,;z,1-s).
\end{align}
}

\begin{Rem}
According to the definition of the Euler factor, the dependence of $C(\alpha(s),\beta(s);z)$ on the parameter $z$ 
is only a dependence on $y=\Imm(z)$ rather than on $z$ itself.
\end{Rem}
\begin{Rem}
The functional equation \eqref{tiry} may be viewed as a natural generalization of the functional equation,
which appears in the work of Asai (see \cite{Asai}), at least, in the case where $K$ is totally real.
\end{Rem}

\begin{Rem}
The reason for which we decided to formulate the functional equation in terms of the completed 
Eisenstein series $\wh{G}_{(\mk{m},\mk{n})}^0(U,p\,;z,s)$ rather than the uncompleted one $G_{(\mk{m},\mk{n})}^0(U,p\,;z,s)$,
is mainly aesthetic. See Section \ref{der}, where we give a functional equation for the uncompleted Eisenstein
series $G_{(\mk{m},\mk{n})}^w(U,p\,;z,s)$, where $w$ is not necessarily equal to $0$.
\end{Rem}

For some applications of the functional equation \eqref{tiry},  
see \cite{Ch10} and \cite{Ch12}. Let us explain briefly the idea behind the two proofs
of Theorem \ref{nice_fn_eq} which are given in Section \ref{fun_eqn}. Both proofs have a common core which uses, in an essential way, 
the explicit description of the constant term of $\wh{G}_{\mk{m},\mk{n}}^0(U,p\,;z,s)$. To simplify the notation, let
\begin{align*}
F(z,s):=\wh{G}_{(\mk{m},\mk{n})}^0(U,p\,;z,s)-(-1)^{\Tr(p)}\cdot e^{2\pi\ii\Tr(\ell_U)}\cdot\frac{\cov{(\mk{n}^*)}}{\cov(\mk{m})}\cdot\wh{G}_{(\mk{n}^*,\mk{m}^*)}^0(U^*,p\,;z,1-s).
\end{align*}
The functional equation is equivalent to show that $[(z,s)\mapsto F(z,s)]\equiv 0$.
The first proof is based on the fact that $F(z,s)$ is {\it square-integrable} at all cusps. 
The proof of the square-integrability is accomplished by combining the explicit description of the 
$\spart{\sigma}{s} \wh{G}_{(\mk{m},\mk{n})}^0(U,p\,;z,s)$, 
at all $\sigma\in\PP^1(K)$, with the functional equation for partial zeta functions which was proved in \cite{Ch4}. 
Then, using our analytic characterization theorem (Theorem \ref{clop}), one may
conclude that $[(z,s)\mapsto F(z,s)]\equiv 0$. Our second proof is based on a direct calculation.
In particular, the functional equation proved in \cite{Ch4} implies that $\spart{\infty}{s}F(z,s)$ and
 $\spart{\infty}{1-s}F(z,s)$ are identically equal to $0$. Finally, the vanishing of the non-constant Fourier coefficients
of $[z\mapsto F(z,s)]$, at the cusp $\infty$, relies essentially on some adroit calculations, and the Euler reflection type formulas proved in Section \ref{ref_for}. 

It is fair to say that our second proof of the analytic continuation of 
$[s\mapsto G_{(\mk{m},\mk{n})}^0(U,p\,;z,s)]$ and of Theorem \ref{nice_fn_eq} is classical and relies ultimately on brute force calculations.
Though, one advantage of this very explicit approach is that it provides, in one stroke, a Kronecker limit formula for the function
$G_{(\mk{m},\mk{n})}^0(U,p\,;z,s)$, i.e., an explicit description the constant term of the Taylor series development
of $[s\mapsto G_{(\mk{m},\mk{n})}^0(U,p\,;z,s)]$ at $s=1$, see \cite{Ch12}. It seems desirable to find a proof of the  meromorphic continuation 
of $\wh{G}_{(\mk{m},\mk{n})}^0(U,p\,;z,s)$ and of its functional equation \eqref{tiry} which {\it does not assume} the following two facts:
\begin{enumerate}[(1)]
 \item  an a priori knowledge 
of the meromorphic continuation in $s$ of $\spart{\infty}{s}G_{(\mk{m},\mk{n})}^0(U,p\,;z,s)$ and $\spart{\infty}{1-s}G_{(\mk{m},\mk{n})}^0(U,p\,;z,s)$, and
 \item an explicit description of the non-constant Fourier series coefficients of $[z\mapsto G_{(\mk{m},\mk{n})}^0(U,p\,;z,s)]$.
\end{enumerate}
Such a proof is worked out in \cite{Ch-Kl} and is 
essentially based on the ideas of Colin de Verdi\`ere introduced in \cite{Ver1}).
In \cite{Ch-Kl}, the first step consists in proving the meromorphic continuation of 
$[s\mapsto G_{(\mk{m},\mk{n})}^0(U,p\,;z,s)]$ by 
applying the Fredholm analyticity theorem to a suitable family of compact operators obtained from a truncation of a graded Laplacian. In a second step, 
the functional equation \eqref{tiry} is deduced from the {\it analytic characterization} of $G_{(\mk{m},\mk{n})}^0(U,p\,;z,s)$ which 
is provided by Theorem \ref{clop} of this monograph. Furthermore, it is also explained in \cite{Ch-Kl} how this approach leads to a 
different proof of the functional equation \eqref{lor} which was proved for the first time in \cite{Ch4}.

{\bf D.} Finally, using an idea of Colmez, we explain in Section \ref{Colmez_trick}, how the knowledge of (2) is sufficient to prove (1).
In essence the idea is simple: From the explicit description of the {\it non-zero} Fourier coefficients $\{a_{\xi}(y,s)\}_{\xi\neq 0}$ of 
$[x\mapsto G_{(\mk{m},\mk{n})}^0(U,p\,;x+\ii y,s)]$, we know that functions $[s\mapsto a_{\xi}(y,s)]$ must satisfy a certain functional equation.
Then, using some linear algebra, one transfers this functional to the constant term $a_0(y,s)$ which involves the partial zeta functions which we
are interested in.

\subsection{Relationships between $G_{(\mk{m},\mk{n})}^w(U,p\, ;z,s)$ with
some $GL_2$-Eisenstein series appearing in the literature}\label{rel}

We would like now to discuss briefly 
the relationships between the $GL_2$-Eisenstein series of this manuscript and the $GL_2$-Eisenstein series appearing in the literature.
For $w\geq 3$, one easily sees that
\begin{align}\label{bord}
F_w(U;z):=G_{(\mk{m},\mk{n})}^w(U,\bz ;z,s)|_{s=0},
\end{align}
is a holomorphic Eisenstein series of parallel weight $w$. The restriction $w\geq 3$ (so that $\frac{w}{2}>1$) is imposed so that the series,
on the right-hand side of \eqref{bord}, converges absolutely (however see Remark \ref{tipi} for the case $w=2$). 
In Section \ref{sec_h}, the Fourier series expansion of $[z\mapsto F_w(U;z)]$ is given. In the special case
where the parameter matrix $U$ has the form $\M{0}{v_1}{0}{v_2}$, the Eisenstein series $F_w(U;z)$
are essentially the holomorphic Eisenstein constructed by Deligne and Ribet in \cite{D-R}, see p. 269. 
From the perspective of \cite{D-R}, the usefulness of these Eisenstein series  comes from the fact 
that, for a finite sum of these (of possibly different weights), 
one may apply, under some assumptions, the so-called \lq\lq $q$-expansion principle\rq\rq\s in order to obtain non-trivial $p$-adic congruences
for the constant term of its Fourier series at the various cusps. From these $p$-adic congruences, it is then possible to
construct $p$-adic $L$-functions.

When $p\in \ZZ^{g}$ and $u_1=u_2=0$, $w\in\ZZ$, and $\mk{m}=\mk{n}=\mk{a}\mk{x}$ the (uncompleted) real analytic Eisenstein series 
\begin{align*}
\frac{G_{(\mk{a}\mk{x},\mk{a}\mk{x})}^w(U,2p;z,\frac{s}{2})}{(2\pi\ii)^{wg+\Tr(p)}\cdot(z-\ov{z})^{p\cdot\bu}
\cdot|y|^{\frac{s}{2}\cdot\bu}}= E_{w,p}(z,s;v_1,v_2;\mk{x},\mk{a}),
\end{align*}
were considered by Shimura (see (2.2) of \cite{Shim4}). Here, $\mk{x}\subseteq K$ is a fractional ideal and
$\mk{a}\subseteq K$ is an integral ideal. From the perspective of \cite{Shim4}, the interest for
such Eisenstein series comes from the fact that their special values, when $s=0$ and $z\in\mk{h}^g$ is a CM point, 
are algebraic numbers multiplied by some transcendental period which depends only on 
$k$ and $p$ and not on the CM point $z$ and the parameter matrix $U$. As it is explained in \cite{Shim4}, 
such a result may be used to show that special values of 
$L$-functions associated to a Gr\"ossencharacter $\psi$ of type $A_0$ (in the sense of Weil)
are equal to an algebraic number times a period which {\it depends only} on the infinity type of $\psi$.

\section{Some formulas}\label{formu}
In this section we introduce various functions which are required in order to write down explicit formulas for 
the non-zero Fourier coefficients of our Eisenstein series. It will also be important for us to study the
meromorphic continuation of these functions on certain domains and also give some growth estimates when the parameters vary
in certain regions.

Let $\alpha$ and $\beta$ be complex numbers and view $z=x+\ii y$ as a variable in $\mk{h}^{\pm}$. We define
the formal sum 
\begin{align}\label{gor0}
f(\alpha,\beta;z):=\sum_{n\in\ZZ}(z+n)^{-\alpha}(\ov{z}+n)^{-\beta},
\end{align}
When $\Ree(\alpha+\beta)>1$, the sum above is absolutely convergent. Note that the continuous function 
$[x\mapsto f(\alpha,\beta;x+\ii y)]$ is invariant under translation by $1$. 
Therefore, one may consider its Fourier series expansion, namely
\begin{align}\label{gor}
\wt{f}(\alpha,\beta;z):=\sum_{n\in\ZZ}\rho_n(y,\alpha,\beta)e^{{2\pi\ii } nx},
\end{align}
where the Fourier coefficients are given by the formulas
\begin{align}\label{timm}\notag
\rho_n(\alpha,\beta;y)&=\int_{0}^{1}f(\alpha,\beta;z)e^{-{2\pi\ii }nx}dx\\
&=\int_{-\infty}^{\infty} z^{-\alpha}(\ov{z})^{-\beta}e^{-{2\pi\ii } nx}dx.
\end{align}
The second equality is justified from the absolute convergence of \eqref{gor0} which allows one to interchange
the order of summation with the integral. We note that the 
integral in \eqref{timm} converges absolutely, since $\Ree(\alpha+\beta)>1$. 
Because $x\mapsto f(\alpha,\beta,x+\ii y)$ is a real analytic function on the one-dimensional torus $\RR/\ZZ$,
we know, from (1) of Theorem \ref{keyy} in Section \ref{Four}, that the Fourier series \eqref{gor} converges absolutely, and that for
all $z\in\mk{h}^{\pm}$, $f(\alpha,\beta;z)=\wt{f}(\alpha,\beta;z)$. In particular,
the Fourier series computes the value of $f(\alpha,\beta;z)$. Moreover, from (3) of Theorem \ref{keyy}, 
we also know that there exists a positive constant $D\in\RR_{>0}$ (depending on $y$),
such that $|\rho_n(\alpha,\beta,y)|\leq e^{-D|n|}$ for all $n\in\ZZ$.

Following Shimura (see p. 84 of \cite{Shim}), for $z=x+\ii y\in \mk{h}^{\pm}$, $t\in\RR$ and $\alpha,\beta\in\CC$ with $\Ree(\alpha+\beta)>1$, we define
\begin{align}\label{nezz}
\tau(\alpha,\beta;t,y)&:=\int_{-\infty}^{\infty} z^{-\alpha}(\ov{z})^{-\beta}e^{-{2\pi\ii } t x}dx.
\end{align}
Note that this integral is absolutely convergent. 
Later on, a more precise inequality of the form $|\tau(\alpha,\beta;t,y)|\leq C_s|t|^{\Ree(\alpha)+\Ree(\beta)-1}e^{-2\pi|ty|}$, as $|ty|\rightarrow\infty$
will be proved (see Proposition \ref{tiopp}). The integral \eqref{nezz} satisfies some obvious symmetries: 
For all $t,y,a\in\RR^{\times}$, we have the following rules:
\begin{enumerate}[(i)]
 \item $\tau(\alpha,\beta; t,-y)=\tau(\beta,\alpha;t,y)$
 \item $\ov{\tau(\alpha,\beta;t,y)}=\tau(\ov{\beta},\ov{\alpha};-t,y)$ (here the \lq\lq bar\rq\rq\s denotes the complex conjugation)
 \item If $y\in\RR^{\times}$ and $a>0$ then $\tau(\alpha,\beta;at,y)=|a|^{\alpha+\beta-1}\tau(\alpha,\beta;t,ay)$.
 \item If $y\in\RR_{>0}$ and $a<0$ then $\tau(\alpha,\beta;at,y)=e^{\pi\ii(\beta-\alpha)}|a|^{\alpha+\beta-1}\tau(\alpha,\beta;t,ay)$.
 \item If $y\in\RR_{<0}$ and $a<0$ then $\tau(\alpha,\beta;at,y)=e^{\pi\ii(\alpha-\beta)}|a|^{\alpha+\beta-1}\tau(\alpha,\beta;t,ay)$.
 \end{enumerate}
All the proofs of these formulas are straightforward. Let us just point out that the
rules (iv) and (v) follow from the formulas $(az)^{\alpha}=|a|^{\alpha}z^{\alpha}e^{-\alpha\pi\ii}$ when $z\in\mk{h}^+$, $a\in\RR_{<0}$; 
and $(az)^{\alpha}=|a|^{\alpha}z^{\alpha}e^{\alpha\pi\ii}$ when $z\in\mk{h}^{-}$, $a\in\RR_{<0}$. These
last two formulas follow from the product rule \eqref{chum}.

It will be important for us to give explicit formulas for $\tau(\alpha,\beta;t,y)$ in the following two cases:
\begin{enumerate}
 \item[(a)] $\Ree(\alpha),\Ree(\beta)>0$ and $\Ree(\alpha+\beta)>1$.
 \item[(b)] $\alpha\in\ZZ_{\geq 2}$ and $\beta=0$ (from the rule (ii) one automatically obtains the formula for the case where $\beta\in\ZZ_{\geq 2}$ and $\alpha=0$), 
\end{enumerate}
\begin{Rem}
Note that these two cases do not overlap, but (b) may be viewed as a limiting case of (a), as $\beta\rightarrow 0^+$ and $\alpha\in\ZZ_{\geq 2}$. 
From the perspective of Eisenstein series, this later observation means that holomorphic Eisenstein series can be viewed as a limit 
(in the $s$ variable) of real analytic Eisenstein series. Explicit formulas in the case (b) (see below) will be used in Section \ref{sec_h} 
to compute the Fourier coefficients of 
holomorphic Eisenstein series.
\end{Rem}
We first focus on obtaining explicit formulas for the case (a). Explicit formulas for the case (b) will
be given in Proposition \ref{posi}.

On page 366 of \cite{Sie5}, Siegel introduced the {\it confluent hypergeometric function} $\sigma(y;\alpha,\beta)$, given by
\begin{align}\label{tornade}
\sigma(y;\alpha,\beta):=\int_{0}^{\infty}(u+1)^{\alpha-1}u^{\beta-1}e^{-yu}du,
\end{align} 
where $y\in\RR_{>0}$, $\Ree(\beta)>0$ and $\alpha$ is any complex number. 
This allowed Siegel to obtain an explicit formula for $\tau(\alpha,\beta;n,y)$ ($n\in\ZZ$) in terms of $\sigma(y;\alpha,\beta)$.
On p. 84 of \cite{Shim}, Siegel states (without a proof) the following explicit formulas for 
the $\tau(\_,\_;\_,\_)$ in terms of the $\sigma(\_;\_,\_)$ function.

\begin{Lemma}\label{Shim_lem}
Let $\alpha,\beta\in\CC$ be such that $\Ree(\alpha)>0,\Ree(\beta)>0$ and 
$\Ree(\alpha+\beta)>1$. Let $t\in\RR$ and $y\in\RR^{\times}$. Then we have:
\begin{align*}
&\tau(\alpha,\beta;t,y)\\[3mm]
&=\left\{
\begin{array}{lll}
(2\pi)^{\alpha+\beta}\ii^{\beta-\alpha}
|t|^{\alpha+\beta-1}\Gamma(\alpha)^{-1}\Gamma(\beta)^{-1}e^{-2\pi |ty|}\sigma(4\pi |ty|;\alpha,\beta) & if & y>0,t>0,  \\ [2mm]
(2\pi)^{\alpha+\beta}\ii^{\alpha-\beta}
|t|^{\alpha+\beta-1}\Gamma(\alpha)^{-1}\Gamma(\beta)^{-1}e^{-2\pi |ty|}\sigma(4\pi |ty|;\beta,\alpha) & if & y<0,t>0,  \\ [2mm]
 (2\pi)^{\alpha+\beta}\ii^{\beta-\alpha}
|t|^{\alpha+\beta-1}\Gamma(\alpha)^{-1}\Gamma(\beta)^{-1}e^{-2\pi |ty|}\sigma(4\pi |ty|;\beta,\alpha) & if & y>0,t<0, \\  [2mm]
(2\pi)^{\alpha+\beta}\ii^{\alpha-\beta}
|t|^{\alpha+\beta-1}\Gamma(\alpha)^{-1}\Gamma(\beta)^{-1}e^{-2\pi |ty|}\sigma(4\pi |ty|;\alpha,\beta) & if & y<0,t<0, \\  [2mm]
(2\pi)\cdot \ii^{\beta-\alpha}\frac{\Gamma(\alpha+\beta-1)}{\Gamma(\alpha)\Gamma(\beta)}(2y)^{1-\alpha-\beta} & if & \mbox{$t=0$, $y>0$},\\ [2mm]
(2\pi)\cdot \ii^{\alpha-\beta}\frac{\Gamma(\alpha+\beta-1)}{\Gamma(\alpha)\Gamma(\beta)}(2|y|)^{1-\alpha-\beta} & if & \mbox{$t=0$, $y<0$}, 
\end{array}\right. 
\end{align*}
where $\Gamma(x)$ stands for the gamma function evaluated at $x$.
\end{Lemma}
{\bf Proof} For a proof of Lemma \ref{Shim_lem}, see the proof of Lemma 1 of \cite{Shim}. \fin
\begin{Rem}
Note that the rules (i), (ii), (iii), (iv) and (v) for the function $\tau(\_,\_;\_,\_)$ can all be verified directly from the above list
of formulas. 
\end{Rem}
\begin{Rem}
Later on, it will be explained (see Corollary \ref{triop0}) that the function $(y,\alpha,\beta)\mapsto\sigma(y;\alpha,\beta)$
admits a single-valued meromorphic continuation when the triple $(y,\alpha,\beta)$ varies over the domain $\CC\bs\RR_{\leq 0}\times \CC\times\CC$.
\end{Rem}

It also will be convenient for us to introduce a modified version of the function $\tau(\alpha,\beta;t,y)$ that 
satisfies additional symmetry relations.
\begin{Def}\label{nop}
Let $\alpha,\beta\in\CC\bs\ZZ_{\leq 0}$, $t\in\RR^{\times}$ and $z=x+\ii y\in\mk{h}^{\pm}$. We set
\begin{align}\label{cyr}
p(\alpha,\beta;t,y):=
\left\{\begin{array}{ccc}
\ii^{\epsilon(y)(\alpha-\beta)}\Gamma(\alpha)(4\pi|ty|)^{\beta}
& \mbox{if} & ty>0, \\
\ii^{\epsilon(y)(\alpha-\beta)}\Gamma(\beta)(4\pi|ty|)^{\alpha}
& \mbox{if} & ty<0,
\end{array}\right.
\end{align}
where $\epsilon(y)=\sign(y)$. We defined the \lq\lq completed $\tau$-function\rq\rq\s as
\begin{align}\label{key_def}
\wh{\tau}(\alpha,\beta;t,y):=
|t|^{1-\alpha-\beta}(2\pi)^{-\alpha-\beta}\cdot
p(\alpha,\beta;t,y)\cdot\tau(\alpha,\beta;t,y).
\end{align}
\end{Def}
\begin{Rem}
We note that the functions $p(\alpha,\beta;t,y)$ and $\wh{\tau}(\alpha,\beta;t,y)$ 
depend only on the triple $(\alpha,\beta,ty)$ rather than the quadruple $(\alpha,\beta,t,y)$ itself. Moreover, the definition of $p(\alpha,\beta;t,y)$
is such that the positions of $\alpha$ and $\beta$, on the right-hand side of \eqref{cyr},
are interchanged, according to the sign of $ty$.
\end{Rem}
The next formulas are easy to prove.
\begin{Prop}\label{magma}
Let $\alpha,\beta\in \CC\bs\ZZ_{\leq 0}$ and let $t,y,a\in\RR^{\times}$. 
Then the following identities hold true:
\begin{enumerate}
 \item If $a>0$, $p(\alpha,\beta;at,y)=p(\alpha,\beta;t,ay)$.
 \item Assume that $\alpha-\beta\in\ZZ$. Then, if $a<0$, we have $p(\alpha,\beta;at,y)=(-1)^{\alpha-\beta}\cdot p(\alpha,\beta;t,ay)$.
 \item Assume that $\alpha-\beta\in\ZZ$, then 
\begin{enumerate}[(a)]
\item $p(\alpha,\beta;-t,y)=(-1)^{\alpha-\beta}\cdot p(\beta,\alpha;t,y)$,
\item $p(\alpha,\beta;t,-y)=(-1)^{\alpha-\beta}\cdot p(\beta,\alpha;t,y)$.
\end{enumerate}
 \item If $t\in\RR^{\times}$ and $y\in\RR_{>0}$ then $p(\alpha,\beta;t,y)=|y|^{\beta}p(\alpha,\beta;t,1)$.
 \item If $t\in\RR^{\times}$ and $y\in\RR_{<0}$ then $p(\alpha,\beta;t,y)=|y|^{\alpha}p(\beta,\alpha;t,1)$.
 \item $p(\alpha,\beta;t,y)p(\beta,\alpha;t,y)=\Gamma(\alpha)\Gamma(\beta)|4\pi ty|^{\alpha+\beta}$.
 \item $\wh{\tau}(\alpha,\beta;at,y)=\wh{\tau}(\alpha,\beta;t,ay)$.
 \item $\wh{\tau}(\alpha,\beta;-t,y)=\wh{\tau}(\beta,\alpha;t,y)$.
 \item $\wh{\tau}(\alpha,\beta;t,-y)=\wh{\tau}(\beta,\alpha;t,y)$.
\end{enumerate}
\end{Prop}

{\bf Proof} The proof is left to the reader. \fin

\subsection{Tricomi's confluent hypergeometric function}\label{frip}
A priori, the function $y\mapsto\tau(\alpha,\beta;1,y)$ (see \eqref{nezz}) is only defined for $y\in\RR^{\times}$,
$\alpha,\beta\in\CC$ with $\Ree(\alpha)>0$, $\Ree(\beta)>0$ and $\Ree(\alpha+\beta)>1$. In order to
to give a meaning to $\tau(\alpha,\beta;1,y)$, when $(\alpha,\beta,y)$ lies outside this region, 
we will express Siegel's confluent hypergeometric function ($\sigma(\_;\_,\_)$) in terms of {\it Tricomi's confluent 
hypergeometric function} ($U(\_,\_;\_)$). In this way, we can take advantage of the extensive study
of $U(\_,\_;\_)$ which was done in \cite{Slat}. In particular, this allows us to transfer the various analytic properties of $U(\_,\_;\_)$ to
$\sigma(\_;\_,\_)$ and to $\tau(\alpha,\beta;1,y)$.

We first define $U(\_,\_;\_)$ in terms of $\sigma(\_;\_,\_)$. Then, we present some of the properties  
of $U(\_,\_;\_)$ that will be needed in the sequel. Let $\alpha,\beta\in\CC$ with $\Ree(\beta)>0$ and set $a:=\beta$, $b:=\alpha+\beta$. 
For any $t\in\RR_{>0}$, Tricomi's confluent hypergeometric function (cf. 3.1.19 of \cite{Slat}) is defined as
\begin{align}\label{trico}
U(a,b;t):=\frac{1}{\Gamma(a)}\sigma(t;b-a,a),
\end{align}
or equivalently,
\begin{align}\label{trico2}
U(\beta,\alpha+\beta;t):=\frac{1}{\Gamma(\beta)}\sigma(t;\alpha,\beta).
\end{align}
One may verify that $U(a,b;t)$ satisfies the Kummer's hypergeometric differential equation
\begin{align}\label{time}
t\frac{d^2w}{dt^2} + (b-t)\frac{dw}{dt} - a\cdot w = 0,
\end{align}
which has a regular singular point at $t=0$ and an irregular singular point at $t=\infty$. Note that
the differential equation \eqref{time} is a degenerate form of the usual Gauss' hypergeometric 
equation, given by
\begin{align}\label{tele}
t(1-t)\frac {d^2 w}{d t^2} + \left[b-(a+c+1) t \right] \frac{d w}{d t} - ac\cdot w = 0.
\end{align}
If one replaces $t$ by $\frac{t}{c}$ in \eqref{tele} (so $\frac{dw}{dt}$ goes to $c\frac{dw}{dt}$ and $\frac{d^2w}{dt^2}$ goes to $c^2\frac{d^2w}{d^2t}$) 
and if one lets $c\rightarrow \infty$,
we obtain the ODE \eqref{time}. 

The next proposition describes some properties of the function $U(a,b;t)$.
\begin{Prop}\label{wee}
\begin{enumerate}
 \item The function $(a,b,z)\mapsto U(a,b;z)$ admits a multi-valued analytic continuation to all of 
 $\CC\times\times\CC\times\CC\bs\{0\}$.
 \item Let $z\in\CC\bs\{0\}$ be fixed. Then,
 function $(a,b)\mapsto U(a,b;z)$ admits a single-valued analytic continuation to all of 
 $\CC\times\CC$.
\item Let $L\subseteq\CC$ be a half-line joining $0$ to $\infty$. Then, the function
$(a,b,z)\mapsto U(a,b;z)$ admits an holomorphic extension to $\CC\times\CC\times(\CC\bs L)$.
\item  Let $L\subseteq\CC$ be a half-line joining $0$ to $\infty$. Let 
$a,b\in\CC$, with $\Ree(a),\Ree(b)>0$ be fixed. Then, 
the one variable holomorphic function $z\mapsto U(b,a+b;z)$ is non-constant.  
\item For $a,b\in\CC$ fixed, we have $U(a,b;z)\sim z^{-a}(1+O(|z|^{-1}))$ as $\Ree(z)\rightarrow \infty$. 
In particular, if $\Ree(a)>0$, we have the asymptotic $U(a,b;z)\sim z^{-a}$ as $\Ree(z) \rightarrow\infty$.
Moreover, if $\Ree(a)>0$, $U(a,b;z)$ is the unique solution of the linear ODE \eqref{tele}, which is, up to a 
multiplicative scalar, bounded as $\Ree(z)\rightarrow\infty$.
\item For all $\beta\in\RR_{>0}$, $\alpha\in\RR$ and $y\in\RR_{>0}$, we have $U(\beta,\alpha+\beta;y)>0$.
\item For all $b\in\CC$ and $z\in\CC\bs\{0\}$, we have $U(0,b;z)=1$.
\item $U(a,a+1;z)=\frac{1}{z^a}$.
\item Let $\alpha,\beta\in\CC$ and $z\in\CC\bs\{0\}$. Then, the quantity $z^\beta\cdot U(\beta,\alpha+\beta;z)$
is invariant under the substitution $\alpha\mapsto 1-\beta$ and $\beta\mapsto 1-\alpha$.
\item $\frac{d}{dz}U(a,b;z)=-a U(a+1,b+1,z)$.
\end{enumerate}
\end{Prop}

\begin{Rem}
Note that in the last statement of (5), we could have replaced \lq\lq which is bounded as $\Ree(z)\rightarrow\infty$\rq\rq\s by
\lq\lq which is bounded by $|\Norm(y)|^N$ for some positive integer $N\in\ZZ_{\geq 1}$, as $\Ree(z)\rightarrow\infty$\rq\rq.
\end{Rem}

{\bf Proof}  Let  $\frac{{}_{1}F_1(a,b;z)}{\Gamma(b)}$ be the {\it normalized Kummer series}.
It follows from the equation (1.3.5) on page 5 of \cite{Slat} that $(a,b,z)\mapsto \frac{{}_{1}F_1(a,b;z)}{\Gamma(b)}$
is a holomorphic function on {\it all} of $\CC^3$. Equation (1.3.1) in \cite{Slat} reads as
\begin{align}\label{1.3.1}
U(a,b;z)=\frac{\Gamma(1-b)}{\Gamma(1+a-b)}{}_{1}F_{1}(a,b;z)+\frac{\Gamma(b-1)}{\Gamma(a)}x^{1-b}{}_{1}F_1(1+a-b,2-b;z).
\end{align}
It is valid, a priori, for all $(a,b,z)\in\CC^3$ with $b\notin\ZZ_{\leq 0}$. In order to handle
the case where $b\in\ZZ_{\leq 0}$, one may look at Equation (1.3.5) of \cite{Slat} which reads as
\begin{align}\label{1.3.5}
U(a,b;z)=\frac{\pi}{\sin(\pi b)}\left(\frac{{}_{1}F_{1}(a,b;z)}{\Gamma(1+a-b)\Gamma(b)}-
\frac{x^{1-b}{}_{1}F_{1}(1+a-b,2-b;z)}{\Gamma(a)\Gamma(2-b)}\right).
\end{align}
Note that \eqref{1.3.5} makes sense when $n\in\ZZ_{\leq 0}$ and $b\rightarrow n$.
From these two equations, we see that $(a,b,z)\mapsto U(a,b;z)$ becomes a multi-valued function on $\CC\times\times\CC\times\CC\bs\{0\}$. The
fact the $U(a,b;z)$ is multi-valued comes from the {\it need of fixing a branch} of $z\mapsto z^{1-b}$. This proves (1).
In particular, if $z\in\CC\bs\{0\}$ is fixed, the function $(a,b)\mapsto U(a,b;z)$ is single-valued on all of $\CC^2$. This proves (2).

The proof of (3) follows directly from Equations \eqref{1.3.5} and \eqref{1.3.1}.

Let us prove (4). Let $a,b\in\CC$ with $\Ree(a),\Ree(b)>0$. We do a proof by contradiction. Assume that 
for all $z\in\CC\bs L$, $z\mapsto U(b,a+b;z)=0$. Then 
\begin{align}\label{maki1}
\lim_{z\rightarrow 0} z^{1-a-b}\cdot\left(z^{a+b-1}\cdot U(b,a+b;z)\right)=0.  
\end{align}
But from (2) of Proposition \ref{geno} (see the proposition below), the limit on the left hand side of \eqref{maki1} is equal to
\begin{align}\label{maki2}
\lim_{z\rightarrow 0} z^{1-a-b}\cdot\frac{\Gamma(a+b-1)}{\Gamma(b)}.
\end{align}
From the infinite product of $\frac{1}{\Gamma(z)}$, we know that $\Gamma(z)\neq 0$ for all $z\in\CC$. Therefore, since
$\Ree(a),\Ree(b)>0$, the quotient $\frac{\Gamma(a+b-1)}{\Gamma(b)}$ is either equal to $\infty$ (if $a+b-1=0$) or different from $0$.
But this is absurd, since the limit in \eqref{maki1} is equal to $0$ and the limit in \eqref{maki2} is equal to $\infty$. 
This concludes the proof of (4).

For a proof of (5), see pages 58-60 of \cite{Slat}. 

The proof of (6) follows directly from the integral \eqref{tornade}
and the identity \eqref{trico}. 

Let us prove (7). When $b\notin\ZZ_{\leq 0}$, 
(7) follows from equation \eqref{1.3.1} and the observation that for all $c\in\CC\bs\ZZ_{\leq 0}$, ${}_{1}F_1(0,c;z)=1$. 
If $b\in\ZZ_{\leq 0}$, then one uses instead equation \eqref{1.3.5}. This proves (6).

The proof of (8) follows from \eqref{tornade} and \eqref{trico}.

The proof of (9) follows from equation (1.4.9) on page 6 of \cite{Slat}. 

For a proof of (10), one may simply differentiate the expression in \eqref{tornade} 
under the integral sign
to find the relation $\frac{d}{dz}\sigma(z;\beta-\alpha,\alpha)=\sigma(z;\beta-\alpha,\alpha+1)$; afterwards,
one uses the identity \eqref{trico}. \fin
\begin{Rem}
It follows from the proof of (1) of Proposition \ref{wee} that for any half-line $L$, in the complex plane, 
joining $0$ to $\infty$, the function $(a,b,z)\mapsto U(a,b;z)$ admits a single-valued holomorphic 
continuation to all of $\CC\times\CC\times(\CC\bs L)$.
\end{Rem}

The next proposition gives the limit values of the function $z\mapsto U(\beta,\alpha+\beta;z)$ when 
$z\rightarrow 0$, and when the complex parameters $(\alpha,\beta)$ lie in certain regions of $\CC^2$.
These two limit formulas will be used later on.
\begin{Prop}\label{geno}
\begin{enumerate}
 \item Assume that $\Ree(\alpha+\beta)<1$, then 
 \begin{align*}
 \lim_{z\rightarrow 0} U(\beta,\alpha+\beta;z)= \frac{\Gamma(1-(\alpha+\beta))}{\Gamma(1-\alpha)}.
 \end{align*}
 \item Assume that $\Ree(\alpha+\beta)>1$, then 
\begin{align*}
\lim_{z\rightarrow 0} z^{\alpha+\beta-1}\cdot U(\beta,\alpha+\beta;z)=\frac{\Gamma(\alpha+\beta-1)}{\Gamma(\beta)}.
\end{align*}
\end{enumerate}
\end{Prop}

{\bf Proof} The proof of (1) and (2) follow from the equations \eqref{1.3.1} and \eqref{1.3.5} and the observation that, for
all $(\alpha,\beta)\in\CC^2$,
\begin{align*}
\lim_{z\rightarrow 0}\frac{{}_{1}F_1(\alpha,\beta;z)}{\Gamma(\beta)}=1.
\end{align*} \fin

\begin{Rem}
Note that the two sets $\{(\alpha,\beta)\in\CC^g:\Ree(\alpha+\beta)<1\}$ and $\{(\alpha,\beta)\in\CC^g:\Ree(\alpha+\beta)>1\}$.
don't overlap.
\end{Rem}

\subsubsection{Relationship between $U(a,b;z)$ and $K_{s}(z)$}

In many papers related to the computation of the Fourier series expansion
$GL_2$-real analytic Eisenstein series, the $K$-Bessel function (or sometimes also called 
the modified Bessel function of the second kind) appears. In this short section, we explain the 
relationship between the $K$-Bessel function and Tricomi's confluent hypergeometric function $U(\_;\_,\_)$.

For $z\in\CC$ with $\Ree(z)>0$, one may define the $K$-Bessel function using its Schl\"afi's integral representation as
\begin{align}\label{rote}\notag
K_s(z)&:=\int_{0}^{\infty} e^{-z\cosh(t)}\cosh(st) dt\\
        &=\frac{1}{2}\int_{0}^{\infty} e^{-\frac{1}{2}z(u+1/u)}(u^s+u^{-s})\frac{du}{u},
\end{align}
where $u=e^{t}$. From equation (1.8.7) of \cite{Slat}, we have
\begin{align}\label{tiol}
U(s,2s;2z)=\frac{(2z)^{\frac{1}{2}-s}}{\sqrt{\pi}} e^{z}K_{s-\frac{1}{2}}(z).
\end{align}

\begin{Rem}
If we define
\begin{align*}
\wt{K}_{s}(z):=\frac{1}{2}\int_{0}^{\infty} e^{-\frac{1}{2}z(u+1/u)}u^s \frac{du}{u},
\end{align*}
then $\wt{K}_{s}(z)$ is invariant 
under the change of variable $u\mapsto \frac{1}{u}$. Therefore, $\wt{K}_{s}(z)=\wt{K}_{-s}(z)$.  
Looking at \eqref{rote}, we  obtain the 
relation $K_s(z)=2\wt{K}_s(z)$. For this reason, some authors also define the $K$-Bessel function by $\wt{K}_{s}(z)$
rather than $K_s(z)$.
\end{Rem}

\subsection{Meromorphic continuation of some functions}\label{mero_cont}
In this section we prove the meromorphic continuation of various functions which have been introduced before.
Ultimately, the existence of such meromorphic continuations
boil down to the analytic properties of the Tricomi's confluent hypergeometric functions $[(a,b,z)\mapsto U(a,b;z)]$ and of
the gamma function $\Gamma(z)$. 

Let $z\in\CC\bs\{0\}$ be fixed. From Proposition \ref{wee}, we know that 
the function $(a,b)\mapsto U(a,b;z)$ admits a holomorphic continuation to all of $\CC^2$.
Looking at the identity \eqref{trico} and using the well-known properties of the gamma function,
we obtain immediately the following corollary:
\begin{Cor}\label{triop0}
Let $t\in\RR\bs\{0\}$ be fixed. Then, function $(\alpha,\beta)\mapsto\sigma(t;\alpha,\beta)$ admits a single-valued meromorphic continuation 
to all of $\CC^2$. Moreover, it is analytic on the domain $\CC\times \CC\bs\ZZ_{\leq 0}$. For a 
fixed pair $(t,\alpha)\in (\RR\bs\{0\})\times \CC$, the one variable meromorphic function 
$\beta\mapsto \sigma(t,\alpha,\beta)$ is holomorphic on all of $\CC\bs\ZZ_{\leq 0}$, with possible 
poles of order at most one at the elements in $\ZZ_{\leq 0}$.
\end{Cor}

Having obtained the meromorphic continuation of $\sigma(t;\alpha,\beta)$, we now focus 
on the $\tau$-function, i.e., $\tau(\alpha,\beta;t,y)$ (see \eqref{nezz}).
Let $t,y\in\RR\bs\{0\}$ be fixed. To fix the idea, let us assume that $t,y>0$. 
From Lemma \ref{Shim_lem} and equation \eqref{trico}, we have
\begin{align}\label{cou}
\tau(\alpha,\beta;t,y)=(2\pi)^{\alpha+\beta}(\ii)^{\beta-\alpha}|t|^{\alpha+\beta-1}
e^{-2\pi |ty|}\cdot \Gamma(\alpha)^{-1}\cdot U(\beta,\alpha+\beta;|4\pi ty|).
\end{align}
Recall that $\alpha\mapsto\Gamma(\alpha)$ is a holomorphic function on all of $\CC$.
We decide to extend analytically the function $(\alpha,\beta)\mapsto\tau(\alpha,\beta;t,y)$ to all 
$(\alpha,\beta)\in\CC\times\CC$, using the right-hand side of \eqref{cou}. In a similar way, we may also define an analytic
continuation of $(\alpha,\beta)\mapsto\tau(\alpha,\beta;t,y)$, to all of $\CC^2$, when 
the pair $(t,y)$ is such that $(\sign(t),\sign(y)))\in\{(-1,+1),(+1,-1),(-1,-1)\}$. 

We thus have
\begin{Prop}\label{tiopp}
Let $t,y\in\RR\bs\{0\}$ be fixed. Then the function $[(\alpha,\beta)\mapsto \tau(\alpha,\beta;t,y)]$ admits an 
analytic continuation to all of $\CC^2$. If $ty>0$, and $\alpha\in\ZZ_{\leq 0}$ then  
$[(t,y)\mapsto\tau(\alpha,\beta;t,y)]\equiv 0$. If $ty<0$, and $\beta\in\ZZ_{\leq 0}$ then  
$[(t,y)\mapsto\tau(\alpha,\beta;t,y)]\equiv 0$. Assume that $\Ree(\alpha)>0$ and $\Ree(\beta)>0$. Then, 
there exists a positive constant $C_s>0$, such that for all $t,y\in\RR\bs\{0\}$ with
$(\sign(t),\sign(y))$ fixed, we have that
\begin{align}\label{cart}
|\tau(\alpha,\beta;t,y)|\leq C_s\cdot |t|^{\Ree(\alpha)+\Ree(\beta)-1}\cdot e^{-2\pi|ty|}. 
\end{align}
\end{Prop}

{\bf Proof} The first three claims follow from the definition of $\tau(\alpha,\beta;t,y)$ and (2) of Proposition 
\ref{wee}. The proof of \eqref{cart} follows from the asymptotic formula for $U(\beta,\alpha+\beta;z)$ 
(and $U(\alpha,\alpha+\beta;z)$) given in (5) of Proposition \ref{wee} and the explicit formula 
of $\tau(\alpha,\beta;t,y)$ in terms of $U(\_,\_;\_)$. \fin

We may now give the explicit formula for $\tau(\alpha,\beta;t,y)$ when $\alpha\in\ZZ_{\geq 2}$ and $\beta=0$ (case (b) in Section \ref{formu}).
\begin{Prop}\label{posi}
Let $t\in\RR^{\times}$ and $z=x+\ii y\in\mk{h}^{\pm}$. Let $\alpha\in\ZZ_{\geq 2}$ and set $\epsilon=\sign(y)$. Then
\begin{align*}
\tau(\alpha,0;t,y)=
\left\{
\begin{array}{ccc}
(-\epsilon\cdot 2\pi\ii)^{\alpha}\frac{|t|^{\alpha-1}}{(\alpha-1)!}e^{-2\pi |t y|} &   \mbox{if}  &  \sign( ty)>0\\
0 & \mbox{if} & \sign(ty)<0
\end{array}
\right.
\end{align*}
and 
\begin{align*}
\tau(0,\alpha;t,y)=
\left\{
\begin{array}{ccc}
0 & \mbox{if} & \sign(ty)>0,\\
(\epsilon\cdot 2\pi\ii)^{\alpha}\frac{|t|^{\alpha-1}}{(\alpha-1)!}e^{-2\pi |t y|} &   \mbox{if}  &  \sign(ty)<0.
\end{array}
\right.
\end{align*}
\end{Prop}

{\bf Proof} Let $t,y\in\RR^{\times}$ be fixed. In virtue of the principle of identity for holomorphic functions
in many variables, it follows from Proposition \ref{tiopp} that the first 4 formulas given in Lemma \ref{Shim_lem},
hold true for all $\alpha,\beta\in\CC$. Finally, a direct computation involving each of these four
cases, combined with (4) of Proposition \ref{wee}, give the formulas above. \fin

In a similar way, for a fixed pair $t,y\in\RR\bs\{0\}$, 
we may extend the definition of $(\alpha,\beta)\mapsto\wh{\tau}(\alpha,\beta;t,y)$ (see Definition \ref{nop}) to all $\CC^2$.
To fix the idea, let us assume that $t,y>0$. Then, from the definition of $\wh{\tau}(\_,\_;\_,\_)$, we have
\begin{align*}
\wh{\tau}(\alpha,\beta;t,y)=e^{-2\pi |ty|}\cdot U(\beta,\alpha+\beta;|4\pi ty|)\cdot|4\pi ty|^{\beta}.
\end{align*}
In particular, $(\alpha,\beta)\mapsto\wh{\tau}(\alpha,\beta;t,y)$ admits an holomorphic
continuation to all of $\CC^2$.
Similarly, one may check that for each of the remaining three possibilities of $(\sign(t),\sign(y))$,
the function $(\alpha,\beta)\mapsto\wh{\tau}(\alpha,\beta;t,y)$ is again holomorphic on all of $\CC^2$.
We thus obtain:
\begin{Prop}\label{ele}
Let $t,y\in\RR\bs\{0\}$ be fixed. Then, function $(\alpha,\beta,z)\mapsto \wh{\tau}(\alpha,\beta;t,y)$ admits
an holomorphic continuation to all of $\CC^2$.
\end{Prop}

\subsection{Reflection formulas}\label{ref_for}

The next two lemmas are the key to show that the non-constant terms of the Fourier series expansion of 
the completed Eisenstein series $\wh{G}_{\mk{m},\mk{n}}^0(U,p\,;z,s)$ (see Definition \ref{tom}) are invariant under the substitution $s\mapsto 1-s$.
\begin{Lemma}\label{key_lemma}
Let $z\in\CC\bs\{0\}$ be fixed. Then, the expression
\begin{align}\label{topi}
z^{\beta}\cdot U(\beta,\alpha+\beta;z),
\end{align}
is invariant under the transformation
$\alpha\mapsto 1-\beta$ and $\beta\mapsto 1-\alpha$.
\end{Lemma}
{\bf Proof} This is a direct consequence of (9) of Proposition \ref{wee}. 
See also Lemma 2 of \cite{Shim}. \fin

\begin{Cor}\label{shimy}
Let $t,y\in\RR^{\times}$. Then the function
$\wh{\tau}(\alpha,\beta;t,y)=\wh{\tau}(\alpha,\beta;1,ty)$
is invariant under the transformation $\alpha\mapsto 1-\beta$ and $\beta\mapsto 1-\alpha$.
\end{Cor}

{\bf Proof} This follows from Lemma \ref{Shim_lem}, Lemma \ref{key_lemma} and
the definition of $\wh{\tau}(\alpha,\beta,t,y)$. \fin 

We will also need the following reflection formula which involves the $p(\_,\_;\_,\_)$ function which appears in Definition \ref{nop}.
\begin{Prop}\label{key_prop2}
Let $\alpha,\beta\in\CC\bs\ZZ$. For $t_1,t_2,y_1,y_2\in\RR^{\times}$, we define
\begin{align}\label{mall}
q(\alpha,\beta;t_1,y_1;t_2,y_2):=\frac{p(\beta,\alpha;t_1,y_1)}{p(\alpha,\beta;t_2,y_2)}.
\end{align}
Then the expression $q(\alpha,\beta;t_1,y_1;t_2,y_2)$ satisfies the following functional equation:
\begin{align}\label{tur}
q(\alpha,\beta;t_1,y_1;t_2,y_2)=q(1-\beta,1-\alpha;t_2,y_2;t_1,y_1)\cdot \left(\frac{\sin(\pi\beta)}{\sin(\pi\alpha)}\right)^{\epsilon}
\cdot\frac{|t_2 y_2|}{|t_1 y_1|},
\end{align}
where
\begin{align*}
\epsilon=
\left\{
\begin{array}{ccc}
1 & \mbox{if} & t_1y_1>0 \mbox{\s\s and\s\s } t_2y_2>0 \\
-1 & \mbox{if} & t_1y_1<0 \mbox{\s\s and\s\s} t_2y_2<0 \\
0  & \mbox{otherwise}
\end{array}
\right.
\end{align*}

Furthermore, assume that $\alpha-\beta\in\ZZ$. Then for any $d\in \RR^{\times}$, we have
\begin{align}\label{turfo}
q(\alpha,\beta;t_1,dy_1;t_2,y_2)=q(\alpha,\beta;t_1,y_1;t_2,dy_2)\cdot|d|^{\alpha+\beta}.
\end{align}
\end{Prop}
Note the swap between the pairs $(t_1,y_1)$ and $(t_2,y_2)$ on the right-hand side of \eqref{tur}.

{\bf Proof} The proof of \eqref{tur} follows directly from the definition of $p(\alpha,\beta;t,y)$
and the Euler's reflection formula for the gamma function. The
proof of \eqref{turfo} follows from the rules (4) and (5) of Proposition \ref{magma}. \fin 

\begin{Cor}\label{forel}
Let $\alpha,\beta\in\CC\bs\ZZ$ with $\alpha-\beta\in\ZZ$. 
Let $d,y\in\RR^{\times}$, then the expression $q(\alpha,\beta;1,dy;1,y)$ satisfies the following functional equation:
\begin{align}
q(\alpha,\beta;1,dy;1,y)=
\left\{
\begin{array}{cc}
(-1)^{\alpha-\beta}\cdot q(1-\beta,1-\alpha;1,dy;1,y)\cdot|d|^{\alpha+\beta-1} & \mbox{if $d>0$}\\[2mm]
q(1-\beta,1-\alpha;1,dy;1,y)\cdot|d|^{\alpha+\beta-1} & \mbox{if $d<0$}\\ 
\end{array}
\right.
\end{align}
\end{Cor}

{\bf Proof} This follows directly from Proposition \ref{key_prop2}. \fin

\section{Lattices in number fields and Fourier series expansion}\label{lat_nu}

\subsection{Lattices in number fields}\label{lat_nu2}
In this section, we decided to gather, for the convenience of the reader, some results about lattices in number fields
which are not easily found in the literature. The main goal of this section is to clarify the relations between two 
operations on lattices which arise in the setting of functional equations in number theory: 
the {\it multiplicative inverse operation} $\mathcal{L}\mapsto\mathcal{L}^{-1}$ and the {\it dual operation} $\mathcal{L}\mapsto\mathcal{L}^*$
(with respect to the trace pairing). These two operations arise naturally when one writes down the 
functional equations of partial zeta functions and $GL_2$-Eisenstein series. In the literature on $GL_2$-Eisenstein series, 
some authors preferred to use the former 
operation while others, the latter. In most papers on $GL_2$-Eisenstein series and partial zeta functions, the authors usually 
restrict themselves to working with lattices of $K$, which are also fractional $\ca{O}_K$-ideals. In this setting, it is easy 
to move from one point of view to the other,
since, if $\mk{a}\subseteq K$ is a fractional $\ca{O}_K$-ideal, then one has the well-known relation $\mk{a}^*=\mk{a}^{-1}\mk{d}_K^{-1}$. Here,
$\mk{d}_K^{-1}$ stands for the inverse of the different ideal of $K$. 
In the present work, we decided to work systematically and almost exclusively with the dual operation. 
Except for this section
and the proof of Proposition \ref{talon}, the multiplicative inverse operation will not be used. Although the role played
by the multiplicative inverse operation in our work is very small,  we think that it is worthwhile
to spend some time clarifying its relationship with the dual operation in order to have a theory which
simultaneously encompasses both points of view and which works for all lattices of $K$, not just those ones which are fractional $\ca{O}_K$-ideals. 
Since all the results in this section are valid for arbitrary number fields, for this section only, 
we let $K$ be an {\it arbitrary} number field of degree $g$ over $\QQ$.

Let $\mathcal{L}\subseteq K$ be a $\ZZ$-lattice, i.e., a free $\ZZ$-module of rank $g$.
For two lattices $\mathcal{L}_1,\mathcal{L}_2\subseteq K$, we define their product as
\begin{align*}
\mathcal{L}_1\mathcal{L}_2=\left\{\sum_{i=1}^n \ell_{1,i}\ell_{2,i}:
\ell_{1,i}\in\mathcal{L}_1,\ell_{2,i}\in\mathcal{L}_2,n\in\ZZ_{\geq 1}\right\}.
\end{align*}
One may verify that $\ca{L}_1\ca{L}_2$ is again a lattice. Moreover, the product operation on lattices is associative. 

\subsubsection{$\ca{O}$-properness}\label{multiply}
Let $\mathcal{L}\subseteq K$ be a lattice. We define
\begin{align*}
\ca{O}_{\mathcal{L}}:=\{\lambda\in K:\lambda \mathcal{L}\subseteq \mathcal{L}\},
\end{align*}
and call $\ca{O}_{\mathcal{L}}$
the {\it multiplier ring} of $\mathcal{L}$ (or the {\it endomorphism ring of $\mathcal{L}$}). If $\ca{O}$ is an order of $K$, such that
$\ca{O}=\ca{O}_\mathcal{L}$, then we say that $\mathcal{L}$ is {\bf $\ca{O}$-proper}. 
So, by definition, for any lattice $\mathcal{L}$, we always have that
$\mathcal{L}$ is $\ca{O}_\mathcal{L}$-proper. 

\subsubsection{The multiplicative inverse operation}\label{mult_inv}
We define the {\it multiplicative inverse of $\mathcal{L}$} to be
\begin{align*}
\mathcal{L}^{-1}:=\{\lambda\in K:\lambda\mathcal{L}\subseteq\ca{O}_{\mathcal{L}}\}. 
\end{align*}
Note that, by definition, $\mathcal{L}^{-1}$ is always an $\ca{O}_{\mathcal{L}}$-module and that
$\mathcal{L}\mathcal{L}^{-1}\subseteq\ca{O}_{\mathcal{L}}$.
Since $\mathcal{L}^{-1}$ is an $\ca{O}_{\mathcal{L}}$-module, it follows that that $\ca{O}_{\ca{L}}\subseteq\ca{O}_{\ca{L}^{-1}}$.
From the previous inclusion, it follows that $\ca{L}\subseteq(\ca{L}^{-1})^{-1}$. As Example \ref{filu} below
shows, the three inclusions above could be strict in general. The example below is inspired from an email exchange
with Keith Conrad, who kindly pointed out to me Exercise 18 on page 94 of \cite{Bor-Sha}.
\begin{Exa}\label{filu}
Let $\theta\in\ov{\QQ}$ be such that $\theta^3=2$ and consider the cubic field $K:=\QQ(\theta)$.
We have that $\ca{O}_K=\ZZ[\theta]$. Let $\mathcal{R}:=\ZZ+2\theta\ZZ+2\theta^2\ZZ$. One may verify that $\mathcal{R}$
is an order of index four in $\ca{O}_K$. Consider the lattice $\mathcal{M}:=4\ZZ+\theta\ZZ+\theta^2\ZZ\subseteq K$. Then 
direct computations (which we leave to the reader) show that
\begin{enumerate}[(i)]
 \item $\ca{O}_{\mathcal{M}}=\mathcal{R}$, 
 \item $\mathcal{M}^2=2\ZZ+2\theta\ZZ+\theta^2\ZZ$,
 \item $\ca{O}_{\mathcal{M}^2}=\ca{O}_K$,
 \item $\mathcal{M}^{-1}=2\ZZ+2\theta\ZZ+\theta^2\ZZ=2\ca{O}_K+\theta^2\ca{O}_K$,
 \item $\ca{O}_{\mathcal{M}^{-1}}=\ca{O}_K$ and $(\mathcal{M}^{-1})^{-1}=\frac{1}{2}(2\ca{O}_K+\theta\ca{O}_K)=
 \frac{1}{2}(2\ZZ+\theta\ZZ+\theta^2\ZZ)\supsetneqq\mathcal{M}$,
 \item $\mathcal{M}\mathcal{M}^{-1}\subseteq 2\ca{O}_K\subsetneqq\mathcal{R}$.
\end{enumerate}
\end{Exa}
Example \ref{filu} is instructive since it shows that if $\mathcal{L}_1,\mathcal{L}_2\subseteq K$ are $\ca{O}$-proper, then 
$\mathcal{L}_1\mathcal{L}_2$ is not necessarily $\ca{O}$-proper. Moreover, the lattice $\mathcal{M}$ above
is such that $\ca{O}_{\mathcal{M}}\neq\ca{O}_{\mathcal{M}^{-1}}$ and
$\mathcal{M}\subsetneqq({\mathcal{M}}^{-1})^{-1}$. In particular, if $Latt_K$ denotes the set of all lattices in $K$,
then the application $[-1]:Latt_K\rightarrow Latt_K$, given by $\mathcal{L}\mapsto\mathcal{L}^{-1}$ is not necessarily involutive.

\subsubsection{$\ca{O}$-invertibility}
Let $\mathcal{L}\subseteq K$ be an $\ca{O}$-module. By definition
of $\ca{O}_{\mathcal{L}}$, we have $\ca{O}\subseteq\ca{O}_{\mathcal{L}}$. We say that $\mathcal{L}$ is {\bf $\ca{O}$-invertible}, if there exists
an $\ca{O}$-module $\mathcal{L}'\subseteq K$ such that $\mathcal{L}\mathcal{L}'=\ca{O}$. 
Since $\ca{O}_{\mathcal{L}}\cdot\mathcal{L}\mathcal{L}'\subseteq\mathcal{L}\mathcal{L}'=\ca{O}$ and $1\in\mathcal{L}\mathcal{L}'=\ca{O}$,
this implies that $\ca{O}_{\mathcal{L}}\subseteq\ca{O}$. 
Therefore, $\ca{O}=\ca{O}_{\mathcal{L}}$. 
Thus, if $\mathcal{L}$ is an $\ca{O}$-invertible module, it is automatically $\ca{O}$-proper. The converse is not true
in general. Indeed, looking at Example \ref{filu}, we have $\ca{O}_{\mathcal{M}}=\mathcal{R}$ 
and $\mathcal{M}\mathcal{M}^{-1}\subsetneqq\mathcal{R}=\ca{O}_{\mathcal{M}}$.
Therefore, the $\ca{O}_{\mathcal{M}}$-proper lattice $\mathcal{M}$ is not $\ca{O}_{\mathcal{M}}$-invertible.
For a further discussion on the discrepancy between the $\ca{O}$-invertibility and $\ca{O}$-properness,
see Remark \ref{st_err}. Let $\mathcal{L}\subseteq K$ be a lattice and assume that there exists an 
$\ca{O}_{\mathcal{L}}$-module $\mathcal{L}'$, such that $\mathcal{L}\mathcal{L}'=\ca{O}_{\mathcal{L}}$. Then such a
lattice $\mathcal{L}'$ is necessarily unique; indeed, this follows directly from the commutativity and the associativity of the product operation
on lattices.
Moreover, one may verify that $\mathcal{L}'$ must be equal to $\mathcal{L}^{-1}$. It follows from the discussion above that 
\begin{align*}
\mathcal{L}\mathcal{L}^{-1}=\ca{O}_{\mathcal{L}}\Longleftrightarrow
1\in\mathcal{L}\mathcal{L}^{-1}\Longleftrightarrow \mbox{$\mathcal{L}$ is $\ca{O}_{\mathcal{L}}$-invertible}.
\end{align*}
Finally, if $\mathcal{L}$
is $\ca{O}_{\mathcal{L}}$-invertible, then one may also check that $(\mathcal{L}^{-1})^{-1}=\mathcal{L}$ and that
$\ca{O}_{\mathcal{L}}=\ca{O}_{\mathcal{L}^{-1}}$.

It is worthwhile to remind the reader of the following set of equivalences for $\ca{O}$-invertibility which will be used
later on in the proof of Corollary \ref{coro}.
\begin{Th}
Let $\ca{O}\subseteq\ca{O}_K$ be an order and let $\mathcal{L}\subseteq K$ be a lattice, which is also an $\ca{O}$-module. Then 
the following statements are equivalent:
\begin{enumerate}
 \item $\mathcal{L}$ is $\ca{O}$-invertible. 
 \item $\mathcal{L}$ is a projective $\ca{O}$-module.
 \item $\mathcal{L}$ is a locally free $\ca{O}$-module, i.e., for each non-zero prime ideal $\mk{p}\subseteq\ca{O}$, one has
 that $\mathcal{L}_{\mk{p}}$ is a free $\ca{O}_{\mk{p}}$-module. 
\end{enumerate}
\end{Th}

{\bf Proof} For a proof of this fact, the author may consult for example \cite{Rot}. \fin

\subsubsection{Conductor of an order and invertibility of prime ideals}\label{Keith}

\begin{Def}
The conductor of an order $\ca{O}\subseteq K$ is defined as
\begin{align*}
\cond(\ca{O})=\mk{c}_{\ca{O}}=\{x\in K: x\ca{O}_K\subseteq\ca{O}\}. 
\end{align*}
\end{Def}
One may check that $\mk{c}_{\ca{O}}$ is the largest integral $\ca{O}_K$-ideal which is included in $\ca{O}$.
In particular, $\mk{c}_{\ca{O}}$ is an integral $\ca{O}$-ideal. 

Let $\ca{O}\subseteq\ca{O}_K$ be a fixed order. 
An ideal $\mk{a}\subseteq\ca{O}$ (i.e. $\mk{a}$ is an $\ca{O}$-module) is called invertible (relative to $\ca{O}$) 
if $\ca{O}_{\mk{a}}=\ca{O}$. The next proposition gives a complete characterization of the invertible {\it prime} ideals of $\ca{O}$.
\begin{Th}\label{pram}
A non-zero prime ideal $\mk{p}\subseteq\ca{O}$ is $\ca{O}$-invertible if and only if $\mk{p}$ is relatively
prime to $\mk{c}_{\ca{O}}$, i.e., $\mk{p}+\mk{c}_{\ca{O}}=\ca{O}$.
\end{Th}

{\bf Proof} See Theorem 6.1 of \cite{Keith_C1}. \fin

\begin{Rem}
Note that the \lq\lq only if direction\rq\rq\; in Theorem \ref{pram} is no longer true if $\mk{p}$ is not prime.     
For example, assume that $\mk{c}_{\ca{O}}\subsetneqq\ca{O}_K$ and let
$c>1$ be the smallest integer inside $\mk{c}_{\ca{O}}$. Then the $\ca{O}$-ideal $c\cdot\ca{O}$ is invertible (since it is
principal) but it is not coprime to $\mk{c}_{\ca{O}}$.
\end{Rem}

\subsubsection{The dual operation}\label{du_op}
We would like now to recall some classical results about dual lattices with respect to the trace pairing. 
For a lattice $\mathcal{L}\subseteq K$, we define
the {\it dual lattice of $\mathcal{L}$} (with respect to the pairing $(x,y)\mapsto \Tr_{K/\QQ}(xy)$) by
\begin{align*}
\mathcal{L}^*:=\{x\in K: \Tr_{K/\QQ}(x\ell)\in\ZZ\s\mbox{for all}\s \ell\in \mathcal{L}\}.
\end{align*}
Note that the $*$ operation is contravariant on the partially ordered set of lattices, i.e., 
if $\mathcal{L}_1\subseteq\mathcal{L}_2$,
then $\mathcal{L}_1^*\supseteq\mathcal{L}_2^*$. Using the notion of the dual $\ZZ$-basis of a given $\ZZ$-basis of $\mathcal{L}$, 
one easily proves that $\mathcal{L}^{**}=\mathcal{L}$. It follows that the map
$\mathcal{L}\mapsto\mathcal{L}^*$ is an involution on the set $Latt_K$.
Let $\mathcal{L}$ be a lattice. We claim that $\mathcal{L}^*$ is $\ca{O}_\mathcal{L}$-proper,
and, therefore, $\ca{O}_\mathcal{L}=\ca{O}_\mathcal{L^*}$.
Indeed, from the definition of $\mathcal{L}^*$, we see that $\ca{O}_\mathcal{L}\cdot 
\mathcal{L}^*\subseteq \mathcal{L}^*$ and, therefore, $\ca{O}_\mathcal{L}\subseteq \ca{O}_{\mathcal{L}^*}$; conversely, 
substituting $\ca{L}$ by $\ca{L}^*$ in the previous inclusion, and using the identity $\mathcal{L}^{**}=\mathcal{L}$,  we deduce that
$\ca{O}_{\mathcal{L}^*}\subseteq \ca{O}_{\mathcal{L}}$. We may summarize the previous observation by saying that a
lattice $\ca{L}$ is $\ca{O}$-proper if and only if $\ca{L}^{*}$ is $\ca{O}$-proper. Let us point out one 
subtle point regarding the dual operation. If $\mathcal{L}$ is $\ca{O}$-invertible, it does not necessarily follow that $\mathcal{L}^*$ is $\ca{O}$-invertible.
Indeed, from Remark \ref{st_err} (see below), this turns out to be false, exactly when $\ca{O}^*$ is not $\ca{O}$-invertible.

\subsubsection{Dual of an order and the different ideal}

Let $\ca{O}\subseteq\ca{O}_K$ be an order. By definition, we have
\begin{align}\label{head}
\ca{O}^*=\{x\in K:\Tr_{K/\QQ}(xy)\in\ZZ\s\mbox{for all $y\in\ca{O}$}\}.
\end{align}
It follows from \eqref{head} that $\ca{O}^*$ is the {\it largest} $\ca{O}$-module $\ca{L}\subseteq K$ 
such that for all $x\in\ca{L}$, $\Tr_{K/\QQ}(x)\in\ZZ$.

When $\ca{O}=\ca{O}_K$ is the maximal order, every fractional ideal of $K$ is 
$\ca{O}_K$-invertible. In particular, it makes sense to define 
\begin{align*}
\mk{d}_K:=\left((\ca{O}_K)^{*}\right)^{-1}.
\end{align*}
The $\ca{O}_K$-fractional ideal $\mk{d}_K$ is called the {\it different ideal of $K$}. Note
that since $\ca{O}\subseteq\ca{O}_K$, we always have $\mk{d}_K^{-1}\subseteq \ca{O}^*$. 

\subsubsection{Relationship between $\mathcal{L}^{-1}$ and $\mathcal{L}^*$}\label{relat}
We would like now to describe some relationships between the lattices $\mathcal{L}^{-1}$ and $\mathcal{L}^*$.

Let $\ca{L}\subseteq K$ be a lattice. Then $\ca{M}:=\ca{L}\ca{L}^*$ is an $\ca{O}_{\ca{L}}$-module
such that for all $x\in\ca{M}$, $\Tr_{K/\QQ}(x)\in\ZZ$. It thus follows that
\begin{align}\label{head2}
\ca{L}\ca{L}^{*}\subseteq (\ca{O}_{\ca{L}})^*. 
\end{align}
Moreover, for every $x\in\ca{L}^{-1}(\ca{O}_{\ca{L}})^{*}$, a direct computation shows that for
all $y\in\ca{L}$, we have $\Tr_{K/\QQ}(xy)\in\ZZ$. It thus follows from the definition of $\ca{L}^*$ that
\begin{align}\label{head3}
\ca{L}^{-1}(\ca{O}_{\ca{L}})^{*}\subseteq \ca{L}^*.
\end{align}
In particular, in the special case where $\ca{L}$ is $\ca{O}_{\ca{L}}$-invertible, we deduce from
\eqref{head2} and \eqref{head3} that $\ca{L}\ca{L}^*=(\ca{O}_{\ca{L}})^*$, and therefore,
\begin{align}\label{head4}
\ca{L}^{*}=\ca{L}^{-1}(\ca{O}_{\ca{L}})^*.
\end{align}

\begin{Rem}\label{st_err}
When we wrote the papers \cite{Ch4} and \cite{Ch8},
we wrongly thought that $\ca{O}$-properness was equivalent to $\ca{O}$-invertibility. 
Fortunately, this does not affect any of the results of the aforementioned papers, since
this fictive equivalence was never used in any of the proofs. 
However, even though these two notions are not equivalent in general, there is a criterion (may be not so well-known), 
which says exactly when they agree. The following two statements are equivalent:
\begin{enumerate}[(i)]
\item For all lattices $\mathcal{L}\subseteq K$ such that $\ca{O}\subseteq\End(\mathcal{L})$, $\mathcal{L}$ is $\ca{O}$-invertible if and only if it is $\ca{O}$-proper.
\item The $\ZZ$-dual $\ca{O}^*$ of $\ca{O}$, with respect to the trace pairing, is $\ca{O}$-invertible.
\end{enumerate}
It follows from this equivalence, that \eqref{head4} necessarily holds true if $(\ca{O}_{\ca{L}})^*$ is $\ca{O}_{\ca{L}}$-invertible.
We note that condition (ii) is always satisfied, if the order $\ca{O}$ is {\it monogenic}, i.e., if $\ca{O}=\ZZ[\mu]$ for some $\mu\in\ca{O}$.
Therefore, if $\ca{O}$ is monogenic, the notions of $\ca{O}$-invertibility and $\ca{O}$-properness agree.
In particular, $\ca{O}$-properness and $\ca{O}$-invertibility are equivalent, when $K$ is a quadratic field,
since any order of a given quadratic field is monogenic. 
For a proof of the equivalence between (i) and (ii), see, for example, Theorem 4.1 of \cite{Keith_C1}. 
\end{Rem}

\subsubsection{Index and covolume}\label{ind_cov}
Let $\mathcal{L}_1,\mathcal{L}_2\subseteq K$ be two lattices. We define the rational index 
\begin{align*}
[\mathcal{L}_1:\mathcal{L}_2]
\end{align*}
as the absolute value of 
the determinant of any $g$-by-$g$ matrix with rational entries which takes a $\ZZ$-basis of $\mathcal{L}_1$ to a $\ZZ$-basis of $\mathcal{L}_2$. 
So we always have $[\ca{L}_1:\ca{L}_2]\in\QQ_{>0}$. The rational index satisfies the transitivity formula $[\ca{L}_1:\ca{L}_2][\ca{L}_2:\ca{L}_3]=[\ca{L}_1:\ca{L}_3]$ for all lattices
$\ca{L}_1,\ca{L}_2,\ca{L}_3\subseteq K$. We also define
\begin{align*}
\Norm(\ca{L}):=[\ca{O}_K:\ca{L}]\in\QQ_{>0}.
\end{align*}
Let $\{\rho_1,\rho_2,\ldots,\rho_r,\sigma_1,\ov{\sigma}_1,\ldots,\sigma_s,\ov{\sigma}_s\}$ ($r+2s=g$), 
be the set of all embeddings of $K$ into $\CC$. Let 
$\theta:K \hookrightarrow \RR^g$ be the embedding given by
\begin{align*}
x \mapsto (\rho_1(x),\ldots,\rho_r(x),2\Ree(\sigma_1(x)),2\Imm(\sigma_1(x)),\ldots,2\Ree(\sigma_s(x)),2\Imm(\sigma_s(x))).
\end{align*}
It is well-known that
\begin{align}\label{oups}
\cov(\mathcal{L}):=\sqrt{|d_K|}\cdot[\ca{O}_K:\mathcal{L}],
\end{align}
where $\cov(\mathcal{L})$ corresponds to the covolume of $\theta(\mathcal{L})$ inside $\RR^g$ with respect
to the Lebesgue measure. The following lemma relates $\cov(\mathcal{L})$ to $\cov(\mathcal{L}^*)$.
\begin{Lemma}\label{bouh}
We have $\cov(\mathcal{L})\cov(\mathcal{L}^*)=1$. 
\end{Lemma}

{\bf Proof} Let $\{e_i\}_{i=1}^{g}$ be a $\ZZ$-basis of $\ca{O}_K$. Let $B=(b_{ij})_{i,j}$ be the $g$ by $g$
square matrix where $b_{ij}=\Tr_{K/\QQ}(e_ie_j)$. Recall that $|d_K|=|\det(B)|$. Let $\{v_i\}_{i=1}^g$
be a $\ZZ$-basis of $\mathcal{L}$. Viewing $B$ as a $\QQ$-bilinear form on $\QQ^g$, we may take the $\ZZ$-dual basis of $\{v_i\}_{i=1}^g$
with respect to $B$, which we denote by $\{v_{i}^{*}\}$. It follows from the definition of $\mathcal{L}^*$ that
$\{v_{i}^{*}\}$ is a $\ZZ$-basis of $\mathcal{L}^*$. Let $T$ (resp. $T^*$) be a matrix (with entries in $\QQ$) such that
$T$ sends $\{e_i\}_{i=1}^{g}$ to $\{v_{i}\}_{i=1}^g$ (resp. on $\{v_{i}^*\}$). We thus have the matrix identity 
$(T^{t})BT^*=I_g$, where $I_g$ is the $g$ by $g$ identity matrix. Taking the determinant, we find 
$[\ca{O}_K:\mathcal{L}]\cdot|d_K|\cdot[\ca{O}_K:\mathcal{L}^*]=1$. The result follows from \eqref{oups}. \fin
\begin{Cor}\label{coro}
For all lattices $\mathcal{L}\subseteq K$, one has 
$[\ca{O}_K:\mathcal{L}]\cdot[\ca{O}_K,\mathcal{L}^*]=\Norm(\ca{L})\cdot\Norm(\ca{L}^*)=\frac{1}{|d_K|}$.
\end{Cor}

\subsubsection{An index formula for the product of two lattices}
We would like now to prove a useful index formula.
\begin{Prop}\label{orte}
Let $\mathcal{O}\subseteq \ca{O}_K$ be an order. Let $\mathcal{L},\mathcal{M}\subseteq K$ be two lattices which are also $\ca{O}$-modules.
Assume that $\mathcal{L}$ or $\mathcal{M}$ is $\ca{O}$-invertible. Then $[\ca{O}:\mathcal{L}]\cdot[\ca{O}:\mathcal{M}]=[\ca{O}:\mathcal{L}\mathcal{M}]$.
\end{Prop}

{\bf Proof} Without loss of generality, one may assume that $\mathcal{L},\mathcal{M}\subseteq\ca{O}$ and that $\mathcal{L}$ is
$\ca{O}$-invertible. We have $[\ca{O}:\mathcal{L}\mathcal{M}]=[\ca{O}:\mathcal{L}][\mathcal{L}:\mathcal{L}\mathcal{M}]$. We claim
that there exists an abelian group isomorphism (non-canonical!) between the two finite $\ZZ$-modules $\mathcal{L}/\mathcal{L}\mathcal{M}$ and $\ca{O}/\mathcal{M}$.
Let $p\in\ZZ_{\geq 2}$ be a prime and set $S_p:=\ca{O}\bs(\mk{p}_1\cup\mk{p_2}\cup\ldots\cup\mk{p}_e)$, where $\mk{p}_i$'s are the distinct 
prime ideals of $\ca{O}$ above $p\ZZ$. The set $S_p$ is multiplicatively closed and the localized ring $\ca{O}_p:=S_p^{-1}\ca{O}$ is a semi-local
ring. Since $\mathcal{L}$ is $\ca{O}$-invertible, it follows that $\mathcal{L}_p:=S_p^{-1}\mathcal{L}$ is a principal $\ca{O}_p$-module.
Therefore, there exists $\pi_p\in\ca{O}_p$, such that $\mathcal{L}_p=\pi_p\ca{O}_p$. Since $\ca{O}_p$ is a flat $\ca{O}$-module, we have
$\mathcal{L}_p/\mathcal{L}_p\mathcal{M}_p\simeq \left(\mathcal{L}/\mathcal{L}\mathcal{M}\right)_p$ and 
$\ca{O}_p/\mathcal{M}_p\simeq \left(\ca{O}/\mathcal{M}\right)_p$. Now consider the map $\varphi:\ca{O}_p/\mathcal{M}_p\rightarrow 
\mathcal{L}_p/\mathcal{L}_p\mathcal{M}_p$,
given by $x+\mathcal{M}_p\mapsto x\pi_p+\mathcal{L}_p\mathcal{M}_p$. It follows that $\varphi$ is an $\ca{O}_p$-module isomorphism. In particular,
we have
\begin{align*}
\mbox{($p$-primary subgroup of $\mathcal{L}/\mathcal{L}\mathcal{M}$) $\simeq\left(\mathcal{L}/\mathcal{L}\mathcal{M}\right)_p$
$\simeq$ ($p$-primary subgroup of $\ca{O}/\mathcal{M}$ )$\simeq \left(\ca{O}/\mathcal{M}\right)_p$} 
\end{align*}
Finally, since $p$ was arbitrary, it follows that $\mathcal{L}/\mathcal{L}\mathcal{M}$, as a finite $\ZZ$-module, 
is isomorphic to $\ca{O}/\mathcal{M}$. \fin

\subsection{Fourier series on the standard $g$-dimensional real torus}\label{Four}

Let $\mathbb{T}^g:=\RR^g/\ZZ^g$ be the standard $g$-dimensional real torus. We think of $\mathbb{T}^g$ as a real analytic manifold
via the natural projection $\pi:\RR^g\rightarrow\mathbb{T}^g$. We would like to record some general results about the convergence 
of Fourier series and the growth of Fourier coefficients of a given function $f:\mathbb{T}^g\rightarrow\CC$. It will be important
for us to impose certain conditions on $f$, so that its Fourier series {\it computes} the value of $f$ at {\it all points of $\mathbb{T}^g$}.
Even though all the functions which appear in this work
are real analytic, we decided to state more general results which can be applied to not necessarily real analytic functions. 
We let $C^k(\mathbb{T}^g)$ 
be the space of functions $f:\mathbb{T}^g\rightarrow\CC$ such that $f$ is of class $C^k$, i.e., all the 
derivatives of $f$ of order less than or equal to $k$ exist and are continuous. We say that a function 
$f:\mathbb{T}^g\rightarrow\CC$ is {\it smooth} if it is of class $C^k$ for all $k\geq 0$. In particular, a real
analytic function $f:\mathbb{T}^g\rightarrow\CC$ is always smooth. We also let
$L^1(\mathbb{T}^g)$  (resp. $L^2(\mathbb{T}^g)$) be the space of integrable functions (resp. square integrable functions)
on $\mathbb{T}^g$ with respect to the Haar measure $dx$ on $\mathbb{T}^g$ normalized, so that
$\int_{\mathcal{T}^g} dx=(2\pi)^g$.
For any $f\in L^1(\mathbb{T}^g)$ and $n\in\ZZ^g$, we let
\begin{align*}
a_{n}(f):=\int_{\mathbb{T}^g} e^{-{2\pi\ii } \langle n, x\rangle} f(x)dx,
\end{align*}
be the $n$-th Fourier coefficient of $f$. Here $\langle\;,\;\rangle$ corresponds to
the standard inner product on $\RR^g$. We also let $||n||_2=\sqrt{\laa n,n\raa}$ be 
the $\ell_2$-norm of the vector $n\in\RR^g$. We may now state a collection of key results regarding Fourier series
of functions on $\mathbb{T}^g$.
\begin{Th}\label{keyy}
\begin{enumerate}
 \item  Let $f\in C^k(\mathbb{T}^g)$ with $k>\frac{g}{2}$. Then $\sum\limits_{n\in\ZZ^g} |a_{n}(f)|<\infty$. Moreover, for all $x\in\mathbb{T}^g$, one has
 $\sum\limits_{n\in\ZZ^g} a_{n}(f)e^{2\pi\ii \langle n, x\rangle}=f(x)$.
 \item (Parseval's theorem) Let $f\in L^2(\mathbb{T}^g)$. Then $\int\limits_{\mathbb{T}^g}  |f(x)|^2dx=\sum\limits_{n\in\ZZ^g} |a_{n}(f)|^2$.
 \item (Paley-Wiener principle) Let $f\in C^k(\mathbb{T}^g)$ for some $k\geq 0$. Then $\lim\limits_{||n||_2\rightarrow\infty} 
 \frac{|a_n(f)|}{1+||n||_2^k}\rightarrow 0$. Furthermore, if $f$ is a real analytic, 
 then there exists a positive constant $D>0$, such that for all $n\geq 0$, 
 $|a_n(f)|\leq e^{-D ||n||_2}$.
\end{enumerate}
\end{Th}

{\bf Proof} For a proof of (1) see
Corollaries 1.8 and 1.9 on p. 249 of \cite{Stein-Weiss}. The proof of (2) is classical and may be found
in most textbooks on Fourier series. For
a proof of the first part of (3), see for example Theorem 3.2.9 of \cite{Gra}. For a proof of the
second part of (3), see for example Proposition 5.4.1 of \cite{Krantz-Parks}. \fin 

The following technical lemma will also be needed. Note that the hypotheses are far from being optimal.
\begin{Lemma}\label{dub}
Let $U\subseteq \RR^r$ be an non-empty open subset. Let 
\begin{align*}
F:\mathbb{T}^g\times U &\rightarrow\CC\\
                 (x,y)      & \mapsto F(x,y)
\end{align*}
be a smooth function, where $x\in\mathbb{T}^g$ and $y\in U$. Let 
\begin{align*}
F(x,y)=\sum_{n\in\ZZ^g} a_{n}(y)e^{2\pi\ii\laa n, x\raa},
\end{align*}
be the Fourier series expansion of $[x\mapsto F(x,y)]$ which exists and computes the value of $F(x,y)$
by (1) of Theorem \ref{keyy}. Then, for each $n\in\ZZ^g$, the
function $[y\mapsto a_{n}(y)]$ is smooth on $U$. 

Let $D$ the linear differential operator $D:=g(y)\frac{\partial^{i_1}}{\partial y_1^{i_1}}
\frac{\partial^{i_2}}{\partial y_2^{i_2}}\ldots \frac{\partial^{i_e}}{\partial y_r^{i_r}}$, where $i_1,\ldots,i_r\geq 0$,
and $g(y)$ is a polynomial in $y$. Then
\begin{align}\label{ore2}
D F(x,y)=
\sum_{n\in\ZZ^g} \left(D a_{n}(y)\right)e^{2\pi\ii\laa n, x\raa}.
\end{align}
Moreover, if we let
$D'=\frac{\partial^{i_1}}{\partial x_1^{i_1}}\frac{\partial^{i_2}}{\partial x_2^{i_2}}
\ldots \frac{\partial^{i_e}}{\partial x_g^{i_g}}$, where $i_1,\ldots,i_g\geq 0$, then
\begin{align}\label{ore3}
D' F(x,y)=\sum_{n\in\ZZ^g} a_{n}(y) \left(D'e^{2\pi\ii\laa n, x\raa}\right).
\end{align}
\end{Lemma}

{\bf Proof} Let us start by proving that the functions $[y\mapsto a_n(y)]$ are smooth.
By definition, we have $a_{n}(y)=\int_{\mathbb{T}^g}F(x,y)e^{-2\pi\ii\laa n, x\raa}dx$. 
Let $D:=\frac{\partial}{\partial y_1}$.
For each compact set $K'\subseteq U$, the function
$[(x,y)\mapsto D F(x,y)]$ is continuous on $K:=\mathbb{T}^g\times K'$. Since $K$ is compact, it follows that 
$[(x,y)\mapsto D F(x,y)]$ is uniformly continuous on $K$ and, therefore, is bounded by a constant.
Moreover, for a fixed $y\in U$, the
function $[x\mapsto D F(x,y)]$ is integrable on $\mathbb{T}^g$. It follows
from all of our assumptions (see for example Lemma 2.2 on page 226 of \cite{La4}), that 
\begin{align}\label{halp}
D a_{n}(y)=
\int_{\mathbb{T}^g} DF(x,y)e^{-2\pi\ii\laa n, x\raa}dx,
\end{align}
i.e., one is allowed to differentiate under the integral sign. 
In particular, it follows from \eqref{halp}, that $[y\mapsto a_{n}(y)]$ is continuous on $K'$. By 
induction on the of $D$, it follows that $[y\mapsto a_{n}(y)]$ is smooth.
Finally, since $K'\subseteq U$ was an arbitrary compact set, this implies that the 
function $[y\mapsto a_{n}(y)]$ is smooth on all of $U$. 

Let us  show \eqref{ore2}. Let $D:=g(y)\frac{\partial^{i_1}}{\partial y_1^{i_1}}$ where $g(y)$
is a polynomial in $y=(y_j)_{j=1}^r$. Let $y\in U$ be fixed. The function $[x\mapsto D F(x,y)]$
is smooth and invariant under a translation by an element in $\mathbb{T}^g$. Therefore, from (1) of Theorem \ref{keyy}, we have
that for all $x\in\mathbb{T}^g$,
\begin{align*}
D F(x,y)=\sum_{n\in\ZZ^g} b_{n}(y)e^{2\pi\ii\laa n,x\raa},
\end{align*}
where $b_{n}(y)\stackrel{def}{=}\int_{\mathbb{T}^g} D F(x,y)e^{-2\pi\ii\laa n, x\raa}dx$. Similarly to the proof
of \eqref{halp}, we have
\begin{align*}
b_{n}(y)=\int_{\mathbb{T}^g} D F(x,y)e^{-2\pi\ii\laa n, x\raa}dx=
D\int_{\mathbb{T}^g} F(x,y)e^{-2\pi\ii\laa n, x\raa}dx=
D a_n(y).
\end{align*}
By induction on the order of $D$ this proves \eqref{ore2}.

It remains to prove \eqref{ore3}. Let $D'=\partial_{x_i}:=\frac{\partial}{\partial x_i}$. Note that
$[(x,y)\mapsto\partial_{x_i} F(x,y)]$ is again a smooth function of $\mathbb{T}^g\times U$. It particular, it admits
a Fourier expansion which computes its value. We have
\begin{align*}
\partial_{x_j} F(x,y)=\sum_{n\in\ZZ^g} c_{n}(y)e^{2\pi\ii\laa n,x\raa},
\end{align*}
where $c_{n}(y)=\int_{\mathbb{T}^g}\left(\partial_{x_j}F(x,y)\right)e^{-2\pi\ii\laa n, x\raa}dx$. Integrating 
by parts, we find that
\begin{align*}
\int_{\mathbb{T}^g}(\partial_{x_i} F(x,y))e^{-2\pi\ii\laa n, x\raa}dx &=
\int_{\mathbb{T}^g}\partial_{x_i} \left(F(x,y)e^{-2\pi\ii\laa n, x\raa}\right)dx-
\int_{\mathbb{T}^g} F(x,y) \left(\partial_{x_i} e^{-2\pi\ii\laa n, x\raa}\right)dx \\
&=2\pi\ii n_j\int_{\mathbb{T}^g} F(x,y) e^{-2\pi\ii\laa n, x\raa}dx\\
&=2\pi\ii n_j a_{n}(y),
\end{align*}
where $n_j$ is the $j$-th coordinate of $n$. It thus follows that
\begin{align*}
\partial_{x_j} F(x,y)= \sum_{n\in\ZZ^g} a_{n}(y) \left(\partial_{x_j}e^{2\pi\ii\laa n,x\raa}\right).
\end{align*}
The result follows by induction on the order of $D'$. This concludes the proof. \fin

\subsection{The monomial $P(\alpha,\beta;z)$ and the product convention}\label{mon_con}
In this section we introduce the basic monomial $P(\alpha,\beta;z)$ and the product convention which gives a rule
on how one should compute $P(\alpha,\beta;z)$.

For $\alpha,\beta\in\CC^{g}$ and $z\in K_\CC^{\pm}$, we define
\begin{align}\label{expo}
P(\alpha,\beta;z):=
z^{\alpha}(\ov{z})^{\beta}=\left(\prod_{i=1}^g(z_i)^{\alpha_i}\prod_{i=1}^{g}(\ov{z_i})^{\beta_i}\right).
\end{align} 
It is a \lq\lq monomial\rq\rq\; in the variables $z$ and $\ov{z}$ where the exponents are allowed to be
complex numbers (not just positive integers). Recall that $z_i^{\alpha_i}:=e^{\alpha_i\log z_i}$ where $\log z_i$ is the principal branch of the logarithm. 
Let us note that
\begin{align*}
P(-\alpha,-\beta;z)=\frac{1}{P(\alpha,\beta;z)}.
\end{align*}
If the weights $\alpha,\beta\in\CC$ are such that $\beta-\alpha=p\in\ZZ^g$, then a 
direct computation shows that, for any $z\in K_{\CC}^{\pm}$, the following identity hold true: 
\begin{align}\label{hibb_00}
P(\alpha,\beta;z)=\frac{|\Norm(z)|^{\alpha+\beta}}{\omega_p(z)}.
\end{align}
Note that \eqref{hibb_00} is still valid if $z\in(\RR_{>0})^g$. However, if $z\in K_\CC^{\times}$
and if $z$ admits a coordinate $z_j$ such that $z_j\in\RR_{<0}$, then \eqref{hibb_00} may fail to be true.
It will be important for the whole work to be able to evaluate the monomial $P(\alpha,\beta;z)$
when $z\in(\RR^{\times})^{g}$, in a way which is compatible with formula \eqref{hibb_00}. Therefore,
we make the following convention:
\begin{Convention}\label{conv}
Let $\alpha_1,\ldots,\alpha_n,\beta_1,\ldots\beta_n,\in\CC$ and let 
$z_1,\ldots,z_n,w_1,\ldots,w_n \in\CC$ be complex variables. Consider the following 
formal products of two types
\begin{enumerate}
 \item  $F(z_1,\ldots,z_n;w_1,\ldots,w_n):=\prod\limits_{j=1}^n z_j^{\alpha_j}\cdot\prod\limits_{j=1}^n(\ov{w}_j)^{\beta_j}$. 
 \item $G(z_1,\ldots,z_n;w_1,\ldots,w_n):=\prod\limits_{j=1}^n (z_j w_j)^{\alpha_j}\prod\limits_{j=1}^n(\ov{z}_j\cdot\ov{w}_j)^{\beta_j}$. 
\end{enumerate}
The \lq\lq bar\rq\rq\; on the $w_i$'s should be thought of as a \lq\lq formal conjugate\rq\rq. 
Let $(a_1,\ldots,a_n;b_1,\ldots,b_m)\in(\CC^{\times})^{g}$. Then we define
\begin{enumerate}
 \item  $F(a_1,\ldots,a_n;b_1,\ldots,b_m):=\lim\limits_{\substack{z_i\rightarrow a_i\\ w_i\rightarrow b_i}} F(a_1,\ldots,a_n;b_1,\ldots,b_m)$
 \item $G(a_1,\ldots,a_n;b_1,\ldots,b_m):=\lim\limits_{\substack{z_i\rightarrow a_i\\ w_i\rightarrow b_i}} G(a_1,\ldots,a_n;b_1,\ldots,b_m)$
\end{enumerate}
where the $z_i,w_i\in \CC^{\times}$ are subjected to the condition that 
\begin{enumerate}
 \item If $a_i\in\RR_{<0}$ then the sequence $\{z_{ik}\}_{k\geq 1}$ converging to $a_i$ is such that
 $0<\arg(z_{ik})<\pi$. 
 \item If $b_i\in\RR_{<0}$ then the sequence $\{w_{ik}\}_{k\geq 1}$ converging to $b_i$ is such that
 $0<\arg(w_{ik})<\pi$. 
\end{enumerate}
In other words, when approaching a negative real number (in the euclidean topology) 
by non-real complex numbers, we do so by using complex numbers in the upper half plane.
\end{Convention}
More generally, we could use the same convention for any algebraic expression in a set of holomorphic variables, 
conjugate variables and different choices of parentheses. For example, we could consider the expression 
$(z_1w_1z_2\ov{w}_2)^{\alpha_1}\cdot (\ov{z}_1 z_2)^{\alpha_2}$, etc. However, for the whole work, only
the formal products (1) and (2) will be needed.

Using this convention, we have the following important proposition:
\begin{Prop}\label{koy_prop}
Let $\alpha,\beta\in\CC^g$ be two weight vectors, such that $\beta-\alpha=p\in\ZZ^g$. Then, for
any $z,w\in K_\CC^{\times}$, the following identities, computed according to Convention \ref{conv}, hold true:
\begin{enumerate}
 \item $P(\alpha,\beta;z)=\frac{|\Norm(z)|^{\alpha+\beta}}{\omega_p(z)}$,
 \item $(z w)^{\alpha}\cdot(\ov{z}\cdot \ov{w})^{\beta}=z^{\alpha}\cdot w^{\alpha}\cdot \ov{z}^{\beta}\cdot\ov{w}^{\beta}$.
\end{enumerate}
\end{Prop}
\begin{Rem}
If $z\in K_\CC^{\pm}\cup\RR_{>0}^g$, then the identities (1) and (2) are valid when 
computed (normally, i.e., without Convention \ref{conv}) in terms of the principal branch of the logarithm on
$\CC\bs\RR_{\leq 0}$.
\end{Rem}

\subsection{The Fourier series expansion of $R_{\mathcal{L}}(\alpha,\beta;z)$}\label{dos}
In this Section, we introduce certain basic infinite sums for which we compute their Fourier series. These infinite sums will be the building
blocks of the real analytic Eisenstein series which appear in Definition \ref{difo}. 

Let $\mathcal{L}\subseteq K$ be a fixed lattice. For a fixed pair $\alpha=(\alpha_i)_{i=1}^g,\beta=(\beta)_{i=1}^g\in\CC^g$ such
that, $\Ree(\alpha_j+\beta_j)>1$ for all $j\in\{1,2,\ldots,g\}$, and $z\in K_\CC^{\pm}$, we define
\begin{align}\label{pap}
R_{\mathcal{L}}(\alpha,\beta;z)&:=\sum_{v\in \mathcal{L}}
(z+v)^{-\alpha}\cdot(\ov{z}+v)^{-\beta}\\  \notag
&=\sum_{v\in \mathcal{L}}P(-\alpha,-\beta;z+v).
\end{align}
Note that the convergence of the summation in \eqref{pap} is absolute.

We would like now to give the Fourier series expansion of $[x\mapsto R_{\mathcal{L}}(\alpha,\beta;x+\ii y)]$. From (1) of Theorem \ref{keyy}, we may deduce that 
\begin{align}\label{subi}
R_{\mathcal{L}}(\alpha,\beta;z)=\sum_{\xi\in \mathcal{L}^*}a_{\xi}(\alpha,\beta;y)e^{{2\pi\ii } \Tr(\xi x)},
\end{align}
for all $z=x+\ii y\in K_\CC^{\pm}$, where $\xi x=(\xi^{(i)}x_i)_{i=1}^g$ and
\begin{align*}
a_{\xi}(\alpha,\beta;y):=\frac{1}{\cov(\mathcal{L})}\int_{\RR^g/\mathcal{L}} e^{-{2\pi\ii } \Tr(\xi x)}
R_{\mathcal{L}}(\alpha,\beta;z)dx.
\end{align*}
Here, $\cov(\ca{L})$ corresponds to the covolume of the lattice $\ca{L}$ (see Section \ref{ind_cov}). 
For a fix element $\xi\in \mathcal{L}^*$ in the dual lattice, we have
\begin{align*}
\int_{\RR^g/\mathcal{L}} e^{-{2\pi\ii } \Tr(\xi x)}R_{\mathcal{L}}(\alpha,\beta;z)dx
&=\int_{\RR^g/\mathcal{L}} e^{-{2\pi\ii } \Tr(\xi x)}\sum_{v\in \mathcal{L}}P(-\alpha,-\beta;z+v)dx\\
&=\int_{\RR^g/\mathcal{L}}e^{{2\pi\ii } \Tr(\xi v)}\sum_{v\in \mathcal{L}}e^{-{2\pi\ii } \Tr(\xi(x+v))}P(-\alpha,-\beta;x+iy+v)dx\\
&=\int_{\RR^g}e^{-{2\pi\ii } \Tr(\xi x)}P(-\alpha,-\beta;x+iy)dx.
\end{align*}
Moreover,
\begin{align*}
\int_{\RR^g}e^{-{2\pi\ii } \Tr(\xi x)}P(-\alpha,-\beta;x+\ii y)dx&=
\prod_{j=1}^g\int_{-\infty}^{\infty}e^{-{2\pi\ii }\xi^{(j)}x_j}z_j^{-\alpha_j}(\ov{z_j})^{-\beta_j}dx_j\\
&=\prod_{j=1}^g\tau(\alpha_j,\beta_j;\xi^{(j)},y_j).
\end{align*}
Therefore, the $\xi$-th Fourier coefficient of $R_{\mathcal{L}}(\alpha,\beta;z)$ is given explicitly by
\begin{align}\label{coeff}
a_{\xi}(\alpha,\beta;y)=\cov(\mathcal{L})^{-1}\prod_{j=1}^g\tau(\alpha_j,\beta_j;\xi^{(j)},y_j).
\end{align}

\subsubsection{A limit formula of $R_{\mathcal{L}}(\alpha,\beta;z)$ when $K=\QQ$}\label{lim_for}
Let us give an nice application which makes use of the explicit Fourier series expansion \eqref{subi}. 
For $x\in\RR\bs\ZZ$, and $s\in\Pi_1$, let 
\begin{align}\label{psi_0}
\zeta(x,s):=\sum\limits_{n\in\ZZ}\frac{1}{|n+x|^s}.
\end{align}
Note that the convergence of the series \eqref{psi_0} is absolute, 
$\zeta(x+1,s)=\zeta(x,s)$ and $\zeta(-x,s)=\zeta(x,s)$. For $x\in\RR$ and $s\in\Pi_1$, let
\begin{align}\label{psi}
\Psi(x,s)=\sum_{k\in\ZZ\bs\{0\}}\frac{e^{2\pi\ii kx}}{|k|^s} \s\s\s.
\end{align}
Similarly to \eqref{psi_0}, for $s\in\Pi_1$ ,the convergence of the series \eqref{psi} is absolute, 
$\Psi(x+1,s)=\Psi(x,s)$ and $\Psi(-x,s)=\Psi(x,s)$. 

For $x\in\RR\bs\ZZ$ and $s$ with $0<\Ree(s)\leq 1$, the series \eqref{psi} converges {\it conditionally} (using the summation
by parts technique) since $\sum_{k=0}^N e^{2\pi\ii kx}$ is bounded as $N\rightarrow\infty$. In this case, the infinite sum 
$\sum\limits_{k\in\ZZ\bs\{0\}}\frac{e^{2\pi\ii kx}}{|k|^s}$ is taken to mean the following limit:
\begin{align}\label{pora}
\lim\limits_{N\rightarrow\infty}\sum\limits_{k=-N}^N\frac{e^{2\pi\ii kx}}{|k|^s},
\end{align}
i.e., the limit of symmetric partial sums. It also follows from the summation by parts technique that
$s\mapsto \Psi(x,s)$ is holomorphic on all of $\Pi_0$. In particular, using the usual Taylor expansion of $\log(1-x)$
at $x=0$, we may deduce from \eqref{pora} that 
\begin{align*}
\Psi(x,1)=-\log(1-e^{2\pi\ii x})-\log(1-e^{-2\pi\ii x})=-2\log |1-e^{2\pi\ii x}|.
\end{align*}
Since $x\in\RR\bs\ZZ$, it follows that $\Psi(x,1)\neq 0$. 

It was proved by Lipschitz (see Epstein's comment 
mentioned above Equation (7) on p. 618 of \cite{Eps}), that $s\mapsto \zeta(x,s)$ and $s\mapsto \Psi(x,s)$ admit a meromorphic continuation to all of $\CC$
and that both functions are interrelated by the following functional equation:
\begin{align}\label{func}
\pi^{-\frac{s}{2}}\Gamma\left(\frac{s}{2}\right)\zeta(x,s)=
\pi^{-\frac{1-s}{2}}\Gamma\left(\frac{1-s}{2}\right)\Psi(-x,1-s).
\end{align}
The functional equation above also follows from the main result proved in \cite{Ch4}.
With the help of the Euler's reflection formula and of the duplication formula 
for the gamma function, we can rewrite \eqref{func} as:
\begin{align}\label{etir}
\zeta(x,s)&=2(2\pi)^{s-1}\Gamma(1-s)\sin\left(\pi s/2\right)\Psi(x,1-s).
\end{align}
In particular, if $s\in\CC$ is fixed and is such that $\Ree(1-s)>0$, then \eqref{etir} can be written as
\begin{align}\label{etir2}
\zeta(x,s)=2(2\pi)^{s-1}\Gamma(1-s)\sin\left(\pi s/2\right)\sum_{k\in\ZZ\bs\{0\}}\frac{e^{2\pi\ii kx}}{|k|^{1-s}}.
\end{align}
In particular, the right-hand side of \eqref{etir2} may be interpreted as the Fourier series expansion
of $x\mapsto\zeta(x,s)$

For $z=x+\ii y\in\CC\bs\ZZ$ and $s\in\Pi_1$, let 
\begin{align}\label{psif}
\wt{\zeta}(z,s):=\sum\limits_{n\in\ZZ}\frac{1}{|n+z|^s}.
\end{align}
Note that $\wt{\zeta}(z+1,s)=\wt{\zeta}(z,s)$ and $\wt{\zeta}(-z,s)=\wt{\zeta}(z,s)$. Obviously, for
$x\in\RR\bs\ZZ$, we have $\zeta(x,s)=\wt{\zeta}(x,s)$. We would like now to explain how \eqref{etir2} may be interpreted 
as the limit of the Fourier series expansion of $x\mapsto\wt{\zeta}(x+\ii y,s)$ when $z=x+\ii y\in\mk{h}^{\pm}$ and
$y\rightarrow 0$. In particular, this provides a \lq\lq new proof\rq\rq\s of \eqref{etir} which does not use 
the Poisson summation formula (even though both proofs rely ultimately on Fourier series expansion associated to functions on the circle). 
From \eqref{subi}, for a fixed $y\in\RR^{\times}$, the Fourier series of $[x+\ii y\mapsto\wt{\zeta}(x+\ii y,s)]$ is given by
\begin{align*}
\wt{\zeta}(z,s)&=\sum_{n\in\ZZ}\tau\left(\frac{s}{2},\frac{s}{2},n,y\right)e^{2\pi\ii nx}\\
      &=2\pi\frac{\Gamma(s-1)}{\Gamma(s/2)^2}(2|y|)^{1-s}+\frac{(2\pi)^s}{\Gamma(s/2)} \sum_{n\neq 0}
      n^{s-1} U\left(\frac{s}{2},s,4\pi|ny|\right)e^{-2\pi|ny|}e^{2\pi\ii nx},\\
\end{align*} 
where the second equality follows from Lemma \ref{Shim_lem}. Note that if $y=0$, then the \lq\lq function\rq\rq\; $[x\mapsto\wt{\zeta}(x,s)]$, for
$x\in\RR$, has singularities at each element of $\ZZ$. This prevents us, a priori, to define its Fourier coefficients.

From Corollary \ref{triop0}, we see that for a fixed $z=x+\ii y\in\mk{h}^{\pm}$, the function $s\mapsto \wt{\zeta}(z,s)$ admits a single-valued
holomorphic continuation to all of $\CC$ which we still denote by $\wt{\zeta}(z,s)$. (In particular, note that for a fixed even integer $k\in\ZZ_{\leq 0}$,
the function $[z\mapsto \wt{\zeta}(z,k)]$ is identically equal to zero). Now fixing $s$ with $\Ree(s)<1$, we may deduce that
\begin{align*}
\zeta(x,s)=\lim_{y\rightarrow 0} \wt{\zeta}(z,s)&=\lim_{y\rightarrow 0} \left[ 2\pi\frac{\Gamma(s-1)}{\Gamma(s/2)^2}(2|y|)^{1-s}+
\frac{(2\pi)^s}{\Gamma(s/2)} \sum_{n\neq 0}
      n^{s-1}U\left(\frac{s}{2},s,4\pi|ny|\right)e^{-2\pi|ny|}e^{2\pi\ii nx} \right]\\
    &=\lim_{y\rightarrow 0}\left[\frac{(2\pi)^s}{\Gamma(s/2)} \sum_{n\neq 0}
      n^{s-1}U\left(\frac{s}{2},s,4\pi|ny|\right)e^{-2\pi|ny|}e^{2\pi\ii nx} \right]\\
      &=2(2\pi)^{s-1}\Gamma(1-s)\sin\left(\pi s/2\right)\Psi(x,1-s).
\end{align*}
The last equality used the fact that the  limit commutes with the summation (which can be justified since 
$x\notin\ZZ$), (1) of Proposition \ref{geno} and the Euler's reflection formula for the gamma function.

\section{Basic properties of real analytic Eisenstein series}\label{chiv}
\subsection{Distinguished subgroups of $GL_1(K)$ and $GL_2(K)$}\label{hubo}
In this section, we introduce various subgroups of $GL_1(K)$ and $GL_2(K)$ which 
naturally intervene in the study of the Eisenstein series 
$G_{(\mk{m},\mk{n})}^w(U,p;z,s)$. 

By definition, recall that a lattice $\ca{L}\subseteq K$ is always supposed to be 
a $\ZZ$-module of {\it maximal} rank, i.e., 
$\ca{L}\simeq\ZZ^g$.
\begin{Def}\label{pilo}
For each row vector $(m,n)\in K^2$, we let $GL_2(K)$ act on the right by the usual matrix multiplication. 
Let $\mk{n},\mk{m}\subseteq K$ be two lattices. We define the sets
\begin{enumerate}
 \item $GL(\mk{m}\oplus\mk{n}):=\{\gamma\in GL_2(K):(\mk{m}\oplus\mk{n})\gamma= \mk{m}\oplus\mk{n}\}$,
 \item $GL^+(\mk{m}\oplus\mk{n}):=\{\gamma\in GL(\mk{m}\oplus\mk{n}):\sg(\det(\gamma))=\ov{\bz}\}$. 
 \item $SL(\mk{m}\oplus\mk{n}):=\{\gamma\in GL_2(\mk{m}\oplus\mk{n}):\det(\gamma)=\bu\}$. 
 \item If $\ca{O}\subseteq\ca{O}_K$ is an order and $\mk{n}\subseteq \ca{O}$ is an integral 
 $\ca{O}$-ideal, we define 
\begin{align*}
\Gamma_{\ca{O}}(\mk{n})=\left\{\gamma\in GL_2(\ca{O}):\gamma\equiv I_2\pmod{\mk{n}}\right\},
\end{align*}
and $\Gamma_{\ca{O}}^+(\mk{n})=\Gamma_{\ca{O}}(\mk{n})\cap GL_2^+(K)$. Here $I_2=\M{1}{0}{0}{1}$ corresponds
to the identity matrix.
 \end{enumerate}
\end{Def}
One may check that all these sets are in fact subgroups of $GL_2(K)$. It is also straight forward from the definition
of $GL(\mk{m},\mk{n})$ that, if $(\mk{m},\mk{n})\gamma=(\mk{m}',\mk{n}')$ for some $\gamma\in GL_2(K)$, then
\begin{align*}
\gamma GL(\mk{m}',\mk{n}')\gamma^{-1}=GL(\mk{m},\mk{n}). 
\end{align*}
In particular, note that for all $\lambda\in K^{\times}$, $GL(\mk{m},\mk{n})=GL(\lambda\mk{m},\lambda\mk{n})$.
\begin{Rem}
In general, if $(\mk{m},\mk{n}),(\mk{m}',\mk{n}')$ are two arbitrary pairs of lattices, there does not necessarily 
exist a $\gamma\in GL_2(K)$, such that $(\mk{m},\mk{n})\gamma=(\mk{m}',\mk{n}')$. However, in the special case where
$\mk{m},\mk{n},\mk{m}',\mk{n}'$ are all fractional $\ca{O}_K$-ideals, then one may show that such a $\gamma$ always exist. 
\end{Rem}

Having fixed an ordering of the embeddings of $K$ into $\RR$, we may view all the groups above as subgroups 
of $\ca{G}(\RR)\simeq GL_2(\RR)^g$.
\begin{Prop}\label{coff}
The group $GL(\mk{m},\mk{n})$, when viewed as a subgroup of $\ca{G}(\RR)$, is discrete.  
Moreover, if $\gamma=\M{a}{b}{c}{d}\in GL(\mk{m},\mk{n})$, then
\begin{enumerate}[(i)]
 \item $\det(\gamma)\in\ca{O}_K^{\times}$,
 \item $a\in\ca{O}_{\mk{m}}$, $d\in\ca{O}_{\mk{n}}$, $b\mk{m}\subseteq\mk{n}$ and $c\mk{n}\subseteq\mk{m}$.
\end{enumerate}
\end{Prop}
Recall here that $\ca{O}_{\mk{m}}$ is the ring of multipliers of $\mk{m}$ (see Section \ref{multiply}).

{\bf Proof} Let us first show that $GL(\mk{m},\mk{n})$ is discrete. Consider the embedding $\iota:K\rightarrow\RR^g$ induced 
by our choice of an ordering of the real embeddings of $K$ into $\RR$.
Note that every lattice $\ca{L}\subseteq K$ is such that $\iota(\ca{L})$ is discrete. 
Let $j:=\iota\times\iota:K^2\rightarrow(\RR^g)^2$. Then $\mk{L}:=j(\mk{m},\mk{n})$ is a discrete $\ZZ$-lattice inside $\RR^{2g}$ of rank $\RR^{2g}$.
Choosing a $\ZZ$-basis $\mk{B}$ of $\mk{L}$, one may construct the following commutative diagram:
\begin{align*}
\xymatrix{
        & GL_{2g}(\RR)                       &              \\
        &GL(\mk{m},\mk{n}) \ar[r]_{\wt{\iota}} \ar[u]_{\wt{\rho}} &\ca{G}(\RR) \ar[ul]_{\rho}
}
\end{align*}
where $\wt{\iota}$ is induced from $\iota$, $\rho$ and $\wt{\rho}$ are injective continuous group homomorphisms such that 
$\Imm(\wt{\rho})\subseteq GL_\ZZ(\mk{L})\simeq GL_{2g}(\ZZ)$.
Since $GL_{\ZZ}(\mk{L})$ is a discrete subgroup of $GL_{2g}(\RR)$ it follows that $\wt{\rho}(GL(\mk{m},\mk{n}))$ is discrete in $GL_{2g}(\RR)$, and, therefore, 
$\iota(GL(\mk{m},\mk{n}))$ is also discrete in $\ca{G}(\RR)$. 

Let us prove (i). Let $\gamma\in GL(\mk{m},\mk{n})$. Then the image
of the $\gamma$ under the $i$-th embedding, i.e., $\gamma^{(i)}\in GL_2(\RR)$ gives rise to two 
complex eigenvalues (counting multiplicities), $\lambda_{i1},\lambda_{i2}\in\CC$. It follows from the commutativity of the diagram, 
that the collection of eigenvalues $\{\lambda_{i1},\lambda_{i2}\}_{i=1}^g$ corresponds to the eigenvalues of $\wt{\rho}(\gamma)$.
But the matrix $\wt{\rho}(\gamma)$ is conjugate (inside $GL_{2g}(\RR)$) to a matrix in $GL_{2g}(\ZZ)$. In particular, all the eigenvalues of $\wt{\rho}(\gamma)$
are algebraic integers, such that their product is equal to $\pm 1$. Therefore, $\{\lambda_{i1},\lambda_{i2}\}_{i=1}^g\subseteq\ca{O}_K^{\times}$.  

Let us prove (ii). Consider the subgroup $(\mk{m},0)\leq (\mk{m},\mk{n})$. Since $(\mk{m},0)\gamma\subseteq (\mk{m},\mk{n})$, we deduce
that $a\mk{m}\subseteq\mk{m}$ and $b\mk{m}\subseteq\mk{n}$. In particular, $a\in\ca{O}_{\mk{m}}$. Similarly, we have  
$(0,\mk{n})\gamma\subseteq (0,\mk{n})$ and thus $d\mk{n}\subseteq\mk{d}$ and $c\mk{n}\subseteq\mk{m}$. \fin
\begin{Rem}
Note that in general, the set $GL(\mk{m},\mk{n})$ is not necessarily a subset of $GL_2(\ca{O}_K)$, even though the diagonal entries
of each of its matrix are elements of $\ca{O}_K$.
\end{Rem}

The next proposition shows that the discrete groups $GL(\mk{m},\mk{n})$ are big. 
\begin{Prop}\label{congru}
Let $\mk{m},\mk{n}\subseteq K$ be two lattices. Then there exists an integer $n\in\ZZ_{\geq 1}$,
such that $\Gamma_{\ca{O}_K}(n\ca{O}_K)\subseteq GL(\mk{m},\mk{n})$.
\end{Prop}

{\bf Proof} Let $\ca{O}=\ca{O}_{\mk{m}}\cap\ca{O}_{\mk{n}}$. In particular, $\mk{m}$ and $\mk{n}$ are $\ca{O}$-module. Define 
\begin{enumerate}
 \item $\mk{b}:=(\mk{m}^{-1}\mk{n})\cap\ca{O}$,
 \item $\mk{c}:=(\mk{n}^{-1}\mk{m})\cap\ca{O}$.
\end{enumerate}
Recall here that if $\ca{L}_1$ and $\ca{L}_2$ are lattices then $\ca{L}_1\ca{L}_2$ means the product of the two lattices
and $\mk{n}^{-1}$ (resp. $\mk{m}^{-1}$) corresponds to the multiplicative inverse of $\mk{n}$ (resp. multiplicative inverse of $\mk{m}$),
see Section \ref{mult_inv}. Let $\mk{d}=\mk{b}\mk{c}\subseteq\ca{O}$. Since $\mk{d}\subseteq\ca{O}_K$ is a lattice, there exists $n\in\ZZ_{\geq 1}$ such that
$[\ca{O}_K:\mk{d}]=n$. In particular $n\ca{O}_K\subseteq\mk{d}$ and $1+n\ca{O}_K\subseteq\ca{O}$. Finally, a direct computation shows that
for all $\gamma\in \Gamma_{\ca{O}_K}(n\ca{O}_K)$, $(\mk{m},\mk{n})\gamma\subseteq (\mk{m},\mk{n})$ and 
$(\mk{m},\mk{n})\gamma^{-1}\subseteq (\mk{m},\mk{n})$. Therefore, $(\mk{m},\mk{n})\gamma=(\mk{m},\mk{n})$, and hence
$\Gamma_{\ca{O}_K}(n\ca{O}_K)\subseteq GL(\mk{m},\mk{n})$. \fin

Recall that a subgroup $\Gamma\leq GL_2(\ca{O}_K)$ is called a {\it congruence subgroup}, if there exists $n\in\ZZ_{\geq 1}$
such that $\Gamma_{\ca{O}_K}(n\ca{O}_K)\leq \Gamma$. It follows from Proposition \ref{congru} that
all the discrete subgroups $GL(\mk{m},\mk{n})\cap GL_2(\ca{O}_K)$ are congruence subgroups.

Let $(G,\cdot)$ be a group and $G_1,G_2\leq G$ be subgroups. Recall that $G_1$ and $G_2$ are said to be {\it commensurable} if
$[G_1:G_1\cap G_2]<\infty$ and $[G_2: G_2\cap G_1]<\infty$. One may easily check that commensurability gives rise to 
an equivalence relation on the set of all subgroups of $G$. Consider now the natural projection map
\begin{align*}
\pi:GL_2(\ca{O}_K)\rightarrow GL_2(\ca{O}_K/n\ca{O}_K).
\end{align*}
Since $\Gamma_{\ca{O}_K}(n\ca{O}_K)=\ker(GL_2(\ca{O}_K)\rightarrow GL_2(\ca{O}_K/n\ca{O}_K))$ and $\# GL_2(\ca{O}_K/n\ca{O}_K)<\infty$, it 
follows that 
\begin{align}\label{ind_fin}
[GL_2(\ca{O}_K):\Gamma_{\ca{O}_K}(n\ca{O}_K)]<\infty.
\end{align}
We have the following useful proposition:
\begin{Prop}\label{annoy}
The group  $GL(\mk{m},\mk{n})$ is commensurable to $GL_2(\ca{O}_K)$.
\end{Prop}

{\bf Proof} In order to prove Proposition \ref{annoy}, one could simply invoke Theorem \ref{sel} in Appendix \ref{app_5},
and then conclude. However, the proof of Theorem \ref{sel} is considerably deeper that the content of Proposition \ref{annoy}.
For this reason, we also give a direct and elementary proof of  Proposition \ref{annoy}.
We already know from Proposition \ref{congru} that there exists $n\in\ZZ_{\geq 1}$ such that $\Gamma_{\ca{O}_K}(n\ca{O}_K)\leq GL(\mk{m},\mk{n})$,
and, therefore, from \eqref{ind_fin} $[GL_2(\ca{O}_K):\Gamma_{\ca{O}_K}(n\ca{O}_K)\cap GL(\mk{m},\mk{n})]<\infty$. 
Therefore, commensurability will follow if we can show that 
\begin{align}\label{tough}
[GL(\mk{m},\mk{n}):\Gamma_{\ca{O}_K}(n\ca{O}_K)]<\infty.
\end{align}
To simplify the notation, we let $\Gamma=GL(\mk{m},\mk{n})$ and $\Gamma'=\Gamma_{\ca{O}_K}(n\ca{O}_K)$.
We do a proof by contradiction. Let us assume that $[\Gamma:\Gamma']=\infty$. We will then derive a contradiction. Define 
\begin{align*}
H:=\{\lambda\cdot I_2\in GL_2(\ca{O}_K):(\lambda\mk{m},\lambda\mk{n})=(\mk{m},\mk{n})\}.
\end{align*}
In particular, $H\leq \Gamma$. Note that for each $h\in H$, there exists an $n\in\ZZ_{\geq 1}$ such that $h^n\in\Gamma'$. Since
$H$ is finitely generated, and that each element of $H$ commutes with the elements of $\Gamma'$, it
follows that $[\Gamma'H:\Gamma']<\infty$. Since by assumption $[\Gamma:\Gamma']<\infty$, we must therefore have 
\begin{align}\label{tough2}
[\Gamma:H\Gamma']=\infty.
\end{align}
Let $\pi$ denote the composition of the two maps $GL_2(K)\rightarrow\ca{G}(\RR)\rightarrow PSL_2(\RR)^g$. Note that 
$\ker(\pi)=\{\lambda\cdot I_2:\lambda\in K^{\times}\}$. 
Let $\ov{\Gamma}=\pi(\Gamma)$ and $\ov{\Gamma'}=\pi(\Gamma'H)=\pi(\Gamma')$. Recall that $PSL_2(\RR)^g$ acts faithfully and isometrically, by M\"obius transformations
on the symmetric space $\mk{h}^g$. Let $\mu$ be the measure on $\mk{h}^g$ induced by integrating the Gau\ss-Bonnet form on the measurable sets of $\mk{h}^g$. By definition,
$\mu$ is left invariant by $PSL_2(\RR)^g$. Since the group $\ov{\Gamma}$ is discrete, it acts faithfully and properly discontinuously on 
$\mk{h}^g$. Moreover, since $\ca{G}(\RR)$ is second-countable, the discrete group $\ov{\Gamma}$ admits a {\it fundamental set} in the sense of \cite{Sie6}, i.e., 
a set $\ca{F}\subseteq\mk{h}^g$ such that
\begin{enumerate}[(i)]
  \item $\bigcup_{\gamma\in\ov{\Gamma}}\;\gamma\ca{F}=\mk{h}^g$,
 \item For all $\gamma\in\ov{\Gamma}\bs\{1\}$, $\gamma\ca{F}\cap\ca{F}=\emptyset$,
 \item $\ca{F}$ is a Borel set. 
\end{enumerate}
Moreover, in \cite{Sie6}, it is shown that the $\mu$-volume of $\ca{F}$, i.e. $\mu(\ca{F})$, is independent of the choice of the fundamental set.  
Since $\ker(\pi)=\{\lambda\cdot I_2:\lambda\in K^{\times}\}$, the map $\pi$ induces a bijection between the following two sets of {\it right} cosets:
$\biglslant{H\Gamma'}{\Gamma}$ and $\biglslant{\ov{\Gamma'}}{\ov{\Gamma}}$ which is explicitly given by $H\Gamma'\cdot\gamma\mapsto \ov{\Gamma'}\cdot\ov{\gamma}$.
Let $\{g_i\in\ov{\Gamma}\}_{i\in I}$ be a complete set of representatives of the right cosets $\biglslant{\ov{\Gamma'}}{\ov{\Gamma}}$. 
Because of \eqref{tough2}, $\#I=\infty$ (we have used here the fact that there is a natural bijection between the left cosets and the 
right cosets which is provided by the inversion).
We may thus view $\ca{F}':=\{g_i\ca{F}:g_i\in I\}$ as a fundamental set of $\ov{\Gamma'}$. 
However, it is well-known that $PSL_2(\ca{O}_K)$ has a fundamental set with a finite strictly positive $\mu$-volume. Therefore, the same is true
for $\ov{\Gamma'}$. Since the measure $\mu$ is countably additive on disjoint Borel sets, we must have $\#I\cdot\mu(\ca{F})=\mu(\ca{F}')>0$. But this is absurd
since $\#I=\infty$. This concludes the proof. \fin

\begin{Def}\label{chan}
For each matrix $U\in M_2(K)$ and $\gamma\in GL_2(K)$, we define
\begin{align}\label{toit0}
U^\gamma:=\left( \gamma^{-1}\V{u_1}{u_2}, \gamma^t\V{v_1}{v_2}\right),
\end{align}
where $\gamma^{-1}$ corresponds to the inverse and $\gamma^t$ to the transpose.
One may check that $U\mapsto U^{\gamma}$ gives a right action of $GL_2(K)$ on $M_2(K)$ which we call 
the upper right action.
\end{Def}
The next proposition gathers some compatibility properties between the cartan involution on $M_2(K)$ (see \eqref{ulu})
and the upper right action defined above.
\begin{Prop}\label{compa}
For all $U\in M_2(K)$ and $\gamma\in SL_2(K)$ one has that
\begin{align}\label{commut}
(U^*)^{\gamma}=(U^{\gamma})^{*}.
\end{align}
For a diagonal matrix $D=\M{d_1}{0}{0}{d_2}\in GL_2(K)$, let us define $D^s:=\M{d_2}{0}{0}{d_1}$. Then 
for all $U\in M_2(K)$, and all diagonal matrices $D$ as above, we have
\begin{align}\label{commut2}
(U^*)^{D}=\left(U^{(D^{-1})^s}\right)^*.
\end{align}
In particular, if $D=\lambda I_2$ is a scalar matrix, then the $*$ action {\it inverse-commutes} with the upper right action
in the sense that $(U^*)^{D}=(U^{D^{-1}})^{*}$.
\end{Prop}

{\bf Proof} The proofs of \eqref{commut} and \eqref{commut2} are straightforward computations. \fin

\begin{Def}\label{unit_def}
Let $V$ be a lattice of $K$. For a pair of elements $a,b\in K$, we define 
\begin{align}\label{m_group}
\mathcal{V}_{a,b,V}:=&\{\epsilon\in\ca{O}_V^{\times}:(\epsilon-1)a\in V,
(\epsilon-1)b\in V^{*},(\epsilon-1)ab\in(\ca{O}_V)^*\}.
\end{align}
Recall here that $\ca{O}_V$ denotes the ring of multipliers of $V$, and $V^*$ denotes the dual lattice with respect to the trace pairing 
(see Section \ref{lat_nu}).
We also define $\mathcal{V}_{a,b,V}^+:=\mathcal{V}_{a,b,V}\cap\ca{O}_K^{\times}(\infty)$, where
$\ca{O}_K^{\times}(\infty)$ denotes the group of totally positive units of $\ca{O}_K$.
\end{Def}
\begin{Rem}\label{turf}
In \cite{Ch4}, a  different group than  $\mathcal{V}_{a,b,V}$ was used. It was 
denoted instead by $\Gamma_{a,b,V}$ and its definition was the same as in \eqref{m_group}, except that the last
condition in \eqref{m_group}, which reads as \lq\lq$ (\epsilon-1)ab\in(\ca{O}_V)^*$\rq\rq\s should be replaced
by \lq\lq$ (\epsilon-1)ab\in\mk{d}_K^{-1}$\rq\rq, where $\mk{d}_K$ is the different ideal of $K$. Note that since
$\mk{d}_K^{-1}\subseteq (\ca{O}_V)^*$, we always have that $\Gamma_{a,b,V}\leq \mathcal{V}_{a,b,V}$. The reason
we preferred to work with $\mathcal{V}_{a,b,V}$ rather than $\Gamma_{a,b,V}$ is because we want property (3) below
to hold true.
\end{Rem}

The next proposition  makes precise the exact dependence of $\ca{V}_{a,b,V}$ (resp. $\ca{V}_{a,b,V}^+$) on the triple $(a,b,V)$.
\begin{Prop}\label{dor}
We have
\begin{enumerate}
 \item $\mathcal{V}_{a,b,V}=\mathcal{V}_{-a,b,V}=\mathcal{V}_{-a,-b,V}=\mathcal{V}_{-b,a,V^*}$.
 \item For all $\lambda\in K\bs\{0\}$, $\mathcal{V}_{a,b,V}=\mathcal{V}_{\lambda a,\frac{b}{\lambda},\lambda V}$ (resp. 
 $\mathcal{V}_{a,b,V}^+=\mathcal{V}_{\lambda a,\frac{b}{\lambda},\lambda V}^+$).
 \item if $a\equiv a'\pmod{V}$ and $b\equiv b'\pmod{V^*}$ then $\mathcal{V}_{a,b,V}=\mathcal{V}_{a',b',V}$
 (resp. $\mathcal{V}_{a,b,V}^+=\mathcal{V}_{a',b',V}^+$).
 \item The subgroup $\mathcal{V}_{a,b,V}^+$ has finite index in $\ca{O}_K^{\times}$.
\end{enumerate}
\end{Prop}

{\bf Proof} The proofs of the first two equalities in (1) follow directly from the definition of $\mathcal{V}_{a,b,V}$. The proof of the
third equality in (1) follows from the fact that $V^{**}=V$ and $\ca{O}_{V}=\ca{O}_{V^*}$ (see Section \ref{du_op}). The proof
of (2) follows from the observations that $(\lambda V)^*=\frac{1}{\lambda}V^*$ and $\ca{O}_{V}=\ca{O}_{\lambda V}$. 
The proof of (3) follows from the observations that
$\ca{O}_{V}=\ca{O}_{V^*}$, $VV^*\subseteq\ca{O}^*$ and the following equivalence: for $\epsilon\in\ca{O}_V$, we have
\begin{align*}
(\epsilon-1)ab\in(\ca{O}_{V})^*\Longleftrightarrow  (\epsilon-1)a'b'\in(\ca{O}_{V})^*.
\end{align*}
It remains to prove (4). Let $\{\epsilon_i:i=1,\ldots,g-1\}\subseteq
\ca{O}_K^{\times}(\infty)$
be a $\ZZ$-basis of $\ca{O}_K^{\times}(\infty)$. From Dirichlet's unit theorem, this is a free $\ZZ$-module of rank $g-1$.
The algebraic number $a$ in the triple $(a,b,V)$ may be written as: $\frac{s}{t}$ for $s,t\in\ca{O}_K$ and $t\neq 0$, where  $s$
is chosen, so that $sV\subseteq V$. 
Choose $N\in\ZZ_{\geq 1}$, such that $N\ca{O}_K\subseteq V$. Let $m=\#(\ca{O}_K/tN\ca{O}_K)^{\times}$. 
Then for all $i\in\{1,\ldots,g-1\}$, we have
$$
\epsilon_i^m\in 1+tN\ca{O}_K\subseteq 1+tV.
$$ 
Therefore, $(\epsilon_i^{m}-1)a=(\epsilon_i^{m}-1)\frac{s}{t}\in sV\subseteq V$. Similarly, we may find an integer
$m'\geq 1$, such that for each $i\in\{1,\ldots,g-1\}$, $(\epsilon_i^{m'}-1)b\in V^*$. Let $m''=\lcm(m,m')$. Finally,
since $(\ca{O}_K^{\times}(\infty))^{m'}\leq\ca{V}_{a,b,V}$ it follows that $[\ca{O}_K^{\times}(\infty):\mathcal{V}_{a,b,V}]\leq (m'')^{g-1}$. \fin

\begin{Def}\label{fromage}
Let $\mk{m}$ and $\mk{n}$ be two lattices of $K$.
Let $U=\M{u_1}{v_1}{u_2}{v_2}\in M_2(K)$ be a parameter matrix. We define 
\begin{align*}
\mathcal{V}_U(\mk{m}\oplus\mk{n}):=\ca{V}_{v_1,u_1,\mk{m}}\cap\ca{V}_{v_2,u_2,\mk{n}},
\end{align*}
and $\mathcal{V}_U^+(\mk{m}\oplus\mk{n}):=\mathcal{V}_U(\mk{m}\oplus\mk{n})\cap\ca{O}_K^+(\infty)$.
\end{Def}
It follows from Proposition \ref{dor} 
that the groups $\mathcal{V}_U(\mk{m}\oplus\mk{n})$ and $\mathcal{V}_U^+(\mk{m}\oplus\mk{n})$  have finite index in $\ca{O}_K^{\times}$.

\begin{Prop}\label{Bren}
We have  
\begin{enumerate}[(i)]
 \item $\mathcal{V}_U(\mk{m}\oplus\mk{n})=\mathcal{V}_{U^*}(\mk{n}^*\oplus\mk{m}^*)$,
where $U\mapsto U^*$ corresponds to the Cartan involution which was defined in \eqref{ulu}. 
 \item Let $D=\M{d_1}{0}{0}{d_2}\in GL_2(\RR)$ be a diagonal matrix. Then 
 \begin{align*}
\mathcal{V}_U(\mk{m}\oplus\mk{n})=
\mathcal{V}_{U^D}(d_1\mk{m}\oplus d_2\mk{n}),
\end{align*}
\item Let $\gamma=\M{0}{-1}{1}{0}$. Then $\mathcal{V}_U(\mk{m}\oplus\mk{n})=\mathcal{V}_{U^{\gamma}}((\mk{m}\oplus\mk{n})\gamma)$.
\end{enumerate}
\end{Prop}

{\bf Proof} This follows directly from the definition of $\mathcal{V}_U(\mk{m}\oplus\mk{n})$ and Proposition \ref{dor}. \fin

\begin{Def}\label{grp}
Let $(\mk{m},\mk{n},U)$ be a triple, where $\mk{m},\mk{n}\subseteq K$ are lattices and $U\in M_2(K)$
is a parameter matrix. We define 
\begin{align*}
\Gamma_U(\mk{m},\mk{n}):=\left\{\gamma\in GL(\mk{m}\oplus\mk{n}):U^{\gamma}-U\in\M{\mk{m}^*}{\mk{m}}{\mk{n}^*}{\mk{n}}\right\}.
\end{align*}
We also let $\Gamma_U^+(\mk{m}\oplus\mk{n}):=\Gamma_U(\mk{m}\oplus\mk{n})\cap GL_2^+(K)$.

Let us assume that $U\in\frac{1}{N}\M{\mk{m}^*}{\mk{m}}{\mk{n}^*}{\mk{n}}$ for some $N\in\ZZ_{\geq 1}$. Then we define 
\begin{align*}
\Gamma_U(\mk{m},\mk{n};N):=\left\{\gamma\in\Gamma_U(\mk{m},\mk{n}):U^{\gamma}-U\in \M{N\mk{m}^*}{N\mk{m}}{N\mk{n}^*}{N\mk{n}}\right\}.
\end{align*}
\end{Def}
If $U\in\frac{1}{N}\M{\mk{m}^*}{\mk{m}}{\mk{n}^*}{\mk{n}}$, then, by definition, we have: $\Gamma_U(\mk{m},\mk{n};N)\leq \Gamma_U(\mk{m},\mk{n})$ and 
$\Gamma_U^+(\mk{m},\mk{n};N)\leq \Gamma_U^+(\mk{m},\mk{n})$.

\begin{Prop}\label{talon}
The groups $\Gamma_U(\mk{m},\mk{n}),\Gamma_U(\mk{m},\mk{n},N), \Gamma_U^+(\mk{m},\mk{n})$ and 
$\Gamma_U^+(\mk{m},\mk{n},N)$ are commensurable to $GL_2(\ca{O}_K)$. 
\end{Prop}

{\bf Proof} Let $\ca{O}=\ca{O}_{\mk{m}}\cap\ca{O}_{\mk{n}}$. Then $\ca{O}$ is an order of $\ca{O}_K$, such that
$\ca{O}\subseteq \ca{O}_{\mk{m}^*}\cap\ca{O}_{\mk{n}^*}$ (in fact it follows from the discussion in Section \ref{du_op} that the
previous inclusion is an equality). Since $\Gamma_U^+(\mk{m},\mk{n};N)$ is contained in all the other subgroups, 
it is enough to show that $\Gamma_U^+(\mk{m},\mk{n};N)\cap GL_2(\ca{O})$ contains a subgroup $\wt{\Gamma}$ (which will be constructed below) 
of finite index in $GL_2(\ca{O}_K)$. Indeed, since $[GL_2(\ca{O}_K):\wt{\Gamma}]<\infty$ and $GL_2(\ca{O}_K)$ is commensurable
to $GL(\mk{m},\mk{n})$ (Proposition \ref{annoy}), we will also have $[\Gamma_U^+(\mk{m},\mk{n};N):\wt{\Gamma}]<\infty$. 

From the short exact sequence 
\begin{align*}
1\rightarrow GL_2^+(\ca{O}_K)\rightarrow GL_2(\ca{O}_K)\stackrel{\det}{\rightarrow} \ca{O}_K^{\times}/\ca{O}_K^{\times}(\infty)\rightarrow 1,
\end{align*}
and the finiteness of $\ca{O}_K^{\times}/\ca{O}_K^{\times}(\infty)$ (since $(\ca{O}_K^{\times})^2\leq \ca{O}_K^{\times}(\infty)$)
it follows that $\wt{\Gamma}$ is of finite index in $GL_2^+(\ca{O}_K)$ if and only if it is of finite index in $GL_2(\ca{O}_K)$.

Define
\begin{enumerate}
 \item $\mk{b}:=\left((\mk{n}^*)^{-1}\mk{m}^*\right)\cap \left(\mk{n}\mk{m}^{-1}\right)\cap\ca{O}$,
 \item $\mk{c}:=\left((\mk{m}^*)^{-1}\mk{n}^*\right)\cap\left(\mk{n}^{-1}\mk{m}\right)\cap\ca{O}$.
\end{enumerate}
One may check that $\mk{b}$ and $\mk{c}$ are $\ca{O}$-modules. Let $\mk{d}:=\mk{b}\mk{c}$. Since 
$\mk{d}\subseteq\ca{O}_K$ is a lattice, there exists an integer $n\in\ZZ_{\geq 1}$ such that 
$[\ca{O}_K:\mk{d}]=n$. In particular, $n\ca{O}_K\leq\mk{d}$ and
$1+n\ca{O}_K\subseteq\ca{O}$. Define $\Gamma':=\Gamma_{\ca{O}_K}^+(n\ca{O}_K)$
Note that $\Gamma'\leq GL_2^+(\ca{O})$. 

By definition of $\Gamma'$, one may check that for all $\gamma\in\Gamma'$,
$(\mk{m}\oplus\mk{n})\gamma\subseteq\mk{m}\oplus\mk{n}$. Since $\gamma$ is also invertible as a linear map of $\ca{O}$-modules, 
it follows also that $(\mk{m}\oplus\mk{n})\gamma^{-1}=\mk{m}\oplus\mk{n}$, and thus 
$(\mk{m}\oplus\mk{n})\gamma=\mk{m}\oplus\mk{n}$.

Consider now the set 
\begin{align*}
W:=\M{\frac{1}{N}\mk{m}^*/(N\mk{m}^*)}{\frac{1}{N}\mk{m}/(N\mk{m})}{\frac{1}{N}\mk{n}^*/(N\mk{n}^*)}{\frac{1}{N}\mk{n}/(N\mk{n})}.
\end{align*}
By definition, $W$ is a finite abelian group (using the usual addition of matrices) of cardinality $N^{8g}$. 
It also follows from the definition of $\Gamma'$ that $W$ is $\Gamma'$-stable, as a set, under the {\it upper right action}. In this way, we
may view $W$ as a right $\Gamma'$-module. Consider now the class
\begin{align*}
 [U]:=U \modd \M{N\mk{m}^*}{N\mk{m}}{N\mk{n}^*}{N\mk{n}}\in W.
\end{align*}
Since $W$ has finite cardinality, it follows that $\wt{\Gamma}:=\Stab_{\Gamma'}([U])$ is
a finite index subgroup of $\Gamma'$ and, therefore, a finite index subgroup of $GL_2^+(\ca{O}_K)$. 
Finally, it follows from the definition of $\wt{\Gamma}$ that for all $\gamma\in\wt{\Gamma}$,
\begin{enumerate}
 \item $\gamma\in GL^+(\mk{m}\oplus\mk{n})\cap \Gamma'\subseteq GL^+(\mk{m}\oplus\mk{n})\cap GL_2(\ca{O})$, 
 \item $U^{\gamma}\equiv U\pmod{\M{N\mk{m}^*}{N\mk{m}}{N\mk{n}^*}{N\mk{n}}}$ (since $\gamma\in \Stab_{\Gamma'}([U])$).
\end{enumerate}
Hence, $\wt{\Gamma}\leq \Gamma_U(\mk{m},\mk{n};N)\cap GL_2(\ca{O})$. The result follows. \fin

\subsubsection{The indices $e_1,e_2,f_1$ and $f_2$}\label{involu}
In this section, we introduce various indices which measure the discrepancy between
various subgroups of $\ca{O}_K^{\times}(\infty)$.

\begin{Def}\label{chien}
Let $(\mk{m},\mk{n},U)$ be a triple and let $\gamma\in GL_2(K)$.
We define the following integers:
\begin{enumerate}
 \item $e_1((\mk{m},\mk{n});U):=[\ca{V}_{v_2,u_2,\mk{n}}^+:\mathcal{V}_U^+(\mk{m},\mk{n})]$,
 \item $e_2((\mk{m},\mk{n});U):=[\ca{V}_{v_1,u_1,\mk{m}}^+:\mathcal{V}_U^+(\mk{m},\mk{n})]$,
 \item $f_1((\mk{m},\mk{n});U;\gamma):=[\mathcal{V}_{U}^+(\mk{m}\oplus\mk{n}):\mathcal{V}_{U^{\gamma}}^+((\mk{m}\oplus\mk{n})\gamma)
 \cap \mathcal{V}_{U}^+(\mk{m}\oplus\mk{n})]$,
 \item $f_2((\mk{m},\mk{n});U;\gamma):=[\mathcal{V}_{U^{\gamma}}^+((\mk{m}\oplus\mk{n})\gamma):
 \mathcal{V}_{U^{\gamma}}^+((\mk{m}\oplus\mk{n})\gamma)\cap \mathcal{V}_{U}^+(\mk{m}\oplus\mk{n})]$.
\end{enumerate}
\end{Def}
When the data $((\mk{m},\mk{n}),U;\gamma)$ is clear from the context, we simply write $e_1$, $e_2$, $f_1$ and $f_2$. 
\begin{Def}
Let $\mk{m},\mk{n}\subseteq K$ be two lattices. We define $(\mk{m},\mk{n})^{s_*}:=(\mk{n}^*,\mk{m}^*)$. 
\end{Def}
Note that $s_{*}$ is an involution on the set of all ordered pairs of lattices of $K$. 
The next proposition describes some symmetries that the invariants $e_1,e_2,f_1$ and $f_2$ satisfy under the Cartan 
involution.
\begin{Prop}\label{todi}
We have:
\begin{enumerate}
\item  $e_1((\mk{m},\mk{n})^{s_*};U^*)=e_2((\mk{m},\mk{n});U)$. \\
   
For all $\gamma\in SL_2(\mk{m}\oplus\mk{n})$, we have 
\item $f_1((\mk{m},\mk{n})^{s_*};U^*;\gamma)=f_1((\mk{m},\mk{n});U;\gamma)$,
\item $f_2((\mk{m},\mk{n})^{s_*};U^*;\gamma)=f_2((\mk{m},\mk{n});U;\gamma)$.
\end{enumerate}
\end{Prop}

{\bf Proof} The proof of (1) follows from Proposition \ref{Bren}. The proofs of (2) and
(3) follows from Proposition \ref{Bren} and from \eqref{commut} of Proposition \ref{compa}. \fin

\subsubsection{The group $\Upsilon(\mk{m},\mk{n})$}

\begin{Def}\label{upsi}
Let $\mk{m},\mk{n}\subseteq K$ be two lattices. We define $\Upsilon(\mk{m},\mk{n})$ to be the subgroup
of $SL_2(K)$ generated by the following three types of elements:
\begin{enumerate}
 \item $\M{0}{-1}{1}{0}$, 
 \item $\M{\lambda}{0}{0}{\lambda^{-1}}$ such that $\lambda\in K^{\times}$,
 \item $\M{1}{\mu}{0}{1}$ such that $\mu\mk{m}\subseteq\mk{n}$.
\end{enumerate}
\end{Def}
The group $\Upsilon(\mk{m},\mk{n})$ will play a key role in the first proof of Theorem \ref{nice_fn_eq}.

\begin{Prop}\label{noot}
The group $\Upsilon(\mk{m},\mk{n})$ has the following three key properties:
\begin{enumerate}
 \item The group  $\Upsilon(\mk{m},\mk{n})$ acts transitively on $\PP^1(K)$,
 \item For all $\gamma\in \Upsilon(\mk{m},\mk{n})$, we have 
\begin{align*}
((\mk{m},\mk{n})\gamma)^{s_*}=((\mk{m},\mk{n})^{s_*})\gamma.
\end{align*}
\item Let $\gamma\in \Upsilon(\mk{m},\mk{n})$ and set $(\mk{m}',\mk{n}'):=(\mk{m},\mk{n})\gamma$. Then 
\begin{align*}
 \cov(\mk{m})\cov(\mk{n})=\cov(\mk{m}')\cov(\mk{n}').
\end{align*}
\end{enumerate}
\end{Prop}

{\bf Proof} Let us prove (1). Let $\sigma=\frac{a}{b}\in K$. Without loss of generality we may assume that
$b$ is such that $b\mk{m}\subseteq\mk{n}$. A direct computation shows that
\begin{align*}
\eta:=\M{-a}{0}{0}{-a^{-1}}\M{0}{-1}{0}{1}\M{1}{b}{0}{1}\M{0}{-1}{0}{1}=\M{a}{0}{b}{a^{-1}}.
\end{align*}
By definition, we have $\eta\in\Upsilon(\mk{m},\mk{n})$ and $\eta\infty=\sigma$. Since $\sigma$
was arbitrary, it follows that $\eta\in\Upsilon(\mk{m},\mk{n})$ acts transitively on $\PP^1(K)$.

The proofs of (2) and (3) follow from a direct calculation for each of the three types of generators
of $\Upsilon(\mk{m},\mk{n})$. \fin

\subsection{The real analytic Eisenstein series $G_{(\mk{m},\mk{n})}^{\alpha(s),\beta(s)}(U;z)$}\label{real_an}
In this section, we rewrite the Eisenstein series $G_{(\mk{m},\mk{n})}^w(U,p\, ;z,s)$
in terms of a family of bi-weights $[\alpha(s),\beta(s)]\in\CC^{g}\times\CC^g$.
\begin{Def}\label{difo}
For each $p\in\ZZ^g$, $w\in\ZZ$ and $s\in\CC$, we associate the following two $g$-tuples of complex numbers (which we also call weights):
\begin{enumerate} 
\item $\alpha=\alpha(s):=(s+w)\cdot\bu-\frac{p}{2}\in\CC^g$,
\item $\beta=\beta(s):=s\cdot\bu+\frac{p}{2}\in\CC^g$.
\end{enumerate}
\end{Def}
In particular, for $1\leq i\leq g$, we have $\alpha(s)_i=s+w-\frac{p_i}{2}$ and  
$\beta(s)_i=s+\frac{p_i}{2}$. When the dependence of $\alpha(s)$ and $\beta(s)$ on $s$ is clear, we may drop the
symbol $s$ from the notation and write instead $\alpha$ and $\beta$. If we want to speak about the $j$-th coordinate of 
$\alpha=\alpha(s)$, we may also write $\alpha_j$ rather than $\alpha(s)_j$.

\begin{Def}\label{difo2}
Let $\ca{Q}=((\mk{m},\mk{n}),U,p,w)$ be a standard quadruple.
For $s\in\CC$, we let $\alpha(s)$ and $\beta(s)$ be the two weights defined as in Definition \ref{difo}. In particular,
$\alpha(s)$ and $\beta(s)$ depend on the parameters $p$ and $k$. For 
$z\in K_\CC^{\pm}$ and $s\in\Pi_{1-\frac{w}{2}}$, we define
\begin{align}\label{tonu}
& G_{\ca{Q}}(z,s)=G_{(\mk{m},\mk{n})}^{\alpha(s),\beta(s)}(U\,;z):=\\ \notag
&\sum_{\mathcal{V}^+\bs\{(0,0)\neq(m+v_1,n+v_2)\in(\mk{m}+v_1,\mk{n}+v_2)\}}
\frac{e^{{2\pi\ii } \Tr(u_1(m+v_1)+u_2(n+v_2))}\cdot|y|^{\bu\cdot s}}{P(\alpha(s),\beta(s);(m+v_1) z+(n+v_2))}.
\end{align}
\end{Def}
We recall here that $P(\alpha,\beta;z)=z^{\alpha}\cdot(\ov{z})^{\beta}$ for $z\in K_{\CC}^{\times}$
where $P(\alpha,\beta;z)$ is computed according to Convention \ref{conv}. The summation in \eqref{tonu} is understood to be taken over 
a complete set of representatives of $(\mk{m}+v_1,\mk{n}+v_2)\bs\{(0,0)\}$ under the left diagonal action of 
$\mathcal{V}^+:=\mathcal{V}_U^+(\mk{m},\mk{n})$. Since $s\in\Pi_{1-\frac{w}{2}}$, the summation on the right-hand side of \eqref{tonu} converges
absolutely (this fact will be proved in Theorem \ref{gr_est}) and, therefore, 
the definition of $G_{(\mk{m},\mk{n})}^{\alpha(s),\beta(s)}(U\,;z)$ makes sense. 

Using (1) of Proposition \ref{koy_prop} and the identity
\begin{align*}
\frac{\omega_{p-w\cdot\bu}(z)}{|\Norm(z)|^{2s+w}}=\frac{\omega_{p}(z)}{\Norm(z)^w|\Norm(z)|^{2s}},
\end{align*}
one may deduce that
\begin{align}\label{oign}
G_{(\mk{m},\mk{n})}^{\alpha(s),\beta(s)}(U\,;z)=G_{(\mk{m},\mk{n})}^{w}(U,p\,;z,s),
\end{align}
where $G_{(\mk{m},\mk{n})}^{w}(U,p\,;z,s)$ is the function which appears in \eqref{mou1}.

%
%

\subsubsection{A transformation formula for $G_{(\mk{m},\mk{n})}^{\alpha(s),\beta(s)}(U;z)$}\label{real_an1}

We may now state a general transformation formula for $G_{(\mk{m},\mk{n})}^{\alpha(s),\beta(s)}(U\,;z)$.
\begin{Prop}\label{pleut}
Let $\ca{Q}=((\mk{m},\mk{n}),U,p,w)$ be a standard quadruple and let 
$G_{(\mk{m},\mk{n}}^{\alpha(s),\beta(s)}(U;z)$ be its associated Eisenstein series.
For all $z\in K_\CC^{\times}=K_{\CC}^{\pm}\cup(\RR^{\times})^g$, $s\in\CC$ and $\gamma\in GL_2(K)$, we have:
\begin{align}\label{chat}
j(\gamma,z)^{-\alpha(s)}\cdot j(\gamma,\ov{z})^{-\beta(s)}\cdot |j(\gamma,z)|^{\bu\cdot 2s}\cdot|\det(\gamma)|^{-\bu\cdot s}
\cdot G_{(\mk{m},\mk{n})}^{\alpha(s),\beta(s)}(U,\gamma z)=f_{\gamma}\cdot G_{(\mk{m},\mk{n})\gamma}^{\alpha(s),\beta(s)}(U^\gamma,z),
\end{align}
where 
\begin{align}\label{indexi}
f_{\gamma}:=\frac{f_1((\mk{m},\mk{n}),U;\gamma)}{f_2((\mk{m},\mk{n}),U;\gamma)}\in\QQ_{>0}.
\end{align}
\end{Prop}
It is understood here that the product $j(\gamma,z)^{-\alpha(s)}\cdot j(\gamma,\ov{z})^{-\beta(s)}$ 
is computed according to Convention \ref{conv}. Note that from the identity (2) of Proposition \ref{koy_prop}, we have
\begin{align*}
j(\gamma,z)^{-\alpha(s)}\cdot j(\gamma,\ov{z})^{-\beta(s)}&=P(-\alpha(s),-\beta(s);j(\gamma,z))
=\frac{\omega_{p-w\cdot\bu}(j(\gamma,z))}{|j(\gamma,z)|^{\bu\cdot (2s+w)}}.
\end{align*}
In particular, the product
\begin{align}\label{ouaq0}
j(\gamma,z)^{-\alpha(s)}\cdot j(\gamma,\ov{z})^{-\beta(s)}\cdot |j(\gamma,z)|^{\bu\cdot 2s}&=\omega_{w\cdot\bu-p}(j(\gamma,z))^{-1}\cdot|j(\gamma,z)|^{-w\cdot\bu},
\end{align}
is {\it independent} of $s$.
\begin{Rem}\label{ouaq}
It follows from \eqref{ouaq0} that the renormalized function 
\begin{align}\label{pouf}
[z\mapsto G_{(\mk{m},\mk{n})}^{\alpha(s),\beta(s)}(U,z)\cdot |y|^{\frac{w}{2}\cdot\bu}],
\end{align}
has unitary weight $\{w\cdot\bu-p;s\}$ in the sense of Section \ref{unit_weight}. In virtue of the identities 
\eqref{tram} and \eqref{oign}, one readily sees that the function in \eqref{pouf} is in fact equal to $G_{(\mk{m},\mk{n})}^{0}(U,p-w\cdot\bu;z,s+\frac{w}{2})$.
\end{Rem}

{\bf Proof of Proposition \ref{pleut}} We need to show that
\begin{align}\label{plast}
&E:=j(\gamma,z)^{-\alpha(s)}\cdot j(\gamma,\ov{z})^{-\beta(s)}\cdot |j(\gamma,z)|^{\bu\cdot 2s}\cdot|\det(\gamma)|^{-\bu\cdot s}\cdot\\[2mm] \notag
 &\hspace{2cm}\sum_{\ca{V}_{U}^+(\mk{m},\mk{n})\bs\{(0,0)\neq(m+v_1,n+v_2)\in(\mk{m}+v_1,\mk{n}+v_2)\}}
\frac{e^{{2\pi\ii } \Tr(u_1(m+v_1)+u_2(n+v_2))}\cdot|\Imm(\gamma z)|^{\bu\cdot s}}{P(\alpha(s),\beta(s);(m+v_1) \gamma z+(n+v_2))},
\end{align}
equals to the right-hand side of \eqref{chat}. First, observe that the following two identities hold true
\begin{enumerate}[(a)]
 \item For all $z\in K_{\CC}^{\pm}=K_{\CC}^{\times}\bs(\RR^{\times})^g$, and $\gamma\in \ca{G}(\RR)$, we have
\begin{align}\label{plast2}
|j(\gamma,z)|^{\bu\cdot 2s}\cdot |\det(\gamma)|^{-\bu\cdot s}|\Imm(\gamma z)|^{\bu\cdot s}=|\Imm(z)|^{\bu\cdot s}.
\end{align}
 \item  For all $\gamma=\M{a}{b}{c}{d}\in\ca{G}(\RR)$, $m,n\in K$ and $z\in K_{\CC}^{\times}$, we have:
\begin{align}\label{off}
\hspace{-1cm}j(\gamma,z)^{-\alpha(s)}j(\gamma,z)^{-\beta(s)}
P\big(-\alpha(s),-\beta(s);m(\gamma z)+n\big)=P\big(-\alpha(s),-\beta(s);m\cdot 
j(\gamma,z) (\gamma z)+n\cdot j(\gamma,z)\big).
\end{align} 
\end{enumerate}
Note that the identity \eqref{off} follows from (2) of Proposition \ref{koy_prop}. Substituting (a) and (b) in
the right-hand side of \eqref{plast}, we obtain:
\begin{align}\label{plast3}
E=\sum_{\mathcal{V}_U^+(\mk{m},\mk{n})\bs\{(0,0)\neq(m+v_1,n+v_2)\in(\mk{m}+v_1,\mk{n}+v_2)\}}
\frac{e^{{2\pi\ii } \Tr(u_1(m+v_1)+u_2(n+v_2))}\cdot|\Imm(\gamma z)|^{\bu\cdot s}}{P(\alpha(s),\beta(s);(m+v_1)j(\gamma,z) (\gamma z)+(n+v_2)j(\gamma,z))}.
\end{align}
We wish now to rewrite the expression
\begin{align*}
P(\alpha(s),\beta(s);(m+v_1)j(\gamma,z) (\gamma z)+(n+v_2)j(\gamma,z))
\end{align*}
differently. For any
$\gamma=\M{a}{b}{c}{d}\in GL_2(K)$ and $(m+v_1,n+v_2)\in(\mk{m}+v_1,\mk{n}+v_2)$, we have the two identities
\begin{align*}
(m+v_1)(cz+d)(\gamma z)+(n+v_2)(cz+d)&=\left(\R{m+v_1}{n+v_2}\right)\cdot\left(\gamma\V{z}{1}\right)\\
&=\left(\R{m+v_1}{n+v_2}\gamma\right)\cdot\left(\V{z}{1}\right),
\end{align*}
and
\begin{align*}
u_1(m+v_1)+u_2(m+v_2)=\left(\R{u_1}{u_2}(\gamma^{t})^{-1}\right)\cdot\left(\gamma^t\V{m+v_1}{n+v_2}\right).
\end{align*}
Therefore, if we consider the column vectors $\V{u_1}{u_2}$ and $\V{v_1}{v_2}$ and the row vector $(m,n)\in\mk{m}\oplus\mk{n}$, 
then the matrix $\gamma$ acts on them as:
\begin{enumerate}
 \item $\V{u_1}{u_2}\mapsto \gamma^{-1}\V{u_1}{u_2}$
 \item $\V{v_1}{v_2}\mapsto \gamma^t\V{v_1}{v_2}$.
 \item $(m,n)\mapsto (m,n)\gamma$.
 \end{enumerate}
Finally, combining all the previous observations, we deduce that
\begin{align*}
E=f_{\gamma}\cdot G_{(\mk{m},\mk{n})\gamma}^{\alpha(s),\beta(s)}(U^\gamma,z).
\end{align*}
Note that the presence of the factor $f_{\gamma}$ above (see \eqref{indexi}) takes care of the difference between the two unit groups $\ca{V}_U^+(\mk{m},\mk{n})$
and $\ca{V}_{U^{\gamma}}^+((\mk{m},\mk{n})\gamma)$. This concludes the proof. \fin

Let us apply Proposition \ref{pleut} to the special case where the matrix $\gamma=\M{d_1}{0}{0}{d_2}$ is a diagonal.
In this special case, it follows that from (2) of Proposition \ref{Bren} that $f_{D}=1$. We have thus obtained the following
corollary:
\begin{Cor}\label{chol}
For any diagonal matrix $D=\M{d_1}{0}{0}{d_2}\in GL_2(\RR)$, we have
\begin{align}\label{chat2}
G_{(\mk{m},\mk{n})}^{\alpha(s),\beta(s)}\left(U;\frac{d_1}{d_2}z\right)= |\Norm(d_1d_2)|^{s}\cdot\omega_{w\cdot\bu-p}(d_2)\cdot|\Norm(d_2)|^{w}\cdot 
G_{(d_1\mk{m},d_2\mk{n})}^{\alpha(s),\beta(s)}(U^{D};z).
\end{align}
\end{Cor}

\subsection{Real analytic modular forms of bi-weight $\{[\alpha,\beta];\mu\}$}\label{real_an2}
In this subsection, we introduce the notion of modular forms of bi-weight $\{[\alpha,\beta];\mu\}$
where $\alpha,\beta,\mu\in\CC^g$. This notion generalizes the notion of modular form of unitary weight
$\{p;s\}$ that was introduced in Section \ref{unit_weight}. Here one should view the weight $\alpha$ (resp. $\beta$) 
as the holomorphic weight (resp. as the anti-holomorphic weight). 

We let $\mathcal{M}_{\mbox{\tiny{cont}}}:=\Maps_{\mbox{\tiny{cont}}}(K_{\CC}^{\pm},\CC)$ be the space of continuous $\CC$-valued
functions on $K_\CC^{\pm}$. Similarly, we let $\mathcal{M}_{\mbox{\tiny{an}}}$ be the subspace of functions of 
$\mathcal{M}_{\mbox{\tiny{cont}}}$ which are real analytic.
\begin{Def}\label{difo3}
Let $\alpha,\beta,\mu\in \CC^g$ and assume that $\alpha-\beta\in\ZZ^g$. Let $F\in\mathcal{M}_{\mbox{\tiny{cont}}}$, $\gamma\in \ca{G}(\RR)$,
$z\in K_\CC^{\pm}$. We define
\begin{align*}
F\big|_{\{[\alpha,\beta];\mu\},\gamma}(z):=j(\gamma,z)^{-\alpha}\cdot j(\gamma,\ov{z})^{-\beta}\cdot|\det(\gamma)|^{-\mu}\cdot F(\gamma z).
\end{align*}
Here, the product $j(\gamma,z)^{-\alpha}\cdot j(\gamma,\ov{z})^{-\beta}$ is computed, according to Convention \ref{conv}.

Let $\Gamma\leq \ca{G}(\RR)$ be a subgroup and let $F\in\mathcal{M}_{\mbox{\tiny{cont}}}$. We say that 
\begin{enumerate}
\item $F(z)$ is {\it almost of bi-weight $\{[\alpha,\beta];\mu\}$ relative to $\Gamma$,} if, for all $\gamma\in\Gamma$, one has that
\begin{align*}
F\big|_{\{[\alpha,\beta];\mu\},\gamma}(z)=\zeta_{\gamma,z}\cdot F(z),
\end{align*}
for some root of unity $\zeta_{\gamma,z}$ which depends only on $\gamma$ and on the connected component of $z$.
\item $F(z)$ is of {\it bi-weight $\{[\alpha,\beta];\mu\}$ relative to $\Gamma$,} if, for all $\gamma\in\Gamma$, one has that
\begin{align*}
F\big|_{\{[\alpha,\beta];\mu\},\gamma}(z)=F(z).
\end{align*}
\end{enumerate}
A real analytic modular form of almost bi-weight $\{[\alpha,\beta];\mu\}$ (of bi-weight $\{[\alpha,\beta];\mu\}$), relative to $\Gamma$, is a 
function $f\in\mathcal{M}_{\mbox{\tiny{an}}}$ which satisfies (1) (resp. which satisfies (2)). When $\Gamma\leq \ca{G}_1(\RR)$, we simply 
speak of modular forms of almost bi-weight $\{[\alpha,\beta]\}$ (resp. of bi-weight $\{[\alpha,\beta]\}$).
\end{Def}
One may check that if $\alpha,\beta\in\CC^g$ are such that $\alpha-\beta\in\ZZ^g$, 
then, for all $\gamma_1,\gamma_2\in GL_2(K)$ and $z\in K_\CC^{\times}$,
one has the identity:
\begin{align}\label{plut}
 j(\gamma_1\gamma_2,z)^{\alpha}j(\gamma_1\gamma_2,\ov{z})^{\beta}=j(\gamma_1,\gamma_2 z)^{\alpha}j(\gamma_2,z)^{\alpha}
 j(\gamma_1,\gamma_2 \ov{z})^{\beta}j(\gamma_2,\ov{z})^{\beta}.
\end{align}
Here the product on the right-hand side is computed according to Convention \ref{conv}.
This identity follows directly from (2) of Proposition \ref{koy_prop}. 
We would like to emphasize that such an identity does not necessarily hold true for arbitrary $\alpha,\beta\in\CC^g$.
It follows from \eqref{plut}, that
$f\mapsto f\big|_{\{[\alpha,\beta];\mu\},\gamma}$ gives rise to a right action of $GL_2(K)$ on $\mathcal{M}_{\mbox{\tiny{cont}}}$.
\begin{Rem}
If $p\in\ZZ^g$, one may check that a function $F\in \mathcal{M}_{\mbox{\tiny{cont}}}$ 
is of bi-weight $\{[p/2,-p/2];s\cdot\bu\}$ relative to $\Gamma$, if and 
only if it is of unitary weight $\{p;s\}$ relative to $\Gamma$. 
Therefore, the notion of a modular form of bi-weight $\{[\alpha,\beta];\mu\}$
encompasses the notion of modular forms of unitary weight $\{p;s\}$. 
\end{Rem}
\begin{Rem}
Let $F\in \mathcal{M}_{\mbox{\tiny{cont}}}$ (resp. $G\in \mathcal{M}_{\mbox{\tiny{cont}}}$) 
be a modular form of bi-weight $\{[\alpha,\beta];\mu\}$ (resp. of bi-weight $\{[\alpha',\beta'];\mu'\}$).
Then $F\cdot G$ is a modular form of bi-weight $\{[\alpha+\alpha',\beta+\beta'];\mu+\mu'\}$.
\end{Rem}

\begin{Exa}
Let $z\in\CC\bs\RR$ and $F(z):=|\Imm(z)|$. Then, the real analytic function $F(z)$ is of bi-weight $\{[-1,-1];1\}$ relative to the group $GL_2(\RR)$.
More generally, let $z\in K_{\CC}^{\pm}$, $s\in\CC$ and $a\in\ZZ^g$. Then the function 
$z\mapsto F_a(z,s):=|\Imm(z)|^{as}$ is a real analytic modular form of bi-weight $\{[-as,-as];\Tr(a)s\}$ relative to the continuous group
$\ca{G}(\RR)$.
\end{Exa}
\begin{Exa}
Let $\alpha(s)$ and $\beta(s)$ be defined as in Definition \ref{difo}. Then, the function
\begin{align*}
z\mapsto G_{(\mk{m},\mk{n})}^{\alpha(s),\beta(s)}(U;z),
\end{align*} 
is a real analytic modular forms of bi-weight $\{w\cdot\bu-\frac{p}{2},\frac{p}{2}];s\cdot\bu\}$ (cf. with Remark \ref{ouaq}).
\end{Exa}

\begin{Rem}
Let $\Gamma\leq \ca{G}_1(\RR)$ be a congruence subgroup. Explicit examples of \lq\lq almost\rq\rq\; 
holomorphic Eisenstein series (only the constant term of their Fourier series expansion may fail to be holomorphic)
of bi-weight $[\alpha,\beta]$, where $\alpha,\beta\in\ZZ^g$ and $\alpha+\beta=r\cdot\bu$ with $r\in\{0,2\}$,
appear naturally in the setting of Eisenstein cohomology of Hilbert modular varieties, see Chapter III.3 of \cite{Frei1}.
\end{Rem}

\subsubsection{Classical weight constraints on holomorphic modular forms of bi-weight $[\alpha,\beta]$}\label{constr}
Let $\alpha,\beta\in\CC^g$ be such that $\alpha-\beta\in\ZZ^g$. Let
$F$ be a real analytic modular form of almost bi-weight $[\alpha,\beta]$ relative
to a discrete subgroup $\Gamma\leq SL_2(K)$ commensurable to $SL_2(\ca{O}_K)$. Since
$x\mapsto F(x+\ii y)$ is translation invariant under a certain lattice $\mathcal{L}$ of $K$ and real analytic in $z$,
it admits a Fourier series expansion (see Section \ref{Fourier}) of the form 
\begin{align}\label{chiv0}
F(z)=a_0(y)+\sum_{\xi\in\mathcal{L}^*} a_{\xi}(y)e^{\Tr(\xi x)}.
\end{align}
The subgroup $\Gamma$ contains a subgroup $\left\{\M{\epsilon}{0}{0}{\epsilon^{-1}}:\epsilon\in\mathcal{V}\right\}$ where
$\mathcal{V}$ is a finite index subgroup of $\ca{O}_K^{\times}(\infty)$. Since $F(z)$ is modular of bi-weight $[\alpha,\beta]$,
the Fourier coefficients of \eqref{chiv0} must satisfy the relation 
\begin{align}\label{nazy}
\epsilon^{\alpha+\beta}a_{\xi\epsilon^{-2}}(\epsilon^2 y)=a_{\xi}(y), 
\end{align}
for all $\epsilon\in\mathcal{V}$ and all $\xi\in\mathcal{L}^*$. Under additional assumptions on the function $F(z)$ and the bi-weight
$[\alpha,\beta]$, it is sometimes possible to rule out the existence of such modular forms. Let us give two such examples which
are taken from \cite{Frei1}.
\begin{enumerate}
 \item Assume that $F(z)$ is holomorphic. In particular, $\beta=\bz$ and $\alpha\in\ZZ^g$. Since
$F(z)$ is holomorphic, the constant term $a_0(y)$ in \eqref{chiv0} is independent of $y$. If $F(z)$ is not cuspidal at
 $\infty$, i.e., if $a_0\neq 0$, then, looking at \eqref{nazy} for $\xi=0$, we readily see that one must have $\alpha\in\ZZ\cdot\bu$, i.e.,
$F(z)$ must have parallel holomorphic weight. Furthermore, if $\alpha\neq\bz$, we may 
restrict $F(z)$ to the diagonal $\Delta\subseteq\mk{h}^g$. Note that $f(z_1):=F(z_1,\ldots,z_1)$ is again a non-zero
modular form in one variable since the image of $\Delta$ under $\Gamma$ is dense. 
Finally, using the well-known fact that there are no holomorphic modular forms on $\mk{h}$ of strictly negative weight 
(for any congruence subgroup of $SL_2(\ZZ)$), 
we deduce that $\alpha>\bz$.  
\item Set $I=\{1,2,\ldots,g\}$. Assume that $F(z)$ is of bi-weight $[\alpha,\beta]$ with $\alpha,\beta\in\ZZ^g$.
Assume, furthermore, that the following two conditions are satisfied:
\begin{enumerate}
\item There exists subsets $A,B\subseteq I$, such that $A\cap B=\emptyset$, $A\cup B=I$;
$F(z)$ is holomorphic in the $z_i$'s for $i\in A$, and $F(z)$ is anti-holomorphic in the $z_i$'s for $i\in B$.
\item There exists two indices $i,j\in\{1,2\ldots,g\}$, such that $\alpha_i+\beta_i\leq 0$ and $\alpha_j+\beta_j>0$.
\end{enumerate}
Then $F$ is identically equal to zero. Note that from (a) and (b), the bi-weight $[\alpha,\beta]$ necessarily satisfies the
conditions $\alpha_j=0$, if $j\in B$ and $\beta_j=0$, if $j\in A$.
\end{enumerate}
The proof of (2) follows from Proposition 4.11 on p. 53 of \cite{Frei1} and the discussion on p. 139-140 of \cite{Frei1}.

\begin{Rem}
In Section \ref{unicity}, we will show that certain non-trivial families of real analytic Eisenstein series $\{E(z,s)\}_{s\in\Pi_1}$, of
a constant bi-weight $[\frac{p}{2},-\frac{p}{2}]$, don't exist.  
\end{Rem}

\subsection{Symmetries and modularity of $G_{(\mk{m},\mk{n})}^{\alpha(s),\beta(s)}(U\,;z)$}\label{symm_a}

Let $\ca{Q}=((\mk{m},\mk{n}),U,p,w)$ be a standard quadruple and let $[\alpha(s),\beta(s)]$ be its associated
bi-weight, where
\begin{enumerate}
 \item $\alpha(s)=(s+w)\cdot\bu-\frac{p}{2}$, 
 \item $\beta(s)=s\cdot\bu+\frac{p}{2}$.
\end{enumerate}
Recall that $G_{\ca{Q}}(z,s)=G_{(\mk{m},\mk{n})}^{\alpha(s),\beta(s)}(U;z)=G_{(\mk{m},\mk{n})}^{w}(U,p\,;z,s)$ is the Eisenstein series which appears in Definition \ref{difo2}.

\subsubsection{Symmetries induced by sign changes of the entries of $U$}
Let us start by pointing out some obvious symmetries of $G_{(\mk{m},\mk{n})}^{\alpha(s),\beta(s)}(U;z)$, when the
parameter matrix $U$ is subjected to involutions induced from sign changes of its entries. Let $\iota_{11},\iota_{12},\iota_{21}$ and
$\iota_{22}$ be the following (set) 
involutions of $M_2(K)$:
\begin{enumerate}
 \item[(j)] $\iota_{11}\left(\M{u_1}{v_1}{u_2}{v_2}\right)=\M{-u_1}{v_1}{u_2}{v_2}$,
 \item[(k)] $\iota_{12}\left(\M{u_1}{v_1}{u_2}{v_2}\right)=\M{u_1}{-v_1}{u_2}{v_2}$,
 \item[(l)] $\iota_{21}\left(\M{u_1}{v_1}{u_2}{v_2}\right)=\M{u_1}{v_1}{-u_2}{v_2}$,
 \item[(m)] $\iota_{22}\left(\M{u_1}{v_1}{u_2}{-v_2}\right)=\M{u_1}{v_1}{u_2}{-v_2}$.
\end{enumerate}
A direct computation reveals that:
\begin{enumerate}
 \item[(n)]  $G_{(\mk{m},\mk{n})}^{\alpha(s),\beta(s)}(\iota_{11}(U);z)=G_{(\mk{m},\mk{n})}^{\alpha(s),\beta(s)}(\iota_{12}(U);-z)$,
 \item[(o)]  $G_{(\mk{m},\mk{n})}^{\alpha(s),\beta(s)}(\iota_{21}(U);z)=(-1)^{\Tr(p-w\cdot\bu)}
 G_{(\mk{m},\mk{n})}^{\alpha(s),\beta(s)}(\iota_{22}(U);-z)$.
 \end{enumerate}
Since complex conjugation is a non-trivial continuous automorphism of $\CC$, it also induces some
symmetries on the real analytic Eisenstein series. Let us give two such examples:
\begin{enumerate}
 \item[(p)] If $w=0$, we have $G_{(\mk{m},\mk{n})}^{0}(U,p\,;\ov{z},s)=G_{(\mk{m},\mk{n})}^{0}(U,-p\,;z,s)$.
 \item[(q)] If $w=0$ and $U=\M{0}{*}{0}{*}$, then 
 \begin{align*}
\ov{G_{(\mk{m},\mk{n})}^{0}(U,p\,;z,s)}=G_{(\mk{m},\mk{n})}^{0}(U,-p\,;z,s).
\end{align*}
\end{enumerate}
The reader will have no difficulty in deriving some more examples. 

\subsubsection{Dependence of $G_{(\mk{m},\mk{n})}^{\alpha(s),\beta(s)}(U\,;z)$ with respect to $U$}

We would like now to point out some symmetries of the expression $G_{(\mk{m},\mk{n})}^{\alpha(s),\beta(s)}(U\,;z)$ which 
make precise its dependence on the entries of the parameter matrix $U$. 
In particular, these symmetries will imply that $[z\mapsto G_{(\mk{m},\mk{n})}^{\alpha(s),\beta(s)}(U\,;z)]$
is almost modular of bi-weight $\{[w\cdot\bu-\frac{p}{2},\frac{p}{2}];s\cdot\bu\}$ relative to the
group $\Gamma_U(\mk{m},\mk{n})$ (see Proposition \ref{mod_for} below).

Let $U=\M{u_1}{v_1}{u_2}{v_2},U'=\M{u_1'}{v_1'}{u_2'}{v_2'}\in M_2(K)$
be two parameter matrices, such that
\begin{align*}
\M{u_1}{v_1}{u_2}{v_1}\equiv \M{u_1'}{v_1'}{u_2'}{v_2'} \pmod{\M{\mk{m}^*}{\mk{m}}{\mk{n}^*}{\mk{n}}}.
\end{align*}
It follows from (3) of Proposition \ref{dor} that 
$\mathcal{V}_U^+(\mk{m}\oplus\mk{n})=\mathcal{V}_{U'}^+(\mk{m}\oplus\mk{n})$. Moreover, a direct computation involving the right-hand
side of \eqref{tonu} shows that:
\begin{align}\label{tiff}
\zeta_{U,U'}\cdot G_{(\mk{m},\mk{n})}^{\alpha(s),\beta(s)}(U;z)=G_{(\mk{m},\mk{n})}^{\alpha(s),\beta(s)}(U';z),
\end{align}
where $U'':=U-U'=\M{u_1''}{v_1''}{u_2''}{v_2''}\in\M{\mk{m}^*}{\mk{m}}{\mk{n}^*}{\mk{n}}$, and
\begin{align}\label{tiff2}
\zeta_{U,U'}=e^{-2\pi\ii\Tr(u_1v_1''+u_1''v_1+u_2v_2''+u_2''v_2)}. 
\end{align}
Therefore, the value $G_{(\mk{m},\mk{n})}^{\alpha(s),\beta(s)}(U;z)$, up to a root of unity, {\it depends only} on $U$ 
modulo $\M{\mk{m}^*}{\mk{m}}{\mk{n}^*}{\mk{n}}$.

We may now state the following important proposition:
\begin{Prop}\label{mod_for}
Let $\ca{Q}=((\mk{m},\mk{n}),U,p,w)$ be a standard quadruple and let $[\alpha(s),\beta(s)]$ be its associated
bi-weight. Set $[\alpha_0,\beta_0]:=[w\cdot\bu-\frac{p}{2}\;,\;\frac{p}{2}]\in(\frac{1}{2}\ZZ)^g\times(\frac{1}{2}\ZZ)^g$.
Let $z\in K_\CC^{\pm}$, $s\in\CC$ and assume that $U\in\frac{1}{N}\M{\mk{m}^*}{\mk{m}}{\mk{n}^*}{\mk{n}}$ for some $N\in\ZZ_{\geq 1}$.
Recall that $\Gamma_U(\mk{m},\mk{n})$ and $\Gamma_U(\mk{m},\mk{n};N)$ are the discrete subgroups of $GL_2(K)$ which appear in Definition \ref{grp}.
Then, for all $\gamma\in\Gamma_U(\mk{m},\mk{n})$, one has that
\begin{align}\label{musc}
G_{(\mk{m},\mk{n})}^{\alpha(s),\beta(s)}\big|_{\{[\alpha_0,\beta_0],s\cdot\bu\},\gamma}(U\,;z)=
\zeta_{U,U^{\gamma}}\cdot G_{(\mk{m},\mk{n})}^{\alpha(s),\beta(s)}(U\,;z),
\end{align}
where $\zeta_{U,U^{\gamma}}$ is a root of unity which is defined as in \eqref{tiff2}. In particular, it follows from \eqref{musc}, that 
the function $z\mapsto G_{(\mk{m},\mk{n})}^{\alpha(s),\beta(s)}(U;z)$ 
\begin{enumerate}
 \item  is of almost bi-weight $\{[\alpha_0,\beta_0];s\cdot\bu\}$ relative
to $\Gamma_U(\mk{m},\mk{n})\leq GL_2(K)$,
 \item is of bi-weight $\{[\alpha_0,\beta_0];s\cdot\bu\}$ relative to $\Gamma_U(\mk{m},\mk{n};N)$.
\end{enumerate}
\end{Prop}

{\bf Proof} This follows from the definitions of $\Gamma_U(\mk{m},\mk{n})$, $\Gamma_U(\mk{m},\mk{n};N)$, 
Proposition \ref{pleut} and the identity \eqref{tiff}.  \fin

\subsubsection{The real matrix torus}
For a lattice $\ca{L}\subseteq \RR^g$, we define its {\it associated real torus} to be $\mathcal{T}_{\ca{L}}:=\RR^g/\ca{L}$.
For an integer $N\in\ZZ_{\geq 1}$, the $N$-torsion of $\mathcal{T}_{\ca{L}}$ is denoted by $\mathcal{T}_{\ca{L}}[N]:=\frac{1}{N}\ca{L}/\ca{L}$.
It is also convenient to define $\ca{T}_{\ca{L}}[\infty]:=\left(\bigcup_{N\geq 1}\ca{T}_{\ca{L}}[N]\right)$, i.e., the subgroup
of $K$-rational points of $\ca{T}_{\ca{L}}$ (which corresponds also to the torsion part of $\ca{T}_{\ca{L}}$).

We define the {\it real matrix torus} associated to an ordered pair of lattices $(\mk{m},\mk{n})$ as:
\begin{align*}
\mathcal{T}_{\mk{m},\mk{n}}:=\M{\RR^g/\mk{m}^*}{\RR^g/\mk{m}}{\RR^g/\mk{n}^*}{\RR^g/\mk{n}}.
\end{align*}
It is understood here that the lattices $\mk{m},\mk{m}^*,\mk{n}$ and $\mk{n}^*$ are embedded in 
$\RR^g$ through a fixed ordering of the $g$ embeddings of $K$ into $\RR$.
For $N,M\in\ZZ_{\geq 1}\cup\{\infty\}$, we define the $(N;M)$-torsion of the real matrix torus $\mathcal{T}_{\mk{m},\mk{n}}$ as
\begin{align*}
\mathcal{T}_{\mk{m},\mk{n}}[N;M]:=
\M{\ca{T}_{\mk{m}^*}[N]}{\ca{T}_{\mk{m}}[M]}{\ca{T}_{\mk{n}^*}[N]}{\ca{T}_{\mk{n}}[M]  }.
\end{align*}

We keep the same notation as in Section \ref{symm_a}. It is also convenient to define the following normalized Eisenstein series:
\begin{align*}
\wt{G}_{(\mk{m},\mk{n})}^{\alpha(s),\beta(s)}(U\,;z):=\frac{G_{(\mk{m},\mk{n})}^{\alpha(s),\beta(s)}(U\,;z)}{e^{2\pi\ii\ell_U}}.
\end{align*}
With the help of the terminology introduced above, we would like to point out two consequences which can be derived from 
the transformation formula \eqref{tiff}:

\begin{enumerate}[(i)]
 \item Let $\V{v_1}{v_2}\in \ca{T}_{\mk{m}}[\infty]\times \ca{T}_{\mk{n}}[\infty]$ be fixed. Then, the map 
\begin{align*}
\left[\left(z,s;\V{u_1}{u_2}\right)\mapsto  \wt{G}_{(\mk{m},\mk{n})}^{\alpha(s),\beta(s)}\left(\M{u_1}{v_1}{u_2}{v_2}\,;z\right)\right],
\end{align*}
descends to a function on $K_\CC^{\pm}\times\Pi_{1-\frac{w}{2}}\times\ca{T}_{\mk{m}^*}[\infty]\times \ca{T}_{\mk{n}^*}[\infty]$.
\item Let $\V{u_1}{u_2}\in \ca{T}_{\mk{m}^*}[\infty]\times \ca{T}_{\mk{n}^*}[\infty]$ be fixed. Then, the map
\begin{align*}
\left[\left(z,s;\V{v_1}{v_2}\right)\mapsto  G_{(\mk{m},\mk{n})}^{\alpha(s),\beta(s)}\left(\M{u_1}{v_1}{u_2}{v_2}\,;z\right)\right],
\end{align*}
descends to a function on $K_\CC^{\pm}\times\Pi_{1-\frac{w}{2}}\times\ca{T}_{\mk{m}}[\infty]\times \ca{T}_{\mk{n}}[\infty]$.
\end{enumerate}
Notice that in (i), it is the normalized function $\wt{G}$ which appears, while in (ii) it is the function $G$. 
\begin{Rem}
Unfortunately, when $g>1$, it is not possible to replace in (i) (resp. (ii)) the (discrete) $K$-rational torus 
$\ca{T}_{\mk{m}^*}[\infty]\times \ca{T}_{\mk{n}^*}[\infty]$ (resp. $\ca{T}_{\mk{m}}[\infty]\times \ca{T}_{\mk{n}}[\infty]$)
by the (continuous) real torus $\ca{T}_{\mk{m}^*}\times \ca{T}_{\mk{n}^*}$ (resp. $\ca{T}_{\mk{m}}\times \ca{T}_{\mk{n}}$). 
Let us first explain the obstruction for (ii). It is crucial in our definition 
\begin{align*}
G_{(\mk{m},\mk{n})}^{\alpha(s),\beta(s)}\left(\M{u_1}{v_1}{u_2}{v_2}\,;z\right)
\end{align*}
to quotient the summation indexing set 
\begin{align*}
S_{v_1,v_2}:=(\mk{m}+v_1,\mk{n}+v_2)\bs\{(0,0)\},
\end{align*}
by a finite index subgroup $\ca{V}$ of $\ca{O}_K^{\times}(\infty)$, in order to have a convergent series.
In particular, the set $S_{v_1,v_2}$ must be {\it stable} under the $\ca{V}$-action. If we allow $v_1,v_2\in\RR^g\bs\iota(K)$, then in general, the corresponding set 
$S_{v_1,v_2}$ will not be stable under {\it any} finite subgroup of $\ca{V}$. Similarly, there is also an obstruction for (i).
From the equation \eqref{obstruc}, we see that if $\epsilon\in\ca{V}$, $m\in\mk{m}$, $n\in\mk{n}$, 
it is crucial that the four elements $(\epsilon-1)u_1v_1,(\epsilon-1)u_2v_2,(\epsilon-1)u_1m,(\epsilon-1)u_2n$ have each an
absolute trace which lies in $\ZZ$; and this won't be case, in general, if $u_1,u_2\in\RR^g\bs\iota(K)$.
However, in the special case where $g=1$, i.e. when $K=\QQ$, the group $\Gamma_U(\mk{m},\mk{n})=\{1\}$. In particular, it is {\it independent} of the parameter matrix 
$U$! In this unique case, the definition of $G_{(\mk{m},\mk{n})}^{\alpha(s),\beta(s)}(U;z)$ is valid for any {\bf real} 
parameter matrix $U\in M_2(\RR)$. Moreover, the map (i) extends to a map on the real torus $\ca{T}_{\mk{m}^*}\times \ca{T}_{\mk{n}^*}$.
Similarly, the map (ii) extends to a map on the real torus $\ca{T}_{\mk{m}}\times \ca{T}_{\mk{n}}$.
\end{Rem}

\subsubsection{Product of two Eisenstein series}

We would like now to explain how one can get rid of the root of unity in the transformation formula \eqref{tiff} by 
considering a judicious product of two Eisenstein series. Let  $\ca{Q}=((\mk{m},\mk{n}),U,p,w)$ and $\wt{\ca{Q}}=((\mk{m},\mk{n}),U,\wt{p},\wt{w})$ be two standard quadruples.
For $s,\wt{s}\in\CC$, we let $[\alpha(s),\beta(s)]$ (resp.$[\wt{\alpha}(\wt{s}),\wt{\beta}(\wt{s})]$ ) be the associated bi-weight of $\ca{Q}$ (resp. $\wt{\ca{Q}}$).
Let $[\alpha_0,\beta_0]:=[w\cdot\bu-\frac{p}{2},\frac{p}{2}]$ and   $[\wt{\alpha}_0,\wt{\beta}_0]:=[\wt{w}\cdot\bu-\frac{\wt{p}}{2},\frac{\wt{p}}{2}]$.
It follows from  Proposition \ref{mod_for} that
\begin{enumerate}[(i)]
 \item $[z\mapsto G_{(\mk{m},\mk{n})}^{\alpha(s),\beta(s)}(U\,;z)]$  is modular of almost bi-weight $\{[\alpha_0,\beta_0];s\cdot\bu\}$
 relative to $\Gamma_U(\mk{m},\mk{n})$.
 \item $[z\mapsto G_{(\mk{m},\mk{n})}^{\wt{\alpha}(\wt{s}),\wt{\beta}(\wt{s})}(\iota_{11}\iota_{12}(U)\,;z)]$  is modular of almost bi-weight 
 $\{[\wt{\alpha}_0,\wt{\beta}_0];\wt{s}\cdot\bu\}$
 relative to $\Gamma_U(\mk{m},\mk{n})$.
\end{enumerate}
Let
\begin{align*}
H_{(\mk{m},\mk{n})}^w(U;z,s,\wt{s}):=G_{(\mk{m},\mk{n})}^{\alpha(s),\beta(s)}(U;z)\cdot G_{(\mk{m},\mk{n})}^{\tilde{\alpha}(\wt{s}),\tilde{\beta}(\wt{s})}(\iota_{11}\iota_{12}(U);z).
\end{align*}
Since $\zeta_{U,U^{\gamma}}\cdot\zeta_{\iota_{11}\iota_{12}(U),\iota_{11}\iota_{12}(U^{\gamma})}=1$, 
it follows from the transformation formula \eqref{tiff} that
\begin{align*}
[(z,s;U)\mapsto H_{(\mk{m},\mk{n})}^w(U;z,s,s')],
\end{align*}
descends to a function on $K_{\CC}^{\pm}\times\Pi_{w_0}\times\mathcal{T}_{\mk{m},\mk{n}}[\infty;\infty]$,
where $w_0=\max\{1-\frac{w}{2},1-\frac{\wt{w}}{2}\}$.
Moreover, if $U\in\frac{1}{N}\M{\mk{m}^*}{\mk{m}}{\mk{n}^*}{\mk{n}}$, it follows from the transformation formula \eqref{musc} that the function
\begin{align*}
[z\mapsto H_{(\mk{m},\mk{n})}^w(U;z,s,\wt{s})],
\end{align*} 
is modular (rather than just almost modular) of bi-weight $\{[\alpha_0+\wt{\alpha}_0,\beta_0+\wt{\beta}_0];(s+\wt{s})\cdot\bu\}$, 
relative to the group $\Gamma_U(\mk{m},\mk{n})$.

\subsection{Definition of the Fourier series expansion of a modular form of bi-weight $[\alpha,\beta]$ at an arbitrary cusp}\label{Fourier}
In this section, we give a precise definition of the {\it Fourier series expansion at a cusp $c$}
of a real analytic modular form $G(z)$. We also try to make precise the dependence of the Fourier series coefficients
on the choice of the local chart which is needed in our definition. 

Let $\Gamma\leq GL_2(K)$ be a subgroup commensurable to $GL_2(\ca{O}_K)$, and let $\alpha,\beta,\mu\in\CC^g$ be weights
which are chosen as in Definition \ref{difo3}. Let $G:K_\CC^{\pm}\rightarrow\CC$ be a real analytic modular form of bi-weight
$\{[\alpha,\beta];\mu\}$ relative to $\Gamma$.
Since $K_\CC^{\pm}=\bigcup\limits_{\ov{p}\in\ov{\mk{S}}}\mk{h}^{\ov{p}}$, every $\Gamma$-modular form may be viewed as 
a collection of $2^g$ modular forms (relative to suitable conjugates of $\Gamma$) on the connected space $\mk{h}^g$. Recall that
\begin{align*}
\ca{U}(K)=\M{1}{K}{0}{1}\leq GL_2(K),
\end{align*}
corresponds the maximal unipotent subgroup of $GL_2(K)$ of upper triangular matrices. We denote by $\pi:\ca{U}(K)\rightarrow K$
the natural projection on the upper-right entry. It will follow from the definition of $\ca{L}_{\eta}$ (the dual lattice of the \lq\lq indexing lattice\rq\rq\; 
associated to the local chart $\eta\in GL_2(K)$), that $\ca{L}_{\eta}$ only depends on $\Gamma^+=\Gamma\cap GL_2^+(K)$ rather than on $\Gamma$ itself. 
Because of the previous observation, there is no loss of generality if we assume from the outset that $\Gamma\leq GL_2^+(K)$. 
In this case, each connected component of $K_{\CC}^{\pm}$ will be stable under $\Gamma$. Moreover, since each component
of $K_\CC^{\pm}$ is real analytically isomorphic to $\mk{h}^g$, it will
be enough to define the Fourier series expansion for modular forms on the (connected) space $\mk{h}^g$. Therefore, in this 
section, we choose to view $G(z)$ as a function on $\mk{h}^g$ which is modular of bi-weight
$\{[\alpha,\beta];\mu\}$, relative to a discrete group $\Gamma\leq GL_2^+(K)$.

We let $Y_{\Gamma}:=\biglslant{\Gamma}{\mk{h}^g}$. As it is well-known, the space  $Y_{\Gamma}$ is 
an open Riemannian orbifold of finite volume, which admits a finite number
of relative $\Gamma$-cusps (see Appendix \ref{app_5} for some background on the notion of cusps). 
We define
\begin{align*}
\Latt(\Gamma):=\pi(\Gamma\cap \ca{U}(K))\leq K.
\end{align*}
By definition, the set $\Latt(\Gamma)$ is a lattice of $K$ of maximal rank.
Let $c\in\PP^1(K)$ be a fixed $K$-rational cusp. If $c=\infty$,
we set $\eta=I_2$. If $c\neq\infty$, we choose, arbitrarily, a \lq\lq local chart\rq\rq\; at the cusp $c$, i.e., a matrix $\eta\in GL_2^+(K)$ such that $\eta(\infty)=c$. 
Note that the group $GL_2^+(K)$ acts transitively on $\PP^1(K)$ (in fact, the group $SL_2(K)$ acts already transitively on $\PP^1(K)$). One can check that the shifted function 
\begin{align*}
G^{\eta}(z):=G\big|_{\{[\alpha,\beta];\mu\},\eta}(z),
\end{align*}
is invariant under the discrete subgroup $H_{\eta}:=(\eta^{-1}\Gamma\eta)\cap\Gamma$. Note that the the function $G^{\eta}(z)$ and
the group $H_{\eta}$ only depend on the left coset $\Gamma\eta$ rather than the matrix $\eta$ itself. We define
\begin{align}\label{heart3}
\mathcal{L}_{\eta}:=\pi\left(H_{\eta}\cap\mathcal{U}(K)\right),
\end{align}
which is  a lattice of maximal rank inside $K$. Let $\theta:K\hookrightarrow \RR^g$ be the embedding 
induced from a chosen ordering of the field embeddings of $K$ into $\RR$. For a fixed value of $y\in\RR_{>0}^g$, the function
$[x\mapsto G^{\eta}(x+\ii y)]$ descends to a function on the real torus $\RR^g/\theta (\ca{L}_{\eta})$. We may
now define the Fourier series expansion of $G(z)$ at the cusp $c$.
\begin{Def}\label{q_expan}
The \lq\lq {\it Fourier
series expansion of $G(z)$ at the cusp $c$, with respect to the local chart $\eta$}\rq\rq\s, is defined as
\begin{align}\label{toq}\notag
\mbox{Fourier series expansion of}\;[\;x\mapsto & G^{\eta}(x+\ii y)\;]:=\\[3mm]
&\hspace{-3cm}a_0^{\eta}(y)+\sum_{\xi\in\mathcal{L}_{\eta}^*\bs\{0\}} 
a_{\xi}^{\eta}(y)e^{2\pi\ii\Tr(\xi x)},
\end{align}
where $z=x+\ii y\in\mk{h}^g$. Here, $\mathcal{L}_{\eta}^*$ denotes the dual lattice of $\mathcal{L}_{\eta}$ with respect to the trace pairing. 
\end{Def}
We note that the Fourier series expansion \eqref{toq} exists and that it computes the value of $G^{\eta}(z)$ in virtue of Theorem \ref{keyy}.

\begin{Rem}
If $c=\infty$, then our definition of \lq\lq the Fourier series expansion at $\infty$\rq\rq\; corresponds to the usual
Fourier series of a function on a torus. Moreover, in this special case, no choice is required and therefore it is uniquely determined.  
For a general cusp $c$ and a choice of a local chart $\eta$ at $c$, 
one may check that the Fourier series expansion given in \eqref{toq} only depends on the left coset $\Gamma\eta$ and not
on the matrix $\eta$ itself. From the previous observation, one may be lead to think that the Fourier series coefficients $\{a_{\xi}^{\eta}(y)\}_{\xi\in\ca{L}_{\eta}^*}$ 
will only depend on the relative cusp $[c]_{\Gamma}$ rather than on the cusp $c$ itself (and, therefore, a fortiori, not on the choice of the local chart $\eta$). 
Unfortunately, this is not the case. For example, let $\eta':=\eta\gamma$,
where $\gamma=\M{1}{\mu}{0}{1}$ and $\mu\in\Latt(\Gamma)$ (so that $\gamma\in\Gamma\cap\ca{U}_K$). 
From the commutativity of the group $\ca{U}(K)$, we readily see that
\begin{align*}
(\gamma^{-1}\eta^{-1}\Gamma\eta\gamma)\cap\Gamma\cap\mathcal{U}(K)=(\eta^{-1}\Gamma\eta)\cap\Gamma\cap\mathcal{U}(K),
\end{align*}
and, therefore, $\ca{L}_{\eta}=\ca{L}_{\eta'}$. In particular, the two Fourier series expansions will have the same indexing set.  
However, in general, given an element $\xi\in\ca{L}_{\eta}=\ca{L}_{\eta'}$, one may check that the coefficient $a_{\xi}^{\eta}(y)$ will
differ from the coefficient $a_{\xi}^{\eta'}(y)$ by the (not necessarily trivial) root of unity $e^{2\pi\ii\Tr(\xi\mu)}$.
Despite this dependence of the Fourier coefficients on the choice of the local chart $\eta$ at the cusp $c$, the
vanishing (or non-vanishing) of $a_{\xi}^{\eta}(y)$ is independent of the choice of $\eta$, and therefore, will {\it only depend} on the relative
cusp $[c]_{\Gamma}$.
Indeed, let $\eta,\eta'\in GL_2^+(K)$ be two local charts at the cusp $c\in\PP^1(K)$. Since $\eta(\infty)=\eta'(\infty)=c$, there exists a matrix
\begin{align}\label{heart2}
\gamma\in \Stab_{GL_2^+(K)}(\infty):=\left\{\M{a}{b}{0}{d}:ad\in K^{\times},ad\gg 0,b\in K\right\},
\end{align}
such that $\eta'=\eta\gamma$, so that
\begin{align}\label{heart5}
G\big|_{\{[\alpha,\beta];\mu\},\eta'}(z)=(G\big|_{\{[\alpha,\beta];\mu\},\eta})\big|_{\{[\alpha,\beta];\mu\},\gamma}(z)=
\frac{\omega_p(d)}{|\Norm(d)|^{2s}}\cdot |ad|^{-\mu}\cdot G|_{\{[\alpha,\beta];\mu\},\eta}\left(\frac{az}{d}+\frac{b}{d}\right).
\end{align}
Finally, it follows from \eqref{heart5} that for all $\xi\in\ca{L}_{\eta}$ (so that $\frac{d\xi}{a}\in\ca{L}_{\eta'}$) and all $y\in\RR_{>0}^g$, 
that
\begin{align*}
a_{\xi}^{\eta'}(y)=0\Longleftrightarrow a_{\frac{d\xi}{a}}^{\eta}\left(\frac{ay}{d}\right)=0.
\end{align*}
\end{Rem}

\begin{Rem}
Note that if  the function $G:\mk{h}^g\rightarrow\CC$ is only assumed to be a real analytic modular form
of {\it almost} bi-weight $\{[\alpha,\beta];\mu\}$ (relative to $\Gamma$ in the sense of Definition \ref{difo3}) then 
Definition \ref{q_expan} still makes perfect sense. 
\end{Rem}

\subsubsection{Cuspidality and square-integrability at a cusp}\label{cusp}
We keep the same notation as in the previous section, except that we now only suppose that 
$G:\mk{h}^g\rightarrow\CC$ is a real analytic modular form of {\it almost} bi-weight $\{[\alpha,\beta];\mu\}$ relative to $\Gamma$.
Let $c\in\PP^1(K)$ be a fixed cusp and let $\eta\in GL_2(K)$ be a local chart at $c$, so that
$\eta\infty=c$. Consider the Fourier series expansion of $[x\mapsto G^{\eta}(x+\ii y)]$, (in the sense of Definition \ref{q_expan})
\begin{align*}
a_0^{\eta}(y)+\sum_{\xi\in\mathcal{L}_{\eta}^*\bs\{0\}} 
a_{\xi}^{\eta}(y)e^{2\pi\ii\Tr(\xi x)}.
\end{align*}
\begin{Def}\label{tok}
We say that $G(z)$ is cuspidal at $c$ if the map $[y\mapsto a_0^{\eta}(y)]\equiv 0$.
We say that $G(z)$ is square-integrable at $c$, if the integral 
\begin{align*}
\int_{B_{[\infty]_\Gamma}} \Big|G^{\eta}(z)|^2 dV(z)<\infty,
\end{align*} 
where $B_{[\infty]_\Gamma}$ corresponds to the neighborhood of the relative cusp $[\infty]_{\Gamma}$
of the orbifold $Y_{\Gamma}$ as defined in Appendix \ref{app_5a}.
\end{Def}
We claim that notion of being cuspidal (or square-integrable) at a relative cusp $[c]_{\Gamma}$ is well-defined, i.e. that it does not depend
on the particular choice of the cusp $c\in [c]_{\Gamma}$, and on the particular choice of 
the local chart $\eta$ at $c$. 

Let $c\in[c]_{\Gamma}$ be fixed. We first show the independence on the choice of the local chart $\eta$ at $c$. 
Let $\eta'\in GL_2^+(K)$ be another local chart at $c$. We thus have $\eta'\infty=c$. In particular, $\eta'=\eta\gamma$ for some
$\gamma\in\M{a}{b}{0}{d}\in\Stab_{GL_2^+(K)}(\infty)$. By definition of the $\{[\alpha,\beta];\mu\}$-slash action, we have
\begin{align}\label{pok}
G\big|_{\{[\alpha,\beta];\mu\},\eta'}(z)\stackrel{\lcdot}{=}(G\big|_{\{[\alpha,\beta];\mu\},\eta})\big|_{\{[\alpha,\beta];\mu\},\gamma}(z)\stackrel{\lcdot}{=}
d^{-\alpha}d^{-\beta}\cdot |ad|^{-\mu}\cdot G|_{\{[\alpha,\beta];\mu\},\eta}\left(\frac{az}{d}+\frac{b}{d}\right),
\end{align}
where the symbol \lq\lq$\stackrel{\lcdot}{=}$\rq\rq\s means equal, up to some root of unity. It follows from \eqref{pok},
that $[y\mapsto a_0^{\eta}(y)]\equiv 0$ (reps.  $G\big|_{\{[\alpha,\beta];\mu\},\eta'}(z)$ is square-integrable) 
if and only if $[y\mapsto a_0^{\eta'}(y)]\equiv 0$ (resp. $G\big|_{\{[\alpha,\beta];\mu\},\eta}(z)$ is square-integrable). 

Let us now prove that the cuspidality at $c$ (or the square-integrability at $c$) only depend on the relative cusp $[c]_{\Gamma}$ on not on cusp $c$ itself.
Let $c'\in [c]_{\Gamma}$. Then, there exists $\gamma\in\Gamma$ such that $\gamma c=c'$.
Let $\eta$ be a local chart at $c$. Then the matrix $\eta':=\gamma\eta$ is a local chart at $c'$. Since $G(z)$ is modular of bi-weight $\{[\alpha,\beta];\mu\}$
relative to $\Gamma$, we have
\begin{align}\label{carti}
(G\big|_{\{[\alpha,\beta];\mu\},\gamma})\big|_{\{[\alpha,\beta];\mu\},\eta}(z) \stackrel{\lcdot}{=}  G\big|_{[\alpha,\beta],\eta'}(z).
\end{align}
It follows from \eqref{carti} that $G(z)$ is cuspidal (resp. square-integrable) at $c$ if and only 
if it is cuspidal (resp. square-integrable) at $c'$.

\begin{Rem}\label{ex_sq}
Note that, in general, being cuspidal at $c$ is a stronger condition than being square integrable at $c$.
Let us illustrate this difference on a concrete example. 
Let $G(z,s):=G_{(\mk{m},\mk{n})}^0(U,p\, ;z,s)$. It follows from Theorem \ref{staff} that constant
term of $G(z,s)$ at a cusp $c\in\PP^1(K)$ is of the form
\begin{align*}
\phi_1^{c}(s)\cdot|\Norm(y)|^s+\phi_2^{c}(s)\cdot|\Norm(y)|^{1-s}, 
\end{align*}
where $\phi_1^{c}(s)$ and $\phi_2^{c}(s)$ are some holomorphic functions on $\Pi_1$. 
For a fixed $s_0\in\Pi_1$, one easily checks that
\begin{enumerate}
 \item $G(z,s_0)$ is cuspidal at $c$ $\Longleftrightarrow$ $\phi_1^{c}(s_0)=\phi_2^{c}(s_0)=0$.
 \item $G(z,s_0)$ is square integrable at $c$ $\Longleftrightarrow$ $\phi_1^{c}(s_0)=0$.
\end{enumerate}
\end{Rem}

\subsection{Growth estimates of $G_{(\mk{m},\mk{n})}^{0}(U,p\,;z,s)$ in the right half-plane $\Ree(s)>1$}\label{groo}

Let $\ca{Q}=((\mk{m},\mk{n}),U,p,w)$ be a standard quadruple and let
$G_{(\mk{m},\mk{n})}^{\alpha(s),\beta(s)}(U;z)=G_{\mk{m},\mk{n}}^w(U,p\,;z,s)$
be the Eisenstein series which appears in Definition \ref{difo2} (or in Definition \ref{eis_ser}). In Section \ref{def_eis},
we claimed that the defining series of $G_{(\mk{m},\mk{n})}^{\alpha(s),\beta(s)}(U;z)$ was absolutely convergent when $s\in\Pi_{1-\frac{w}{2}}$. The 
goal of this section is to prove this result by comparing it to a finite sum of {\it real analytic Poincar\'e-Eisenstein series
of weight $0$}. In fact, we don't only prove the absolute convergence of 
$G_{(\mk{m},\mk{n})}^{\alpha(s),\beta(s)}(U;z)$,  but we also 
provide some non-trivial growth estimates of $[z\mapsto G_{(\mk{m},\mk{n})}^{\alpha(s),\beta(s)}(U;z)]$,
when $z$ tends to an arbitrary cusp $c\in\PP^1(K)$. All these estimates boil down to growth
estimates for the classical real analytic Poincar\'e-Eisenstein series $E_i(z,s)$ (of weight $0$), for $i\in\{1,\ldots,h\}$, 
where $h$ is the number of relative $\Gamma_K$-cups ($h$ also corresponds to the wide ideal class group of $K$). 
Here $\Gamma_K:=SL_2(\ca{O}_K)$ is the classical {\it Hilbert modular group}.
For the sake of completeness and for the benefit of the reader, 
we included in Appendix \ref{app_6} a proof of the key Proposition \ref{poids} which
proof is based on the notion of a {\it point-pair invariant kernel}; a notion, which appeared prominently in the work of Selberg, see \cite{Sel1}. 
\begin{Rem}
In Section \ref{lun} (see Theorem \ref{lune}), using the explicit formula
for the Fourier series expansion of $[z\mapsto G_{(\mk{m},\mk{n})}^{\alpha(s),\beta(s)}(U;z)]$, we provide sharper growth estimates that the ones given below, 
which hold true for for {\it all} $s\in\CC$ (away from the poles); 
in particular these estimates remain valid when $s$ approaches the boundary of convergence of the series in \eqref{tonu}. Note though that the proof  
of Theorem \ref{lune} lies much deeper that the proof of Proposition \ref{gr_est}, since it requires the explicit description of 
the Fourier series expansion of $[z\mapsto G_{(\mk{m},\mk{n})}^{\alpha(s),\beta(s)}(U;z)]$ and the fact that
$[s\mapsto G_{(\mk{m},\mk{n})}^{\alpha(s),\beta(s)}(U;z)]$ admits a meromorphic continuation to all of $\CC$.
\end{Rem}
\begin{Prop}\label{gr_est}
The series on the right-hand side of \eqref{mou1}
converges absolutely for $s>1-\frac{w}{2}$. Moreover, for each $s>1-\frac{w}{2}$ fixed,
there exists positive constants $C_{1,s}$, such that
for all $z\in K_\CC^{\times}$, we have
\begin{align}\label{fat0}
\Big|G_{(\mk{m},\mk{n})}^{\alpha(s),\beta(s)}(U\,;z)-e_1\cdot \delta_{\mk{m}}(v_1)\cdot 
\wt{Z}_{\mk{n}}(v_2,u_2,\omega_{\ov{p}}\omega_{\epsilon},2s+w)|\Norm(y)|^s\Big|\leq C_{1,s}\cdot|\Norm(y)|^{1-s-w}.
\end{align}
Here, $\epsilon=\ov{\bz}$, if $w$ is even, $\epsilon=\ov{\bu}$, if $w$ is odd and
$e_1:=e_1((\mk{m},\mk{n}),U)$ is the positive integer which appears in Definition \ref{chien}.

Let $\mk{c}=\frac{a'}{c'}\in K\subseteq\PP^1(K)$ be a finite cusp. Let
$\mathcal{R}$ be defined as in Definition \ref{set_def} and let
\begin{align}\label{hot}
\mathcal{R}[\mk{c}]:=\{(c,d)\in\mathcal{R}:ca'+dc'=0\}. 
\end{align}
Note that the pairs in $\mathcal{R}$ are in one-to-one correspondence with the terms in the series \eqref{tonu} 
(or the series \eqref{mou1}), which tend to infinity when $z\rightarrow\mk{c}$. Moreover, for each $s>1$ fixed,
there exists a positive constant $C_{2,s}\in\RR_{>0}$ such that,
for all $z\in K_\CC^{\pm}$,
\begin{align}\label{fat}
\Big|G_{(\mk{m},\mk{n})}^{\alpha(s),\beta(s)}(U\,;z)-\sum_{(c,d)\in\mathcal{R}[\mk{c}]}
\frac{e^{{2\pi\ii } \Tr(u_1 c+u_2 d)}\cdot|y|^{\bu\cdot s}}{P(\alpha(s),\beta(s);cz+d)}\Big|
\leq C_{2,s}\cdot |\Norm(y)|^{s-1},
\end{align}
In particular, in the special case where $\mathcal{R}[\mk{c}]=\emptyset$, for any fixed
value of $s>1$, we find that $G_{(\mk{m},\mk{n})}^{\alpha(s),\beta(s)}(U\,;z)=O(|\Norm(y)|^{s-1})$
as $z\rightarrow\mk{c}$.
\end{Prop}
\begin{Rem}
Note that the right-hand side of \eqref{fat}
{\it does not} depend on $w$ while \eqref{fat0} does. 
\end{Rem}
\begin{Rem}\label{couch}
Using the general transformation formula \eqref{chat} for $G_{(\mk{m},\mk{n})}^{\alpha(s),\beta(s)}(U\,;z)$, one can show that the
inequalities \eqref{fat0} and \eqref{fat} are in fact equivalent. It is also
possible to prove \eqref{fat} {\it directly} using the same set of ideas (the point-pair invariant kernel method) as in the proof of Proposition \ref{poids}.
\end{Rem}

{\bf Proof} Thanks to the formula \eqref{tram}, in order to show \eqref{fat0} and \eqref{fat}, it is enough 
to prove it in the special case, where $w=0$. From now on, we suppose that $w=0$. From
Remark \ref{couch}, it is also enough to prove the growth estimate \eqref{fat0}, i.e., when the cusp $\mk{c}=\infty$.

Let $z=x+\ii y\in K_{\CC}^{\pm}$ and $s>1$ be fixed. We let
$\mathcal{T}$ be defined as in Definition \ref{set_def}, except that 
we {\it further impose} that each matrix element in $\mathcal{T}$ lies in $SL_2(K)$ rather than just in $GL_2(K)$. 
This additional restriction will eliminate the presence of some determinant factors. We also let
\begin{enumerate}
\item $\mathcal{R}:=\left\{(c,d):\M{*}{*}{c}{d}\in\mathcal{T}\right\}$,
\item $\wt{\mathcal{T}}=\left\{\M{a}{b}{c}{d}\in\mathcal{T}:c\neq 0\right\}$.
\item $\wt{\mathcal{R}}=\{(c,d)\in\mathcal{R}:c\neq 0\}$.
\end{enumerate}
Note that the pairs in $\ca{R}\bs\wt{\ca{R}}$ are in one-to-one correspondence with the terms of the series \eqref{tonu} which 
tend to infinity as $z\rightarrow\infty$. Here $\infty$ stands for the cusp $[1,0]\in\PP^1(K)$.
Note also that $\mathcal{T}\bs\wt{\mathcal{T}}\neq\emptyset\Longleftrightarrow v_1\in\mk{m}\Longleftrightarrow \delta_{\mk{m}}(v_1)=1$. 
In the case where $\delta_{\mk{m}}(v_1)=1$, a direct and easy computation shows (see for example the first part of the proof of Theorem \ref{key_thmm}) that
\begin{align}\label{church0}
\sum_{(0,d)\in\ca{R}\bs\wt{\ca{R}}} 
\frac{e^{2\pi\ii\Tr(u_1(m+v_1)+u_2(m+v_2))}|y|^{\bu\cdot s}}{P(\alpha(s),\beta(s),n+v_2)}=e_1\cdot \wt{Z}_{\mk{n}}(v_2,u_2,\omega_{\ov{p}},2s)\cdot|y|^{\bu\cdot s},
\end{align}
where $e_1:=e_1((\mk{m},\mk{n}),U)$ is the positive integer which appears in Definition \ref{chien} and
$\wt{Z}_{\mk{n}}(v_2,u_2,\omega_{\ov{p}};s)$ is the zeta function which appears in Definition \ref{nor_zeta}. 

For any $c,d\in K$ and $z\in K_{\CC}^{\pm}$, we have the inequality
\begin{align}\label{church}
\left|\frac{e^{2\pi\ii\Tr(u_1(m+v_1)+u_2(m+v_2))}\cdot|y|^{\bu\cdot s}}{P(\alpha(s),\beta(s),(m+v_1)z+(n+v_2))}
\right|\leq \frac{|y|^{\bu\cdot s}}{|\Norm(cz+d)|^{2s}}.
\end{align}
Combining \eqref{church0} with \eqref{church} we thus obtain
\begin{align*}
|G_{(\mk{m},\mk{n})}^{\alpha(s),\beta(s)}(U\,;z)-e_1\cdot \delta_{\mk{m}}(v_1)\cdot Z_{\mk{n}}(v_2,u_2;\omega_{\ov{p}};2s)|\Norm(y)|^s|
& \leq \sum_{(c,d)\in\wt{\mathcal{R}}}
\frac{|y|^{\bu\cdot s}}{|\Norm(cz+d)|^{2s}}\\[3mm]
&=\sum_{\gamma\in\wt{\mathcal{T}}}|\Imm(\gamma z)|^{\bu\cdot s}.
\end{align*}
In order to prove the inequality \eqref{fat0}, it is enough to show that, for each $s>1$ fixed,
there exists a positive constant $C_{1,s}\in\RR_{>0}$,
such that for all $z\in K_\CC^{\pm}$,
\begin{align}\label{hug}
\sum_{\gamma\in\wt{\mathcal{T}}}|\Imm(\gamma z)|^{\bu\cdot s}\leq C_{1,s}\cdot|\Norm(y)|^{1-s}.
\end{align}

Since $(\mk{m}+v_1,\mk{n}+v_2)\subseteq K^2$ is a discrete subset,
there exists a rational integer $M_1\in\ZZ_{\geq 1}$, such that 
\begin{align}\label{fas}
M_1\cdot(\mk{m}+v_1,\mk{n}+v_2)\subseteq\ca{O}_K\times\ca{O}_K.
\end{align}
Recall that $\ca{V}^+:=\ca{V}_U^+(\mk{m},\mk{n})$ is a subgroup of totally positive units 
which acts diagonally on the set $(\mk{m}+v_1,\mk{n}+v_2)$, and which is indexing the summation in \eqref{tonu}.
The group $\ca{V}^+$ acts as well on the set $(\ca{O}_K\bs\{0\})\times \ca{O}_K$. 
We let $\mathcal{R}'$ be a complete set of representatives of $(\ca{O}_K\bs\{0\})\times \ca{O}_K$ under the diagonal action of 
$\ca{V}^+$. It follows from \eqref{fas} that 
\begin{align}\label{fas2}
M_1\cdot\wt{\mathcal{R}}\subseteq \mathcal{R}'\pmod{\mathcal{V}^+}.
\end{align}
From \eqref{fas2}, we get
\begin{align*}
\frac{1}{M_1^{2s}}\sum_{\gamma\in\wt{\mathcal{T}}}|\Imm(\gamma z)|^{\bu\cdot s}\leq \sum_{(c,d)\in\ca{R}'}\frac{|y|^{\bu\cdot s}}{|\Norm(cz+d)|^{2s}}.
\end{align*}
Therefore, in order to prove \eqref{hug}, it is enough to show that for each
$s\in>1$ fixed, there exists a constant $C_{1,s}'$, such that for all $z\in K_\CC^{\pm}$,
\begin{align}\label{fas3}
 \sum_{(c,d)\in\mathcal{R}'}
\frac{|y|^{\bu\cdot s}}{|\Norm(cz+d)|^{2s}}\leq C_{1,s}'\cdot|\Norm(y)|^{1-s}.
\end{align}
If we let $\ca{R}''$ be a complete set of representatives of  
$(\ca{O}_K\bs\{0\})\times \ca{O}_K$ under the diagonal action of $\ca{O}_K^{\times}$ (these units are not necessarily totally positive), then
\begin{align}\label{chap}
\sum_{(c,d)\in\mathcal{R}'}
\frac{|y|^{\bu\cdot s}}{|\Norm(cz+d)|^{2s}}=[\ca{V}^+:\ca{O}_K^{\times}]\cdot
\sum_{(c,d)\in\mathcal{R}''}\frac{|y|^{\bu\cdot s}}{|\Norm(cz+d)|^{2s}}.
\end{align}
From \eqref{chap}, we see that the inequality \eqref{fas3} is equivalent to
show that for each $s\in\PP_1$ fixed, there exists a positive constant $C_{1,s}''\in\RR_{>0}$, such that for all $z\in K_\CC^{\pm}$,
\begin{align}\label{fas4}
\sum_{(c,d)\in\mathcal{R}''}
\frac{|y|^{\bu\cdot s}}{|\Norm(cz+d)|^{2s}}\leq C_{1,s}''\cdot|\Norm(y)|^{1-s},
\end{align}
Finally, the proof of the inequality \eqref{fas4} follows 
from the following key proposition:

\begin{Prop}\label{poids}
$\ca{R}''$ be a complete set of representatives of  
$(\ca{O}_K\bs\{0\})\times \ca{O}_K$ under the diagonal action of $\ca{O}_K^{\times}$. Then, if
$s>1$ and $z\in K_{\CC}^{\pm}$, the series 
\begin{align}\label{chag}
\ca{E}(z,s):=\sum_{(c,d)\in\mathcal{R}''}
\frac{|y|^{\bu\cdot s}}{|\Norm(cz+d)|^{2s}},
\end{align}
converges (and therefore necessarily absolutely). Moreover, for each $s>1$ fixed, 
there exists a positive constant $C_s>0$, such that,
for all all $z\in K_\CC^{\pm}$,
\begin{align}\label{chag2}
\sum_{(c,d)\in\mathcal{R}''}\frac{|y|^{\bu\cdot s}}{|\Norm(cz+d)|^{2s}}\leq C_s\cdot|\Norm(y)|^{1-s}.
\end{align}
\end{Prop}

{\bf Proof} See Appendix \ref{app_6} for a proof of \eqref{chag2} which uses the notion of a point-pair invariant kernel.
\begin{Rem}
Assuming only the absolute convergence of $E(z,s)$ (which can be proved in an elementary way), 
one may also obtain the inequality \eqref{chag2} as a consequence of
Theorem \ref{key_thmm}, which provides the Fourier series expansion of $E(z,s)$ at
the cusp $\infty$.
\end{Rem}

\section{Maa\ss\;-graded operators and the Hilbert space $L^2(\mk{h};\Gamma;p)$}\label{car0}

\subsection{Mass\ss\;-graded operators and partial-graded Laplacians}\label{def_mas}

Let us now introduce an important class of graded differential operators on the space of smooth
functions $\mathcal{M}_{\infty}=\Maps_{\mbox{\tiny{smooth}}}(K_{\CC}^{\pm},\CC)$. For some motivations 
and historical references for these operators when $g=1$, see Section 2.1 of \cite{Bump}.

Let $s\in\CC$ and $p=(p_j)_{j=1}^{g}\in\ZZ^g$ be fixed. For $f\in\ca{M}_{\infty}$, $z\in K_\CC^{\pm}$ and $\gamma\in\ca{G}(\RR)$, recall that
the $\{p;s\}$-slash action is defined as 
\begin{align}\label{boute}
f\big|_{\{p;s\},\gamma}(z)=\omega_p(j(\gamma,z))^{-1}\cdot|\det(\gamma)|^{-\bu\cdot s}\cdot f(\gamma z).
\end{align}
For $j\in\{1,\ldots,g\}$, we let
\begin{align*}
R_{j,p_j}=(z_j-\ov{z}_j)\frac{\partial}{\partial z_j}+\frac{p_j}{2},
\end{align*}
be the  {\it Maa\ss\; raising operator} of weight $p_j$ with respect to the $j$-th coordinate. Similarly, we let
\begin{align*}
L_{j,p_j}=-(z_j-\ov{z}_j)\frac{\partial}{\partial \ov{z}_i}-\frac{p_j}{2},
\end{align*}
be the {\it Ma\ss\; lowering operator} of weight $p_j$ with respect to the $j$-th coordinate of $z$. We also define
\begin{align}\label{grad_lap}
\Delta_{j,p_j}:=-L_{j,p_j+2}\circ R_{j,p_j} -\frac{p_j}{2}\left(1+\frac{p_j}{2}\right)
&=-R_{j,p_j-2}\circ L_{j,p_j} +\frac{p_j}{2}\left(1-\frac{p_j}{2}\right)\\
&=-y_j^2\left(\frac{\partial^2}{\partial x_j^2}+ \frac{\partial^2}{\partial y_j^2}\right)+\ii p_j\cdot y_j\frac{\partial}{\partial x_j},
\end{align}
be the {\it weight $p_j$ hyperbolic Laplacian} with respect to the $j$-th coordinate of $z$. More concisely,
we sometimes simply say that $\Delta_{j,p_j}$ is a partial-graded Laplacian (of weight $p_j$ with respect to the $j$-th coordinate of $z$).

The differential operators
$L_{j,*}, R_{j,*}$ and $\Delta_{j,*}$ are compatible with the $\{p;s\}$-slash action in the sense
that the following associative formulas hold true: Let $f\in\mathcal{M}_{\infty}$, then, for all $\gamma\in\ca{G}(\RR)$, we have
\begin{enumerate}[(a)]
 \item $(R_{j,p_j}f)\big|_{\{p+2e_j;s\},\gamma}=R_{j,p_j}\left(f\big|_{\{p;s\},\gamma}\right)$, 
 \item $(L_{j,p_j}f)\big|_{\{p-2 e_j;s\},\gamma}=L_{j,p_j}\left(f\big|_{\{p;s\},\gamma}\right)$,
 \item $(\Delta_{j,p_j}f)\big|_{\{p;s\},\gamma}=\Delta_{j,p_j}\left(f\big|_{\{p;s\},\gamma}\right)$.
\end{enumerate}
Recall here that $e_j\in\ZZ^g$ is the row vector with $1$ place in the $j$-th coordinate and zero elsewhere.
For a proof of these facts, see, for example, Lemma 2.1.1 of \cite{Bump}.
These associative laws are sometimes summarized by saying that the operators $\Delta_{j,p_j}$, $R_{j,p_j}$ and $L_{j,p_j}$
commute with the slash operation $\big|_{\{*;s\},\gamma}$. The previous symbol $*$ represents the weight which changes
according to the rules (a),(b) or (c).

\subsubsection{Behavior of $G_{(\mk{m},\mk{n})}^0(U,p\, ;z,s))$ under Maa\ss\,-graded operators}\label{eigen1}

In this section, we show that the real analytic Eisenstein series 
$G_{(\mk{m},\mk{n})}^w(U,p\, ;z,s)$, in the special case where $w=0$,
is stable under the partial-graded Laplacians. Note that, if $w\neq 0$, 
this is no longer true.

Let $s\in\CC$ and $p\in\ZZ^g$ be fixed. 
For $z\in(\mk{h}^{\pm})^g$, the function
\begin{align*}
z\mapsto f_s(z):=|\Imm(z)|^{s\cdot\bu}
\end{align*} 
has unitary weight $\{[-s,-s],s\cdot\bu\}$ relative to the continuous group $\ca{G}(\RR)$. 
Moreover, writing $\Imm(z)$ as $\Imm(z)=\frac{(z-\ov{z})}{2\ii}$, a direct computation
shows that, for all $z\in K_{\CC}^{\pm}$ and $j\in\{1,\ldots,g\}$ that
\begin{enumerate}
 \item[(d)]  $(R_{j,p_j} f_s(z))=\left(s+\frac{p_j}{2}\right)\cdot f_s(z)$,
 \item[(e)]  $(L_{j,p_j} f_s(z))=\left(s-\frac{p_j}{2}\right)\cdot f_s(z)$,
 \item[(f)]  $(\Delta_{j,p_j} f_s(z))=s(1-s)\cdot f_s(z)$.
\end{enumerate}

Each term appearing in the series in \eqref{mou1} is of the form
\begin{align}\label{oop}
e^{2\pi\ii \Tr(\theta)}\cdot f_s\big|_{\{-p,s\},\gamma}(z)
=\omega_{-p}(j(\gamma,z))^{-1}\cdot|\det(\gamma)|^{-\bu\cdot s}\cdot|\Imm(\gamma z)|^{\bu\cdot s}\cdot e^{2\pi\ii \Tr(\theta)},
\end{align}
for suitable elements $\gamma\in GL_2(K)$ and $\theta\in K$. It follows from \eqref{oop} that for all $\gamma\in\ca{G}(\RR)$,
\begin{align}\label{mou2}
R_{j,-p_j}\left(f_s\big|_{\{-p;s\},\gamma}(z)\right)= (R_{j,-p_j}f_s)\Big|_{\{-p+2e_j;s\},\gamma}(z)=
\left(s-\frac{p_j}{2}\right)\cdot f_s\big|_{\{-p+2e_j;s\},\gamma}(z),
\end{align}
where the first equality follows from (a) in Section \ref{def_mas}, and the second equality follows
from (d) in Section \ref{eigen1}. More explicitly, the identity 
\eqref{mou2} can be written as:
\begin{align}\label{nerff}
R_{j,-p_j} \Big(\omega_{p}(j(\gamma,z))\cdot|\Imm(\gamma z)|^{\bu\cdot s}\Big)=
(s-p_j/2)\cdot\omega_{p-2e_j}(j(\gamma,z))\cdot|\Imm(\gamma z)|^{\bu\cdot s}. 
\end{align}
From \eqref{nerff}, we readily deduce that
\begin{align}\label{nerf1}
R_{j,-p_j} (G_{(\mk{m},\mk{n})}(U,p\, ;z,s))=(s-p_j/2)\cdot G_{(\mk{m},\mk{n})}(U,p'\,; z,s), 
\end{align}
where $p'=p-2e_j$. A similar computation shows that
\begin{align}\label{nerf2}
L_{j,-p_j} (G_{(\mk{m},\mk{n})}^0(U,p\, ;z,s))=(s+p_j/2)\cdot G_{(\mk{m},\mk{n})}^0(U,{}'p, z,s), 
\end{align}
where ${}'p=p+2e_j$. In particular, it follows from \eqref{nerf1}, \eqref{nerf2} and the definition
of $\Delta_{j,-p_j}$ that
\begin{align}\label{tim7}
\Delta_{j,-p_j} (G_{(\mk{m},\mk{n})}^0(U,p\, ;z,s))=s(1-s)\cdot G_{(\mk{m},\mk{n})}^0(U,p\,; z,s).
\end{align}
\begin{Rem}
The identity \eqref{tim7} implies that the coefficients of the Fourier series expansion
of $[z\mapsto G_{(\mk{m},\mk{n})}^0(U,p\, ;z,s)]$ satisfy a certain linear system of ODEs of order $2$ in $g$ variables (see Appendix \ref{app_3}), 
and that the coefficients of the Taylor series expansion of $[s\mapsto G_{(\mk{m},\mk{n})}^0(U,p\, ;z,s)]$, around $s=1$, 
satisfy a recurrence relation of order two (see Appendix \ref{app_4}).
\end{Rem}
\subsection{The Hilbert space $L^2(\mk{h};\Gamma;p)$}\label{Hil}
Let $\Gamma\leq SL_2(K)$ be a discrete subgroup commensurable to the Hilbert modular group $SL_2(\ca{O}_K)$ and let 
$p\in\ZZ^g$. We define
\begin{align*}
C^{\infty}(\mk{h}^g;\Gamma;p):=\left\{f:\mk{h}^g\rightarrow\CC|\;\mbox{$f$ is a smooth function such that $f\big|_{\{p\},\gamma}=f$}\right\},
\end{align*}
and
\begin{align*}
C_c^{\infty}(\mk{h}^g;\Gamma;p):=\left\{f\in C^{\infty}(\mk{h}^g;\Gamma;p)|\;\mbox{$f$ has compact support}\right\}.
\end{align*}
Both spaces are $\CC$-vector spaces with the obvious inclusion $C_c^{\infty}(\mk{h}^g;\Gamma;p)\subseteq
C^{\infty}(\mk{h}^g;\Gamma;p)$. 

Consider the space $Y_{\Gamma}:=\mk{h}^g/\Gamma$ endowed with the Poincar\'e metric (see Appendix \ref{app_5}).
The space $Y_{\Gamma}$ is a Riemannian orbifold. The $\CC$-vector $C^{\infty}(\mk{h}^g;\Gamma;p)$ may be viewed naturally as the set of global smooth sections of 
an appropriate complex line bundle $\ca{L}_p$ over $Y_{\Gamma}$. However, we won't use this point of view here.
For each $p\in\ZZ$, let $A_p:=C^{\infty}(\mk{h}^g;\Gamma;p)$. Note that vector space $A_{\bz}(Y_{\Gamma}):=C^{\infty}(\mk{h}^g;\Gamma;\bz)$ 
admits a ring structure with a unit (where the multiplication is the pointwise product of two functions), and the $\CC$-vector space
$A_p$ (for any $p\in\ZZ$) admits the structure of an $A_{\bz}(Y_{\Gamma})$-module. Note that one may work simultaneously with all integral weights $p\in\ZZ^g$
by considering the $\ZZ^g$-graded $A_{\bz}(Y_{\Gamma})$-module
\begin{align*}
A(Y_{\Gamma}):=\bigoplus_{p\in\ZZ^g} A_{p}(Y_{\Gamma}).
\end{align*}


We would like now to define an inner product on the vector spaces 
$C_c^{\infty}(\mk{h}^g;\Gamma;p)$.
\begin{Def}
For $f,g\in C_c^{\infty}(\mk{h}^g;\Gamma;p)$, we define the Petersson inner product of $f$ and $g$ as 
\begin{align}\label{doof2}
\laa f,g\raa:=\int_{Y_{\Gamma}} f(z)\ov{g(z)}dV<\infty,
\end{align}
where $dV$ is the Poincar\'e volume form normalized as in \eqref{dof} of Appendix \ref{app_5}.
\end{Def}
Note that the integral above is well-defined (and bounded), since $f(z)\ov{g(z)}$ may be viewed as 
a {\it smooth function} on $Y_{\Gamma}$ with compact support.
\begin{Def}
We let $L^2(\mk{h}^g;\Gamma;p)$ be the Hilbert space associated to the pre-Hilbert space $(C_c^{\infty}(\mk{h}^g;\Gamma;p),\laa\;,\;\raa)$.
\end{Def}

\subsubsection{Duality of Maa\ss\;-graded operators}\label{duali}

The next proposition is an easy generalization of Proposition 2.1.3 on p. 135 of \cite{Bump}.
\begin{Prop}
Let $f\in L^2(\mk{h}^g;\Gamma;p)\cap C^{\infty}(\mk{h}^g;\Gamma;p)$ and $g\in 
L^2(\mk{h}^g;\Gamma;p')\cap C^{\infty}(\mk{h}^g;\Gamma;p')$, 
where $p'=p+2e_j$. Here $e_j$ is the standard vector with $1$ in the $j$-th coordinate and $0$
elsewhere. Assume that $R_{j,p}(f)\in L^2(\mk{h}^g;\Gamma;p')$ and that $L_{j,p'}(g)\in L^2(\mk{h}^g;\Gamma;p)$. Then
\begin{align}\label{pho}
\laa R_{j,p_j}\;f,g\raa=\laa f,-L_{j,p_j+2}\;g\raa.
\end{align}
\end{Prop}
Note the presence of the $-$ sign on the right-hand side of the above equality. If we forget momentarily about the index and the weight, 
this result simply means that the left adjoint of $R$ (with respect to the Petersson inner product) is $-L$.

{\bf Proof} The proof is essentially the same as the proof of Proposition 2.1.3 of \cite{Bump}, except that one
needs to replace the one-dimensional complex space $\mk{h}$ by $\mk{h}^g$.
Let $j\in\{1,\ldots,g\}$ be the index which appears in the statement of the proposition.
Consider the degree $2g-1$ smooth differential form on $Y_{\Gamma}$
\begin{align*}
\eta:=\frac{2}{\ii}\cdot y_j^{-1}f(z)\ov{g(z)}\;\left(\frac{\ii}{2}\frac{dz_1\wedge d\ov{z}_1}{y_1^2}\right)\wedge \cdots \wedge
\left(\frac{\ii}{2}\frac{\wh{dz_j}\wedge d\ov{z}_j}{y_j}\right)\wedge\cdots \wedge
\left(\frac{\ii}{2}\frac{dz_g\wedge d\ov{z}_g}{y_g^2}\right),
\end{align*}
where the $\wh{dz_j}$ means that we omit this term. Then, using the analogue of Stokes' theorem for square-integrable
functions, we have
\begin{align*}
0=\int_{\partial Y_{\Gamma}} \eta=\int_{Y_{\Gamma}} d\eta =-\int_{Y_{\Gamma}} 
\left( (R_{j,p_j}f(z))\cdot\ov{g(z)}+f\cdot(\ov{L_{j,p_j+2}\;g(z)})\right)\; dV,
\end{align*}
where the symbol $d$ in $d\eta$ corresponds to the de Rham exterior derivative of $\eta$, and $dV$ corresponds
to the Poincar\'e volume form (see Appendix \ref{app_4}). But the equality above is equivalent to \eqref{pho}. \fin

\begin{Cor}\label{nof}
Let $f\in L^2(\mk{h}^g;\Gamma;p)\cap C^{\infty}(\mk{h}^g;\Gamma;p)$ and assume that there exists an index $j\in\{1,\ldots,g\}$ such that
$L_{j,p_j}f\in L^2(\mk{h}^g;\Gamma;{'}p)$, where ${'}p=p-2e_j$, and $\Delta_{j,p_j}f\in L^2(\mk{h}^g;\Gamma;p)$. Then
\begin{align*}
\laa \Delta_{j,p_j} f,f\raa\geq\;\frac{p_j}{2}\left(1-\frac{p_j}{2}\right)\cdot\laa f,f\raa. 
\end{align*}
\end{Cor}

{\bf Proof} We have
\begin{align*}
 \laa \Delta_{j,p_j} f,f\raa &=\laa-R_{j,p_j-2}L_{j,p_j}f,f\raa +\frac{p_j}{2}\left(1-\frac{p_j}{2}\right)\laa f,f\raa\\
 &=\laa L_{j,p_j}f,L_{j,p_j} f\raa+\frac{p_j}{2}\left(1-\frac{p_j}{2}\right)\laa f,f\raa\\
 &\geq \frac{p_j}{2}\left(1-\frac{p_j}{2}\right)\laa f,f\raa.
\end{align*}
For the second equality we have used the previously proven fact that $-L$ is the left adjoint of $R$. The result follows. \fin

\section{Families of real analytic modular forms of unitary weight $\{p;s\}$}\label{fam_eis}

\subsection{A qualitative description of the Fourier series coefficients}\label{qual}
Let $\ca{Q}=((\mk{m},\mk{n}),U,p,0)$ be a standard quadruple and let  $G_{\ca{Q}}(z,s)=G_{(\mk{m},\mk{n})}^0(U,p\,;z,s)$ be its associated
Eisenstein series. So far, we have proved that the function $[(z,s)\mapsto G_{(\mk{m},\mk{n})}^0(U,p\,;z,s)]$ satisfies the following 
three main properties:
\begin{enumerate}
 \item $[z\mapsto G_{(\mk{m},\mk{n})}^0(U,p\,;z,s)]$ is modular of unitary weight $\{-p;s\}$ with respect to 
 a discrete subgroup of $GL_2(K)$ (see Proposition \ref{mod_for}).
 \item $[z\mapsto G_{(\mk{m},\mk{n})}^0(U,p\,;z,s)]$ is an eigenvector with eigenvalue $s(1-s)$ for each of the partial
 graded Laplacians $\Delta_{j,-p_j}$, $j\in\{1,\ldots,g\}$ (see Section \ref{eigen1}).
 \item For $\Ree(s)>1$ fixed, the function $[z\mapsto G_{(\mk{m},\mk{n})}^0(U,p\,;z,s)]$ satisfies some growth conditions when 
 $z$ tends to a cusp (see Proposition \ref{gr_est}).
\end{enumerate}
For a real number $a\in\RR$, recall that
\begin{align*}
\Pi_a:=\{s\in\CC:\Ree(s)>a\}.
\end{align*}
As $s$ varies in $\Pi_1$, we view $\{[z\mapsto G_{(\mk{m},\mk{n})}^0(U,p\,;z,s)]\}_{s\in\Pi_1}$ as a family of real 
analytic Eisenstein series of unitary weight $\{-p;s\}$ with respect to a fixed discrete subgroup $\Gamma\leq GL_2(K)$. 

Let $s\in\Pi_1$ be fixed and let
\begin{align}\label{Four_ser}
a_0(y,s)+\sum_{\xi\in\mathcal{L}^*\bs\{0\}}  a_{\xi}(y,s)e^{2\pi\ii\Tr(\xi x)},
\end{align}
be the Fourier series expansion of the Eisenstein series $[z\mapsto G_{(\mk{m},\mk{n})}^0(U,p\,;z,s)]$ at the cusp $\infty$.
For $\xi\in\ca{L}^*\bs\{0\}$ fixed, we know, from Theorem \ref{keyy}, that $|a_{\xi}(y,s)|$ decays exponentially to zero as $\Tr(|\xi|)\rightarrow\infty$.
In this subsection, we explain how the growth estimate \eqref{fat0}
implies that the Fourier coefficient $a_{\xi}(y,s)$ also decays exponentially to zero, as $\Tr(|y|)\rightarrow\infty$. 
The last fact is proved by showing that the non-constant Fourier coefficient $a_{\xi}(y,s)$ must be of a very 
particular shape. Moreover, we also show that the constant Fourier coefficient $a_0(y,s)$ has a
very precise shape. We would like to emphasize that all these results are proved {\it indirectly},
in the sense that they don't use the explicit description of the Fourier coefficient $a_{\xi}(y,s)$, but use only
the properties (1), (2) and (3) stated above. 

The proposition below provides a {\it qualitative description} of the Fourier series coefficients
of  certain families of real analytic modular forms $\{F(z,s)\}_{s\in\Pi_1}$ which satisfy the properties (1), (2) and (3) above. In fact, on closer
inspection, only the properties (2), (3) and the existence of a Fourier series expansion as in \eqref{Four_ser} are required in the
Proof of Theorem \ref{staff}.
\begin{Th}\label{staff}
Let $T\subseteq \CC$ be a discrete subset. 
Let $F(z,s):\mk{h}^g\times(\Pi_1\bs T)\rightarrow \CC$ be a function,
such that $[(z,s)\mapsto F(z,s)]$ is real analytic, and $[s\mapsto F(z,s)]$ is holomorphic on $\Pi_1\bs T$. Furthermore, we assume that 
$F(z,s)$ satisfies the following 3 properties:
\begin{enumerate}[(i)]
 \item There exists a lattice $\mathcal{L}\subseteq\RR^g$, such that, for all $(z,s)\in \mk{h}^g\times(\Pi_1\bs T)$ and
 all $\ell\in\mathcal{L}$, $F(z+\ell,s)=F(z,s)$. 
 \item There exists a function $\phi_1:(\Pi_1\bs T)\rightarrow\CC$, such that, for any fixed value
 of $s\in\Pi_1$, there exists a constant $C_s>0$, such that for all $z\in \mk{h}^g$,
 \begin{align*}
 |F(z,s)-\phi_1(s)\cdot|\Norm(y)|^s|\leq C_s\cdot|\Norm(y)|^{1-s},
 \end{align*}
 \item There exists an integral weight $p=(p_j)_{j=1}^g\in\ZZ^g$, such that for $j\in\{1,\ldots,g\}$, and all
 $(z,s)\in\mk{h}^g\times(\Pi_1\bs T)$, one has that $\Delta_{j,p_j} F(z,s)=s(1-s) F(z,s)$.
\end{enumerate}
Let 
\begin{align*}
F(z,s)=a_0(y,s)+\sum_{\xi\in\mathcal{L}^*\bs\{0\}}  a_{\xi}(y,s)e^{2\pi\ii\Tr(\xi x)}
\end{align*}
be the Fourier series expansion of $F(z,s)$, which exists by (i). Then, 
the constant Fourier coefficient $a_0(y,s)$ is of the form 
\begin{align}\label{nos1}
a_0(y,s)=\phi_1(s)\cdot|\Norm(y)|^s+\phi_2(s)\cdot|\Norm(y)|^{1-s},
\end{align}
for some holomorphic functions $\phi_1(s)$ and $\phi_2(s)$ on $\Pi_1\bs T$. 

Let $M_p:=m_p+1$ where $m_p:=\max\left(\{\pm \frac{p_j}{2}:1\leq j\leq g\}\right)$. Note that by definition $M_p\geq 1$.
Then, for $\xi\in\ca{L}^*\bs\{0\}$, there exists a function  
$[s\mapsto c_{\xi}(s)]$, holomorphic on $\Pi_{M_p}\bs T$, such that
\begin{align}\label{nos2}
a_{\xi}(y,s)=c_{\xi}(s)\cdot B_{\xi}(y;p;s),
\end{align}
where $B_{\xi}(y;p;s)$ is the function which appears in equation \eqref{pers} of Appendix \ref{app_3}. 
\end{Th}

{\bf Proof} Let us start by proving \eqref{nos1}. By assumption, the function $[z\mapsto F(z,s)]$ is $\ca{L}$-invariant. 
Moreover, it is an eigenvector with respect to the graded Laplacians $\Delta_{j,p_j}$ (for $j\in\{1,\ldots,g\}$)
with eigenvalue $s(1-s)$. It follows from the discussion on the first page of Appendix \ref{app_3} that
\begin{align*}
a_0(y,s)=\sum_{\mu\in\{s,1-s\}^{g}} b_{\mu}(s)y^{\mu}, 
\end{align*}
where $\mu=(\mu_j)_{j=1}^g$ is a vector of length $g$ with $\mu_j\in\{s,1-s\}$ and $b_{\mu}(s)\in\CC$. Since
\begin{align*}
a_0(y,s)-\phi_1(s)|\Norm(y)|^s,
\end{align*}
is bounded as $|\Norm(y)|\rightarrow\infty$, this forces $b_{\mu}(s)\equiv 0$, if $\mu\notin\{s\cdot\bu,(1-s)\cdot\bu\}$
and $b_{s\cdot\bu}(s)=\phi_1(s)$. We set $\phi_2(s):=b_{(1-s)\cdot\bu}(s)$. It thus follows that
\begin{align*}
a_0(y,s)=\phi_1(s)\cdot|\Norm(y)|^s+\phi_2(s)\cdot|\Norm(y)|^{1-s}.
\end{align*}
Moreover, since
$[s\mapsto F(z,s)]$ is holomorphic on $\Pi_1\bs T$, it follows from Lemma \ref{dub} that
$\phi_1(s)$ and $\phi_2(s)$ are also holomorphic on $\Pi_1\bs T$.

Let us prove \eqref{nos2}. If we apply Parseval's theorem (see (2) of Theorem \ref{keyy}) to the periodic function
\begin{align*}
[x\mapsto (F(x+\ii y,s)-\phi_1(s)(\Norm y)^s)],
\end{align*}
and combine it with the the growth estimate given in (ii) of Theorem \ref{staff}, we may deduce that
for each $s\in\Pi_1\bs T$ fixed, there exists a constant $D_s>0$, such that, for all $z\in \mk{h}^g$,
\begin{align}\label{tol}
\sum_{\xi\in\mathcal{L}^*\bs\{0\}}  |a_{\xi}(y,s)|^2\leq D_s\cdot |\Norm(y)|^{2-2s}.
\end{align}
In particular, from \eqref{tol}, we deduce that, for a fixed $\xi\in\ca{L}^*\bs\{0\}$, 
the function $[y\mapsto |a_{\xi}(y,s)|]$ is bounded as $|\Norm(y)|\rightarrow\infty$.
Moreover, by (iii), $F(z,s)$ is an eigenvector with respect to the partial-graded Laplacians $\Delta_{j,p_j}$ 
(for $j\in\{1,\ldots,g\}$), with eigenvalue $s(1-s)$. Under these assumptions, we may apply Proposition \ref{congre} in Appendix \ref{app_3}
to obtain that 
\begin{align}\label{doss}
a_{\xi}(y,s)=c_{\xi}(s)\cdot B_{\xi}(y;p;s),
\end{align}
where $B_{\xi}(y;p;s)$ is the function which appears in equation \eqref{pers} of Appendix \ref{app_3}.
If $[y\mapsto a_{\xi}(y,s)]\equiv 0$, we take $c_{\xi}(s)$ to be identically equal to zero.
In this way, the function $s\mapsto c_{\xi}(s)$ is a holomorphic function on $\Pi_{M_p}\bs T$, and \eqref{doss} is satisfied.

Let us assume now that $[(y,s)\mapsto a_{\xi}(y,s)]\not\equiv 0$. We intend to show that there exists a
holomorphic function $c_{\xi}(s)$ on $\Pi_{M_p}\bs T$ such that $a_{\xi}(y,s)=c_{\xi}(s)\cdot B_{\xi}(y;p;s)$.
Note that such a function  $c_{\xi}(s)$ is necessarily not identically equal to zero.
Since $[s\mapsto F(z,s)]$ is holomorphic on $\Pi_1\bs T$, it follows from Lemma \ref{dub}
that, for each $y\in(\RR^{\times})^g$ fixed, the function $[s\mapsto a_{\xi}(y,s)]$ is holomorphic on $\Pi_1\bs T$.
For each $s_0\in \Pi_{M_p}\bs T$ fixed, let us choose $y_{s_0}=(y_{s_0,j})_{j=1}^g\in\mk{h}^g$, such that $\min_{j=1}^g y_{s_0,j}$ is large enough.
Then, from (5) of Proposition \ref{wee}, we know that $B_{\xi}(y_{s_0};p;s_0)\neq 0$, and, therefore, we may define
$c_{\xi}(s_0)=\frac{a_{\xi}(y_{s_0},s_0)}{B_{\xi}(y_{s_0};p;s_0)}$. Finally, for each fixed
$y\in\RR_{>0}^g$, the functions $[s\mapsto a_{\xi}(y,s)]$ and $[s\mapsto B_{\xi}(y;p;s)]$ are holomorphic on $\Pi_{M_p}\bs T$.
It thus follows that $[s\mapsto c_{\xi}(s)]$ is a holomorphic on $\Pi_{M_p}\bs T$. 
This concludes the proof. \fin

The next corollary provides some strong restrictions that the Fourier coefficients
of $G_{(\mk{m},\mk{n})}^w(U,p\,;z,s)$ must satisfy (even if $w\neq 0$).
\begin{Cor}\label{kcor00}
Let $\ca{Q}=((\mk{m},\mk{n}),U,p,w)$ be a standard quadruple. Consider the Eisenstein series 
$G_{(\mk{m},\mk{n})}^w(U,p\,;z,s)$ associated to $\ca{Q}$. We view $G_{(\mk{m},\mk{n})}^w(U,p\,;z,s)$ as
a function on $\mk{h}^g\times\Pi_{1-\frac{w}{2}}$, where $z\in \mk{h}^g$ and $s\in\Pi_{1-\frac{w}{2}}$.
The constant term of the Fourier series expansion of $G_{(\mk{m},\mk{n})}^w(U,p\,;z,s)$ at $\infty$ is of the form:
\begin{align}\label{coeur}
\phi_1(s)\cdot|\Norm(y)|^s+\phi_2(s)\cdot|\Norm(y)|^{1-s-w},
\end{align}
where $\phi_1(s)$ and $\phi_2(s)$ are holomorphic functions on $\Pi_{1-\frac{w}{2}}$. Let
$\ca{L}\subseteq K$ be a lattice, such that $[x\mapsto G_{(\mk{m},\mk{n})}^w(U,p\,;z,s)]$ is
invariant under $\ca{L}$. Let $\xi\in\mathcal{L}^*\bs\{0\}$ and let $M_p=m_p+1$ where 
$m_p:=\max_{j=1}^g\left\{\pm \frac{p_j}{2}\right\}$. Then, for each $s\in\Pi_{M_p}$, there
exists a positive constant $C_s>0$, such that 
\begin{align}\label{coeur2}
|a_{\xi}(y,s)|\leq C_s\cdot \frac{|\Norm(\xi)|^{2\Ree(s)+w-1}\cdot e^{-6||\xi y||_{\infty}}}{|\Norm(y)|^{w/2}},
\end{align}
for all $y\in(\RR^{\times})^g$. Here $||\_||_{\infty}$ corresponds to the $\ell_{\infty}$-norm of a vector in $\RR^g$.
\end{Cor}

{\bf Proof} First note that by making use of the identity \eqref{tram}, we easily see that Corollary \ref{kcor00} follows from the special case when 
$w=0$. Since $w=0$, we may thus apply Theorem \ref{staff} to the family of real analytic modular forms 
$\{G_{(\mk{m},\mk{n})}^0(U,p\,;z,s)\}_{s\in\Pi_1}$. Now, the equality \eqref{coeur} follows from Theorem \ref{staff}, and the
inequality \eqref{coeur2} follows Theorem \ref{staff} and Corollary \ref{tab} in Appendix \ref{app_3}.  \fin


\subsection{Existence of a non-zero Fourier coefficient}\label{exis}
In this section, we show that a non-zero modular function $F:\mk{h}^g\times(\Pi_1\bs T)\rightarrow\CC$
of unitary weight $\{p;s\}$, which satisfies all the assumptions of
Theorem \ref{staff}, must have at least one non-constant Fourier coefficient which is not identically equal to zero, i.e., 
there must exist $\xi\in\ca{L}^*\bs\{0\}$, such that $[(y,s)\mapsto a_{\xi}(y,s)]\not\equiv 0$. Moreover,
we show that if $a_{\xi}(y_0,s_0)\neq 0$, for some $(y_0,s_0)\in \mk{h}^g\times(\Pi_1\bs T)$ and $\xi\in\ca{L}^*\bs\{0\}$,
then $a_{\xi}(y,s)$ is not vanishing on a tiny vertical strip of $\mk{h}^g\times(\Pi_1\bs T)$ which contains a small neighborhood
of $s_0$.

\begin{Prop}\label{cliff}
Let $F(z,s)$ be a function as in Theorem \ref{staff}. Assume, furthermore, that $[z\mapsto F(z,s)]$
has unitary weight $p$ with respect to a discrete subgroup $\Gamma\leq SL_2(K)$ commensurable to
$SL_2(\ca{O}_K)$. Assume that $[(z,s)\mapsto F(z,s)]\not\equiv 0$ and let
\begin{align*}
a_0(y,s)+\sum_{\xi\in\mathcal{L}^*\bs\{0\}}  a_{\xi}(y,s)e^{2\pi\ii\Tr(\xi x)},
\end{align*}
be the Fourier series expansion of $[z\mapsto F(z,s)]$ at the cusp $\infty$. Let $M_p=m_p+1$ where $m_p:=\max_{j}\{\pm \frac{p_j}{2}\}$.
Then, given a real number $a\geq M_p$, 
there exists an $s_0\in \RR_{>a}\bs T$, such that the set 
\begin{align*}
\ca{T}_{s_0}:=\{\xi\in \ca{L}^{*}\bs\{0\}:\mbox{there exists $y\in(\RR^{\times})^g$ such that $a_{\xi}(y,s_0)\neq 0$}\},
\end{align*}
is non-empty. Moreover,  for each $\xi\in\ca{T}_{s_0}$, there exists a small open disc $D(s_0,\epsilon)=\{s\in\CC:|s-s_0|<\epsilon\}$, with $\epsilon>0$ 
(depending on $\xi$), such that, for all $y\in(\RR^{\times})^g$ and all $s\in D(s_0,\epsilon)$, $a_{\xi}(y,s_0)\neq 0$.
\end{Prop}

\begin{Rem}
In the case where $g>1$, it follows directly from the identity \eqref{nazy} that if $\ca{T}_{s_0}\neq \emptyset$, then 
necessarily $\#\ca{T}_{s_0}=\infty$.
\end{Rem}

{\bf Proof}\; We do a proof by contradiction. Assume that there exists a real number $a\geq M_p$,
such that for all pairs $(\xi,s)\in(\ca{L}^{*}\bs\{0\})\times(\RR_{>a}\bs T)$, 
the function $[y\mapsto a_{\xi}(y,s)]\equiv 0$. From \eqref{nos1} of Theorem \ref{staff}, we know that, for each 
$s\in\RR_{>a}\bs T$,
\begin{align*}
F(z,s)=a_0(y,s)=\phi_1(s)\cdot|\Norm(y)|^{s}+\phi_2(s)\cdot|\Norm(y)|^{1-s},
\end{align*}
where  $\phi_1(s)$ and $\phi_2(s)$ are holomorphic functions on $\Pi_a\bs T$.
By assumption, $F(z,s)$ has integral unitary weight $p$ relative to a suitable
congruence subgroup $\Gamma\leq SL_2(K)$. The function $|\Norm(y)|^s$ (resp. $|\Norm(y)|^{1-s}$ ) 
is of unitary weight $-s\cdot\bu$ (resp. is of unitary weight $-(1-s))\cdot\bu$.
Since for all $s\in\RR_{>a}$, we have $s\neq (1-s)$, it follows that $[s\mapsto\phi_1(s)]\big|_{\RR_{>a}\bs T}\equiv 0$ 
and $[s\mapsto\phi_2(s)]\big|_{\RR_{>a}\bs T}\equiv 0$. Since $\phi_1(s)$ and $\phi_2(s)$ are holomorphic on $\Pi_{1}\bs T$,
it follows, from the identity principle, that $[s\mapsto\phi_1(s)]\equiv 0$ and $[s\mapsto \phi_2(s)]\equiv 0$ on
$\Pi_{1}\bs T$. But then $[(z,s)\rightarrow F(z,s)]\equiv 0$, which is a contradiction.

From the previous paragraph, given $a\geq M_p$, there must exist $s_0\in (\RR_{>a}\bs T)$, 
such that $[y\mapsto a_{\xi}(y,s_0)]\not\equiv 0$. Now choose arbitrarily $\xi_0\in\ca{T}_{s_0}$. Since 
$[y\mapsto a_{\xi_0}(y,s_0)]\not\equiv 0$, there exists $y_0\in(\RR^{\times})^g$,
such that $a_{\xi_0}(y_0,s_0)\neq 0$. From \eqref{nos2} of Theorem \ref{staff},
we know that $a_{\xi_0}(y,s)=c_{\xi_0}(s)\cdot B_{\xi_0}(y;p;s)$ for some holomorphic function 
$c_{\xi_0}(s)$ on $\Pi_{a}\bs T$. Since $a_{\xi_0}(y_0,s_0)\neq 0$, it follows that
$c_{\xi_0}(s_0)\neq 0$. Moreover, from (6) of Proposition \ref{wee},
we know that for all $y\in(\RR^{\times})^g$, $B_{\xi_0}(y;p;s_0)\neq 0$. By continuity of $c_{\xi}(s)$,
we may thus find a small open disc $D(s_0,\epsilon)$ (with $\epsilon>0$), such that for all $s\in D(s_0,\epsilon)$,
$c_{\xi}(s)\neq 0$. Therefore, it follows that, for all $s\in  D(s_0,\epsilon)$ and all $y\in(\RR^{\times})^g$,
that $a_{\xi_0}(y,s)\neq 0$. This concludes the proof. \fin

\subsection{Non-existence of certain square-integrable real analytic families of modular forms of unitrary weight $p$}\label{unicity}
The next proposition should be viewed as a complement to the results presented in Section \ref{constr}.
We let $\Gamma\leq SL_2(K)$ be a discrete subgroup commensurable to $SL_2(\ca{O}_K)$ and $Y_{\Gamma}:=\biglslant{\Gamma}{\mk{h}^g}$ be the associated orbifold.
\begin{Prop}\label{bard}
Let  $T\subseteq\RR_{>1}$ be a discrete subset. Let $G: \mk{h}^g\times(\RR_{>1}\bs T)\rightarrow\CC$
be a smooth function in $(z,s)$ which is real analytic in $s$. Let $p\in\ZZ^g$ and $M_p=m_p+1$ where $m_p=\max_{j=1}^g\{\pm\frac{p_j}{2}\}$.
Assume that the following three properties hold true:
\begin{enumerate}
 \item For each $s_0\in\RR_{>1}\bs T$, the function $[z\mapsto G(z,s_0)]$ is modular of unitary weight $p$ with respect to $\Gamma$.
 \item There exists on index $j\in\{1,\ldots,g\}$, such that $\Delta_{j,p} G(z,s)=s(1-s)G(z,s)$.
 \item The set 
\begin{align*}
\ca{D}&:=\Big\{s\in\RR_{>M_p}\bs T:[z\mapsto G(z,s)], [z\mapsto\Delta_{j,{p_j}}G(z,s)]\in L^2(\mk{h}^g;\Gamma;p),\\[2mm] 
&\hspace{2cm}\mbox{\;\;and\;\;} [z\mapsto L_{j,p_j}G(z,s)]\in L^2(\mk{h}^g;\Gamma;{'}_{}p)\Big\},
\end{align*}
is uncountable. Here ${'}_{}p=p-2e_j$ and the index $j$ which appears in the definition of $\ca{D}$ is the same as in (2).
\end{enumerate}
Then $[(z,s)\mapsto G(z,s)]\equiv 0$.
\end{Prop}

{\bf Proof} We do a proof by contradiction. Suppose that $[(z,s)\mapsto G(z,s)]\not\equiv 0$. 
By assumption, for each $s_0\in\ca{D}$, the functions $[z\mapsto G(z,s_0)]\in L^2(\mk{h}^g;\Gamma;p)$
$[z\mapsto \Delta_{j,p_j}G(z,s_0)]\in L^2(\mk{h}^g;\Gamma;p)$ and $[z\mapsto L_{j,p_j}G(z,s_0)]\in  L^2(\mk{h}^g;\Gamma;p)$. 
We claim that that this implies that $[z\mapsto G(z,s_0)]\equiv 0$. Indeed, assume that $s_0\in\ca{D}$ is such that 
$[z\mapsto G(z,s_0)]\not\equiv 0$.  Then, applying Corollary \ref{nof} to $[z\mapsto G(z,s_0)]$, we deduce that 
\begin{align*}
s_0(1-s_0)\geq p_j/2(1-p_j/2). 
\end{align*}
But this is impossible since $s_0\in\ca{D}$. We thus have proved that,
\begin{align}\label{yor}
\mbox{for all $s\in\ca{D}$, $[z\mapsto G(z,s)]\equiv 0$}. 
\end{align}
We claim that \eqref{yor} implies that $[(z,s)\rightarrow G(z,s)]\equiv 0$. Let
$z_0\in\mk{h}^g$ be an arbitrary fixed element. By assumption, the function $[s\mapsto G(z_0,s)]$
is real analytic on $\Pi_1\bs T$ and it vanishes on all of $\ca{D}$. But since
$\ca{D}$ is not discrete, it follows that $[s\mapsto G(z_0,s)]\equiv 0$. Finally, since
$z_0$ was arbitrary, it follows that $[(z,s)\mapsto G(z,s)]\equiv 0$. This concludes the proof. \fin

\subsection{An analytic characterization of certain families of real analytic modular forms of unitary weight $p$}\label{car}

The theorem below provides an analytic characterization of certain families of real analytic
modular forms of a fixed unitary weight $p$. In particular, this theorem may be applied to the real
analytic family $\{G_{(\mk{m},\mk{n})}^0(U,p\, ;z,s))\}_{s\in\Pi_1}$.
\begin{Th}\label{clop}
Let $\Gamma\leq SL_2(K)$ be a subgroup commensurable to $SL_2(\ca{O}_K)$. 
Let $\{\mk{c}_1,\ldots,\mk{c}_h\}$ be a complete set of representatives
of the relative $\Gamma$-cusps. Let $T\subseteq\Pi_1$ be a discrete subset, 
let $\{A_{i}(s)\}_{1\leq i\leq h}$ be a finite collection of holomorphic functions on $\Pi_1\bs T$, and
let $\{\gamma\}_{i=1}^h$ be a collection of matrices in $SL_2(K)$ such that such that for all $1\leq i\leq h$, $\gamma_i\infty=\mk{c}_i$. 
Then, there exists at most one function $G$ which satisfies the following four properties:
\begin{enumerate}
 \item $G:\mk{h}^g\times(\Pi_1\bs T)\rightarrow\CC$, $(z,s)\mapsto G(z,s)$, is smooth in $(z,s)$ and 
 $[s\mapsto G(z,s)]$ is holomorphic on $\Pi_1\bs T$. 
 \item There exists a weight $p=(p_j)_{j=1}^g\in\ZZ^g$, such that $[z\mapsto G(z,s)]$ is modular of unitary weight $p$ relative
 to $\Gamma$.
 \item  For each $j\in\{1,2,\ldots,g\}$, we have that 
 $\Delta_{j,p_j}G(z,s)=s(1-s)G(z,s)$.
 \item  For each $1\leq i\leq h$, and $s\in\Pi_1\bs T$ fixed, we have
\begin{align}\label{cir}
\Big| G\big|_{\{p\},\gamma_{i}}(z,s)- A_{i}(s)\Norm(y)^{s}\Big|=O(\Norm(y)^{1-s}),
\end{align}
as $\Norm(y)\rightarrow\infty$.
\end{enumerate}
In condition (4), note that the constant term of the Fourier series expansion at the cusp $\mk{c}_i$
does not depend on the choice of the local chart $\gamma_{i}\in SL_2(K)$. 
\end{Th}

{\bf Proof} Let $G_1(z,s)$ and $G_2(z,s)$ be two functions which satisfy (1), (2), (3) and (4).
Set $G(z,s):=G_1(z,s)-G_2(z,s)$. We want to show that $[(z,s)\mapsto G(z,s)]\equiv 0$. First notice that,
because of condition (4), for any $s\in\Pi_1\bs T$ fixed, $[z\mapsto G(z,s)]\in L^2(\mk{h}^g;\Gamma;p)$, i.e., it is square-integrable 
on $Y_\Gamma$. Let
\begin{align*}
G(z,s)=\sum_{\xi\in\ca{L}^{*}} a_{\xi}(y,s)e^{2\pi\ii\Tr(\xi x)},
\end{align*}
be the Fourier series expansion of $G(z,s)$ at $\infty$, where $\ca{L}\subseteq K$ is a suitable lattice.
Let $M_p=m_p+1$ where $m_p=\max_{j}\{\pm \frac{p_j}{2}\}$. For each $\xi\in\ca{L}^*\bs\{0\}$, 
it follows from Theorem \ref{staff} that 
\begin{align*}
a_{\xi}(y,s)=c_{\xi}(s)\cdot B_{\xi}(y;p;s),
\end{align*}
where the function $c_{\xi}(s)$ is holomorphic on $\Pi_{M_p}\bs T$. For each $s\in\Pi_{M_p}\bs T$ fixed, 
it follows from (10) of Proposition \ref{wee} that $[y\mapsto\frac{\partial}{\partial y}a_{\xi}(y,s)]$ decays
exponentially to $0$ as $||y||_{\infty}\rightarrow \infty$. Therefore, for all $j\in\{1,\ldots,g\}$,
\begin{align}\label{pech}
\mbox{ $[z\mapsto \frac{\partial}{\partial y_j} G(z,s)]$
decays exponentially to $0$ as $||y||_{\infty}\rightarrow\infty$}.
\end{align}
Combining \eqref{pech} and the identity \eqref{ore3} of Proposition \ref{dub},
we may deduce that, for all $j\in\{1,\ldots,g\}$, $[z\mapsto L_{j,p_j}G(z,s)]\in L^2(Y_{\Gamma},{'}_{}p)$ where
${'}_{}p=p-2e_j$. Let us fix, arbitrarily, an index $j_0\in\{1,\ldots,g\}$. Using Proposition \ref{cliff}, we may deduce that, 
for each $b\geq M_p$, the set
\begin{align*}
\ca{D}_b:=\{s\in\RR_{>b}\bs T:[z\mapsto G(z,s)]\in L^2(\mk{h}^g;\Gamma;p) \mbox{\;\;and\;\;} [z\mapsto L_{j,p_j}G(z,s)]\in L^2(\mk{h}^g;\Gamma;{'}_{}p)\},
\end{align*}
is uncountable. Note that if $s\in\ca{D}_b$, then automatically $[z\mapsto \Delta_{j_0,p_{j_0}}G(z,s)]\in L^2(\mk{h}^g;\Gamma;p)$,
since $\Delta_{j_0,p_{j_0}}G(z,s)=s(1-s)G(z,s)$. Finally, applying Proposition \ref{bard}
to $G(z,s)$, we conclude that $[(z,s)\mapsto G(z,s)]\equiv 0$. \fin

\begin{Rem}
Using the ideas of Colin de Verdi\`ere in \cite{Ver1}, it is explained in \cite{Ch-Kl} 
how Theorem \ref{clop} above may be used to give a  different proof of the meromorphic continuation 
and the functional equation of the Eisenstein $[s\mapsto G_{(\mk{m},\mk{n})}^0(U,p;z,s)]$. As a consequence of this approach, 
we obtain a new proof of the analytic continuation of the zeta function $[s\mapsto Z_{\mk{m}}(u,v;\omega_{\ov{p}};s)]$, 
which circumvents the use of the Poisson summation formula for lattices in euclidean spaces.
The analytic continuation is proved instead using the Fredholm analyticity theorem for a suitable family of compact operators. 
\end{Rem}

\subsection{The $\spartsing{c}{s}$ and the $\spartsing{c}{1-s}$ of a family of real analytic modular forms}\label{party0}
From the results proved in the previous sections, it should be clear now that the collection of 
constant terms of the Fourier series expansion of a real analytic family of modular forms characterizes the family. 
In this section, we formalize this observation by introducing the $\spartsing{c}{s}$ and the $\spartsing{c}{1-s}$
to any family of modular forms $\{G(z,s)\}_{s\in\Pi_1\bs T}$ which satisfy the same assumptions as the ones studied in the previous 
sections.
\begin{Def}\label{party}
Let $T\subseteq\Pi_1$ be a discrete subset and let $\Gamma\leq SL_2(K)$ be a discrete
subgroup commensurable to $SL_2(\ca{O}_K)$. Let $\{G(z,s)\}_{s\in\Pi_1\bs T}$ be a real analytic
family of modular form of bi-weight $[\alpha,\beta]$ where $\alpha,\beta\in\CC^g$ and $\alpha-\beta=\ZZ^g$.
For each cusp $c\in\PP^1(K)$, choose a local chart $\eta\in SL_2(K)$ at $c$. Let
\begin{align*}
a_0^{\eta}(y,s)=\phi_1^{\eta}(s)\cdot|\Norm(y)|^s+\phi_2^{\eta}(s)\cdot|\Norm(y)|^{1-s},
\end{align*}
be the constant term of the Fourier series expansion of $[z\mapsto G(z,s)]$ at the cusp $c$. 
Note that $a_0^{\eta}(y,s)$ really does depend on the choice of the local chart at $\eta$ at $c$ (see Remark \ref{depen} below). We define 
\begin{enumerate}
 \item $\spart{c}{s} G(z,s)=\spart{c}{s} G:=[s\mapsto\phi_1^{\eta}(s)]$,
 \item $\spart{c}{1-s} G(z,s)=\spart{c}{1-s} G:=[s\mapsto\phi_2^{\eta}(s)]$.
\end{enumerate}
\end{Def}
We view the expressions  $\phi_1^{\eta}(s)=\spart{c}{s} G$ and $\phi_2^{\eta}(s)=\spart{c}{1-s} G$ as functions in the variable
$s$ where $s\in\Pi_1\bs T$. In particular, it follows that for $s_0\in\Pi_1\bs T$,
$\spart{c}{s} G\big|_{s_0}=\phi_1^{\eta}(s_0)$ and $\spart{c}{1-s} G\big|_{s_0}=\phi_2^{\eta}(s_0)$.
\begin{Rem}\label{depen}
If $\eta,\wt{\eta}\in SL_2(K)$ are two local charts at $c$, then there must exists $\delta=\M{a}{b}{0}{a^{-1}}\in \Stab_{SL_2(K)}(\infty)$ such that 
$\eta\delta=\wt{\eta}$. From this, it follows that
\begin{enumerate}[(i)]
 \item  $\phi_1^{\eta}(s)\cdot (a^{\alpha}a^{\beta})\cdot|\Norm(a)|^{2s}=\phi_1^{\wt{\eta}}(s)$,
 \item  $\phi_2^{\eta}(s)\cdot (a^{\alpha}a^{\beta})\cdot|\Norm(a)|^{2(1-s)}=\phi_2^{\wt{\eta}}(s)$.
\end{enumerate}
Even though the functions $\spart{c}{s} G$ and $\spart{c}{s} G$ do depend on the choice of the local chart at the cusp $c$,
we see from (i) and (ii) that the dependence is really mild. In practice, the slightly ambiguous notation $\spart{c}{s} G$ creates no difficulty. 
\end{Rem}
\begin{Rem}
Given $s_0\in\Pi_1\bs T$ which is fixed, one should think of the complex number
$\spart{c}{s} G|_{s_0}\in\CC$ as the {\it obstruction} for the modular form $[z\mapsto G(z,s_0)]$ to be square-integrable in a neighborhood 
of the cusp $c$ in the sense of Definition \ref{tok} (cf. with Remark \ref{ex_sq}). In particular, if $\spart{c}{s} G|_{s_0}=0$, 
then the modular form $[z\mapsto G(z,s_0)]$ is square-integrable at $c$.
\end{Rem}
\begin{Rem}\label{dual_fam}
In the special case where the family $\{G(z,s)\}_{s\in\Pi_1\bs T}$ admits a meromorphic continuation (in the variable $s$) 
to all of $\CC$, one may define its {\it dual family} as
\begin{align*}
\{G^*(z,s)\}_{s\in-\Pi_0\bs T^*},
\end{align*}
where
$G^*(z,s):=G(z,1-s)$ and 
\begin{align*}
T^*=\{s\in-\Pi_0:\mbox{there exists $z\in\mk{h}^g$, such that $[s\mapsto G(z,s)]$ has a pole}\}.
\end{align*}
It follows directly from the definition of $G^*(z,s)$ that for all $c\in\PP^1(K)$,
\begin{enumerate}
 \item $\spart{c}{s} G^*=\spart{c}{1-s} G$,
 \item $\spart{c}{1-s} G^*=\spart{c}{s} G$.
\end{enumerate}
\end{Rem}

\subsection{Writing $G_{(\mk{m},\mk{n})}^{0}(U,p\,;z,s)$ as a linear sum of classical real analytic Poincar\'e-Eisenstein series}\label{writing_app_7}
In this section, we explain how Theorem \ref{clop} (the analytic characterization theorem) can be used to write 
the Eisenstein series $G_{\ca{Q}}(z,s)=G_{(\mk{m},\mk{n})}^0(U,p\,;z,s)$ as a sum of {\it classical real analytic Poincar\'e-Eisenstein series}
of unitary weight $-p$. Here the \lq\lq coefficients\rq\rq\; of this sum will be prescribed by $\{\spart{c}{s} G_{\ca{Q}}\}_{c}$, where 
$c$ goes over a complete set of representatives of relative $\Gamma$-cusps, where $\Gamma\leq SL_2(K)$ is a suitable discrete subgroup.

Let $\Gamma\subseteq SL_2(K)$ be a subgroup commensurable to $SL_2(\ca{O}_K)$. Let $Y_{\Gamma}=\biglslant{\Gamma}{\mk{h}^g}$ be 
the associated Riemannian orbifold. For each cusp $c\in\PP^1(K)$, we define
\begin{align*}
\Gamma_c^+:=\{\gamma\in \Gamma:\gamma c=c\;\;\mbox{and the eigenvalues of $\gamma$ are totally positive}\}. 
\end{align*}
One may check that $\Gamma_c^+$ is a subgroup. If $c_1,c_2\in\PP^1(K)$ and $\eta\in SL_2(K)$ is such that $\eta c_1=c_2$, then
one may also check $\eta\Gamma_{c_1}^+\eta^{-1}=\Gamma_{c_2}^+$. 

Let $c\in\PP^1(K)$ and choose arbitrarily a local chart at $c$, i.e., a matrix $\sigma\in SL_2(K)$
such that $\sigma\infty=c$. Let also $p\in\ZZ^g$ be a fixed integral weight.
For $z\in\mk{h}^g$ and $s\in\Pi_1$, we define
\begin{align}\label{banu}
\wt{E}_{\Gamma,\sigma}(p;z,s):=\sum_{\gamma\in\left(\biglslant{\Gamma_{c}^+}{\Gamma}\right)} \omega_p(j(\sigma^{-1}\gamma, z))\cdot\Imm(\sigma^{-1}\gamma z)^s,
\end{align}
where the summation is understood to be taken over a complete set of representatives of the right cosets $\biglslant{\Gamma_{c}^+}{\Gamma}$.
One may check that $\wt{E}_{\Gamma,\sigma}(p;z,s)$ does not depend on the choice of representatives of $\biglslant{\Gamma_{c}^+}{\Gamma}$.
However, it does depend on the choice of the local chart at $c$ in a simple way. If $\sigma'$ is another local chart at $c$, then there must exists
a matrix $\gamma=\M{a}{b}{0}{a^{-1}}\in\Stab_{SL_2(K)}(\infty)$ such that $\sigma\gamma=\sigma'$. A direct calculation shows that
\begin{align}\label{changem}
\omega_p(a)\cdot|\Norm(a)|^{-2s}\cdot\wt{E}_{\Gamma,\sigma}(p;z,s)=\wt{E}_{\Gamma,\sigma'}(p;z,s).
\end{align}
In order to simplify the presentation, we allow ourself to use the slightly ambiguous notation 
$\wt{E}_{\Gamma,c}(p;z,s)$ to mean $\wt{E}_{\Gamma,\sigma}(p;z,s)$. The reader should simply keep in mind that the
expression $\wt{E}_{\Gamma,c}(p;z,s)$ really does depend on a choice of a local chart at the cusp $c$ and not 
just on the cusp $c$ itself. Moreover, a change of local chart at a given cusp simply changes the function as is described in \eqref{changem}.

It follows from \eqref{banu} and \eqref{abeil} that, for all $\gamma\in \Gamma$, one has
\begin{align}\label{cartif}
\wt{E}_{\Gamma,c}\big|_{\{-p\},\gamma}(p;z,s)=\wt{E}_{\Gamma,c}(p;z,s).
\end{align}
Therefore, $[z\mapsto\wt{E}_{\Gamma,c}(p;z,s)]$ is a real analytic modular form of unitary weight $-p$ relative to $\Gamma$.

Let $c,c'\in\PP^1(K)$ be two $\Gamma$-equivalent cusps. In other words, there exists $\eta\in\Gamma$ such that $\eta c=c'$. We claim
that $\wt{E}_{\Gamma,c}(p;z,s)=\wt{E}_{\Gamma,c'}(p;z,s)$. Let us prove it. Let $\{\gamma_j\}_{j\in J}$ is a complete set
of representatives of $\biglslant{\Gamma_{c}^+}{\Gamma}$. Then $\{\eta\gamma_j\eta^{-1}\}_{j\in J}$ is 
a complete set of representatives of $\biglslant{\Gamma_{c'}^+}{\Gamma}$. We have
\begin{align*}
\wt{E}_{\Gamma,c'}(p;z,s)&=\sum_{j\in J} \omega_p(j(\sigma'^{-1}\eta\gamma_j\eta^{-1}, z))\cdot\Imm(\sigma'^{-1}\eta\gamma_j\eta^{-1} z)^s\\
&=\sum_{j\in J} \omega_p(j(\sigma^{-1}\gamma_j\eta^{-1}, z))\cdot\Imm(\sigma^{-1}\gamma_j\eta^{-1} z)^s\\
&=\sum_{j\in J} \omega_p(j(\sigma^{-1}\gamma_j,\eta^{-1} z))\cdot \omega_p(j(\eta^{-1}, z))\cdot\Imm(\sigma^{-1}\gamma_j\eta^{-1} z)^s\\
&=\wt{E}_{\Gamma,c}\big|_{\{-p\},\eta^{-1}}(p;z,s)=\wt{E}_{\Gamma,c}(p;z,s),
\end{align*}
where the third equality follows from \eqref{abeil} and the fifth equality from \eqref{cartif}. 

It follows from Theorem,
\ref{staff} that the constant term of the Fourier series of $\wt{E}_{\Gamma,c}(p;z,s)$ at the cusp $d$ is of the form
\begin{align*}
\phi_{1,c,d}(s)\Norm(y)^s+\phi_{2,c,d}(s)\Norm(y)^{1-s},
\end{align*}
where $\phi_{1,c,d}(s)$ and $\phi_{2,c,d}(s)$ are holomorphic functions on $\Pi_1$. If $[d]_{\Gamma}\neq [c]_{\Gamma}$, 
it follows from the definition of $\wt{E}_{\Gamma,c}(p;z,s)$ and an analogue of inequality \eqref{fat0} that 
the function $\phi_{1,c,d}(s)\equiv 0$. In particular, for all fixed $s\in\Pi_1$ and all cusp $d\in\PP^1(K)$
such that $[d]_{\Gamma}\neq [c]_{\Gamma}$, the function $[z\mapsto \wt{E}_{\Gamma,c}(p;z,s)]$ is {\it square integrable at the cusp
$d$}.
\begin{Rem}
Using the techniques in \cite{Ch-Kl}, it is possible to prove that $[s\mapsto \wt{E}_{\Gamma,c}(p;z,s)]$ admits a meromorphic continuation
to all of $\CC$ and that it satisfies a functional equation.
\end{Rem}

Let $\{c_1,\ldots,c_h\}\subseteq\PP^1(K)$ be a complete set of representatives of the relative $\Gamma$-cusps.
We make the following assumption:
\begin{Assumption}\label{assump}
For all $i\in\{1,\ldots,h\}$, $[z\mapsto\wt{E}_{\Gamma,c_i}(p;z,s)]\not\equiv 0$.
\end{Assumption}
It follows from the above discussion
that $[s\mapsto \phi_{1,c_i,c_i}(s)]\not\equiv 0$. Indeed, if $[s\mapsto \phi_{1,c_i,c_i}(s)]\equiv 0$,
then, for $i\in\{1,\ldots,h\}$ fixed, the two distinct families of Eisenstein series
\begin{enumerate}
 \item $\{0\}_{s\in\Pi_1}$ (the trivial family),
 \item $\{\wt{E}_{\Gamma,c_i}(p;z,s)]\}_{s\in\Pi_1}$, 
\end{enumerate}
satisfy the properties (1), (2), (3) and (4) of Theorem \ref{clop}. But this is absurd since it contradicts  
the uniqueness statement of Theorem \ref{clop}. Therefore, it makes sense to define the {\it characteristic Eisenstein series
at the relative cusp $[c_i]_{\Gamma}$} as
\begin{align*}
E_{\Gamma,c_i}(p;z,s):=\frac{\wt{E}_{\Gamma,c_i}(p;z,s)}{\phi_{1,c_i,c_i}(s)}.
\end{align*}

So, by definition,
\begin{enumerate}[(i)]
 \item $\spart{c_i}{s}E_{\Gamma,c_i}(p;z,s)=1$,
 \item $\spart{c_j}{s}E_{\Gamma,c_i}(p;z,s)=0$ for $j\neq i$.
\end{enumerate}
Now let $\ca{Q}=((\mk{m},\mk{n}),U,p,0)$ be a standard quadruple and let $G_{(\mk{m},\mk{n})}^0(U,p\,;z,s)$
be its associated real analytic Eisenstein series. Choose an integer $N\in\ZZ_{\geq 1}$, such that $U\in\frac{1}{N}\M{\mk{m}^*}{\mk{m}}{\mk{n}^*}{\mk{n}}$.
Let $\Gamma:=\Gamma_U(\mk{m},\mk{n};N)\leq GL_2(K)$ (see Definition \ref{grp}). We know from Proposition \ref{mod_for}
that $G_{(\mk{m},\mk{n})}^0(U,p\,;z,s)$ is a modular form of unitary weight $-p$ relative to $\Gamma$.
Let $\{c_1,\ldots,c_h\}\subseteq\PP^1(K)$ be a complete set of representatives of relative $\Gamma$-cusps. For 
each $i\in\{1,\ldots,h\}$, the constant term of the Fourier series expansion of $[z\mapsto G_{(\mk{m},\mk{n})}^0(U,p\,;z,s)]$ at 
the cusp $c_i$ is of the form
\begin{align*}
A_{c_i}(s)\Norm(y)^s+ B_{c_i}(s)\Norm(y)^{1-s},
\end{align*}
for some holomorphic functions $A_{c_i}(s)$ and $B_{c_i}(s)$ on $\Pi_1$. We claim that
\begin{align}\label{identit}
G_{(\mk{m},\mk{n})}^0(U,p\,;z,s)=\sum_{j=1}^h A_{c_j}(s)\cdot E_{\Gamma,c_j}(p;z,s).
\end{align}
Indeed, if we let
\begin{align*}
F(U,p\;;z,s):=G_{(\mk{m},\mk{n})}^0(U,p\,;z,s)-\sum_{j=1}^h A_{c_j}(s)\cdot E_{\Gamma,c_j}(p;z,s),
\end{align*}
then $F(U,p\; ;,z,s)$ is square-integrable at all the cusps $c_j$'s. Therefore, from the uniqueness statement of 
Theorem \ref{clop}, we deduce that $F(U,p;\,z,s)\equiv 0$.
\begin{Rem}
Note that all the $A_{c_j}(s)$'s appearing in \eqref{identit} are explicitly given by the formulas appearing in Theorem \ref{key_thmm}. Therefore, the right-hand
side of \eqref{identit} may be viewed as an explicit
writing of $G_{(\mk{m},\mk{n})}^0(U,p\,;z,s)$ as a sum of real analytic Poincar\'e-Eisenstein series 
$\{E_{\Gamma,c_j}(p;z,s)\}_{j=1}^h$.
\end{Rem}
\begin{Rem}
Let us give a special situation where Assumption \ref{assump} can be easily verified. 
Let us assume that $h_K=1$ (the class number of $K$) and that  $\Gamma\unlhd SL_2(\ca{O}_K)$. 
In this case, it is easy to see that
\begin{align*}
[z\mapsto \wt{E}_{\Gamma,\infty}(p;z,s)]\not\equiv 0\Longrightarrow \left\{\mbox{ $[z\mapsto \wt{E}_{\Gamma,c_i}(p;z,s)]\not\equiv 0$ for all $i\in\{1,\ldots,h\}$}\right\}.
\end{align*}
Indeed, since $h_K=1$, for each $i\in\{1,\ldots,h\}$, there exists $\sigma_i\in SL_2(\ca{O}_K)$, such that $\sigma_i\infty=c_i$. Finally, by the normality of 
$\Gamma$ in $SL_2(\ca{O}_K)$, one may simply notice that
\begin{align*}
\wt{E}_{\Gamma,\infty}\big|_{\{-p\},\sigma_i^{-1}}(p;z,s)=\wt{E}_{\Gamma,c_i}(p;z,s).
\end{align*}
\end{Rem}

\section{Fourier series expansion and meromorphic continuation of $G_{(\mk{m},\mk{n})}^{0}(U,p\,;z,s)$}

\subsection{Explicit Fourier series expansion at the cusp $\infty$}\label{nico0}
Let $\ca{Q}=((\mk{m},\mk{n}),U,p,s)$ be a standard quadruple and let $G_{\ca{Q}}(z;s)=G_{(\mk{m},\mk{n})}^w(U,p;z,s)=
G_{(\mk{m},\mk{n})}^{\alpha(s),\beta(s)}(U,z)$ its associated Eisenstein series. Recall that
$\alpha(s),\beta(s)\in\CC^g$ are the two weight vectors defined as in Definition \ref{difo}. 
In order to simplify the notation, we shall write sometimes $\alpha$ (resp. $\beta$) rather
than $\alpha(s)$ (resp. $\beta(s)$), and similarly, for $j\in\{1,\ldots,g\}$, $\alpha_j$ (resp $\beta_j$)
rather than $\alpha_j(s)$ (resp. $\beta_j(s)$).

In this section, we compute explicitly the Fourier series expansion of 
\begin{align*}
[x\mapsto G_{(\mk{m},\mk{n})}^{\alpha(s),\beta(s)}(U\,;x+\ii y)]
\end{align*} 
at the cusp $\infty=\frac{1}{0}$. Note that, from the transformation formula \eqref{chat}, one may obtain the Fourier series expansion of 
$G_{(\mk{m},\mk{n})}^{\alpha(s),\beta(s)}(U;z)$ at any other cusp.
\begin{Th}\label{key_thmm}
Let $U=\M{u_1}{v_1}{u_2}{v_2}$ be a parameter matrix and let $z\in K_{\CC}^{\pm}$. Assume that 
$s\in\Pi_{1-\frac{w}{2}}$. Then Fourier series expansion of $[x\mapsto G_{(\mk{m},\mk{n})}^{\alpha,\beta}(U\,;x+\ii y)]$ is given by
\begin{align}\label{eqn_a} 
&G_{(\mk{m},\mk{n})}^{\alpha,\beta}(U\,;z)=(T_1+T_2+T_3)\cdot|y|^{s\cdot\bu},
\end{align}
where 
\begin{align*}
T_1=\delta_{\mk{m}}(v_1)\sum_{\mathcal{V}^+\bs\{0\neq(n+v_2)\in\mk{n}+v_2\}}
\frac{\omega_{\ov{p}+w\cdot\ov{\bu}}(n+v_2)\cdot e^{{2\pi\ii }\Tr(u_2(n+v_2))}}{\prod_{i=1}^g |n^{(i)}+v_2^{(i)}|^{\alpha_i+\beta_i}},
\end{align*}
\begin{align*}
&T_2=\cov(\mk{n})^{-1}\cdot \delta_{\mk{n}^*}(u_2)\cdot e^{{2\pi\ii } \Tr(u_2 v_2)}(2\pi)^g(\ii)^{\Tr(p-w\cdot\bu)}\cdot (-1)^{\Tr(\sg(z)\cdot (p-w\cdot\bu))}\cdot
\prod_{j=1}^{g}(2|y_j|)^{1-\alpha_j-\beta_j}\\ 
&\hspace{2cm}\cdot\prod_{j=1}^g\frac{\Gamma(\alpha_j+\beta_j-1)}{\Gamma(\alpha_j)\Gamma(\beta_j)}
\sum_{\mathcal{V}^+\bs\{0\neq(m+v_1)\in\mk{m}+v_1\}}
\frac{\omega_{\ov{p}+w\cdot\ov{\bu}}(m+v_1)\cdot e^{{2\pi\ii } \Tr(u_1(m+v_1))}}{\prod_{j=1}^g |m^{(j)}+v_1^{(j)}|^{\alpha_j+\beta_j-1}},
\end{align*}
and 
\begin{align*}
&T_3=\cov(\mk{n})^{-1}\cdot e^{{2\pi\ii } \Tr(u_2 v_2)}
\sum_{\mathcal{V}^+\bs\{0\neq\xi_1\in\mk{m}+v_1\}}
e^{{2\pi\ii } \Tr(u_1\xi_1)}\\  \notag
&\hspace{3cm}\cdot\sum_{0\neq\xi_2\in (\mk{n}^*-u_2)}\left(\prod_{j=1}^g
\tau\left(\alpha_j,\beta_j,\xi_2^{(j)},\xi_1^{(j)}y_j\right)\right)
e^{2\pi\ii\Tr(\xi_2 v_2)}\cdot e^{{2\pi\ii } \Tr\left(\xi_2\xi_1 x\right)}.
\end{align*}
Here,
\begin{enumerate}[(i)]
 \item within the expression of $T_1$, $\delta_{\mk{m}}(v_1)=1$ if $v_1\in\mk{m}$ and $0$ otherwise,
 \item within the expression of $T_2$, $\delta_{\mk{n}^*}(u_2)=1$ if $u_2\in\mk{n}^*$ and $0$ otherwise,
 \item $\mathcal{V}^+=\mathcal{V}_U^+(\mk{m},\mk{n})$ is the subgroup of totally positive units of $\ca{O}_K$
which appears in Definition \ref{fromage}.
\end{enumerate}
Note that since $2s+w-1>1$, the series $T_1$ and $T_2$ are absolutely convergent.
\end{Th}
Before starting the proof, let us point out that, since $s\in\Pi_{1-\frac{w}{2}}$, the defining 
series of $G_{(\mk{m},\mk{n})}^{\alpha(s),\beta(s)}(U\,;z)$ 
converges absolutely. Therefore, all sum rearrangements appearing in the proof below are justified. Moreover, since
$[z\mapsto G_{(\mk{m},\mk{n})}^{\alpha,\beta}(U\,;z)]$ is real analytic, it follows from (1) of Theorem \ref{keyy},
that the Fourier series of $[x\mapsto G_{(\mk{m},\mk{n})}^{\alpha,\beta}(U\,;x+\ii y)]$ converges absolutely
to the value $G_{(\mk{m},\mk{n})}^{\alpha,\beta}(U\,;z)$ for all $z\in K_{\CC}^{\pm}$. 
\begin{Rem}
We would like to explain how to see directly that the $T_3$ term is independent of the set of representatives of
$\{0\neq\xi_1\in\mk{m}+v_1\}$ modulo $\mathcal{V}^+$. Let us replace $\xi_1$ by $\epsilon\xi_1$ in the $T_3$ term where 
$\epsilon\in\ca{V}^+$. Note that $e^{{2\pi\ii } \Tr(u_1\xi_1)}=e^{{2\pi\ii } \Tr(u_1\epsilon\xi_1)}$ and 
that
\begin{align*}
\prod_{j=1}^g
\tau\left(\alpha_j,\beta_j,\xi_2^{(j)},\epsilon^{(j)}\xi_1^{(j)}y_j\right)=
\prod_{j=1}^g|\epsilon^{(j)}|^{1-\alpha_j-\beta_j}\cdot\tau\left(\alpha_j,\beta_j,\epsilon^{(j)}\xi_2^{(j)},\xi_1^{(j)}y_j\right).
\end{align*}
Since, for all $j\in\{1,\ldots,g\}$, $\alpha_j+\beta_j=2s+w$ is {\it independent} of $j$, and 
$\Norm(\epsilon)=1$, we thus see that, the $T_3$ term is, indeed, independent from the choice of representatives 
of $\biglslant{\mathcal{V}^+}{\{0\neq\xi_1\in\mk{m}+v_1\}}$.
\end{Rem}

{\bf Proof} We have
\begin{align}\label{first_eq} \notag
\frac{G_{(\mk{m},\mk{n})}^{\alpha,\beta}(U\,;z)}{|y|^{s\cdot\bu}} 
&=\sum_{\mathcal{V}^+\bs\{(0,0)\neq(m+v_1,n+v_2)\in(\mk{m}+v_1,\mk{n}+v_2)\}}
\frac{e^{{2\pi\ii } \Tr(u_1(m+v_1)+u_2(n+v_2))} }{ P(\alpha,\beta;(m+v_1) z+(n+v_2))}\\ 
&=\delta_{\mk{m}}(v_1)\sum_{\mathcal{V}^+\bs\{0\neq(n+v_2)\in\mk{n}+v_2\}}
\frac{e^{{2\pi\ii }\Tr(u_2(n+v_2))}}{P(\alpha,\beta;n+v_2)}\\ \notag
&\hspace{1cm}+\sum_{\mathcal{V}^+\bs\{0\neq(m+v_1)\in\mk{m}+v_1\}}\sum_{n\in\mk{n}}
\frac{e^{{2\pi\ii } \Tr(u_1(m+v_1)+u_2(n+v_2))}}{P(\alpha,\beta;(m+v_1) z+(n+v_2))},
\end{align}
where $P(\alpha,\beta;z)$ is the monomial (in the variables $z$ and $\ov{z}$) which appears in \eqref{expo}.
We note that the summation in the first equality mentioned above has been divided in two parts. The first part
corresponds to the terms indexed by the pairs $(m+v_1,n+v_2)$, such that $m=-v_1$, 
and the second part corresponds to the terms indexed by the pairs $(m+v_1,n+v_2)$, such that $m\neq -v_1$. From Proposition \ref{koy_prop} and
the above computation, we deduce that
\begin{align*}
T_1=\delta_{\mk{m}}(v_1)\sum_{\mathcal{V}^+\bs\{0\neq(n+v_2)\in\mk{n}+v_2\}}
\frac{e^{{2\pi\ii }\Tr(u_2(n+v_2))}}{P(\alpha,\beta;n+v_2)}.
\end{align*}

It remains to prove the formulas for $T_2$ and $T_3$. 
Let us concentrate on the double summation in \eqref{first_eq}. Let us choose arbitrarily a writing
for $u_2=\frac{e}{h}$, where $e\in\mk{n}^*$ and $h\in\ca{O}_K\bs\{0\}$, with $h\gg 0$. We have
\begin{align}\label{noz}\notag
&\sum_{\mathcal{V}^+\bs\{0\neq(m+v_1)\in\mk{m}+v_1\}}\sum_{n\in\mk{n}}
\frac{e^{{2\pi\ii } \Tr(u_1(m+v_1)+u_2(n+v_2))}}{P(\alpha,\beta;(m+v_1) z+(n+v_2))}\\ \notag
&=\sum_{\mathcal{V}^+\bs\{0\neq(m+v_1)\in\mk{m}+v_1\}}\hspace{0.3cm}\sum_{\mu \pmod{\frac{\mk{n}}{h\mk{n}}}}
e^{2\pi\ii \Tr(u_2(\mu+v_2))}
\sum_{\substack{n\in\mk{n}\\n\equiv\mu \pmod{h\mk{n}}}} \frac{e^{{2\pi\ii } \Tr(u_1(m+v_1))}}{P(\alpha,\beta;(m+v_1) z+(n+v_2))}\\
&=\sum_{\mathcal{V}^+\bs\{0\neq(m+v_1)\in\mk{m}+v_1\}}\sum_{\mu \pmod{\frac{\mk{n}}{h\mk{n}}}}e^{2\pi\ii \Tr(u_2(\mu+v_2))}
\sum_{n\in\mk{n}}
\frac{e^{{2\pi\ii } \Tr(u_1(m+v_1))}}{P(\alpha,\beta;(m+v_1) z+(hn+\mk{\mu}+v_2))}.
\end{align}

Now, let us compute the Fourier series expansion of 
\begin{align*}
\Big[z\mapsto \sum\limits_{n\in\mk{n}}\frac{e^{{2\pi\ii } \Tr(u_1(m+v_1))}}{P(\alpha,\beta;(m+v_1) z+(hn+\mk{\mu}+v_2))}\Big].
\end{align*}
We have
\begin{align*}
& \frac{1}{P(\alpha,\beta;(m+v_1)z+(hn+\mu+v_2))}=\left(\prod_{j=1}^g (h^{(j)})^{-\alpha_j-\beta_j}\right)\cdot P\left(-\alpha,-\beta;
\left(\frac{m+v_1}{h}\right)z+n+\frac{(\mu+v_2)}{h}\right).
\end{align*}
Therefore, it follows that
\begin{align}\label{kilo}\notag
&\sum_{n\in\mk{n}}P(-\alpha,-\beta;(m+v_1)z+(hn+\mu+v_2))\\  \notag
&=\left(\prod_{j=1}^g (h^{(j)})^{-\alpha_j-\beta_j}\right)\sum_{n\in\mk{n}} P\left(-\alpha,-\beta;
\left(\frac{m+v_1}{h}\right)z+n+\frac{(\mu+v_2)}{h}\right)\\ \notag
&=\left(\prod_{j=1}^g (h^{(j)})^{-\alpha_j-\beta_j}\right)\cdot R_{\mk{n}}\left(\alpha,\beta;\left(\frac{m+v_1}{h}\right)z+\frac{(\mu+v_2)}{h}\right)\\ \notag
&=\left(\prod_{j=1}^g (h^{(j)})^{-\alpha_j-\beta_j}\right)\sum_{\xi\in\mk{n}^*} a_{\xi}\left(\alpha,\beta,\frac{(m+v_1)y}{h}\right)
e^{{2\pi\ii } \Tr(\xi(\frac{(m+v_1)}{h}x+\frac{\mu+v_2}{h}))}\\ 
&=\cov(\mk{n})^{-1}\left(\prod_{j=1}^g (h^{(j)})^{-\alpha_j-\beta_j}\right)
\sum_{\xi\in\mk{n}^*} \left(\prod_{j=1}^g \tau\left(\alpha_j,\beta_j,\xi^{(j)},\frac{(m+v_1)y_j}{h^{(j)}}\right)\right)
e^{{2\pi\ii } \Tr\left(\xi\left(\frac{(m+v_1)}{h}x+\frac{\mu+v_2}{h}\right)\right)}.
\end{align}
In the second equality, the function $R_{\mk{n}}(\alpha,\beta;z)$ is the one appearing in  \eqref{pap} which corresponds to the summation
of the monomial $P(\alpha,\beta;z)$ over the lattice $\mk{n}$. The third equality stated above follows from \eqref{subi}, and the fourth equality follows from \eqref{coeff}.
The function $\tau(\_,\_,\_,\_)$ is the function which appears in \eqref{nezz}.

Now substituting the expression on the right-hand side of \eqref{kilo} into the right-hand side of \eqref{noz}, and re-indexing the summation, we find
\begin{align*}
&\sum_{\mathcal{V}^+\bs\{0\neq(m+v_1)\in\mk{m}+v_1\}}\sum_{n\in\mk{n}}
\frac{e^{{2\pi\ii } \Tr(u_1(m+v_1)+u_2(n+v_2))}}{P(\alpha,\beta;(m+v_1) z+(n+v_2))}\\
&=\prod_{j=1}^g 
(h^{(j)})^{-\alpha_j-\beta_j}\sum_{\mathcal{V}^+\bs\{0\neq(m+v_1)\in\mk{m}+v_1\}}e^{{2\pi\ii } \Tr(u_1(m+v_1))}
\sum_{\mu\pmod{\frac{\mk{n}}{h\mk{n}}}}e^{{2\pi\ii }\Tr(\frac{e}{h}(\mu+v_2))}\\
&\hspace{1cm}\cdot \cov(\mk{n})^{-1}\sum_{\xi\in\mk{n}^*}\prod_{j=1}^g
\tau\left(\alpha_j,\beta_j,\xi^{(j)},(m^{(j)}+v_1^{(j)})\frac{y_j}{h^{(j)}}\right)
e^{{2\pi\ii } \Tr\left(\xi\left(\frac{(m+v_1)}{h}x+\frac{\mu+v_2}{h}\right)\right)}.
\end{align*}
The right-hand side of the above equality can be broken into two parts: the first part corresponding
to the terms for which $\xi=0$, and the second part corresponding to the terms for which $\xi\neq 0$. 
Therefore, we have
\begin{align}\label{chev}
&\cov(\mk{n})\sum_{\mathcal{V}^+\bs\{0\neq(m+v_1)\in\mk{m}+v_1\}}\sum_{n\in\mk{n}}
\frac{e^{{2\pi\ii } \Tr(u_1(m+v_1)+u_2(n+v_2))}}{P(\alpha,\beta;(m+v_1) z+(n+v_2))}=\wt{T}_2+\wt{T_3},
\end{align}
where
\begin{align*}
\wt{T}_2&=\prod_{j=1}^g (h^{(j)})^{-\alpha_j-\beta_j}\sum_{\mathcal{V}^+\bs\{0\neq(m+v_1)\in\mk{m}+v_1\}}
\Big(\sum_{\mu\pmod{\frac{\mk{n}}{h\mk{n}}}}e^{{2\pi\ii } \Tr(u_1(m+v_1))}e^{{2\pi\ii }\Tr(u_2(\mu+v_2))}\cdot\\
&\hspace{3cm}\prod_{j=1}^g\tau\left(\alpha_j,\beta_j,0,(m^{(j)}+v_1^{(j)})\frac{y_j}{h}\right)\Big),\\ 
\end{align*}
and
\begin{align*}
\wt{T}_3&=\prod_{j=1}^g (h^{(j)})^{-\alpha_j-\beta_j}\sum_{\mathcal{V}^+\bs\{0\neq(m+v_1)\in\mk{m}+v_1\}}
e^{{2\pi\ii } \Tr(u_1(m+v_1))}
\sum_{\mu\pmod{\frac{\mk{n}}{h\mk{n}}}}e^{{2\pi\ii } \Tr(u_2(\mu+v_2))}\\ 
&\hspace{1cm}\cdot\sum_{0\neq\xi\in\mk{n}^*}\prod_{j=1}^g
\tau\left(\alpha_j,\beta_j,\xi^{(j)},(m^{(j)}+v_1^{(j)})\frac{y_j}{h^{(j)}}\right)
e^{{2\pi\ii } \Tr\left(\xi\left(\frac{(m+v_1)}{h}x+\frac{\mu+v_2}{h}\right)\right)}.
\end{align*}
Using the explicit formula for $\tau\left(\alpha_j,\beta_j,0,(m^{(j)}+v_1^{(j)})\frac{y_j}{h}\right)$ appearing
in Lemma \ref{Shim_lem} and the well-known trigonometric identity (see for example Lemma 6.6 on p. 474 of \cite{Neu})
\begin{align}\label{piff}
\sum_{\mu\in\frac{\mk{n}}{h\mk{n}}}e^{{2\pi\ii }\Tr(x\mu)}=
\left\{
\begin{array}{cl}
|\Norm(h)| & \mbox{if}\s\s x\in\mk{n}^*,\\ 
  0      &  \mbox{if}\s\s x\in\frac{1}{h}\mk{n}^*\bs\mk{n}^*,
\end{array}\right.
\end{align}
we see that the $\wt{T}_2$ term may be rewritten as:
\begin{align}\label{chi}
& \wt{T}_2=\delta_{\mk{n}^*}(u_2)\cdot e^{{2\pi\ii } \Tr(u_2 v_2)}(2\pi)^{g}(\ii)^{\Tr(\beta-\alpha)}\cdot\prod_{j=1}^g (h^{(j)})^{-\alpha_j-\beta_j}\cdot
\prod_{j=1}^g\frac{\Gamma(\alpha_j+\beta_j-1)}{\Gamma(\alpha_j)\Gamma(\beta_j)}\\ \notag
&\cdot\sum_{\mathcal{V}^+\bs\{0\neq(m+v_1)\in\mk{m}+v_1\}}
(-1)^{\Tr(\sg((m+v_1)\frac{y}{h})\cdot (\beta-\alpha))}\cdot e^{{2\pi\ii } \Tr(u_1(m+v_1))}
\prod_{j=1}^g \left|2 (m^{(j)}+v_1^{(j)})\frac{y_j}{h^{(j)}}\right|^{1-\alpha_j-\beta_j}.
\end{align}
After some simplifications, one finds that
\begin{align*}
\wt{T}_2 &=\delta_{\mk{n}^*}(u_2)\cdot e^{{2\pi\ii } \Tr(u_2 v_2)}(2\pi)^{g}(\ii)^{\Tr(\beta-\alpha)}\cdot (-1)^{\Tr(\sg(z)\cdot (\beta-\alpha))}\cdot
\prod_{j=1}^g (2|y_j|)^{1-\alpha_j-\beta_j}
\prod_{j=1}^g\frac{\Gamma(\alpha_j+\beta_j-1)}{\Gamma(\alpha_j)\Gamma(\beta_j)}\\ \notag
&\cdot\sum_{\mathcal{V}^+\bs\{0\neq(m+v_1)\in\mk{m}+v_1\}} e^{{2\pi\ii }\Tr(u_1(m+v_1))}\cdot\omega_{\ov{p}+w\cdot\ov{\bu}}(m+v_1)\cdot
\prod_{j=1}^g \left|m^{(j)}+v_1^{(j)}\right|^{1-\alpha_j-\beta_j}\\[2mm]
&=T_2.
\end{align*}
Note that, up to the sign $(-1)^{\Tr(\sg(z)\cdot (\beta-\alpha))}$, the expression for $T_2$ depends only on $|y|$ rather than $y$ itself.

Let us now handle the $\wt{T}_3$ term. Using Lemma \ref{Shim_lem}, the $\wt{T}_3$ term may be rewritten as:
\begin{align}\label{fin_eq0}
\wt{T}_3&=e^{2\pi\ii\Tr(u_2v_2)}\cdot\prod_{j=1}^g (h^{(j)})^{-\alpha_j-\beta_j}\sum_{\mathcal{V}^+\bs\{0\neq(m+v_1)\in\mk{m}+v_1\}}
e^{{2\pi\ii } \Tr(u_1(m+v_1))}\\ \notag
&\hspace{-1cm}\cdot\sum_{0\neq\xi\in\mk{n}^*}
\left(\prod_{j=1}^g
\tau(\alpha_j,\beta_j,\xi^{(j)},(m^{(j)}+v_1^{(j)})\frac{y_j}{h^{(j)}})\right)
e^{2\pi\ii\Tr(\frac{\xi}{h}(m+v_1)x)}\cdot e^{{2\pi\ii } \Tr(\frac{\xi}{h}v_2)}\sum_{\mu\pmod{\frac{\mk{n}}{h\mk{n}}}}
e^{2\pi\ii\Tr(\left(\frac{\xi}{h}+u_2\right)\mu)}.
\end{align}

Note that for $\mu\in\mk{n}$ and $\xi\in\mk{n}^*$, we have $\frac{\xi}{h}+u_2\in\frac{1}{h}\mk{n}^*$.
Using again \eqref{piff}, we see that the last summation of \eqref{fin_eq0} is zero unless $\frac{\xi}{h}+u_2\in\mk{n}^*$, 
which is equivalent to say that $\xi\in (\mk{n}^*h-hu_2)\subseteq\mk{n}^*$. 

Moreover, using (iii) of Section \ref{formu}, we have that
\begin{align*}
(h^{(j)})^{1-\alpha_j-\beta_j}\cdot\tau\left(\alpha_j,\beta_j,\xi^{(j)},(m^{(j)}+v_1^{(j)})\frac{y_j}{h^{(j)}}\right)=
\tau\left(\alpha_j,\beta_j,\frac{\xi^{(j)}}{h^{(j)}},(m^{(j)}+v_1^{(j)})y_j\right).
\end{align*}
Therefore, we find that
\begin{align}\label{fin_eq1}
\wt{T}_3&=e^{2\pi\ii\Tr(u_2v_2)}\sum_{\mathcal{V}^+\bs\{0\neq(m+v_1)\in\mk{m}+v_1\}}
e^{{2\pi\ii } \Tr(u_1(m+v_1))}\\ \notag
&\hspace{1cm}\cdot\sum_{0\neq \xi\in h\cdot\mk{n}^*-h \cdot u_2}
\left(\prod_{j=1}^g
\tau(\alpha_j,\beta_j,\frac{\xi^{(j)}}{h^{(j)}},(m^{(j)}+v_1^{(j)})y_j)\right)
e^{{2\pi\ii\Tr(\frac{\xi}{h}v_2) }}\cdot e^{{2\pi\ii } \Tr\left(\frac{\xi}{h}\left((m+v_1)x\right)\right)}.
\end{align}
Finally, combining \eqref{first_eq}, \eqref{chi} and \eqref{fin_eq1}, we see that $\wt{T}_3=T_3$. \fin

\subsubsection{Rewriting the constant term in terms of $\varphi_p(s)$ and $\psi_p(s)$}\label{nico}

Let $\alpha(s),\beta(s)$ be as in Definition \ref{difo} and assume that $w=0$.
It follows from Theorem \ref{key_thmm} that the constant term of the Fourier series expansion 
of $[x\mapsto G_{(\mk{m},\mk{n})}^{\alpha(s),\beta(s)}(U\,;x+\ii y)]$ is given explicitly by
\begin{align}\label{rue}
 \delta_{\mk{m}}(v_1)\phi_1(s)\Norm(y)^s+\delta_{\mk{n}^*}(u_2)\phi_2(s)\Norm(y)^{1-s},
\end{align}
where
\begin{align}\label{hibou0}
\phi_1(s):=e_1
\cdot \wt{Z}_{\mk{n}}(v_2,u_2,\omega_{\ov{p}};2s),
\end{align}
and 
\begin{align}\label{hibou}
\phi_2(s):= &\cov(\mk{n})^{-1}e^{{2\pi\ii }\Tr(u_2 v_2)}(\ii)^{\Tr(p)}(2\pi)^{g}\cdot 2^{g(1-2s)}\cdot(-1)^{\Tr(\sg(z)\cdot p)} \\[2mm] \notag
& \hspace{3cm} \left(\prod_{j=1}^g\frac{\Gamma(2s-1))}{\Gamma(s-p_j/2)\Gamma(s+p_j/2)}\right)\cdot
e_2\cdot\wt{Z}_{\mk{m}}(v_1,u_1,\omega_{\ov{p}};2s-1).
\end{align}
Note that the functions $\phi_1(s)$ and $\phi_2(s)$ do not depend on the variable $y$. Here:
\begin{enumerate}[(i)]
\item $\wt{Z}_{\mk{n}}(a,b,\omega_{\ov{p}};s):=\frac{Z_{\mk{n}}(a,b;\omega_{\ov{p}};s)}{[\ca{O}_K:\mk{n}]^s}$, where
$Z_{\mk{n}}(a,b;\omega_{\ov{p}};s)$ is the zeta function which appears in Appendix \ref{app_1b},
\item $e_1=e_1((\mk{m},\mk{n});U)$ and $e_2=e_2((\mk{m},\mk{n});U)$ are the indices which appear in Definition \ref{chien},
 \end{enumerate}
Let $\varphi_{p}(s),\psi_{p}(s)$ be the two functions which appear in Definition \ref{key_defi}. A direct
calculation shows that
\begin{align}\label{habou}
\varphi_p(1-s)\cdot\psi_p(s)=\left(\prod_{j=1}^g\frac{\Gamma(2s-1))}{\Gamma(s-p_j/2)\Gamma(s+p_j/2)}\right).
\end{align}
Substituting the left-hand 
side of \eqref{habou} in the right-hand side of \eqref{hibou}, we find that
\begin{align}\label{habo}
\phi_2(s):= &\cov(\mk{n})^{-1}e^{{2\pi\ii }\Tr(u_2 v_2)}(2\pi)^{g}\cdot 2^{-2g(s-\frac{1}{2})}\cdot \\[2mm] \notag
& \hspace{3cm} (\ii)^{\Tr(p)}\cdot(-1)^{\Tr(\sg(z)\cdot p)}\cdot\varphi_p(1-s)\cdot\psi_p(s)\cdot e_2\cdot\wt{Z}_{\mk{m}}(v_1,u_1,\omega_{\ov{p}};2s-1).
\end{align}

\subsubsection{Rewriting the $T_3$ term as a standard Fourier series}\label{nuz}
If one looks at the $T_3$ term of Theorem \ref{key_thmm}, it involves exponential terms of the form $e^{2\pi\ii\Tr(\xi_1\xi_2)}$, where
$\xi_1$ and $\xi_2$ are elements lying in a translate of a certain lattice. In general, given an element $d=\xi_1\xi_2$, 
there will be many pairs $(\xi_1',\xi_2')$ such that $\xi_1'\xi_2'=d$.
The goal of this subsection is to rewrite the $T_3$ term in Theorem \ref{key_thmm} as a sum over such $d$'s, i.e., 
as a {\it standard Fourier series}. 

We recall some notations. We let $\ca{Q}:=((\mk{m},\mk{n}),U,p,w)$ be
a standard quadruple and we let $\alpha(s),\beta(s)\in\CC^g$ be the associated weights, i.e.,
$\alpha(s)=(s+w)\cdot\bu-p/2\in\CC^g$ and $\beta(s)=s+p/2\in\CC^g$ for $s\in\CC$. We also let $\ca{V}^+=\ca{V}_U^+(\mk{m},\mk{n})
\leq\ca{O}_K^+(\infty)$ be the subgroup, which appears in Definition \ref{fromage}. 
We now define two sets that will be indexing the summation of a rewriting of the $T_3$ term. We define
\begin{enumerate}
 \item $D:=\{d\in K\bs\{0\}: \exists\s \xi_1\in(\mk{m}+v_1)\s\mbox{and}\s\xi_2\in (\mk{n}^*-u_2)
\s\mbox{such that}\s d=\xi_1\xi_2\}$.
 \item For each $d\in D$, we let $R_d:=\{(\xi_1,\xi_2)\in (\mk{m}+v_1)\times(\mk{n}^*-u_2):\xi_1\xi_2=d\}$.
\end{enumerate}
Since $\mk{n}^*-u_2$ and $\mk{m}+v_1$ are stable subsets under the action of $\ca{V}^+$, it makes sense to define a
{\it twisted diagonal action} $*$ of $\ca{V}^+$ on the set $R_d$ in the following way: for $\epsilon\in\ca{V}^+$, and
$(\xi_1,\xi_2)\in R_d$, we let $\epsilon*(\xi_1,\xi_2)=(\epsilon\xi_1,\epsilon^{-1}\xi_2)$. From now on, 
whenever we consider the quotient set $R_d/\mathcal{V}^+$, it will always be understood that the quotient will be taken
with respect to the $*$-action.

For $x\in\RR^{\times}$, recall that $\sg(x)=\ov{0}\in\ZZ/2\ZZ$ if $x>0$, and $\sg(x)=\ov{1}\in\ZZ/2\ZZ$ if $x<0$. 
It is convenient for us to make the following definition
\begin{Def}\label{roty}
Let $p\in\ZZ^g$ be a fixed weight. For $\xi\in K^{\times}$, we define 
\begin{align*}
\ca{W}_p(\xi):=-\sum\limits_{\substack{1\leq j\leq g\\ \xi^{(j)}<0  }} p_j.
\end{align*}
\end{Def}
Note that, by definition of $\ca{W}_p$, for any $\xi\in K^{\times}$, we have
\begin{align}\label{tarte}
(-1)^{\ca{W}_p(\xi)}=(-1)^{\Tr(p\cdot\sg(\xi))}.
\end{align}
The next proposition gives the writing of the $T_3$ term as a standard Fourier series.
\begin{Prop}\label{oily}
The  $T_3$ term in Theorem \ref{key_thmm} can be rewritten as
\begin{align}\label{tim00}
T_3&=\cov(\mk{n})^{-1}\cdot e^{2\pi\ii\Tr(u_2 v_2)}\cdot
\sum_{d\in D} B_{\xi}(y;p;s)
\cdot b_{d}(s)\cdot e^{{2\pi\ii } \Tr(dx))}\cdot|\Norm(y)|^{s},
\end{align}
where
\begin{align}\label{coffi0}
B_{d}(y;p;s):=\prod_{j=1}^g \tau\left(\alpha_j(s),\beta_j(s),1,d^{(j)} y_j\right),
\end{align}
and
\begin{align}\label{coffi}
b_{d}(s)= \left(\sum_{(\xi_1,\xi_2)\in R_d/\mathcal{V}^+}(-1)^{\ca{W}_{p+w\cdot\bu}(\xi_2)}|\Norm(\xi_2)|^{2s-1+w}e^{{2\pi\ii } \Tr(u_1\xi_1)}\cdot 
e^{{2\pi\ii } \Tr(\xi_2 v_2)}\right).
\end{align}
\end{Prop}

{\bf Proof} From the equations (iii),(iv) and (v) of Section \ref{formu}, we may deduce that, if $\xi_2\neq 0$, 
the following identity holds true:
\begin{align}\label{gori}
\tau\left(\alpha_j(s),\beta_j(s);\xi_2^{(j)},\xi_1^{(j)} y_j\right)=
(-1)^{\epsilon_j}|\xi_2^{(j)}|^{\alpha_j(s)+\beta_j(s)-1}\tau\left(\alpha_j(s),\beta_j(s);1,d^{(j)}y_j\right),
\end{align}
where $\epsilon_j=\sg(d^{(j)})$. Now, substituting \eqref{gori} in the $T_3$ term of Theorem
\ref{key_thmm} we find that
\begin{align*}
T_3&=\cov(\mk{n})^{-1}\cdot e^{{2\pi\ii } \Tr(u_2 v_2)}
\sum_{\mathcal{V}^+\bs\{0\neq\xi_1\in\mk{m}+v_1\}}
e^{{2\pi\ii } \Tr(u_1\xi_1)}\\  \notag
&\cdot\sum_{0\neq\xi_2\in (\mk{n}^*-u_2)}
(-1)^{\ca{W}_{p+w\cdot\bu}(\xi_2)}|\Norm(\xi_2)|^{2s-1+w}\left(\prod_{j=1}^g
\tau\left(\alpha_j,\beta_j,1,\xi_1^{(j)}\xi_2^{(j)}y_j\right)\right)
e^{2\pi\ii\Tr(\xi_2 v_2)}\cdot e^{{2\pi\ii } \Tr\left(\xi_2\xi_1 x\right)}.
\end{align*}
Let $\{\xi_{1r}\}_{r\in\ca{R}}$ be a complete set of representatives of  $\{0\neq\xi_1\in\mk{m}+v_1\}$ modulo $\ca{V}^+$.
By definition of $D$, for each $d\in D$, there exists at least one index $r_d\in\ca{R}$, such 
$\xi_{1r_d}\xi_{2r_d}=d$ for some $\xi_{2r_d}\in\mk{n}^*-u_2$. Note that $\xi_{2r_d}$ is uniquely determined by
the pair $(d,\xi_{1r_d})$. If the representative $\xi_{1r}$ is replaced by the representative
$\epsilon\xi_{1r}$, for some $\epsilon\in\ca{V}^+$, then this has the effect of replacing 
$\xi_{2r_{d}}$ by $\epsilon^{-1}\xi_{2r_d}$. From these observations, \eqref{tim00} follows. \fin

\begin{Rem}
Note that if there exists a unit $\epsilon\in\ca{V}:=\ca{V}_U(\mk{m},\mk{n})$, such that $\Norm(\epsilon)=-1$, then the coefficient
$b_{d}(s)$ in \eqref{coffi} is identically equal to zero.
\end{Rem}

\subsection{Meromorphic continuation of $[s\mapsto G_{(\mk{m},\mk{n})}^{0}(U,p\,;z,s)]$}\label{lun}
\begin{Th}\label{lune}
The function $[s\mapsto G_{(\mk{m},\mk{n})}^{0}(U,p\,;z,s)]$ admits a
meromorphic continuation to all of $\CC$ with poles of order at most one at
$s\in\{\frac{1}{2},1\}$. More precisely, there exist holomorphic functions $f_1,f_2:\CC\rightarrow\CC$,
and a real analytic function $R:K_{\CC}^{\pm}\times\CC\rightarrow\CC$, holomorphic in $s\in\CC$, 
such that, for all $s\in\CC$ and all $z\in K_{\CC}^{\pm}$, we have
\begin{align}\label{poum0}
&G_{(\mk{m},\mk{n})}^{\alpha(s),\beta(s)}(U\,;z)\\[1mm] \notag
&=\frac{f_1(s)}{(s-\frac{1}{2})^{\alpha_1}}
\cdot|\Norm(y)|^{s}+\frac{f_2(s)}{(s-\frac{1}{2})^{\alpha_2}(s-1)^{\beta}}|\Norm(y)|^{1-s} +R(z,s)\\[1mm] \notag
&=\phi_1(s)|\Norm(y)|^{s}+\phi_2(s)|\Norm(y)|^{1-s} +R(z,s),
\end{align}
where $\alpha_1=\alpha_2=\beta=1$, $\phi_1(s)=\frac{f_1(s)}{(s-\frac{1}{2})^{\alpha_1}}$ and
$\phi_2(s)=\frac{f_2(s)}{(s-\frac{1}{2})^{\alpha_2}(s-1)^{\beta}}$. Here, the function $R$ satisfies the following additional property:
for each  $s\in\CC$ fixed, it decreases exponentially to $0$ as $||y||_\infty\rightarrow\infty$. Furthermore,
\begin{enumerate}
 \item If $v_1\notin\mk{m}$, then one can take $f_1(s)\equiv 0$.
 \item If $u_2\notin\mk{n}^*$, then one can take $f_2(s)\equiv 0$.
 \item If there exists a unit $\epsilon\in\ca{V}_{v_2,u_2,\mk{m}}$, such that $\omega_{\ov{p}}(\epsilon)=-1$, then one can take 
 $f_1(s)\equiv 0$.
 \item If there exists a unit $\epsilon\in\ca{V}_{v_1,u_1,\mk{m}}$, such that $\omega_{\ov{p}}(\epsilon)=-1$, then one can take $f_2(s)\equiv 0$.
 \item If there exists $i\in\{1,\ldots,g\}$, such that $\ov{p}_i\neq\ov{0}$, then one can take $\alpha_1=\alpha_2=\beta=0$.
 \item If $u_2\notin\mk{n}^*$ and $\ov{p}=\ov{\bz}$, then one can take $\alpha_1=0$.
 \item If $u_1\notin\mk{m}^*$ and $\ov{p}=\ov{\bz}$, then one can take $\beta=0$.
 \item If $\ov{p}=\ov{\bz}$, $\mk{n}=\mk{m}^*$, $u_2\in\mk{m}$ and 
 $(-u_2,v_2)=(v_1,u_1)$ (the last equality being equivalent to $U=U^*$), then the pole of $\phi_1(s)$ at $s=\frac{1}{2}$
 cancels the pole of $\phi_2(s)$ at $s=\frac{1}{2}$.
\end{enumerate}
\end{Th}

{\bf Proof} The proof follows by combining together Theorem \ref{key_thmm}, Theorem \ref{main_result_previous}, 
Corollary \ref{ranki}, and Proposition \ref{tiopp}. We leave the details to the reader. \fin

\begin{Exa}
Let $K=\QQ$, $\mk{m}=\mk{n}=\ZZ$, $U=\M{0}{0}{0}{0}$, $p=w=0$. Then, for $z\in\CC\bs\RR$ and $s\in\Pi_1$, we have
\begin{align*}
G(z,s):=G_{(\mk{m},\mk{n})}(U,p;z,s)=\sum_{(0,0)\neq(m,n)\in\ZZ^2}\frac{\Imm(z)^s}{|mz+n|^{2s}}.
\end{align*}
It follows from Theorem \ref{lune} that $[s\mapsto G(z,s)]$ admits a meromorphic continuation to all of $\CC$ with 
a single pole of order one at $s=1$. Indeed, the constant term of the Fourier series expansion of $G(z,s)$ is given by
\begin{align}\label{vide}
2\zeta_\QQ(2s)y^s+\frac{2^{2-2s}\pi\cdot\Gamma(2s-1)}{\Gamma(s)^2}\cdot 2\zeta_\QQ(2s-1)y^{1-s},
\end{align}
where $\zeta_\QQ(s)=\frac{1}{2}\wt{Z}_{\ZZ}((0,0),0;s)$ corresponds the Riemann's zeta function. Note that the poles 
of each of the two terms in \eqref{vide}, at $s=\frac{1}{2}$, cancel out 
(this follows from (8) of Theorem \ref{lune}), and the second term admits a pole of order one at $s=1$.
Moreover, at $s=0$, the term $R(z,s)$ in \eqref{poum0} vanishes identically, and consequently, $G(z,0)=2\zeta(0)=-1$.
\end{Exa}

Let $\ca{Q}=((\mk{m},\mk{n}),U,p,0)$ be a standard quadruple and let 
\begin{align*}
G^0(z,s):=G_{\ca{Q}}(z,s)=G_{(\mk{m},\mk{n})}^{0}(U,p\,;z,s),
\end{align*}
be its associated Eisenstein series. We view $G^0(z,s)$ as a function on $\in\mk{h}^g\times\Pi_1$. 
Using the fact that the family $\{G^0(z,s)\}_{s\in\Pi_1}$ admits a meromorphic continuation to all
of $\CC$ in conjunction with Theorem \ref{clop} we obtain the following corollary (cf. Remark \ref{dual_fam}):
\begin{Cor}\label{partt}
Assume that $[(z,s)\mapsto G^0(z,s)]\not\equiv 0$. Then there exists cusps $c,c'\in\PP^1(K)$, such that 
\begin{enumerate}[(i)]
 \item $\spart{c}{s} G^0\not\equiv 0$,
 \item $\spart{c'}{1-s} G^0\not\equiv 0$.
\end{enumerate}
\end{Cor}

{\bf Proof} We do a proof by contradiction. Assume that $\spart{c}{s} G^0\equiv 0$ for all $c\in\PP^1(K)$. Then,
from Theorem \ref{clop}, we must have $[(z,s)\mapsto G^0(z,s)]|_{\Pi_1}\equiv 0$, which contradicts our initial assumption. Therefore,
there must exists $c\in\PP^1(K)$, such that $\spart{c}{s} G^0\not\equiv 0$. This proves (i). For the proof
of (ii) we consider instead the following family of modular forms: 
\begin{align*}
F(z,s):=G^0(z,1-s),
\end{align*} 
where $(z,s)\in \mk{h}^g\times\Pi_1$. Note that the definition of $F(z,s)$ makes sense,
since $[s\mapsto G^0(z,s)]$ {\it admits a meromorphic continuation} to all of $\CC$. Finally, from the first part of the proof, 
we know that there must exist a cusp $c'\in\PP^1(K)$ such that $\spart{c'}{1-s} G^0(z,s)=\spart{c'}{s} F(z,s)\not\equiv 0$. This concludes the proof. \fin

\subsection{Explicit Fourier series of holomorphic Eisenstein series}\label{sec_h}
We keep the same notation as in Section \ref{nuz}.
\begin{Prop}\label{bateau}
Let $p=\bz\in\ZZ^g$ be the trivial weight and let
$w\in\ZZ_{\geq 3}$ (we impose this restriction in order to have absolute convergence). 
Then, for $z\in K_\CC^{\pm}$, we have
\begin{align}\label{from}
F_w(U;z):=\sum_{\mathcal{V}^+\bs\{(0,0)\neq(m+v_1,n+v_2)\in(\mk{m}+v_1,\mk{n}+v_2)\}}
\frac{e^{{2\pi\ii } \Tr(u_1(m+v_1)+u_2(n+v_2))}}{\Norm((m+v_1)z+(n+v_2))^{w}}=\wt{T}_1+\wt{T}_2+\wt{T}_3,
\end{align}
where
\begin{align*}
\wt{T}_1=\delta_{\mk{m}}(v_1)\cdot e_1\cdot\wt{Z}_{\mk{n}}(u_2,v_2,\omega_{\epsilon(w)},0),
\end{align*}
\begin{align*}
&\wt{T}_2=0,
\end{align*}
and
\begin{align}\label{pouce}
&\wt{T}_3=\cov(\mk{n})^{-1}\cdot e^{2\pi\ii\Tr(u_2 v_2)}
\left(\frac{(-2\pi i)^w}{(w-1)!}\right)^{g}\\[1mm] \notag
&\cdot \sum_{\substack{d\in D\\ d\gg 0}}\left(\sum_{(\xi_1,\xi_2)\in R_d/\mathcal{V}^+}
\big(
e^{{2\pi\ii } \Tr(u_1\xi_1)}e^{{2\pi\ii } \Tr(u_2\xi_2)}\Norm(\xi_2)|^{w-1}\sign(\Norm(\xi_2))^w\big)\right)e^{{2\pi\ii } \Tr(d z)}.
\end{align}
Here,
\begin{enumerate}[(i)]
 \item $R_d$ is the set which appears in the beginning of Section \ref{nuz}.
 \item $\mathcal{V}^+:=\mathcal{V}_U^+(\mk{m},\mk{n})$ is the subgroup of the group
 of totally positive units which appears in Definition \ref{fromage}.
 \item  $\epsilon(w)=\ov{\bu}\in\ov{\mk{S}}$ if $w$ is odd, and $\epsilon(w)=\ov{\bz}\in\ov{\mk{S}}$ 
if $w$ is even.
\item $e_1=e_1((\mk{m},\mk{n}),U)$ is the positive rational number which appears in Definition \ref{chien}.
 \end{enumerate}
\end{Prop}
Note that since $w\geq 3$, the right-hand side of \eqref{from} makes sense since the summation converges absolutely.

{\bf Proof} Since the defining series of $F_{w}(U;z)$ converges absolutely, we have by definition of $G_{(\mk{m},\mk{n})}^w(U,\bz\, ;z,s)$
that 
\begin{align*}
G_{(\mk{m},\mk{n})}^w(U,\bz\, ;z,s)\Big|_{s=0}=F_{w}(U;z).
\end{align*}
Therefore,
\begin{align*}
F_{w}(U\,;z)=T_1+T_2+T_3
\end{align*}
where the $T_i$'s are the terms appearing in  Theorem \ref{key_thmm}. Since
$p=\bz$, we have $\alpha(s)=(s+w)\cdot\bu$ and $\beta(s)=s\cdot\bu$. By inspection, one
sees that $\wt{T}_1=T_1$. The vanishing of the 
$T_2$ term comes from the presence of a term $\Gamma(\beta_j(0))$ (where $\beta_j(0)=0$) 
in the denominator of the expression of $T_2$ of Theorem \ref{key_thmm},
and the fact that zeta function appearing in the $T_2$ term 
does not have a pole when $\alpha_j(s)+\beta_j(s)-1=2s+w-1$ (recall that $w\geq 3$). 

It remains to show
that $T_3=\wt{T}_3$. From Proposition \ref{oily}, the $T_3$ term can be rewritten as the right-hand
side of \eqref{tim00}. Now, let $(\xi_1,\xi_2)\in R_d$, so that $\xi_1\xi_2=d$. Then, for any $w\in\ZZ_{\geq 2}$, a direct
calculation shows that
\begin{align*}
e^{{2\pi }\Tr\left(d y)\right)}
\left(\prod_{j=1}^g\tau(w,0;\xi_2^{(j)},\xi_1^{(j)}y_j)\right)
&=|\Norm(\xi_2)|^{w-1}\sign(\Norm(\xi_2))^w \left(\prod_{j=1}^g\tau(w,0;1,d^{(j)}y_j)\right)\\
&=\left\{
\begin{array}{ccc}
0 & \mbox{if $d$ is not totally positive}\\
\left(\frac{(-2\pi i)^w}{(w-1)!}\right)^{g}|\Norm(\xi_2)|^{w-1}\cdot\sign(\Norm(\xi_2))^w & \mbox{if $d$ is totally positive}
\end{array}
\right.
\end{align*}
The first equality follows from (iii) and (iv) of Section \ref{formu}.
The second equality follows from Proposition \ref{posi}. This concludes the proof. \fin

\begin{Rem}
The fact that the only terms
appearing in the Fourier series are supported on totally positive elements of $K$ is an instance 
of the so-called \lq\lq K\"ocher principle\rq\rq\s for holomorphic automorphic forms, with two or more complex variables.
\end{Rem}

\begin{Rem}\label{tipi}
In the case where $g\geq 1$ and $w=2$, the summation on the right-hand side of \eqref{from} does not
converge uniformly. However, from the Theorem \ref{Sonne}, we know that $[s\mapsto G_{(\mk{m},\mk{n})}^2(U,\bz\, ;z,s)]$
admits a holomorphic continuation in a neighborhood of $s=0$. Therefore, the quantity
\begin{align*}
G_{(\mk{m},\mk{n})}^2(U,\bz\,;z,s)\Big|_{s=0}
\end{align*}
makes sense. We may thus use the above expression in order to give a meaning to the formal sum 
\begin{align*}
\mbox{\raisebox{0.3cm}{\lq\lq}}\sum\limits_{\mathcal{V}^+\bs\{(0,0)\neq(m+v_1,n+v_2)\in(\mk{m}+v_1,\mk{n}+v_2)\}}
\frac{e^{{2\pi\ii } \Tr(u_1(m+v_1)+u_2(n+v_2))}}{\Norm((m+v_1)z+(n+v_2))^{2}}.\;\;\mbox{\raisebox{0.3cm}{\rq\rq}}
\end{align*}
\end{Rem}

\begin{Exa}
For example, assume that $g=1$, $w=2$, $p=0$, $\mk{m}=\mk{n}=\ZZ$, and $U=\M{0}{0}{0}{0}$. It follows from Theorem \ref{key_thmm}
that, for  $z\in\mk{h}$, we have
\begin{align}\label{perd}
&\lim_{s\rightarrow 0}G_{(\mk{m},\mk{n})}^2(U,\bz\,;z,s)\\ \notag
&=\lim_{s\rightarrow  0} \left(2\zeta(2s+2)+(2\pi)(2y)^{-2s-1}\frac{\Gamma(2s+1)}{\Gamma(s+2)\Gamma(s)}2\zeta(2s+1)\right)
+(-2\pi\ii)^2\cdot 2\cdot\sum_{n\geq 1}\sigma_1(n)e^{2\pi\ii nz},\\ \notag
&=\frac{\pi^2}{3}-\frac{\pi}{y}-8\pi^2\cdot \sum_{n\geq 1}\sigma_1(n)e^{2\pi\ii nz},
\end{align}
where $\sigma_1(n)=\sum\limits_{\substack{d|n\\ d>0 }} d$ and $\zeta(s)$ is the Riemann zeta function. In order
to obtain the right-hand side of \eqref{perd}, we have used the following observations: 
\begin{enumerate}
 \item $\wt{Z}_{\ZZ}((0,0),\omega_{\ov{0}},s)=2\zeta(s)$,
 \item $\lim\limits_{s\rightarrow 0} \left(2\zeta(2s+1)\frac{\Gamma(2s+1)}{\Gamma(s+2)\Gamma(s)}\right)=1$,
 \item $\mathcal{V}^+=\{1\}$,
 \item $R_d/\mathcal{V}^+=R_d=\{(d_1,d_2)\in\ZZ^2:d_1d_2=d\}$, so that $\# R_d/\mathcal{V}^+=2\sigma_1(d)$.
\end{enumerate}
When $g>1$ and $w=2$, the limit
\begin{align*}
\lim_{s\rightarrow 0} \left(\frac{\Gamma(2s+1)}{\Gamma(s+2)\Gamma(s)}\right)^g Z_{\mk{m}}(0,0;\ov{p};2s+1), 
\end{align*}
always vanishes. This follows from the observation that 
$s\mapsto Z_{\mk{m}}(0,0;\ov{0};s)$ has a simple pole at $s=1$ (this being a consequence of Theorem \ref{main_result_previous}) 
and, similarly,  $s\mapsto \Gamma(s)$ is also known to have a simple pole at $s=0$. Therefore, in this case 
\begin{align*}
G_2(z):=G_{(\mk{m},\mk{n})}^2(U,\bz\, ;z,s)\Big|_{s=0},
\end{align*}
is a {\it holomorphic} modular form of parallel weight $2$.
\end{Exa}

\begin{Rem}
Since $G_2(\frac{-1}{z})=z^2G_2(z)$ and $\ii$ is a fixed point of the involution $z\mapsto\frac{-1}{z}$,
it follows that  $G_2(\ii)=0$. Therefore, after evaluating \eqref{perd} at $z=\ii$ we find the non-trivial identity
\begin{align}\label{Schloe}
\frac{\pi^2}{6}=\zeta(2)=\frac{\pi}{2}+4\pi^2\sum_{n\geq 1}\sigma_1(n)e^{-2\pi n}. 
\end{align}
This identity seems to have appeared for the first time in \cite{Schloe}. One may also obtain similar identities 
for the value of the Riemann zeta function $\zeta(k)$, when $k>2$ and $k\equiv 2\mod{4}$, by replacing $G_2(z)$ by
$G_k(z):=G_{(\ZZ,\ZZ)}^k(\M{0}{0}{0}{0},\bz\, ;z,0)$ and using again the observation that $G_k(\ii)=0$.
\end{Rem}

\section{The completed Eisenstein series and its functional equation}

\subsection{The Euler factors at $\infty$ and the completed Eisenstein series}\label{pars}

Let $p=(p_i)_{i=1}^g\in\ZZ^g$ be a fixed integral weight. As before, for $s\in \CC$, we let  
\begin{enumerate} 
\item $\alpha(s)=s\cdot\bu-\frac{p}{2}\in\CC^g$,
\item $\beta(s)=s\cdot\bu+\frac{p}{2}\in\CC^g$.
\end{enumerate}
In particular, for $1\leq i\leq g$, we have $\alpha(s)_i=s-\frac{p_i}{2}$ and  
$\beta(s)_i=s+\frac{p_i}{2}$. Notice that under the involution $s\mapsto 1-s$ we have
\begin{enumerate}[(i)]
 \item $\alpha(s)\mapsto\alpha(1-s)=1-\beta(s)$,
 \item $\beta(s)\mapsto\beta(1-s)=1-\alpha(s)$,
 \item $\beta(s)-\alpha(s)\mapsto \beta(1-s)-\alpha(1-s)=\beta(s)-\alpha(s)$. 
\end{enumerate}

\begin{Def}\label{tom}
Let $\ca{Q}=(\mk{m},\mk{n},U,p,0)$ be a standard quadruple and let $\alpha(s),\beta(s)$ be defined as above. 
For each $z=x+\ii y\in K_{\CC}^{\pm}$ and $s\in\CC$, we define 
\begin{align*}
C(\alpha(s),\beta(s);z):=&|\Norm(y)|^s\prod_{j=1}^g \left(2^{\alpha_j(s)+\beta_j(s)}
\frac{\Gamma({\alpha}_j(s))\cdot\Gamma({\beta}_j(s))}{p({\alpha_j(s)},{\beta}_j(s);1,y_j)}\right),
\end{align*}
where $p(\_,\_;\_,\_)$ is the function which appears in Definition \ref{nop}. We call $C(\alpha(s),\beta(s);z)$
the \lq\lq Euler factor at $\infty$\rq\rq\s associated to the Eisenstein series $G_{(\mk{m},\mk{n})}^{\alpha(s),\beta(s)}(U;z)$.
We view $C(\alpha(s),\beta(s);z)$ as a function in $(z,s)$. The {\bf completed Eisenstein series} of 
$G_{(\mk{m},\mk{n})}^{\alpha(s),\beta(s)}(U;z)$ is defined as
\begin{align*}
\wh{G}_{(\mk{m},\mk{n})}^{\alpha(s),\beta(s)}(U\,;z):=C(\alpha(s),\beta(s);z)\cdot G_{(\mk{m},\mk{n})}^{\alpha(s),\beta(s)}(U\,;z).
\end{align*}
\end{Def}

\begin{Exa}
For example, assume that $K=\QQ$ and $\ca{Q}=((\ZZ,\ZZ),U,0,w)$. Let 
\begin{align*}
G(z,s):=G_{(\ZZ,\ZZ)}^{\alpha(s),\beta(s)}(U\,;z)=G_{(\ZZ,\ZZ)}^{0}(U,0;z,s)
\end{align*}
be the associated Eisenstein
series to $\ca{Q}$. Then, in this special case, $C(\alpha(s),\beta(s);z)=\pi^{-s}\Gamma(s)$ is independent of $z$, and
\begin{align*}
\wh{G}(z,s)=\pi^{-s}\cdot\Gamma(s)\cdot G(z,s).
\end{align*}

\end{Exa}

The next proposition gives a precise formula for the quotient of two Euler factors. 
\begin{Prop}\label{tappe}
The exponents of the variables $|y_i|$'s ($i\in\{1,\ldots,g\}$) in $C(\alpha(s),\beta(s);z)$ {\it do not depend on $s$}.
Moreover, we have
\begin{align}\label{toff}
\frac{C(\alpha(s),\beta(s);z)}{C(\alpha(1-s),\beta(1-s);z)}=(-1)^{\Tr(\sg(z)\cdot p)}
\cdot\pi^{-2g(s-1/2)}\cdot(-\ii)^{\Tr(p)}\cdot\varphi_p(s),
\end{align}
where $\varphi_p(s)$ is the function which appears in Definition \ref{key_defi}. Here, the product
$\sg(z)\cdot p$ is computed inside the product ring $\ZZ^g$
\end{Prop}

{\bf Proof} We let $\epsilon_j=\sign(y_j)$. The exponent of $|y_j|$ in $p(\alpha_j(s),\beta_j(s);1,y_j)$ is 
given by $|y_j|^{s+\frac{\epsilon_j p_j}{2}}$. Therefore the exponent of $|y_j|$ in $C(\alpha(s),\beta(s);z)$
is given by $-\frac{\epsilon_j p_j}{2}$ which is indeed independent of $s$. This proves the first assertion. 
Let us prove the second part. To simplify the notation, we write $\alpha_j$ (resp. $\beta_j$)
instead of $\alpha_j(s)$ (resp. $\beta_j(s)$). With this convention, it follows that
$\alpha(1-s)_j=1-\beta_j$ and $\beta(1-s)_j=1-\alpha_j$. By definition, we have
\begin{align}\label{choco1}
\frac{C(\alpha(s),\beta(s);z)}{C(\alpha(1-s),\beta(1-s);z)}&=\\  \notag
&\hspace{-2cm} 2^{4s-2}\cdot|\Norm(y)|^{2s-1}\prod_{j=1}^g \left(\frac{\Gamma(\alpha_j)\Gamma(\beta_j)}{\Gamma(1-\beta_j)\Gamma(1-\alpha_j)}
\cdot \frac{p(1-\beta_j,1-\alpha_j;1,y_j)}{p(\alpha_j,\beta_j;1,y_j)} \right).
\end{align}
From the definition of $p(\_,\_;\_,\_)$, we have
\begin{align}\label{choco2}
\frac{p(1-\beta_j,1-\alpha_j;1,y_j)}{p(\alpha_j,\beta_j;1,y_j)}=
\left\{
\begin{array}{cc}
\frac{\Gamma(1-\beta_j)}{\Gamma(\alpha_j)}\cdot \left|4\pi y_j\right|^{1-\alpha_j-\beta_j} & \mbox{if $y_j>0$} \\
\frac{\Gamma(1-\alpha_j)}{\Gamma(\beta_j)}\cdot \left|4\pi y_j\right|^{1-\alpha_j-\beta_j} & \mbox{if $y_j<0$}
\end{array} 
\right.
\end{align}
Moreover, from the duplication formula for the gamma function, we have 
\begin{align}\label{choco3}
\frac{\Gamma(1-\alpha_j)}{\Gamma(\beta_j)}=(-1)^{\alpha_j-\beta_j}\cdot\frac{\Gamma(1-\beta_j)}{\Gamma(\alpha_j)}.
\end{align}
Substituting the right hand side of \eqref{choco3} into the right hand side of \eqref{choco2} (when $y_j<0$), the equality \eqref{choco1}
can be rewritten as
\begin{align}\label{choco33}
\frac{C(\alpha(s),\beta(s);z)}{C(\alpha(1-s),\beta(1-s);z)}=
\pi^{-2g(s-\frac{1}{2})}\prod_{j=1}^g \left((-1)^{\epsilon_j (\alpha_j-\beta_j)}\cdot\frac{\Gamma(1-\beta_j)}{\Gamma(\alpha_j)}\right).
\end{align}
By definition of $\varphi_p(s)$ (see Definition \ref{key_defi}) we have
\begin{align}\label{choco4}
\varphi_p(1-s)=\prod_{j=1}^g (\ii)^{\beta_j-\alpha_j}\cdot\frac{\Gamma(1-\alpha_j)}{\Gamma(\beta_j)}.
\end{align}
Substituting \eqref{choco4} in the right-hand side of \eqref{choco33}, we finally obtain
\begin{align*}
\frac{C(\alpha(s),\beta(s);z)}{C(\alpha(1-s),\beta(1-s);z)}&=\pi^{-2g(s-\frac{1}{2})}\cdot (-1)^{\Tr(\sg(z)\cdot p)}\cdot(-\ii)^{\Tr(p)}\varphi_p(1-s).
\end{align*}
This concludes the proof. \fin

\begin{Rem}
Note that from the left-hand side
of \eqref{toff} we may deduce that the function $\varphi_p(s)$ satisfies the functional equation 
\begin{align*}
\varphi_p(1-s)=(-1)^{\Tr(p)}\cdot\varphi_p(s)^{-1}.
\end{align*}
This can also be verified directly from the definition of $\varphi_p(s)$.
\end{Rem}

\subsubsection{Meromorphic continuation of the completed Eisenstein series}\label{nozze}

\begin{Th}\label{Sonne}
Let $\ca{Q}=(\mk{m},\mk{n},U,p,0)$ be a standard quadruple and let $\wh{G}_{\mk{m},\mk{n}}^{\alpha(s),\beta(s)}(U;z)$ be the
associated completed Eisenstein series to $\ca{Q}$. Let $m_p=\max_{j=1}^g\{\pm p_j\}$ and
let 
\begin{align*}
S_p:=\left\{\pm\frac{n}{2}\in\frac{1}{2}\ZZ:n\in\ZZ, -m_p\leq n\leq m_p\right\}.
\end{align*}
For each $z\in K_{\CC}^{\pm}$ fixed, we consider the function
\begin{align*}
F_z(s):=[s\mapsto\wh{G}_{\mk{m},\mk{n}}^{\alpha(s),\beta(s)}(U;z)],
\end{align*}
in the variable $s\in\CC\bs S_p$. Then the function $[s\mapsto F_z(s)]$ is analytic on the open set $\CC\bs S_p$. Furthermore, 
for each $s_0\in\CC\bs S_p$, we have $\ord_{s=s_0} F_z(s)\geq -(g+1)$. 
\end{Th}

{\bf Proof} This follows from the first part of Theorem \ref{lune} and the definition of $C(\alpha(s),\beta(s);z)$. \fin

\subsubsection{Fourier series expansion of the completed Eisenstein series}\label{secu}

The next proposition gives the Fourier series expansion of the completed Eisenstein series
as a sum of three terms. 
\begin{Prop}\label{bird}
We keep the same notation as in Section \ref{pars}. In particular $w=0$.
Let $z\in K_{\CC}^{\pm}$ and $s\in\Pi_1$. Then the Fourier series expansion 
of the completed Eisenstein series is given by
\begin{align*}
\wh{G}_{(\mk{m},\mk{n})}^{\alpha(s),\beta(s)}(U\,;z)=S_1+S_2+S_3,
\end{align*}
where
\begin{align}\label{cota}
&S_1=\delta_{\mk{m}}(v_1)\cdot C(\alpha(s),\beta(s);z)\cdot\phi_1(s)\cdot|\Norm(y)|^s,
\end{align} 
\begin{align}\label{cota2}  
S_2&=\delta_{\mk{n}^*}(u_2)\cdot C(\alpha(s),\beta(s);z)\cdot\phi_2(s)\cdot|\Norm(y)|^{1-s}.
\end{align}
Here $\phi_1(s)$ and $\phi_2(s)$ are defined as in Section \ref{nico}. Moreover, the term $S_3$ is given by
\begin{align}\label{gout}
&S_3=\cov(\mk{n})^{-1}\cdot e^{2\pi\ii\Tr(u_2v_2)}
\cdot 2^{2sg}\sum_{\mathcal{V}^+\bs\{0\neq\xi_1\in\mk{m}+v_1\}}
e^{{2\pi\ii } \Tr(u_1\xi_1)}\cdot\\[2mm] \notag
&\hspace{-1cm}\sum_{0\neq\xi_2\in (\mk{n}^*-u_2)}\Bigg(\prod_{j=1}^g \left(\frac{\Gamma(\alpha_j)
\Gamma(\beta_j)}{p(\alpha_j,\beta_j,1,y_j)}\right)
\cdot\tau\left(\alpha_j,\beta_j,\xi_2^{(j)},\xi_1^{(j)})\right)\Bigg)
e^{2\pi\ii\Tr(\xi_2 v_2)}\cdot e^{{2\pi\ii } \Tr(\xi_2\xi_1 x)}\cdot|\Norm(y)|^{2s}.
\end{align}
\end{Prop}
{\bf Proof} This follows directly from Theorem \ref{key_thmm} and the definition of the Euler factor at infinity. \fin 

\begin{Rem}
It is possible to rewrite $S_1$ (resp. $S_2$) in terms of the completed zeta function
$\wh{Z}_{\mk{n}}(v_2,u_2,\omega_{\ov{p}},2s)$ (resp. the completed zeta function $\wh{Z}_{\mk{m}}(v_1,u_1,\omega_{\ov{p}},2s-1)$), although
the analogous formulas to \eqref{cota} and \eqref{cota2} are more complicated. For this reason, we have preferred to leave it as in \eqref{cota} and \eqref{cota2}. 
\end{Rem}

\subsubsection{Rewriting the $S_3$ term}\label{nuz2}
The goal of this subsection is to rewrite the $S_3$ term in Proposition \ref{bird} in a more symmetrical way. This rewriting of the $S_3$ term
will be one of the key ingredients in the second proof of Theorem \ref{nice_fn_eq}. We use the same
notation as in Section \ref{nuz}. 

It follows from the Proposition \ref{oily} and Proposition \ref{bird} that
\begin{align}\label{tim0} \notag
S_3&=\cov(\mk{n})^{-1}\cdot e^{2\pi\ii\Tr(u_2 v_2)}\cdot
2^{2s g}\sum_{d\in D} \Bigg(\prod_{j=1}^g \frac{\Gamma(\alpha_j)
\Gamma(\beta_j)}{p(\alpha_j,\beta_j,1,y_j)}\cdot\tau\left(\alpha_j,\beta_j,1,d^{(j)}\frac{y_j}{h^{(j)}}\right)\Bigg)
\\[2mm]
&\cdot\left(\sum_{(\xi_1,\xi_2)\in R_d/\mathcal{V}^+}(-1)^{\ca{W}_p(\xi_2)}|\Norm(\xi_2)|^{2s-1}e^{{2\pi\ii } \Tr(u_1\xi_1)}\cdot 
e^{{2\pi\ii } \Tr(\xi_2 v_2)}\right)\cdot e^{{2\pi\ii } \Tr(dx))}\cdot|\Norm(y)|^{2s}.
\end{align}

Now, we would like to rewrite the expression on the right-hand side of 
\eqref{tim0} in such a way that the functions $\wh{\tau}(\_,\_,\_,\_)$ (see Definition \ref{nop}) and
$q(\_,\_;\_,\_;\_,\_)$ (see \eqref{mall} in Proposition \ref{key_prop2}) appear together. 

Let us first explain how to make the function $\wh{\tau}(\_,\_,\_,\_)$ appear. A direct calculation which
uses the definition of functions $p(\_,\_,\_,\_)$, $\wh{\tau}(\_,\_,\_,\_)$ and the identity (5) of Proposition \ref{magma}, 
for $d\neq 0$, gives the following identity:
\begin{align}\label{timq}\notag
&\prod_{j=1}^g\Bigg(\frac{p(\beta_j,\alpha_j,1,d^{(j)}y_j)}{p(\alpha_j,\beta_j,1,y_j)}
\cdot\wh{\tau}\left(\alpha_j,\beta_j;1,d^{(j)}y_j\right)\Bigg)\\
&=\Bigg(\prod_{j=1}^g\frac{\Gamma(\alpha_j)\Gamma(\beta_j)}{p(\alpha_j,\beta_j,1,y_j)}
\cdot\tau(\alpha_j,\beta_j,1,d^{(j)}y_j)\Bigg)|\Norm(y)|^{2s}|\Norm(d)|^{2s}\cdot 2^{2s g}.
\end{align}
\begin{Rem}
Note the swap of the arguments $\alpha_j$ and $\beta_j$ in the expression $\frac{p(\beta_j,\alpha_j,*,*)}{p(\alpha_j,\beta_j,*,*)}$ 
of the left-hand side of the above equality.
\end{Rem}

After substituting the left-hand side of \eqref{timq} in the right-hand side of \eqref{tim0}, we obtain
\begin{align}\label{tye}\notag
& S_3=\cov(\mk{n})^{-1}\cdot e^{2\pi\ii\Tr(u_2v_2)}\\ 
&\hspace{2cm}\sum_{d\in D}
|\Norm(d)|^{-2s}\cdot\Bigg(\prod_{j=1}^g 
\frac{p(\beta_j,\alpha_j,1,d^{(j)}y_j)}{p(\alpha_j,\beta_j,1,y_j)}
\cdot\wh{\tau}\left(\alpha_j,\beta_j,1,d^{(j)}y_j\right)\Bigg)e^{{2\pi\ii } \Tr(dx)}\\[2mm] \notag
&\sum_{(\xi_1,\xi_2)\in R_d/\mathcal{V}^+}(-1)^{\ca{W}_p(\xi_2)}\cdot e^{{2\pi\ii } \Tr(u_1 \xi_1)}
e^{{2\pi\ii } \Tr(v_2 \xi_2)}|\Norm(\xi_2)|^{2s-1}\cdot|\Norm(y)|^{2s}.
\end{align}
From the definition of the function $q(u,v;t_1,y_1,t_2,y_2)$ (see \eqref{mall}) we have
\begin{align*}
q(\alpha_j,\beta_j;1,d^{(j)}y_j;1,y_j)=\frac{p(\beta_j,\alpha_j,1,d^{(j)}y_j)}{p(\alpha_j,\beta_j,1,y_j)}.
\end{align*}
Moreover, for $d=\xi_1\xi_2\neq 0$, one has
\begin{align*}
|\Norm(d)|^{-s}\cdot|\Norm(\xi_2)|^{2s-1}=|\Norm(\xi_1)|^{-s}\cdot|\Norm(\xi_2)|^{-(1-s)}.
\end{align*}
Combining the previous two observations with \eqref{tye}, we finally obtain the following:
\begin{Prop}\label{roi}
We have 
\begin{align}\label{roi2} 
& S_3=\cov(\mk{n})^{-1}\cdot e^{2\pi\ii\Tr(u_2 v_2)}\cdot\\ \notag
&\hspace{1cm }\sum_{d\in D}\Bigg(\prod_{j=1}^g \Big(|\Norm(d)|^{-s}\cdot q(\alpha_j,\beta_j;1,d^{(j)}y_j;1,y_j)
\cdot\wh{\tau}\left(\alpha_j,\beta_j,1,d^{(j)}y_j\right)\Big)\\[2mm] \notag
&\hspace{2cm} e^{{2\pi\ii } \Tr(dx)}\sum_{(\xi_1,\xi_2)\in R_d/\mathcal{V}^+}(-1)^{\ca{W}_p(\xi_2)} e^{{2\pi\ii } 
\Tr(u_1\xi_1)}\cdot e^{{2\pi\ii } \Tr(\xi_2 v_2)}\cdot|\Norm(\xi_1)|^{-s}\cdot|\Norm(\xi_2)|^{-(1-s)}\Bigg)\cdot|\Norm(y)|^{2s}.
\end{align}
\end{Prop}

\subsection{Two proofs of the functional equation of the completed Eisenstein series}\label{fun_eqn}
We are now ready to prove the third main theorem of this work.
We let $((\mk{m},\mk{n}),U,p,0)$ be a standard quadruple. We also let $\alpha(s),\beta(s)\in\CC^g$
be the corresponding weights as in Definition \ref{difo} of Section \ref{hubo}.

\begin{Th}\label{nice_fn_eq}
For  any $z\in K_\CC$ and $s\in\CC$, the following functional equation holds true
\begin{align}\label{chemin}
\wh{G}_{(\mk{m},\mk{n})}^{0}(U,p,z,s)=
(-1)^{\Tr(p)}\cdot e^{2\pi\ii\Tr(\ell_U)}\cdot\frac{\cov(\mk{n}^*)}{\cov(\mk{m})}
\cdot\wh{G}_{(\mk{n}^*,\mk{m}^*)}^{0}(U^*,p,z,1-s).
\end{align}
\end{Th}

{\bf Proof} We give two proofs of this result. Both of these proofs use in an essential way the explicit description of the
constant term of the Fourier series expansion of $[z\mapsto\wh{G}_{(\mk{m},\mk{m}^*)}^{0}(U,p,z,s)]$ at the cusp $\infty$.

Using Proposition \ref{bird} and Proposition \ref{oily}, we may break each of the Fourier series expansion
at $\infty$ of $\wh{G}_{(\mk{m},\mk{n})}^{0}(U,p,z,s)$ and $\wh{G}_{(\mk{n}^*,\mk{m}^*)}^{0}(U^*,p,z,s)$ into three parts:
\begin{align*}
\wh{G}_{(\mk{m},\mk{n})}^{0}(U,p,z,s)=S_1(s)+S_2(s)+S_3(s),
\end{align*} 
and 
\begin{align*}
\wh{G}_{(\mk{n}^*,\mk{m}^*)}^{0}(U^*,p,z,s)=S_1^*(s)+S_2^*(s)+S_3^*(s).
\end{align*}
Note that the terms $S_1(s),S_2(s)$ and $S_3(s)$ (resp. $S_1^*(s),S_2^*(s)$ and $S_3^*(s)$)
are associated to the standard quadruple $\mathcal{Q}:=((\mk{m},\mk{n}),U,p,0)$ (resp. $\mathcal{Q}^*:=((\mk{n}^*,\mk{m}^*),U^*,p,0)$)

We claim that 
\begin{enumerate}[(i)]
 \item  $S_1(s)= (-1)^{\Tr(p)}\cdot e^{2\pi\ii\Tr(\ell_U)}\cdot\frac{\cov(\mk{n}^*)}{\cov(\mk{m})}\cdot S_2^*(1-s)$,
 \item  $S_2(s)= (-1)^{\Tr(p)}\cdot e^{2\pi\ii\Tr(\ell_U)}\cdot\frac{\cov(\mk{n}^*)}{\cov(\mk{m})}\cdot S_1^*(1-s)$,
 \item  $S_3(s)= (-1)^{\Tr(p)}\cdot e^{2\pi\ii\Tr(\ell_U)}\cdot\frac{\cov(\mk{n}^*)}{\cov(\mk{m})}\cdot S_3^*(1-s)$.
\end{enumerate}
First, note that if $\delta_{\mk{m}}(v_1)=0$, then (i) is trivially true. Similarly, if $\delta_{\mk{n}^*}(u_2)=0$, then (ii) is trivially true.
So it is enough to prove (i) (resp. (ii)) under the assumption that $\delta_{\mk{m}}(v_1)=1$ (resp.$\delta_{\mk{n}^*}(u_2)=1$).

Let us start by showing (i) and (ii). It is enough to prove (i) since the proof for (ii) is similar. The identity (i) is equivalent to show that
\begin{align}\label{bana}
C(\alpha(s),\beta(s);z)\cdot \phi_1(s)=(-1)^{\Tr(p)}\cdot e^{2\pi\ii\Tr(\ell_U)}
\cdot\frac{\cov(\mk{n}^*)}{\cov(\mk{m})}\cdot C(\alpha(1-s),\beta(1-s);z)\cdot\phi_2^*(1-s),
\end{align}
where $\phi_1(s)$ (resp. $\phi_2^*(s)$) is the expression appearing in \eqref{hibou0} (resp. in \eqref{hibou}) 
which is associated to the quadruple $\ca{Q}$ (resp. to 
the quadruple $\ca{Q}^*$). Using the identity \eqref{toff}, we may rewrite \eqref{bana} as
\begin{align}\label{bana2}
 \phi_1(s)=(-1)^{\Tr(p)}\cdot e^{2\pi\ii\Tr(\ell_U)}\cdot\frac{\cov(\mk{n}^*)}{\cov(\mk{m})}\cdot(-1)^{\Tr(p\cdot\sg(z))}\cdot\pi^{-2g(1/2-s)}
 \cdot(-\ii)^{\Tr(p)}\cdot\varphi_p(1-s)
 \cdot\phi_2^*(1-s).
\end{align}
Now using the observations that $e_1(\mk{m},\mk{n},U)=e_2(\mk{n}^*,\mk{m}^*,U^*)$, $\varphi_p(1-s)=(-1)^{\Tr(p)}\cdot\varphi_p(s)^{-1}$, 
and unfolding the definitions of $\phi_1(s)$ and $\phi_2^*(s)$ in equality \eqref{bana2}, we find
\begin{align}\label{bana3}
&\wt{Z}_{\mk{m}}(v_2,u_2,\omega_{\ov{p}},2s)
=\cov(\mk{m})^{-1}\cdot e^{{2\pi\ii }\Tr(u_2 v_2)}\cdot(2\pi)^{2g s}\cdot\psi_p(1-s)\cdot\wt{Z}_{\mk{m}^*}(-u_2,v_2,\omega_{\ov{p}},1-2s).
\end{align}
But \eqref{bana3} is equivalent to the functional equation \eqref{tofou2} of Corollary \ref{fato}. This proves (i)

We thus have shown that for all pairs of lattices $\mk{m},\mk{n}\subseteq K$, and all parameter matrix $U\in M_2(K)$,
\begin{align*}
F_{(\mk{m},\mk{n})}(U; z,s):=\wh{G}_{(\mk{m},\mk{n})}^{0}(U,p,z,s)-(-1)^{\Tr(p)}\cdot 
e^{2\pi\ii\Tr(\ell_U)}\cdot\frac{\cov(\mk{n}^*)}{\cov(\mk{m})}\cdot\wh{G}_{(\mk{n}^*,\mk{m}^*)}^{0}(U^*,p,z,1-s),
\end{align*}
is square-integrable at the cusp $\infty$ (see Section \ref{cusp}) of the modular variety  
$Y_{\Gamma}$, where $\Gamma=\Gamma_U(\mk{m},\mk{n};N)$. We would like now to show that $F_{(\mk{m},\mk{n})}(U;z,s)$
is square-integrable at all cusps of $\PP^{1}(K)$. Recall that $\Upsilon(\mk{m},\mk{n})$ (see Definition \ref{upsi}) 
is a certain subgroup of $SL_2(K)$ which acts transitively on $\PP^1(K)$. Let $\gamma\in\Upsilon(\mk{m},\mk{n})$ be fixed. 
Combining together the following four facts:
\begin{enumerate}[(a)]
 \item the transformation formula of Proposition \ref{pleut} for $\gamma$ applied to $F_{(\mk{m},\mk{n})}(U; z,s)$,
 \item the identity \eqref{commut} of Proposition \ref{compa} which says that Cartan involution commutes with the upper-right action,
 \item the identities (1.b) and (2.b) of Section \ref{sec_beg},
 \item and Proposition \ref{noot},
\end{enumerate}
we may deduce that
\begin{align}\label{niif}
\omega_p(j(\gamma,z))\cdot F_{(\mk{m},\mk{n})}(U;\gamma z,s)=h_{\gamma}\cdot F_{(\mk{m},\mk{n})\gamma}(U^{\gamma}; z,s),
\end{align}
where $h_\gamma\in\QQ_{>0}$. It follows from \eqref{niif} and what we have just proved before, that $F_{(\mk{m},\mk{n})}(U; z,s)$ is also square-integrable at the cusp $\gamma\infty$.
Since $\Upsilon(\mk{m},\mk{n})$ acts transitively on $\PP^1(K)$ (by (1) of Proposition \ref{noot}),
it follows that $[z\mapsto F_{(\mk{m},\mk{n})}(U;\gamma z,s)]$ is square-integrable on all of $Y_{\Gamma}$.
Finally, from Theorem \ref{clop}, we deduce that $[(z,s)\mapsto F_{(\mk{m},\mk{n})}(U;z,s)]\equiv 0$ which proves \eqref{chemin}. 
In particular, (ii) and (iii) hold true. 

Let us give a second proof of \eqref{chemin} by showing (iii) directly. From Proposition \ref{roi}, 
we may rewrite the $S_3(s)$ term in the following way:
\begin{align}\label{roy2} 
& S_3(s)=\cov(\mk{n})^{-1}\cdot\\ \notag
&\hspace{1cm} e^{2\pi\ii\Tr(u_2 v_2)}
\sum_{d\in D}\Bigg(\prod_{j=1}^g \Big(\left(q(\alpha_j(s),\beta_j(s);1,d^{(j)}y_j;1,y_j)|\cdot d^{(j)}|^{-s}\right)
\cdot\wh{\tau}\left(\alpha_j(s),\beta_j(s),1,d^{(j)}y_j\right)\Big)\\[2mm] \notag
&\hspace{1.5cm}e^{{2\pi\ii } \Tr(dx))}\cdot\Bigg(\sum_{(\xi_1,\xi_2)\in R_d/\mathcal{V}^+}(-1)^{\ca{W}_p(\xi_2)} 
e^{{2\pi\ii }\Tr(u_1\xi_1)}e^{{2\pi\ii }\Tr(v_2\xi_2)}|\Norm(\xi_1)|^{-s}|\Norm(\xi_2)|^{-(1-s)}\Bigg).
\end{align}
Similarly, we have
\begin{align*}
& S_3^*(s)=\cov(\mk{m}^*)^{-1}\cdot\\ \notag
&\hspace{1cm} e^{2\pi\ii\Tr(u_2 v_2)}
\sum_{d\in D}\Bigg(\prod_{j=1}^g \Big(\left(q(\alpha_j(s),\beta_j(s);1,d^{(j)}y_j;1,y_j)|\cdot d^{(j)}|^{-s}\right)
\cdot\wh{\tau}\left(\alpha_j(s),\beta_j(s),1,d^{(j)}y_j\right)\Big)\\[2mm] \notag
&\hspace{1.5cm}e^{{2\pi\ii } \Tr(dx))}\cdot\Bigg(\sum_{(\xi_1^*,\xi_2^*)\in R_d^*/\mathcal{V}^+}(-1)^{\ca{W}_p(\xi_2^*)} 
e^{{2\pi\ii }\Tr(u_1\xi_2^*)}e^{{2\pi\ii }\Tr(v_2\xi_1^*)}|\Norm(\xi_1^*)|^{-s}|\Norm(\xi_2^*)|^{-(1-s)}\Bigg).
\end{align*}
Recall that
\begin{align*}
R_d:=\{(\xi_1,\xi_2)\in (\mk{m}+v_1)\times(\mk{n}^*-u_2):\xi_1\xi_2=d\},
\end{align*}
and
\begin{align*}
R_d^*:=\{(\xi_1^*,\xi_2^*)\in (\mk{n}^*-u_2)\times(\mk{m}+v_1):\xi_1^*\xi_2^*=d\}.
\end{align*}
In particular, the map $\sigma:R_d\rightarrow R_d^*$ given by $(\xi_1,\xi_2)\mapsto (\xi_2,\xi_1)$ induces a bijection 
of sets.

We make four observations from which (iii) will follow:
\begin{enumerate}[(1)]
\item From Corollary \ref{shimy}, the function
\begin{align*}
\wh{\tau}\left(\alpha_j(s),\beta_j(s),1,d^{(j)}y_j\right),
\end{align*}
is invariant under the involution $s\mapsto 1-s$
\item For any $j\in\{1,\ldots,g\}$, it follows from Corollary \ref{forel} that, for $d\neq 0$, 
the following identity hold true:

\vspace{-0.5cm}

\begin{align}\label{rott}
&|d^{(j)}|^{-s}\cdot q(\alpha_j(s),\beta_j(s);1,d^{(j)}y_j;1,y_j)=\\[2mm] \notag
&\hspace{1.5cm}(-1)^{p_j\cdot(\sg(d^{(j)})+\ov{1})}\cdot
q(\alpha_j(1-s),\beta_j(1-s);1,d^{(j)}y_j;1,y_j)\cdot|d^{(j)}|^{-(1-s)},
\end{align}

Here $\sg(d^{(j)})=\ov{0}\in\ZZ/2\ZZ$ if $d^{(j)}>0$ and $\sg(d^{(j)})=\ov{1}\in\ZZ/2\ZZ$ if $d^{(j)}<0$.

\vspace{0.3cm}

To each $d\in D$ and $(\xi_1,\xi_2)\in R_d$, we associate the pair $(\xi_1^*,\xi_2^*):=\sigma(\xi_1,\xi_2)=(\xi_2,\xi_1)\in R_d^*$.
\item It follows from \eqref{tarte} that
\begin{align*}
(-1)^{\ca{W}_p(\xi_2)}(-1)^{\ca{W}_p(\xi_2^*)}=(-1)^{\ca{W}_p(\xi_2)}(-1)^{\ca{W}_p(\xi_1)}=(-1)^{\Tr(p\cdot\sg(d))},
\end{align*}
\item and a direct inspection which uses (3) reveals that
\begin{align*}
&\sum_{(\xi_1,\xi_2)\in R_d/\ca{V}^+} (-1)^{\ca{W}_p(\xi_2)}\cdot|\xi_1|^{-s}
|\xi_2|^{-(1-s)}e^{2\pi\ii\Tr(\xi_2 v_2)}e^{2\pi\ii\Tr(\xi_1 u_1)}\\ 
&=(-1)^{\Tr(p\cdot\sg(d))}\sum_{(\xi_1^*,\xi_2^*)\in R_d^*/\ca{V}^+} (-1)^{\ca{W}_p(\xi_2^*)}\cdot|\xi_2^*|^{-s}
|\xi_1^*|^{-(1-s)}e^{2\pi\ii\Tr(\xi_1^* v_2)}e^{2\pi\ii\Tr(\xi_2^* u_1)}.
\end{align*}
\end{enumerate}
Combining altogether these four observations, we finally obtain the identity
\begin{align}\label{jek}
S_3(s)= (-1)^{\Tr(p)}\cdot e^{2\pi\ii\Tr(\ell_U)}\cdot \frac{\cov(\mk{n}^*)}{\cov(\mk{m})}\cdot S_3^*(1-s).
\end{align}

This concludes the proof. \fin

\subsection{Functional equation for the uncompleted Eisenstein series}\label{der}

In this section, we give the functional equation for the {\it uncompleted} Eisenstein series.
Let $\ca{Q}=((\mk{m},\mk{n}),U,p,w)$ be a standard quadruple and let 
\begin{align*}
G_{\ca{Q}}(z,s)=G_{(\mk{m},\mk{n})}^w(U,p\,;z,s)=G_{(\mk{m},\mk{n})}^{\alpha(s),\beta(s)}(z,s),
\end{align*}
be its associated Eisenstein series. For $s\in\CC$, recall that 
$\alpha(s),\beta(s)\in\CC^g$ are defined as in Section \ref{pars}. For $z\in K_\CC^{\pm}$
and $s\in\CC$, it is convenient to define
\begin{align*}
\Phi_p(s;\sg(z))&:=\frac{C(\alpha(s),\beta(s);z)}{C(\alpha(1-s),\beta(1-s);z)}\\[2mm]
&=(-1)^{\Tr(\sg(z)\cdot p)}
\cdot\pi^{-2g(s-1/2)}\cdot(-\ii)^{\Tr(p)}\cdot\varphi_p(s),
\end{align*}
where the second equality follows from Proposition \ref{tappe}. It follows from the definition of $\varphi_p(s)$ (see Definition \ref{key_defi}) that for $s\in\RR$, 
the value $\Phi_p(s;\sg(z))$ (when defined) is again a real number.

We may now formulate a functional equation for the uncompleted Eisenstein series.
\begin{Th}\label{nice_fn_eq2} 
We have
\begin{align}\label{sexy}
&G_{(\mk{m},\mk{n})}^w(U,p\,;z,s)\\ \notag
&=e^{2\pi\ii\Tr(\ell_U)}\cdot (-1)^{\Tr(p-w\cdot\bu)}\cdot\frac{\cov(\mk{n}^*)}{\cov(\mk{m})}\cdot\Phi_{p-w\cdot\bu}\left(1-s-\frac{w}{2};\sg(z)\right)\cdot 
G_{(\mk{n}^*,\mk{m}^*)}^w(U^*,p,z,1-s-w).
\end{align}
\end{Th}

{\bf Proof} This follows directly from Theorem \ref{nice_fn_eq} and the identity \eqref{toff}. \fin

\begin{Exa}
Let us assume that $\mk{m}=\mk{n}=\ZZ$, $g=1$, $w=0$ and $z\in\mk{h}$. Then, for $p\in\ZZ$ and $U\in M_2(\QQ)$, the functional equation 
\eqref{sexy} reduces to
\begin{align}\label{dop}
G_{(\ZZ,\ZZ)}^0(U,p\,;z,s)=e^{2\pi\ii\Tr(\ell_U)}\cdot\pi^{-2(\frac{1}{2}-s)}\cdot
\frac{\Gamma(1-s+p/2)}{\Gamma(s+p/2)}\cdot G_{(\ZZ,\ZZ)}^0(U^*,p\,;z,1-s).
\end{align}
One may check that \eqref{dop} is equivalent to the functional equation which appears on page 55 of \cite{Sie1}.

\end{Exa}

\subsubsection{A symmetry induced from the complex conjugation}
We would like in this short section to explain how the complex conjugation applied to the functional equation in \eqref{sexy} leads to a 
non-trivial symmetry of the function $\Phi_p(s;\sg(z))$. We choose a standard
quadruple $\ca{Q}$ of the form $\ca{Q}=\left((\mk{m},\mk{n}),U=\M{0}{0}{0}{0},p,0\right)$. Let us assume, furthermore, that
\begin{align}\label{aret}
[(z,s)\mapsto\wh{G}_{(\mk{m},\mk{n})}^{0}(U,p,z,s)]\not\equiv 0.
\end{align}
Since \eqref{aret} holds true, we may choose $z_0\in\mk{h}^g$, such that  
$[s\mapsto\wh{G}_{(\mk{m},\mk{n})}^{0}(U,p,z_0,1-s)]\neq 0$. Let
us fix an element $s\in\CC$ such that $\wh{G}_{(\mk{m},\mk{n})}^{0}(U,p,z_0,1-s)\neq 0$.

From Theorem \ref{nice_fn_eq2}, we deduce that
\begin{align}\label{pass0}
\Phi_p(s;\sg(z_0))=(-1)^{\Tr(p)}\cdot\frac{G_{(\mk{m},\mk{n})}^0(U,p\,;z_0,s)}{G_{(\mk{m},\mk{n})}^0(U^*,p\,;z_0,1-s)}.
\end{align}
\begin{enumerate}
\item Because $u_1=u_2=0$, we have  $\ov{G_{(\mk{m},\mk{n})}^0(U,p\,;z_0,s)}=G_{(\mk{m},\mk{n})}^0(U,-p\,;z_0,s)$.
\item Because $v_1=v_2=0$, we have  $\ov{G_{(\mk{m},\mk{n})}^0(U^*,p\,;z_0,s)}=G_{(\mk{m},\mk{n})}^0(U^*,-p\,;z_0,s)$.
\end{enumerate}
Recall that $[s\mapsto \Phi_{p}(s;\sg(z_0))]$ is a real valued function. Taking the complex conjugate of the right-hand side of \eqref{pass0}, 
we find that
\begin{align}\label{sexy2}
\Phi_{p}(s;\sg(z_0))=\frac{\ov{G_{(\mk{m},\mk{n})}^0(U,p\,;z_0,s)}}{\ov{G_{(\mk{m},\mk{n})}^0(U^*,p,z_0,1-s)}}&=
\frac{G_{(\mk{m},\mk{n})}^0(U,-p\,;z_0,s)}{G_{(\mk{m},\mk{n})}^0(U^*,-p,z_0,1-s)}\\ \notag
&=(-1)^{\Tr(p)}\cdot\Phi_{-p}(s;\sg(z_0)).
\end{align}
From \eqref{sexy2}, we deduce that
\begin{align}\label{aret3}
\varphi_{-p}(s)=\varphi_{p}(s).
\end{align}

In fact, as one may naturally guess at this point, the identity \eqref{aret3} is unconditionally true. 
This follows from the Euler's reflection formula for the gamma function, see (3) of Proposition \ref{kalk}.

\subsection{Another proof of the meromorphic continuation and the functional equation
of $s\mapsto \wt{Z}_V(a,b,\omega_{\ov{p}};s)$}\label{Colmez_trick}

The goal of this section is to give a proof different from the one given in \cite{Ch4} 
of the meromorphic continuation and the functional equation of the partial zeta function 
$s\mapsto\wt{Z}_V(a,b,\omega_{\ov{p}};s)$ (see Definition \ref{nor_zeta}). The method presented here is based on an idea
of Colmez (see Theorem III.4.5 on p. 90 of \cite{Colmez1}).

Consider a standard quadruple $\ca{Q}$ of the form $\ca{Q}=\left((\mk{m},\mk{n}),U=\M{u_1}{0}{u_2}{v_2},p,0\right)$. In order
to simplify the notation, we set
\begin{align*}
G_{\ca{Q}}(z,s):=G_{(\mk{m},\mk{n})}^0(U,p;z,s).
\end{align*}
It follows from Theorem \ref{key_thmm} that the Fourier series expansion of
$[z\mapsto G_{\ca{Q}}(z,s)]$ at the cusp $\infty$ can be written in the form
\begin{align}\label{vonk1}
e_1\cdot \wt{Z}_{\mk{n}}(v_2,u_2,\omega_{\ov{p}},2s)\cdot|\Norm(y)|^s+
\Omega_p(s;\sg(z))\cdot \delta_{\mk{n}^*}(u_2)\cdot e_2\cdot\wt{Z}_{\mk{m}}(0,u_1,\omega_{\ov{p}},2s-1)\cdot|\Norm(y)|^{1-s}+R_{\ca{Q}}(z,s),
\end{align}
where 
\begin{align*}
\Omega_p(s;\sg(z)):=&\cov(\mk{n})^{-1}e^{{2\pi\ii }\Tr(u_2 v_2)}(\ii)^{\Tr(p)}(2\pi)^{g}\cdot 2^{g(1-2s)}\cdot(-1)^{\Tr(\sg(z)\cdot p)}
\cdot\left(\prod_{j=1}^g\frac{\Gamma(2s-1))}{\Gamma(s-p_j/2)\Gamma(s+p_j/2)}\right),
\end{align*}
$e_1,e_2$ are positive rational numbers as defined in Section \ref{nico}, and $R_{\ca{Q}}(z,s)$ corresponds to the $T_3$ term
of Theorem \ref{key_thmm}. It follows from Proposition \ref{tiopp} that for a fixed $z\in K_\CC^{\pm}$, the function
$[s\mapsto R_{\ca{Q}}(z,s)]$ admits a {\it holomorphic continuation to all of $\CC$}. We may thus view $[(z,s)\mapsto R_{\ca{Q}}(z,s)]$
as a function on $K_{\CC}^{\pm}\times\CC$; it is real analytic in $z$ and holomorphic in $s$. Consider now the matrix
$\gamma=\M{0}{-1}{1}{0}\in SL_2(K)$. It follows from Proposition \ref{pleut} that
\begin{align}\label{vonk2}
G_{\ca{Q}}(\gamma z,s)=\omega_{-p}(z)\cdot G_{\ca{Q}'}(z,s),
\end{align}
where $\ca{Q}'=\left((\mk{n},\mk{m}),\wt{U}=\M{u_2}{v_2}{-u_1}{0},p,0\right)$ (notice the swap between the positions of $\mk{m}$ and $\mk{n}$). Note that the
rational factor $f_{\gamma}$ which appears in Proposition \ref{pleut} has disappeared since $f_{\gamma}=1$. Similarly
to \eqref{vonk1}, the Fourier series expansion of
$[z\mapsto G_{\ca{Q'}}(z,s)]$ at the cusp $\infty$ can be written in the following way
\begin{align}\label{vonk3}
&e_1'\cdot \delta_{\mk{n}}(v_2)\cdot \wt{Z}_{\mk{m}}(0,-u_1,\omega_{\ov{p}},2s)\cdot|\Norm(y)|^s+\\ \notag
&\hspace{2cm}\Omega_p(s;\sg(z))\cdot \delta_{\mk{m}^*}(u_1)\cdot e_2'\cdot\wt{Z}_{\mk{n}}(v_2,u_2,\omega_{\ov{p}},2s-1)\cdot|\Norm(y)|^{1-s}+R_{Q'}(z,s).
\end{align}

Fix an element $a\in\RR$ with $0<a<1$ and set $z_a=(\sqrt{a^{-1}-1}+\ii)\cdot\bu\in\mk{h}^g$. By definition of $z_a$, we have
$\Imm(z_a)=\bu$ and $\Imm(\frac{-1}{z_a})=a\cdot\bu$. From \eqref{vonk2}, we deduce that 
\begin{align}\label{vonk4}
w_a\cdot G_{\ca{Q}'}\left(z_a,s\right)=G_{\ca{Q}}\left(\frac{-1}{z_a},s\right),
\end{align} 
where $w_a=\omega_{-p}(z_a)$ is a complex number lying on the unit circle. Substituting \eqref{vonk1} (resp. \eqref{vonk3}) 
in the right hand side of \eqref{vonk4} (resp. in the left-hand side of \eqref{vonk4}),
we may deduce the identity 
\begin{align}\label{vonk5} 
& T_{\ca{Q}',\ca{Q}}(z_a,s)=T(z_a,s):=w_a \cdot R_{\ca{Q}'}\left(z_a,s\right)-R_{\ca{Q}}(-1/z_a,s)\\[2mm] \notag
& =e_1\cdot \wt{Z}_{\mk{n}}(v_2,u_2,\omega_{\ov{p}},2s)\cdot a^{gs}+
\Omega_p(s;\bz)\cdot \delta_{\mk{n}^*}(u_2)\cdot e_2\cdot\wt{Z}_{\mk{m}}(0,u_1,\omega_{\ov{p}},2s-1)\cdot a^{g(1-s)}\\  \notag
&\hspace{2cm}-w_a\left(e_1'\cdot \delta_{\mk{n}}(v_2)\cdot \wt{Z}_{\mk{m}}(0,-u_1,\omega_{\ov{p}},2s)+
\Omega_p(s;\bz)\cdot e_2'\cdot \delta_{\mk{m}^*}(u_1)\cdot\wt{Z}_{\mk{m}^*}(v_2,u_2,\omega_{\ov{p}},2s-1)\right).
\end{align}
For $s\in\Pi_1$ fixed, consider the matrix 
\begin{align*}
M=M_{a,s}:=
\left(
\begin{array}{ccc}
a^{gs} & a^{g(1-s)} & w_a\\
a^{2gs} & a^{2g(1-s)} & w_{a^2}\\
a^{3gs} & a^{3g(1-s)} & w_{a^3}\\
\end{array}
\right).
\end{align*}
It follows from \eqref{vonk5} that we have the following linear system:
\begin{align}\label{kuce}
 M\cdot V=
\left(
\begin{array}{c}
T(z_{a},s)\\
T(z_{a^2},s)\\
T(z_{a^3},s)
\end{array}
\right),
\end{align}
where the coordinates of the column vector 
$V=
\left(
\begin{array}{c}
V_1\\
V_2\\
V_3
\end{array}
\right)
$ 
are explicitly given by 
\begin{enumerate}
 \item $V_1=e_1\cdot \wt{Z}_{\mk{n}}(v_2,u_2,\omega_{\ov{p}},2s)$ 
 \item $V_2=\Omega_p(s;\bz)\cdot \delta_{\mk{n}^*}(u_2)\cdot e_2\cdot\wt{Z}_{\mk{m}}(0,u_1,\omega_{\ov{p}},2s-1)$,
 \item $V_3=e_1'\cdot \delta_{\mk{n}}(v_2)\cdot \wt{Z}_{\mk{m}}(0,-u_1,\omega_{\ov{p}},2s)+
\Omega_p(s;\bz)\cdot e_2'\cdot \delta_{\mk{m}^*}(u_1)\cdot\wt{Z}_{\mk{m}^*}(v_2,u_2,\omega_{\ov{p}},2s-1)$.
\end{enumerate}
If $\Ree(s)$ is large enough, one may check that $\det(M_{a,s})\neq 0$. Therefore, the map $[s\mapsto \det(M_{a,s})]$ is a holomorphic function on $\CC$ which is 
{\it not identically equal to zero}. Using this observation to the inverse of the linear system \eqref{kuce}, we automatically 
obtain that the coordinates of 
the column vector $V$ are meromorphic functions in $s$. In particular, the function $s\mapsto \wt{Z}_{\mk{n}}(v_2,u_2,\omega_{\ov{p}},2s)$ admits a meromorphic 
continuation to all of $\CC$. Since the triple $(\mk{n},v_2,u_2)$ was arbitrary, it follows that, for each standard quadruple
$\ca{Q}''$ and $z\in K_{\CC}^{\pm}$ fixed, the function 
$[s\mapsto G_{\ca{Q}''}(z,s)]$ admits a meromorphic continuation to all of $\CC$. Moreover, from the identity principle for holomorphic functions,
the transformation formula of Proposition \ref{pleut} holds true for all $\gamma\in SL_2(K)$ and all $s\in\CC$ (away from the possible poles in $s$).

Let us now show that $[s\mapsto \wt{Z}_{\mk{n}}(v_2,u_2,\omega_{\ov{p}};s)]$ satisfies the functional equation stated
in Theorem \ref{enig}. As before, we let $\ca{Q}=\left((\mk{m},\mk{n}),U=\M{u_1}{0}{u_2}{v_2},p,0\right)$.
Consider the real analytic Eisenstein series
\begin{align*}
F(z,s):=G_{\ca{Q}}(z,s)-\Lambda_p(s;\sg(z))\cdot G_{\ca{Q}^*}(z,1-s),
\end{align*}
where
\begin{align*}
\Lambda_p(s;\sg(z))=\Phi_p(1-s;\sg(z))\cdot (-1)^{\Tr(p)}\cdot e^{2\pi\ii\Tr(\ell_U)}\cdot \frac{\cov(\mk{n}^*)}{\cov(\mk{m})}.
\end{align*}
From what has been proved before, the function $[(z,s)\mapsto F(z,s)]$ is well defined on 
all of $K_{\CC}^{\pm}\times\CC$ (away from the possible poles in $s$) and is meromorphic in $s$. Moreover, the function
$[z\mapsto F(z,s)]$ has unitary weight $-p$ with respect to a suitable congruence subgroup of $SL_2(K)$.
The Fourier series expansion of $[z\mapsto F(z,s)]$ at the cusp $\infty$ can be written as
\begin{align*}
\phi_1(s)\cdot|\Norm(y)|^s+\phi_2(s)\cdot|\Norm(y)|^{1-s}+R(z,s).
\end{align*}
From the functional equation \eqref{jek} proved during our second proof of Theorem \ref{nice_fn_eq}, we know
that $[(z,s)\mapsto R(z,s)]\equiv 0$. From the particular shape of the constant term of $F(z,s)$, and the fact that $[z\mapsto F(z,s)]$ 
has unitary weight $-p$, we deduce that $\phi_1(s)\equiv 0$ and $\phi_2(s)\equiv 0$. Set
$Z(s):=\wt{Z}_{\mk{n}}(v_2,u_2,\omega_{\ov{p}};s)$ and $Z^*(s):=\wt{Z}_{\mk{n}^*}(-u_2,v_2,\omega_{\ov{p}};s)$.
Unfolding the meaning of $\phi_1(s)\equiv 0$ and using
the fact that $\delta_{\mk{m}}(v_1)=1$, we obtain a non-trivial identity (which relates the value of $Z(2s)$ to the  value of $Z^{*}(1-2s)$ ),
which is, after inspection, seen to be equivalent to the functional equation \eqref{aler0}. This concludes the proof. \fin

\begin{Rem}
By definition, the function 
$$
[s\mapsto \wt{Z}_{\mk{m}^*}(v_2,u_2,\omega_{\ov{p}};s)]
$$ 
is holomorphic on the right half-plane $\Pi_1$. Moreover, from the functional equation \eqref{aler0}, one may also deduce that 
it is holomorphic on the left half-plane $-\Pi_{0}=\{s\in\CC:\Ree(s)<0\}$. But what about the closed 
critical strip $\ca{S}:=\{s\in\CC:0\leq\Ree(s)\leq 1\}$ ?
In principle, it should be also possible, using the method presented in this section, to show that
$[s\mapsto \wt{Z}_{\mk{m}^*}(v_2,u_2,\omega_{\ov{p}};s)$ is holomorphic on $\ca{S}\bs\{0,1\}$.
(Note that when $p=0$ and $s_0\in\{0,1\}$, the matrix $M_{a,s_0}$ is singular). One possible strategy to do so is the following: Let 
$s_0\in\ca{S}\bs\{0,1\}$ be fixed. Show that there exists $0<a_{s_0}<1$ (which may depend potentially on $s_0$), such that
$[s\mapsto \det(M_{a_{s_0},s})]$ does not vanish in a small neighborhood of $s_0$. However, the author of this work did not
try to prove (or disprove) the existence of such an element $a_{s_0}$.
\end{Rem}

\subsection{Functional equation for weighted sums of $G_{(\mk{m}_i,\mk{n}_i)}^{w}(U_i,p;z,s)$'s}\label{weighted}

In this section, we give a generalization of Theorem \ref{nice_fn_eq} for certain weighted
sums of (uncompleted) Eisenstein series of the type $G_{(\mk{m}_i,\mk{n}_i)}^{w}(U_i,p;z,s)$, where the pairs
$((\mk{m}_i,\mk{n}_i),U_i)$ are allowed to vary but where the parameters $(p,w)$ are fixed.

For $r\in\ZZ_{\geq 1}$, it is convenient to define
\begin{align*}
\ca{L}_r(K):=\{(\mk{m}_1,\mk{m}_2,\ldots,\mk{m}_r):\mbox{$\mk{m}_i\subseteq K$ is a lattice}\},
\end{align*}
i.e., the set of ordered $r$-tuples of lattices. For $\ca{M}=(\mk{m}_i)_{i=1}^r\in \ca{L}_r(K)$, we define 
\begin{align*}
\ca{M}^*:=(\mk{m}_1^*,\ldots,\mk{m}_r^*),
\end{align*}
where $\mk{m}_i^*$ is the dual lattice of $\mk{m}_i$ with respect to the trace pairing. Similarly,
for each $r$-tuples of parameter matrices $\mk{U}=(U_i)_{i=1}^r\in (M_2(K))^r$ we also define
$\mk{U}^*=(U_i^*)_{i=1}^r$.

\begin{Def}\label{dual_fa}
Let  $\ca{M}=(\mk{m}_i)_{i=1}^r\in\ca{L}_r(K)$, $\ca{N}=(\mk{n}_i)_{i=1}^r\in\ca{L}_r(K)$ be two $r$-tuples
of lattices and let $w\in\ZZ$. To each $r$-tuple $\ca{F}:=\{f_i(s)\}_{i=1}^r$ of 
meromorphic functions on $\CC$ we define its dual family, relative to the triple $((\ca{M},\ca{N},w))$, 
as the unique $r$-tuple of meromorphic functions on $\CC$, $(f_i^*(s))_{i=1}^r$, such that, for all $i\in\{1,\ldots,r\}$,
\begin{align}\label{dualit}
f_i(s)=\frac{\cov(\mk{m}_i)}{\cov(\mk{n}_i^*)}\cdot f_i^*(1-s-w). 
\end{align}
\end{Def}

The next theorem may be viewed as a direct generalization of Theorem \ref{nice_fn_eq}.
\begin{Th}\label{gener}
Let $p\in\ZZ^g$, $w\in\ZZ$ and let $\mk{U}=(U_i)_{i=1}^r\in M_2(K)$ be an $r$-tuple of parameter matrices
such that $\ell_{U_i}=:\ell_0$ is independent of $i$. Let $\ca{F}:=(f_i(s))_{i=1}^r$ be an $r$-tuple of 
$r$ meromorphic functions on $\CC$ and let 
$\ca{M}=(\mk{m}_i)_{i=1}^r\in\ca{L}_r(K),\ca{N}=(\mk{n}_i)_{i=1}^r\in\ca{L}_r(K)$. To the 
quintuple $(\mk{U},\ca{M},\ca{N},p,w)$ we associate the following Eisenstein series:
\begin{align}\label{pass}
G_{(\ca{M},\ca{N});\ca{F}}^w(\mk{U},p;z,s):=\sum_{i=1}^r f_i(s)\cdot G_{(\mk{m}_i,\mk{n}_i)}^{w}(U_i,p;z,s).
\end{align}
Then
\begin{align}\label{dejar}
\hspace{-1cm} G_{(\ca{M},\ca{N});\ca{F}}^w(\mk{U},p;z,s)=(-1)^{\Tr(p-w\cdot\bu)}\cdot e^{2\pi\ii\Tr(\ell_0)}\cdot 
\Phi_{p-w\cdot\bu}\left(1-s-\frac{w}{2};\sg(z)\right)
\cdot G_{(\ca{N}^*,\ca{M}^*);\ca{F}^*}^w(\mk{U}^*,p;z,1-s-w), 
\end{align}
where $\ca{F}^*$ is the dual family of $\ca{F}$, relative to the triple $(\ca{M},\ca{N},w)$.
\end{Th}

{\bf Proof} This follows directly from Theorem \ref{nice_fn_eq2}. We leave the details to the reader. \fin

\begin{Rem}
Note that the symmetric group $S_r$ \lq\lq acts\rq\rq\; naturally on the left on any $r$-tuple of elements chosen from a given set, 
simply by permuting its coordinates. It
follows from the definition of $G_{(\ca{M},\ca{N});\ca{F}}^w(\mk{U},p;z,s)$ that for any $\sigma\in S_r$, one has
$G_{(\ca{M},\ca{N});\ca{F}}^w(\mk{U},p;z,s)=G_{({}^{\sigma}\ca{M},{}^{\sigma}\ca{N});{}^{\sigma}\ca{F}}^w({}^{\sigma}\mk{U},p;z,s)$
\end{Rem}

\subsubsection{Functional equation for sums of ray class invariants of $K$}\label{weighted2}
In this section, we give an application of Theorem \ref{gener} by proving  functional equations
for certain {\it weighted sums of Eisenstein series}, where each term appearing in the summation {\it depends only} 
on the ray class of an ideal $\mk{m}$, and not on the ideal $\mk{m}$ itself. This section is organized as follows. 
Firstly, we specify the type of weighted sums that
will be considered (using the notation of Section \ref{weighted}). 
Secondly, we explain how one can obtain ray class invariants when the parameter matrix $U$ has a particular shape.
Finally, we end this section by stating a functional equation where each side of the functional equation 
is a sum of ray class invariants.

Let $p\in\ZZ^g$ and $w\in\ZZ$ be fixed. Choose $\ca{M}=(\mk{m}_i)_{i=1}^r\in\ca{L}_r(K)$ and set $\ca{N}=\ca{M}$. Let
$\wt{\Omega}:\{\mk{m}_1,\ldots,\mk{m}_r\}\rightarrow\CC$ be a function and assume that there exists
$\zeta\in\CC$, such that for all $i\in\{1,\ldots,r\}$,
\begin{align}\label{condi}
\wt{\Omega}(\mk{m}_i)=\zeta\cdot\wt{\Omega}(\mk{m}_i^*).
\end{align}
Later on, we will further impose that the function $\wt{\Omega}$ is induced from a function (which is denoted below by $\Omega$) on a certain ray class group.

Let $\ca{F}=(f_i(s))_{i=1}^r$ be the following $r$-tuple of holomorphic functions:
\begin{align*}
f_i(s):=\wt{\Omega}(\mk{m}_i)\cdot|\Norm(\mk{m}_i)|^{2s+w}.
\end{align*}
Let $\mk{U}=(U_i)_{i=1}^r\in (M_2(K))^r$ be an $r$-tuple of parameter matrices such that $\ell_{U_i}=:\ell_{0}$
is independent of $i$. To the previously chosen data, we associate the following weighted sum of Eisenstein series:
\begin{align*}
G_{(\ca{M},\ca{M}),\ca{F}}^w(\mk{U},p;z,s):=\sum_{i=1}^r f_i(s)\cdot G_{(\mk{m}_i,\mk{m}_i)}^{w}(U_i,p;z,s),
\end{align*}
where the left-hand side is defined as in \eqref{pass}. Let $\ca{F}^*:=(f_i(s))_{i=1}^r$ 
be the dual family of meromorphic functions associated to the family $\ca{F}$, relative to the triple $(\ca{M},\ca{M},w)$,
in the sense of Section \ref{weighted}. By definition, we have
\begin{align}\label{dejar2}\notag
f_i^*(s)&:=\wt{\Omega}(\mk{m}_i)\cdot\frac{\cov(\mk{m}_i^*)}{\cov(\mk{m}_i)}\cdot |\Norm(\mk{m}_i)|^{2(1-s-w)+w} \\
      &=\zeta\cdot\wt{\Omega}(\mk{m}_i^*)\cdot |d_K|^{2s+w-1}\cdot|\Norm(\mk{m}_i^*)|^{2s+w}.
\end{align}
It follows directly from the functional equation \eqref{dejar} that
\begin{align}\label{chom0}
&G_{(\ca{M},\ca{M}),\ca{F}}^w(\mk{U},p;z,s)\\ \notag
&=(-1)^{\Tr(p)}\cdot e^{2\pi\ii\Tr(\ell_0)}\cdot\zeta\cdot|d_K|^{1-2s-w}\cdot
\Phi_{p-w\cdot\bu}\left(1-s-\frac{w}{2};\sg(z)\right)\cdot G_{(\ca{M}^*,\ca{M}^*),\ca{G}}^w(\mk{U}^*,p;z,1-s-w),
\end{align}
where $\ca{G}=(g_i(s))_{i=1}^r$ and $g_i(s):=\wt{\Omega}(\mk{m}^*)\cdot|\Norm(\mk{m}_i^*)|^{2s+w}$. In particular, note
that the function $f_i(s)$ has the {\it same shape} as the function $g_i(s)$, except that one needs to replace the lattice 
$\mk{m}_i$ by its dual lattice $\mk{m}_i^*$. Also, note that if $\ca{F}^*=(f_i^*(s))_{i=1}^r$ denotes the dual family of $\ca{F}$ (in the sense of Definition \ref{dual_fa}),
then one may check that $\zeta\cdot|d_K|^{2s+w-1}\cdot g_i(s)=f_i^*(s)$ for $1\leq i\leq r$.

In order to construct Eisenstein series which depend only on a ray class of a lattice (and not on the lattice itself), 
it is important to assume that the parameter matrix $U$ in the expression $G_{(\mk{m}_i,\mk{m}_i)}^{w}(U,p;z,s)$
has a particular shape. Let us assume that $U=\M{u_1}{0}{u_2}{0}$, where $u_1,u_2\in\frac{1}{N}(\ca{O}_K)^*=\frac{1}{N}\mk{d}_K^{-1}$, for some integer
$N\in\ZZ_{\geq 1}$. Recall here that $\mk{d}_K$ corresponds to different ideal of $K$.
Note that $\ell_U=0$ whenever $U$ has the prescribed shape as above. Let us assume that $N$ is coprime to $\mk{d}_K$.
Now, consider the modulus $\varpi:=N\ca{O}_K\infty_1\ldots\infty_g$ of $K$, where $\{\infty_1,\ldots,\infty_g\}$
corresponds to the set of real places of $K$. For two integral ideals $\mk{m},\mk{m}'$ coprime to $N\ca{O}_K$, 
we write $\mk{m}\sim_{\varpi}\mk{m}'$, if there exists $\lambda\in \mk{m}'\mk{m}^{-1}N+1$, such that  $\lambda\mk{m}=\mk{m}'$ and
$\lambda\gg 0$. One may check that $\sim_{\varpi}$ is an equivalence relation.
The equivalence class of $\mk{m}$ will be denoted by $[\mk{m}]$. The set of equivalence classes form 
a group which we denote by $C\ell_{\varpi}(K)$, it the so-called {\it (narrow) ray class group of modulus $\varpi$}.
If $\mk{m}$ is a fractional $\ca{O}_K$-ideal, we know that $\mk{m}^{*}=\mk{d}_K^{-1}\mk{m}^{-1}$ (see Section \ref{relat}).
In particular, if $\mk{m}$ is coprime to $N$, then $\mk{m}^{*}$ is again coprime to $N$. Moreover, we also see that
$\mk{m}\sim_{\varpi}\mk{m}'\Longleftrightarrow \mk{m}^{*}\sim_{\varpi}(\mk{m}')^*$. Therefore, for
$\ca{C}=[\mk{m}]\in C\ell_{\varpi}(K)$, it makes sense to define $\ca{C}^*=[\mk{m}^*]$.

We would like now to explain that the value of the additive character on $(\mk{m},\mk{m})$ given by 
$(m_1,m_2)\mapsto e^{2\pi\ii\Tr(u_1m_1+u_2m_2)}$, only depends on \lq\lq ray class\rq\rq\; of the pair
$(m_1,m_2)$ and not on the pair itself. The element $\lambda$ in the previous paragraph can be written as 
$\lambda=\left(\sum_i a_i b_i N\right)+1$
for some elements $a_i\in\mk{m}'$ and $b_i\in\mk{m}^{-1}$. Let $m_1,m_2\in\mk{m}$. Using the previous 
writing of $\lambda$, one readily sees that
\begin{align}\label{vict0}
(\lambda-1)u_1m_1+(\lambda-1)u_2m_2\in(\ca{O}_K)^*=\mk{d}_K^{-1}.
\end{align}
In particular, from \eqref{vict0}, we deduce that
\begin{align}\label{vict}
e^{2\pi\ii\Tr(u_1m_1+u_2m_2)}= e^{2\pi\ii\Tr(u_1\lambda m_1+u_2\lambda m_2)}.
\end{align}

We keep the same notation (and the same assumptions) as above. 
Recall that $\ca{M}=(\mk{m}_i)_{i=1}^r\in\ca{L}_r(K)$. We further require that the lattices $\mk{m}_i$'s are integral ideals coprime to $N\ca{O}_K$. 
Let us choose a function $\Omega:C\ell_{\varpi}(K)\rightarrow\CC$, for which there exists $\zeta\in\CC^{\times}$ such that
\begin{align*}
\Omega([\mk{m}])=\zeta\cdot\Omega([\mk{m}^*]),
\end{align*}
for all classes $[\mk{m}]\in C\ell_{\varpi}(K)$. It follows from \eqref{vict} and
the definition of $G_{(\mk{m}_i,\mk{m}_i)}^{w}(U,p;z,s)$ that the expression
\begin{align}\label{chom}
H_{\ca{C}_i}^w(U,p\; ; z,s):=\Omega(\ca{C}_i)\cdot|\Norm(\mk{m}_i)|^{2s+w}\cdot G_{(\mk{m}_i,\mk{m}_i)}^{w}(U,p;z,s),
\end{align}
{\it depends only} on the narrow ray class of $\ca{C}_i=[\mk{m}_i]$, and not on the particular choice 
of the ideal $\mk{m}_i$ which represents the class $\ca{C}_i$. Similarly, it follows from Proposition \ref{mod_for} 
(or the functional equation in Theorem \ref{nice_fn_eq2}) that the expression
\begin{align}\label{chom2}
H_{\ca{C}_i^*}^{w}(U^*,p\; ; z,s):=\Omega(\ca{C}_i^*)\cdot|\Norm(\mk{m}_i^*)|^{2s+w}\cdot G_{(\mk{m}_i^*,\mk{m}_i^*)}^{w}(U^*,p;z,s),
\end{align}
depends only on the narrow ray class of $\ca{C}_i=[\mk{m}_i]$ (or the narrow class $\ca{C}_i^*$, since one determines the other). 

Let $(\ca{C}_j)_{j=1}^g\in (C\ell_{\varpi}(K))^r$ and $\mk{U}=(U_j)_{j=1}^r$ where each $U_j$ is either
of the form $\M{*}{0}{*}{0}$ or of the form $\M{0}{*}{0}{*}$. To all the previously chosen 
data we associate the following weighted sum of Eisenstein series:
\begin{align}\label{chom22}
G_{\ca{E}}^w(\mk{U},p;z,s):=\sum_{i=1}^r  H_{\ca{C}_i}^w(U_i,p\; ; z,s),
\end{align}
and
\begin{align*}
G_{\ca{E}^*}^w(\mk{U}^*,p;z,s):=\sum_{i=1}^r  H_{\ca{C}_i^*}^w(U_i^*,p\; ; z,s).
\end{align*}
Note that each of these two expressions is a {\it sum of ray class invariants of modulus $\varpi$}.

Finally, it follows directly from \eqref{chom0} that the following functional equation holds true:
\begin{align}\label{chom3}
G_{\ca{E}}^w(\mk{U},p;z,s)=(-1)^{\Tr(p)}\cdot \zeta\cdot |d_K|^{1-2s-w}\cdot
\Phi_{p-w\cdot\bu}\left(1-s-\frac{w}{2};\sg(z)\right)\cdot G_{\ca{E}^*}^w(\mk{U}^*,p;z,1-s-w).
\end{align}

\subsection{Hecke's real analytic Eisenstein series  and singular moduli}
The goal of this section is to give a concrete example of ray class invariants (in the sense of Section \ref{weighted2})
which has been previously considered by some authors (see \cite{Hecke4} and \cite{Gr-Za2}). Moreover, this same example provides 
an illustration on how one may use the functional equation \eqref{chom3} in order to predict a non-trivial order of vanishing at $s=0$.

Let us start by giving the original motivation for the example considered in this section.
Let $K$ be a real quadratic field. Among other things, in \cite{Hecke4}, Hecke intended to construct a {\it non-zero holomorphic} Hilbert modular
form of parallel weight $(1,1)$ for the {\it full} Hilbert modular group $SL_2(\ca{O}_K)$. He first observed that if $\ca{O}_K$
contains a unit of negative norm then his summation technique (Hecke's trick) only constructs a function
which is identically zero in the two variables $(z,s)$, where $z\in\mk{h}^2$ and $s\in\CC$.
He then tried to show that if no such unit exists, then his summation technique leads indeed to a {\it non-zero} holomorphic Hilbert
modular form of weight $(1,1)$ relative to $SL_2(\ca{O}_K)$. 
He thought he had proved the last result, but as we will see, Hecke made a sign mistake in one of his calculations 
which prevented him to see that, even in this situation (i.e. assuming that no unit of norm $-1$ exists), his construction still 
produces the zero function. However, in this case, the function considered is not identically zero in the $(z,s)$ variables
even though it vanishes identically in the variable $z$ when $s=0$.
In \cite{Gr-Za2}, Gross and Zagier took advantage of this
situation which allowed them to give an explicit formula for the factorization of the absolute norm of a difference of two singular moduli.

Let $d_1,d_2<-4$ be two negative fundamental discriminants which are coprime and let 
$D=d_1d_2$. Note that these assumptions imply that there exists at least one prime $q\equiv 3\pmod{4}$ such that
$q|D$. Let $K=\QQ(\sqrt{D})$ be the real quadratic field of discriminant $D$. Since
$\left(\frac{-1}{q}\right)=-1$, it follows that $\ca{O}_K$ has no unit of norm $-1$. 
Let $\Pic^+(\ca{O}_K)$ be the narrow ideal class group of $K$, 
 $\Pic(\ca{O}_K)$ be the wide ideal class group of $K$, and let 
\begin{align}
\pi:\Pic^+(\ca{O}_K)\rightarrow\Pic(\ca{O}_K), 
\end{align}
be the natural projection. Since $K$ has not unit of norm $-1$, it follows that $\ker(\pi)=\{[\ca{O}_K],[\mk{d}_K]\}$, 
where $(\sqrt{D})=\mk{d}_K=((\ca{O}_K)^*)^{-1}$ corresponds to the different ideal of $K$. Let 
\begin{align*}
\chi:\Pic^+(\ca{O}_K)\rightarrow \{\pm 1\},
\end{align*} 
be the genus character (see p. 60-61 of \cite{Sie1}) associated to the factorization  $D=d_1d_2$. Using the 
reciprocity law for the Kronecker symbol (see for example the calculation done p. 96 of \cite{Sie1}), one may prove that
$\chi([\mk{d}_K])=-1$.

Set $U=\M{0}{0}{0}{0}$, $p=\bz\in\ZZ^2$, $w=1$ and $\Omega=\chi$. As was observed in Section \ref{weighted2}, the function
\begin{align}\label{fiery}
\mk{m}\mapsto H_{\mk{m}}(z,s):=\Omega(\mk{m})\cdot|\Norm\mk{m}|^{2s+1}\cdot G_{(\mk{m},\mk{m})}^{1}(U,p;z,s),
\end{align} 
descends to a function on $\Pic^+(\ca{O}_K)$. In fact, on closer inspection, one sees that it descends
to a function on $\Pic(\ca{O}_K)$. Therefore, if $\ca{C}=[\mk{m}]\in \Pic(\ca{O}_K)$, it makes
sense to define
\begin{align*}
H_{\ca{C}}(z,s):=H_{\mk{m}}(z,s).
\end{align*}
Note that since $K$ has no fundamental unit of norm $-1$, there is {\it a priori no reason}
for the function $[(z,s)\mapsto H_{\ca{C}}(z,s)]$ to be identically equal to zero. By definition,
for each $s\in\CC$, the function $[z\mapsto H_{\ca{C}}(z,s)]$ has bi-weight $[(1,1);(0,0)]$ relative to the Hilbert modular
group $SL_2(\ca{O}_K)$. However, in general, for a given value of $s$ it is not holomorphic in $z$. 
In order to construct a {\it holomorphic} Hilbert modular form,
Hecke's idea was to evaluate $H_{\ca{C}}(z,s)$ at $s=0$, which makes sense by the analytic continuation
of $[s\mapsto H_{\ca{C}}(z,s)]$.

We would like now to show, with the help of the functional equation \eqref{chom3}, that the function
$[z\mapsto H_{\ca{C}}(z,0)]$ is identically equal to zero. First note that for $\ca{C}=[\mk{m}]\in \Pic(\ca{O}_K)$, we have
$\ca{C}^*=[\mk{m}^{*}]=[\mk{m}^{-1}]=\ca{C}^{-1}$. Therefore,  using the functional equation \eqref{chom3}
we deduce that
\begin{align}\label{coex1}\notag
H_{\ca{C}}(z,s)&=-\Phi_{-\bu}\left(\frac{1}{2}-s;\sg(z)\right)\cdot |d_K|^{-2s}\cdot H_{\ca{C}^*}(z,-s)\\
&=-\Phi_{-\bu}\left(\frac{1}{2}-s;\sg(z)\right)\cdot |d_K|^{-2s}\cdot H_{\ca{C}^{-1}}(z,-s).
\end{align}
In the first equality we have used the fact that $U=U^*$ and $\zeta=-1$. Let us assume furthermore that $\ca{C}=\ca{C}^{-1}$, i.e., that
$\ca{C}$ is an ambiguous ideal class. Then, using the fact that
\begin{align*}
-\Phi_{-\bu}\left(\frac{1}{2}-s;\sg(z)\right)\cdot|d_K|^{-2s}\Big|_{s=0}=-1,
\end{align*}
in \eqref{coex1}, we deduce that $[z\mapsto H_{\ca{C}}(z,0)]\equiv 0$. When $\ca{C}=[\ca{O}_K]$ is the trivial class, 
Hecke calculated in \cite{Hecke4} the constant term of the Fourier series expansion of $[z\mapsto H_{\ca{C}}(z,0)]$ and
thought that it was non-vanishing because of a sign mistake he had made. This error  prevented him to see the vanishing.
In general, if $\ca{C}=[\mk{m}]$ is not an ambiguous ideal, the function $[z\mapsto H_{\ca{C}}(z,0)]$ also vanishes.
Indeed, one may always find a matrix $\gamma\in GL_2^+(K)$ such that
$(\ca{O}_K,\ca{O}_K)\gamma=(\mk{m},\mk{m})$. Then from the previous vanishing result and the transformation formula \eqref{chat} of Proposition \ref{pleut}
we again deduce that $[z\mapsto H_{\ca{C}}(z,0)]\equiv 0$. 
\begin{Rem}
It is also possible to prove the previous vanishing result
by showing directly that the term $T_3(s)|_{s=0}$ which appears in Proposition \ref{oily} is equal to zero.
Indeed, let $\ca{C}=[\mk{m}]\in\Pic(\ca{O}_K)$ and let $D$ be the set associated to the pair $(\mk{m},\mk{m})$
which appears in the beginning of Section \ref{nuz}. First note that if $d\in D$ is not totally positive then $B_d(y;p,0)=0$. Therefore,
the only terms of $T_3(0)$ that can contribute are indexed by an element $d\in D$ with $d\gg 0$.
Let us choose $d\in D$ with $d\gg 0$ and consider the map
$\sigma:(\xi_1,\xi_2)\mapsto (\xi_2 \sqrt{D},\frac{\xi_1}{\sqrt{D}})$, where $(\xi_1,\xi_2)\in R_d$. It is straightforward 
to check $\sigma$ induces an automorphism of set of $R_d$. Moreover, since $d\gg 0$ we have $\sign(\xi_2)=-\sign(\frac{\xi_1}{\sqrt{D}})$. 
It follows from the previous observation that $b_d(0)=0$, where $b_d(s)$ is given by \eqref{coffi}. Finally, since $[z\mapsto H_{\ca{C}}(z,0)]$ has holomorphic weight
$(1,1)$ and that all of its non-zero Fourier series coefficients vanish, this forces also its constant Fourier coefficient to vanish,
and therefore, $[z\mapsto H_{\ca{C}}(z,0)]\equiv 0$.
\end{Rem}

However, the story does not end up here. As was previously observed, the function $[(z,s)\mapsto H_{\ca{C}}(z,s)]$ 
is potentially not identically equal to zero. Therefore,
it makes sense to define
\begin{align*}
\wt{H}_{\ca{C}}(z):=\left(\frac{d}{ds}H_{\ca{C}}(z,s)\big|_{s=0}\right).
\end{align*} 
By construction, $[z\mapsto \wt{H}_{\ca{C}}(z)]$ is a real analytic modular form of bi-weight $[(1,1);(0,0)]$ relative
to $SL_2(\ca{O}_K)$. In particular, if we specialize $\wt{H}_{\ca{C}}(z)$ on the diagonal of $\mk{h}^2$, we obtain the
function
\begin{align*}
F_{\ca{C}}(z_1):=\wt{H}_{\ca{C}}(z_1,z_1),
\end{align*} 
which is a  one variable real analytic modular form of weight $2$ with
respect to $SL_2(\ZZ)=SL_2(\ca{O}_K)^{\Gal(K/\QQ)}$. Moreover, one can prove that $F_{\ca{C}}(z_1)$ is non-trivial! 

Now, let $\ca{E}=(\ca{C}_j)_{j=1}^r$ where $\{\ca{C}_1,\ldots,\ca{C}_r\}=\Pic(\ca{O}_K)$ (the {\it wide} ideal class group)
and set $\mk{U}=(U_j)_{j=1}^r$ where $U_j=\M{0}{0}{0}{0}$ for all $j\in\{1,\ldots,r\}$, $p=\bz\in\ZZ^2$, $w=1$ and $\Omega=\chi$. 
Let $G_{\ca{E}}^1(\mk{U},\bz;z,s)$ be the Eisenstein series associated
to the previous data as defined in \eqref{chom22}. By definition, we have
\begin{align*}
G_{\ca{E}}^1(\mk{U},\bz;z,s)=\sum_{k=1}^{r} H_{\ca{C}}(z,s),
\end{align*}
and 
\begin{align*}
GZ(z_1):=\left(\frac{d}{ds} G_{\ca{E}}^1(\mk{U},\bz;z,s)|_{s=0}\right)\Big|_{z_1=z_2}=\sum_{k=1}^{r} F_{\ca{C}_k}(z_1).
\end{align*}
In \cite{Gr-Za2}, Gross and Zagier obtained a precise relationship between the first Fourier coefficient of 
$GZ(z_1)$ and the logarithm of $|\Norm_{L/\QQ}(j(\tau_1)-j(\tau_2))|$. 
Here $j:\mk{h}\rightarrow\CC$
corresponds to the modular $j$-invariant, $\tau_1\in\mk{h}$ (resp $\tau_2$) is a quadratic irrationality
of discriminant $d_1$ (resp. $d_2$) and $L=\QQ(j(\tau_1),j(\tau_2))$.
In particular, their result provides an explicit formula for the prime factorization of  
$|\Norm_{L/\QQ}(j(\tau_1)-j(\tau_2))|$.


\appendix

\section{Appendix}

\subsection{Functional equation of partial zeta functions}\label{app_1a}
Recall that $K$ is a totally real number field of dimension $g$ over $\QQ$.
Let $V\subseteq K$ be a lattice, $a,b\in K$ and $\ov{p}\in\ov{\mk{S}}$ be a signature element. 
For each quadruple $\ca{Q}=(a,b,\ov{p},V)$, we have defined in \cite{Ch4} the 
partial zeta function 
\begin{align}\label{lawn}
Z_V(a,b;\omega_{\ov{p}};s):=[\ca{O}_K:V]^s\sum_{\substack{a+v\in\mathfrak{R}\\ a+v\neq 0}}
\omega_{\ov{p}}(a+v)\cdot\frac{e^{{2\pi\ii } \Tr_{K/\QQ}(b(a+v))}}{|\Norm_{K/\QQ}(a+v)|^{s}},
\end{align}
where the general term of the summation is twisted simultaneously by a (finite order) additive character of $V$ ($v\mapsto e^{2\pi\ii\Tr_{K/\QQ}(bv)}$), 
and by a (multiplicative) sign character of $K^{\times}$ ($\lambda\mapsto \omega_{\ov{p}}(\lambda)$).
The indexing set $\mathfrak{R}$ below the summation is a complete
set of representatives of $\{a+V\}/\mathcal{V}^+_{a,b,V}$ (see Definition \ref{unit_def} for the meaning of $\mathcal{V}^+_{a,b,V}$).
By definition of $\mathcal{V}^+_{a,b,V}$, one may check that the coset $a+V$ is stable under the multiplication by
any element of $\mathcal{V}^+_{a,b,V}$. Therefore, the quotient makes sense.
It is easy to see that \eqref{lawn} does not depend on the set of representatives $\mathfrak{R}$ and that it 
converges absolutely for any complex number $s$, such that $\Ree(s)>1$. With the help of Proposition \ref{dor}, one may verify that
\begin{enumerate} 
 \item[(ii)] If $a\equiv a' \pmod{V}$ and $b\equiv b'\pmod {V^*}$, then 
 \begin{align*}
 Z_V(a,b;\omega_{\ov{p}};s)=e^{2\pi\ii\Tr((b-b')a)}\cdot Z_V(a',b';\omega_{\ov{p}};s).
 \end{align*}
 \item[(iii)] For all $\lambda\in K\bs\{0\}$, one has that 
 \begin{align*}
 Z_{\lambda V}(\lambda a,\frac{b}{\lambda};\omega_{\ov{p}};s)=
 \omega_{\ov{p}}(\lambda)\cdot Z_V(a,b;\omega_{\ov{p}};s).
 \end{align*}
 \end{enumerate}
\begin{Rem}\label{vanish}
The function $Z_V(a,b;\omega_{\ov{p}};s)$ {\it may happen to be identically zero}. For example, if 
there exists a unit $\epsilon\in\ca{V}_{a,b,V}$ (see Definition \ref{unit_def}; these units are no more required to be totally positive),
such that $\omega_{\ov{p}}(\epsilon)=-1$, then it follows from (ii) and (iii) that $Z_V(a,b;\omega_{\ov{p}};s)\equiv 0$.
Note that if such a unit exists, then necessarily $\ov{p}\neq\ov{\bz}$. The author does not know if the non-existence
of such a unit automatically implies that $Z_V(a,b;\omega_{\ov{p}};s)\not\equiv 0$. However, in some very special cases,
it is sometimes possible to prove the converse. Let us give one such example. Let us assume that 
the degree $[K:\QQ]$ is even, $a=b=0$ and $\ov{p}=\ov{\bu}$. (Note that if $[K:\QQ]$ were odd  then automatically
$[s\mapsto Z_V(0,0;\omega_{\ov{\bu}};s)]\equiv 0$ since $-1\in\ca{V}_{0,0,V}$). Furthermore, let us assume
that $V$ is an integral $\ca{O}_K$-ideal, that $K$ is Galois over $\QQ$
with Galois group $\Gal(K/\QQ)$, and that for all $\sigma\in\Gal(K/\QQ)$, $V^{\sigma}=V$. Then, under
all these assumptions, we claim that if $[s\mapsto Z_V(a,b;\omega_{\ov{p}};s)\equiv 0]$, then there exists
$\epsilon\in\ca{O}_K^{\times}=\ca{V}_{0,0,V}$, such that $\omega_{\ov{\bu}}(\epsilon)=-1$. Let us prove it. 
From the Chebotarev's density theorem, there exists $\lambda\in V$,
such that $|\Norm(\lambda V^{-1})|=p$ is a rational prime (in fact
there are infinitely many such pairs $(\lambda,p)$).
From the uniqueness of the writing of Dirichlet series, the vanishing of $[s\mapsto Z_{V}(0,0;\omega_{\ov{p}};s)]\equiv 0$ implies
that there exists $\lambda'\in V$, such that $\omega_{\ov{\bu}}(\lambda)|\Norm_{K/\QQ}(\lambda)|=-\omega_{\ov{\bu}}(\lambda')|\Norm_{K/\QQ}(\lambda')|$.
Note that $\mk{p}=\lambda V^{-1}$ and $\mk{p}'=\lambda' V^{-1}$ are prime ideals of $K$ above $p$.
Since $K/\QQ$ is Galois, there must exists $\sigma\in\Gal(K/\QQ)$ such that 
$\lambda V^{-1}=\mk{p}=\mk{p}'^{\sigma}=\lambda'^{\sigma}(V^{\sigma})^{-1}=\lambda'^{\sigma}V^{-1}$.
Let $\epsilon=\frac{\lambda}{\lambda'^{\sigma}}$. Then $\epsilon\in\ca{V}_{0,0,\ca{O}_K}$ and 
$\omega_{\ov{\bu}}(\epsilon)=\frac{\omega_{\ov{\bu}}(\lambda)}{\omega_{\ov{\bu}}(\lambda'^{\sigma})}=
\frac{\omega_{\ov{\bu}}(\lambda)}{\omega_{\ov{\bu}}(\lambda')}=-1$.
\end{Rem}
We would like now to state a special case of a functional equation that was proved in \cite{Ch4} for the completed
zeta function associated to $Z_V(a,b;\omega_{\ov{p}};s)$. In order to do so, we need to introduce  
the \lq\lq Euler factor at $\infty$\rq\rq\s which is associated to $Z_V(a,b;\omega_{\ov{p}};s)$. We define
\begin{align*}
F_{\ov{p},K}(s/2):=|d_K|^{s/2}\pi^{-gs/2}\prod_{i=1}^{g}\Gamma\left(\frac{s+[\ov{p}_i]}{2}\right),
\end{align*}
where $\Norm(\mk{d}_K)=d_K>0$ is the discriminant of $K$ and $\Gamma(x)$ stands for the usual gamma function evaluated at $x$. 
Recall here that $[\ov{p}_i]=0\in\ZZ$ if $\ov{p}_i=\ov{0}$, and  $[\ov{p}_i]=1\in\ZZ$ if $\ov{p}_i=1$.
As our notation suggests, the Euler factor at infinity $F_{\ov{p},K}(s/2)$ depends on the field $K$, but not on the lattice $V$ itself.

We may now state a special case of the main theorem that was proved in \cite{Ch4}.
\begin{Th}\label{main_result_previous} 
Let
\begin{align}\label{comple}
\wh{Z}_V(a,b,\omega_{\ov{p}};s):=F_{\ov{p},K}(\frac{s}{2})\cdot Z_V(a,b;\omega_{\ov{p}};s)
\end{align}
be the completed zeta function of $Z_V(a,b;\omega_{\ov{p}};s)$.
Firstly, $\wh{Z}_V(a,b,\omega_{\ov{p}};s)$ admits an analytic continuation to $\CC\bs\{0,1\}$ and has at most a pole
of order one at $s\in\{0,1\}$. Secondly, 
$\wh{Z}_V(a,b,\omega_{\ov{p}};s)$ satisfies the following functional equation:
\begin{align}\label{functional equation}
(\ii)^{\Tr([\ov{p}])}\cdot e^{{-2\pi\ii }\Tr_{K/\QQ}(ab)}\cdot\wh{Z}_V(a,b,\omega_{\ov{p}};s)=\wh{Z}_{V^*}(-b,a,\omega_{\ov{p}},1-s),
\end{align}
where  $[\ov{p}]:=([\ov{p}_1],[\ov{p}_2],\ldots,[\ov{p}_g])$. Thirdly, the function $s\mapsto\wh{Z}_V(a,b,\omega_{\ov{p}};s)$ admits
\begin{enumerate}
 \item  a pole of order one at $s=1$, if and only if,
 $p_{i}=0$ for all $i$ and $-b\in V^*$.
 \item a pole of order one at $s=0$, if and only if,
 $p_{i}=0$ for all $i$ and $a\in V$,
\end{enumerate}
\end{Th}

\begin{Rem}
The presence of the fourth root of unity $\ii=\sqrt{-1}$ in \eqref{functional equation} is in accordance with
the observation that the map $(a,b,V)\mapsto (-b,a,V^*)$ is an automorphism of order four of the set of
triples $\{(a,b,V):a,b\in K,\,\mbox{$V\subseteq K$ is a lattice}\}$.
\end{Rem}

\begin{Rem}
Unfortunately, there was one sign mistake in the statement of the functional equation given in \cite{Ch4}.
The factor $(-\ii)^{\Tr(p)}e^{-2\pi\ii \Tr_{K/\QQ}(ab)}$ appearing on the left-hand side of (1.3) of \cite{Ch4}
should be read instead as $(\ii)^{\Tr(p)}e^{-2\pi\ii \Tr_{K/\QQ}(ab)}$. There is also a change of
notations between the paper \cite{Ch4} and the present work. The function
denoted by $\Psi_V(a,b,\omega,s)$ in (1.2) of \cite{Ch4} (resp. $Z_V(a,b;\omega;s)$ in Theorem 1.1 of \cite{Ch4}),
is now denoted by $Z_V(a,b;\omega;s)$ (resp. $\wh{Z}_V(a,b,\omega_{\ov{p}};s)$). Moreover, the underlying unit
group that was used in \cite{Ch4}, which was denoted by $\Gamma_{a,b,V}$, is slightly different from the group
$\ca{V}_{a,b,V}$. See Remark \ref{turf} for an explanation of the reason for this small change of definition. 
\end{Rem}
\begin{Rem}
The main theorem proved in \cite{Ch4} holds true for any number field, and not just for totally real fields. 
\end{Rem}

Here are two non-trivial corollaries of Theorem \ref{main_result_previous}.
\begin{Cor}
Assume that $\wh{Z}_V((a,b),\omega_{\ov{p}};s)\not\equiv 0$. Then the quotient
\begin{align*}
(\ii)^{\Tr([\ov{p}])}\cdot\frac{\wh{Z}_V((a,b),\omega_{\ov{p}};s)}{\wh{Z}_{V^*}((-b,a),\omega_{\ov{p}},1-s)},
\end{align*}
depends only on $\Tr_{K/\QQ}(ab)\pmod{\ZZ}$.
\end{Cor}

The next corollary gives some lower bound for the order of vanishing of the uncompleted zeta function
$Z_{V}(a,b;\omega_{\ov{p}};s)$ at non-positive integers.
\begin{Cor}\label{ranki}
Define
\begin{enumerate}[(a)]
 \item  $r_{\ov{p}}^0:=\#\{1\leq j\leq g:\ov{p}_j=\ov{0}\}$. 
 \item  $r_{\ov{p}}^1:=\#\{1\leq j\leq g:\ov{p}_j=\ov{1}\}$. 
\end{enumerate}
Then for $\ell\in\ZZ_{\leq -1}$, we have
\begin{enumerate}
 \item if $\ell\equiv 0\pmod{2}$, then $\ord_{s=\ell} Z_V(a,b;\omega_{\ov{p}};s)\geq  r_{\ov{p}}^0$. 
 \item if $\ell\equiv 1\pmod{2}$, then $\ord_{s=\ell} Z_V(a,b;\omega_{\ov{p}};s)\geq  r_{\ov{p}}^1$. 
\end{enumerate}
In particular, the (uncompleted) zeta function $[s\mapsto Z_V(a,b;\omega_{\ov{p}};s)]$ is holomorphic on 
the left half-plane $-\Pi_0=\{s\in\CC:\Ree(s)<0\}$.
Furthermore, if we assume that 
\begin{align*}
\ord_{s=1} Z_{V^*}(-b,a;\omega_{\ov{p}};s)\geq 0,
\end{align*}
(i.e. condition (1) of Theorem \ref{main_result_previous} is not
fulfilled), then $\ord_{s=0} Z_V(a,b;\omega_{\ov{p}};s)\geq r_{\ov{p}}$.
\end{Cor}

{\bf Proof} If one rewrites the functional equation \eqref{functional equation} in terms of the uncompleted zeta function then
one obtains
\begin{align}\label{jabl}
(\ii)^{\Tr([\ov{p}])}\cdot e^{{-2\pi\ii }\Tr_{K/\QQ}(ab)}\cdot Z_V(a,b,\omega_{\ov{p}};s)= 
\frac{F_{K,\ov{p}}(\frac{1-s}{2})}{F_{K,\ov{p}}(\frac{s}{2})}\cdot Z_{V^*}(-b,a,\omega_{\ov{p}},1-s).
\end{align}
The result now follows directly from \eqref{jabl} and the well-known fact that $\Gamma(x)$ has a pole of order one 
on each element of the set $\ZZ_{\leq 0}$, and is holomorphic elsewhere. \fin

\subsubsection{The normalized zeta function $\wt{Z}_V(a,b,\omega_{\ov{p}};s)$}\label{app_1b}

It will be convenient to renormalize the uncompleted zeta function $Z_V(a,b;\omega_{\ov{p}};s)$
and rewrite the functional equation \eqref{functional equation} in terms of this renormalization.
\begin{Def}\label{nor_zeta}
We define
\begin{align*}
\wt{Z}_V(a,b,\omega_{\ov{p}};s):=\frac{Z_V(a,b;\omega_{\ov{p}};s)}{[\ca{O}_K:V]^s}.
\end{align*}
For $p\in\ZZ^g$ and a lattice $V$, we also define
\begin{align}\label{culo}
\lambda_{p}(s):=(-\ii)^{\Tr([\ov{p}])}\cdot|d_K|^{s-\frac{1}{2}}\cdot\frac{F_{\ov{p},K}(\frac{s}{2})}{F_{\ov{p},K}(\frac{1-s}{2})}. 
\end{align}
\end{Def}
One may check that $\lambda_p(s)$ {\it does not depend} on $|d_K|$. Moreover,
it is clear from the definition that the dependence of $\lambda_{p}(s)$ 
on $p$ is in fact only a dependence on $\ov{p}\in\ov{\mk{S}}=(\ZZ/2\ZZ)^g$. Finally, a direct calculation shows that
\begin{align}\label{tyrr1}
\lambda_p(1-s)=(-1)^{\Tr(\ov{p})}\cdot\lambda_p(s)^{-1}.
\end{align}

With the notation introduced above, the functional equation in Theorem \eqref{main_result_previous} may be rewritten in the following way:
\begin{Th}\label{enig}
For any $p\in\ZZ^g$, we have the functional equation
\begin{align}\label{aler0}
\wt{Z}_V((a,b),\omega_{\ov{p}};s)=\cov(V^*)\cdot e^{2\pi\ii\Tr_{K/\QQ}(ab)}\cdot
\lambda_{p}(1-s)\cdot \wt{Z}_{V^*}((-b,a),\omega_{\ov{p}},1-s),
\end{align}
where $\cov(V^*)=[\ca{O}_K:V]\cdot|d_K|^{1/2}$ is the covolume of $V^*$ (see Section \ref{lat_nu}).
\end{Th}

{\bf Proof} This follows directly from \eqref{functional equation} and the definitions of 
$\wt{Z}_V((a,b),\omega_{\ov{p}};s)$ and $\lambda_{p}(1-s)$. \fin

\subsubsection{Rewriting the functional equation}\label{app_1d}

In this subsection, we derive various identities which involve certain products of values
of the $\Gamma$ function. These identities play a key role in the proof of Theorem \ref{nice_fn_eq}.

Using \eqref{functional equation}, the fact that $\cov(V)\cov(V^*)=1$ (see Lemma \ref{bouh}),
and replacing $s$ by $2s$ in \eqref{aler0}, we deduce that
\begin{align}\label{fenster}
\wt{Z}_V(a,b,\omega_{\ov{p}},2s)=\cov(V^*)\pi^{2gs-g/2}e^{2\pi\ii
\Tr_{K/\QQ}(ab)}(-\ii)^{\Tr([\ov{p}])}\left(\prod_{j=1}^g \frac{\Gamma(\frac{1}{2}-\alpha_j(s))}{\Gamma(\beta_j(s))}\right)
\wt{Z}_{V^*}(-b,a,\omega_{\ov{p}},1-2s),
\end{align}
where $\alpha_j(s)=s-\frac{[p_j]}{2}$ and $\beta_j(s)=s+\frac{[p_j]}{2}$. Applying the change of variables
$s\mapsto 1-s$ in \eqref{fenster} (so that $\alpha_j(s)\mapsto\alpha_j(1-s)=1-\beta_j(s)$ and $\beta_j(s)\mapsto \beta_j(1-s)=1-\alpha_j(s)$), 
we get
\begin{align}\label{fenster2}
\wt{Z}_V(a,b,\omega_{\ov{p}},2(1-s))=\cov(V^*)\pi^{\frac{3}{2}g-2gs}(-\ii)^{\Tr([\ov{p}])}e^{2\pi\ii\Tr_{K/\QQ}(ab)}
\left(\prod_{j=1}^g \frac{\Gamma(\beta_j(s)-\frac{1}{2})}{\Gamma(1-\alpha_j(s))}\right)
\wt{Z}_{V^*}(-b,a,\omega_{\ov{p}},2s-1).
\end{align}
Now recall that the duplication formula for the gamma function reads as
\begin{align}\label{gam}
\Gamma(s)\Gamma(s+1/2)=2^{1-2s}\cdot\sqrt{\pi}\cdot\Gamma(2s).
\end{align}
Replacing $s$ by $s-\frac{1}{2}$ in \eqref{gam}, we get
\begin{align}\label{gam2}
\Gamma\left(s-\frac{1}{2}\right)\Gamma(s)=2^{2-2s}\sqrt{\pi}\cdot\Gamma\left(2s-1\right).
\end{align}
Since either $(\alpha_j(s),\beta_j(s))=(s,s)$ or $(\alpha_j(s),\beta_j(s))=(s-1/2,s+1/2)$, it thus follows from \eqref{gam2} that 
\begin{align}\label{ong}
2^{(2-2s)g}\pi^{g/2}\prod_{j=1}^g \Gamma(\alpha_j(s)+\beta_j(s)-1)=\prod_{j=1}^g\Gamma(\alpha_j(s))\Gamma\left(\beta_j(s)-\frac{1}{2}\right).
\end{align}
Using \eqref{ong}, we may rewrite \eqref{fenster2} as
\begin{align}\label{toy}\notag
&\wt{Z}_V(a,b,\omega_{\ov{p}},2(1-s))\\
&=\pi^{2g(1-s)}2^{2g(1-s)}\left(\cov(V^*)(-\ii)^{\Tr([\ov{p}])}e^{2\pi\ii\Tr_{K/\QQ}(ab)}
\left(\prod_{j=1}^g \frac{\Gamma(\alpha_j(s)+\beta_j(s)-1)}{\Gamma(\alpha_j(s))\Gamma(1-\alpha_j(s))}\right)
\wt{Z}_{V^*}(-b,a,\omega_{\ov{p}},2s-1)\right).
\end{align}
Here, the key identity \eqref{toy} will be used in Theorem \ref{nice_fn_eq}. 

Now comes a key observation. 
\begin{Prop}\label{fatu}
Let $s\in\CC$ and $p\in\ZZ$, and set $\alpha=s-p/2$ and $\beta=s+p/2$. Then, the quantity
\begin{align}\label{gorge0}
(-\ii)^{\beta-\alpha}\frac{\Gamma(\alpha+\beta-1)}{\Gamma(\alpha)\Gamma(1-\alpha)} 
\end{align}
depends only on the image $\ov{p}\in\ZZ/2\ZZ$ rather than $p$ itself.
More precisely, we have
\begin{align}\label{gorge}
(-\ii)^{\beta-\alpha}\frac{\Gamma(\alpha+\beta-1)}{\Gamma(\alpha)\Gamma(1-\alpha)}=
\left\{
\begin{array}{cc}
\frac{\Gamma(2s-1)}{\pi}\cdot\sin(\pi s) & \mbox{if $p$ is even}\\[2mm]
-\ii\frac{\Gamma(2s-1)}{\pi}\cdot\cos(\pi s) & \mbox{if $p$ is odd}
\end{array}
\right.
\end{align}
In particular, the involution $(\alpha,\beta)\mapsto (\beta,\alpha)$ leaves \eqref{gorge0} invariant. 
\end{Prop}

{\bf Proof} This follows from the Euler's reflection formula for the gamma function. \fin

\subsubsection{The functions $\varphi_p(s)$ and $\psi_p(s)$}\label{app_1e}

In order to take advantage of Proposition \ref{fatu}, it is convenient to make the following definitions:
\begin{Def}\label{key_defi}
Let $p\in\ZZ^g$ and let $\alpha(s),\beta(s)\in\CC^g$ be the associated weights, where $\alpha_j(s)=s-\frac{p_j}{2}$
and $\beta_j(s)=s+\frac{p_j}{2}$ for $j=1,\ldots,g$. Let also $\ov{q}\in\ov{\mk{S}}$. For each, $j\in\{1,\ldots,g\}$, we define
\begin{enumerate}
 \item $\alpha_j^*(s):=\alpha_j(s)$ and $\beta_j^*(s):=\beta_j(s)$ if 
 $\ov{q}_j=\ov{0}$.
 \item $\alpha_j^*(s):=\beta_j(s)$ and $\beta_j^*(s):=\alpha_j(s)$ if 
 $\ov{q}_j=\ov{1}$.
\end{enumerate}
We may now define the functions $\psi_{p}(\ov{q};s)$ and $\varphi_{p}(\ov{q};s)$. We set
\begin{align*}
\psi_{p}(\ov{q};s):=\prod_{j=1}^g(-\ii)^{\beta_j^*(s)-\alpha_j^*(s)}
\left(\frac{\Gamma(\alpha_j^*(s)+\beta_j^*(s)-1)}{\Gamma(\alpha_j^*(s))\Gamma(1-\alpha_j^*(s))}\right).
\end{align*}
and
\begin{align*}
\varphi_{p}(\ov{q};s):=\prod_{j=1}^g\left(\ii^{\beta_j^*(s)-\alpha_j^*(s)}\cdot
\frac{\Gamma(\beta_j^*(s))}{\Gamma(1-\alpha_j^*(s))}\right).
\end{align*}
\end{Def}
For the convenience of the reader, we present five obvious observations 
concerning the functions $\psi_p(\ov{q};s)$ and $\varphi_p(\ov{q};s)$.
\begin{Prop}\label{kalk}
\begin{enumerate}
 \item It follows from Proposition \ref{fatu} that $\psi_p(\ov{q};s)$ depends only on $\ov{p}\in\ov{\mk{S}}$ rather than $p$ itself.
 \item Since the involution $s\mapsto 1-s$ takes $\alpha_j^*(s)\mapsto 1-\beta_j^*(s)$ and $\beta_j^*(s)\mapsto 1-\alpha_j^*(s)$,
 it follows that $\varphi_{p}(\ov{q};1-s)=(-1)^{\Tr(p)}\cdot\varphi_{p}(\ov{q};s)^{-1}$ for all $s\in\CC$ where it is defined.
 \item Euler's reflection formula for the gamma function implies that  $\varphi_{-p}(\ov{q};s)=\varphi_{p}(\ov{q};s)$.
 \item It follows from Proposition \ref{fatu} that for any $\ov{q}\in\mk{S}$, $\psi_p(\ov{q};s)=\psi_p(\ov{\bz};s)$.
 \item It follows from Euler's reflection formula that for any $\ov{q}\in\mk{S}$, $\varphi_p(\ov{q};s)=\varphi_p(\ov{\bz};s)$.
\end{enumerate}
\end{Prop}
Because of (4) and (5), we see that $\psi_p(\ov{q};s)$ $\varphi_p(\ov{q};s)$ don't depend, a posteriori, on the parameter $\ov{q}$ and, therefore,
from now on, we simply write $\psi_p(s)$ and $\varphi_p(s)$.
\begin{Rem}
Contrary to the function $\psi_p(s)$,
the function $\varphi_p(s)$ really depends on $p$ itself rather than just on $\ov{p}$.
\end{Rem}

Comparing the functional equations appearing in \eqref{toy} and \eqref{aler0}, we derive the following relations:
\begin{Cor}\label{fato}
We have 
\begin{align}\label{tofou2}
&\wt{Z}_V(a,b,\omega_{\ov{p}},2(1-s))=
\cov(V^*)\cdot e^{2\pi\ii\Tr_{K/\QQ}(ab)}\cdot (2\pi)^{2g(1-s)}\cdot\psi_p(s)\cdot\wt{Z}_{V^*}(-b,a,\omega_{\ov{p}},2s-1),
\end{align}
and
\begin{align}\label{tofou}
\psi_p(s)\cdot(2\pi)^{2g(1-s)}=\lambda_p(2s-1),
\end{align}
\end{Cor}

\subsection{A certain linear system of ODEs of order $2$ in $g$-variables}\label{app_3}
Let $\Pi_1=\{s\in\CC:\Ree(s)>1\}$ and let $T\subseteq\Pi_1$ be a discrete subset. 
Suppose that 
\begin{align*}
G(z,s):\mk{h}^g\times(\Pi_1\bs T)\rightarrow\CC
\end{align*}
is a real analytic function for which there exists a lattice $\ca{L}\subseteq \RR^g$ such that, 
for each $s\in\Pi_1\bs T$ and $y\in(\RR^{\times})^g$ fixed, the function $[x\mapsto G(x+\ii y,s)]$ is $\ca{L}$-invariant. 
We may thus view $[x\mapsto G(x+\ii y,s)]$ as a function on the real torus $\RR^g/\ca{L}$. In particular, it admits  a
Fourier series expansion of the form
\begin{align}\label{but2}
a_0(y,s)+\sum_{\xi\in\mathcal{L}^*\bs\{0\}} 
a_{\xi}(y,s)e^{2\pi\ii\Tr(\xi x)},
\end{align}
where $\ca{L}^*$ is the dual lattice of $\ca{L}$. We assume, furthermore, that for each $j\in J:=\{1,\ldots,g\}$ and $s\in\Pi_1\bs T$,
\begin{align}\label{but1}
\Delta_{j,p_j} G(z,s)=s(1-s) G(z,s),
\end{align}
where $\Delta_{j,p_j}$ is the partial-graded Laplacian defined in \eqref{grad_lap}.
Then, combining the identities \eqref{but2} and \eqref{but1} for all $j\in J$, we find
that the function $y\mapsto a_{\xi}(y,s)$ ($\xi\in\ca{L}^*$) lies in the common kernels of the following $g$
partial differential operators:
\begin{align*}
D_{j,p_j}[\xi;s]:=-y_j^2 \frac{\partial^2}{\partial y_j^2}-2\pi p_j y_j\xi^{(j)}+4\pi y_j^2 (\xi^{(j)})^2-s(1-s) \hspace{1cm}(j\in J).
\end{align*}
Here, $y_j$ corresponds to the $j$-th coordinate of $y\in K_{\CC}^{\pm}$, and
$\xi^{(j)}$ corresponds to the image of $\xi$ under the $j$-th embedding of $K$ into $\RR$. 
For each $p\in\ZZ^g$, $\xi\in\ca{L}^*$ and $s\in\CC$, we define
\begin{align*}
\ca{V}_{p}(\xi;s):=\bigcap_{j=1}^g\ker(D_{j,p_j}[\xi,s]).
\end{align*}
From the theory of linear ODE's, we know that $\dim_{\CC}(\ca{V}_p(\xi;s))=2^g$. 
\begin{Exa}
When $\xi=0$, the solution 
space to $\ca{V}_{p}(\xi;s)$ is given explicitly by 
\begin{align*}
\sum_{\mu\in\{s,1-s\}^{g}} \CC y^{\mu}, 
\end{align*}
where each element $\mu=(\mu_j)_{j=1}^g$ is a vector of length $g$, where its $j$-th coordinate $\mu_j\in\{s,1-s\}$.
In this special case, the solution space {\it does not depend} on the weight $p$.
\end{Exa}
\begin{Rem}
When $\xi\neq 0$, one may also provide an explicit description of the vector space 
$\ca{V}_{p}(\xi;s)$ with the help of Kummer's hypergeometric function ${}_{1}F_{1}(a,b;z)$ but we will not
need it. 
\end{Rem}

In what follows, we will be particularly interested in the following subspace of $\ca{V}_{p}(\xi;s)$:
\begin{align*}
\ca{V}_{p}^b(\xi;s):=\left\{f(y)\in\ca{V}_{p}(\xi;s):
\mbox{$f(y)$ is bounded as $|\Norm(y)|\rightarrow\infty$}\right\}.
\end{align*}
If fact, it will be proved below that if the parameter  $s$ lies in a certain half plane which depends only
on the weight $p$, then $\ca{V}_{p}^b(\xi;s)$ is one-dimensional. 

For $p\in\ZZ$, let $A(s)=s\cdot\bu+\frac{p}{2}\in\CC^g$ and $B(s)=s\cdot\bu-\frac{p}{2}\in\CC^g$
(note that the $A(s)=\beta(s)$ and $B(s)=\alpha(s)$ where $\alpha(s)$ and $\beta(s)$ are defined
as in Section \ref{pars}). It follows from the computation done in Section \ref{dos} (or from a direct calculation) that
\begin{align}\label{pers}
B_{\xi}(y;p;s):=\left[y\mapsto\prod_{j=1}^g \tau(A_j(s),B_j(s);1,\xi^{(j)}y_j)\right],
\end{align}
is a function in the solution space $\ca{V}_{p}(\xi;s)$. With no further assumption on $s$, it may happen that $[y\mapsto B_{\xi}(y;p;s)\equiv 0]$. 
We now impose further assumptions on $s$ in order to guarantee that \eqref{pers} is not identically equal to $0$ and that it
is bounded as $|\Norm(y)|\rightarrow\infty$.
Let $m_p:=\max\limits_{j\in J}\{\pm\frac{p_j}{2}\}$. Then we have the following proposition:
\begin{Prop}\label{congre}
Let $\xi\in\ca{L}^*\bs\{0\}$ and let $s\in\CC$ with $\Ree(s)>m_p$. Then  $\ca{V}_{p}^b(\xi;s)$ is a one-dimensional complex vector space.
Moreover, $[y\mapsto B_{\xi}(y;p;s)]$ is
a non-zero element of $\ca{V}_{p}^b(\xi;s)$. Furthermore, if $s\in\RR_{>m_p}$, then
for all $y\in(\RR^{\times})^g$, $B_{\xi}(y;p;s)\neq 0$. 
\end{Prop}

{\bf Proof} The fact that $\ca{V}_{p}^b(\xi;s)$ is at most a one-dimensional complex vector space follows from the 
second part of (5) of Proposition \ref{wee}. The fact that $[y\mapsto B_{\xi}(y;p;s)]$ 
lies in  $\ca{V}_{p}(\xi;s)$ follows from the computations done in Section \ref{dos} and
the observation that the real analytic Eisenstein series constructed in this manuscript are eigenvectors (of eigenvalue $s(1-s)$) 
of the partial Laplacians $\Delta_{j,-p_j}$ (for $j\in J$). It is also possible to prove by a direct (but tedious) calculation that  
$[y\mapsto B_{\xi}(y;p;s)]$ is killed by the differential operator $D_{j,p_j}[\xi;s]$.
The fact that $[y\mapsto B_{\xi}(y;p;s)]$ is bounded as $|\Norm(y)|\rightarrow\infty$
follows from the definition of $B_{\xi}(y;p;s)$ and (5) of Proposition \ref{wee}. The fact that, for a fixed
$s\in\RR_{>m_p}$, $[y\mapsto B_{\xi}(y;p;s)\not\equiv 0]$
follows from (4) of Proposition \ref{wee}. Finally, the non-vanishing of 
$[y\mapsto B_{\xi}(y;p;s)$, when $s\in\RR_{>m_p}$ and $y\in(\RR^{\times})^g$,
follows from (6) of Proposition \ref{wee}. This concludes the proof. \fin

For a real number $p\in[1,\infty)$, or $p=\infty$, and a vector $w\in\CC^g$, we let $||w||_p$ denote the $\ell_p$-norm. In particular,
$||w||_{\infty}=\max\limits_{j\in J}|w_j|$.
\begin{Cor}\label{tab}
Let $\xi\in\ca{L}^*\bs\{0\}$ and let $s\in\CC$ with $\Ree(s)>m_p$. Then 
\begin{align}\label{tub}
\ca{V}_{p}^b(\xi;s)=\CC\cdot B_{\xi}(y;p;s).
\end{align}
Moreover, there exists a positive constant
$C_s>0$ (which depends only on $s$ and not on $p$), such that 
for all $\xi\in\ca{L}^{*}\bs\{0\}$, and all $y\in\RR_{>0}^g$, 
\begin{align}\label{tub2}
f_{\xi}(y)\leq C_s\cdot |\Norm(\xi)|^{2s-1}\cdot e^{-6\cdot ||\xi y||_{\infty}}.
\end{align}
\end{Cor}

{\bf Proof} The proof of \eqref{tub} follows directly from Proposition \ref{congre}. The proof of \eqref{tub2}
follows from the inequalities $||\xi y||_1 \geq ||\xi y||_{\infty}$ and the explicit
upper bound \eqref{cart} given in Proposition \ref{tiopp}. \fin

For a real number $a\in\RR$, recall that $\Pi_a:=\{s\in\CC:\Ree(s)>a\}$.
We may now state the main result of this appendix.
\begin{Prop}\label{congre2}
Let $T\subseteq\Pi_1$ be a discrete subset. Let $G(z,s):\mk{h}^g\times (\Pi_1\bs T)\rightarrow\CC$ be a real analytic function in $(z,s)$.
Suppose that there exists a lattice $\ca{L}\subseteq \RR^g$, such that, for each $s\in\Pi_1\bs T$ and $y\in(\RR^{\times})^g$ fixed, 
the function $[x\mapsto G(x+\ii y,s)]$ is $\ca{L}$-invariant. Let $p\in\ZZ^g$, and
assume that for each $j\in J$, $\Delta_{j,p_j}G(z,s)=s(1-s)G(z,s)$. Finally, suppose also that for each 
pair $(x,s)\in \RR^g\times(\Pi_{m_p}\bs T)$ fixed, the function $[y\mapsto |G(x+\ii y,s)|]$ is bounded as $|\Norm(y)|\rightarrow\infty$. Then
the Fourier coefficient $a_{\xi}(y,s)$ lies in the solution space $\ca{V}_{p}^b(\xi;s)$.
\end{Prop}

{\bf Proof} From the previous discussion, we already know that, for each $\xi\in\ca{L}^{*}$, 
the coefficient $a_{\xi}(y,s)$ lies in the solution space $\ca{V}_p(\xi;s)$ of complex dimension $2^g$. 
By assumption, for each pair $(x,s)\in \RR^g\times \left(\Pi_{m_p}\bs T\right)$ fixed,
the function $[y\mapsto |G(x+\ii y,s)|]$ is bounded as $|\Norm(y)|\rightarrow\infty$. It follows
from Parseval's theorem (see (2) of Theorem \ref{keyy}) that $|a_{\xi}(y,s)|$ is bounded as $|\Norm(y)|\rightarrow 0$,
and, therefore, from Proposition \ref{congre}, we deduce that $a_{\xi}(y,s)\in \ca{V}_p^b(y,s)$. \fin

\subsection{Recurrence formula for the Taylor series coefficients around $s=1$}\label{app_4}
It follows from Theorem \ref{key_thmm} that 
the function $[s\mapsto G_{(\mk{m},\mk{n})}^0(U,p\,;z,s)]$ admits at most a pole of order one 
at $s=1$. Let 
\begin{align}\label{eto}
G_{(\mk{m},\mk{n})}^0(U,p\,;z,s)=\sum_{n\geq -1} c_n(z)(s-1)^n,
\end{align}
be its Taylor series development around $s=1$. It is proved in Section \ref{car0},
that for each $j\in J:=\{1,2,\ldots,g\}$, $G_{(\mk{m},\mk{n})}^0(U,p\,;z,s)$ is an eigenvector
with eigenvalue $s(1-s)$ for the partial-graded Laplacian $\Delta_{j,-p_j}$. Therefore, if we apply, $\Delta_{j,-p_j}$ to
\eqref{eto}, we find that the real analytic functions $c_n(z)$ must satisfy, for each $j\in J$, the following
relations:
\begin{enumerate}
 \item $\Delta_{j,-p_j} (c_{-1}(z))=0$, 
 \item $\Delta_{j,-p_j} (c_0(z))=-c_{-1}(z)$,
 \item $\Delta_{j,-p_j} (c_n(z))=-c_{n-1}(z)-c_{n-2}(z)$ for $n\geq 1$.
\end{enumerate}

\subsection{Riemannian metric on $\mk{h}^g$ and the distance to a cusp}\label{app_5}
We view $\mk{h}^g$ as a complex analytic manifold. Since the group $SL_2(\RR)^g$ acts transitively on $\mk{h}^g$
and that $\Stab_{SL_2(\RR)}(\ii,\ii,\ldots,\ii)=SO(2)^g$, it follows that 
\begin{align}\label{sein}
\mk{h}^g\simeq SL_2(\RR)^g/SO(2)^g.
\end{align}
From the isomorphism \eqref{sein} and the compactness of $SO(2)^g$, we may deduce from Proposition 1.6 
of \cite{Shim0} that a subgroup $\Gamma\leq SL_2(\RR)^g$ is discrete if and only if it
acts {\it properly discontinuously} on $\mk{h}^g$.

The group of {\it complex analytic isomorphisms} of $\mk{h}^g$, which we denote
by $\Hol(\mk{h}^g)$, fits into the following short exact sequence:
\begin{align*}
1\rightarrow PSL_2(\RR)^g\rightarrow\Hol(\mk{h}^g)\rightarrow S_g\rightarrow 1.
\end{align*}
Here, $S_g$ denotes the symmetric group of degree $g$, which we let act on $\mk{h}^g$, by permuting the $g$ coordinates of $\mk{h}^g$. 
We endow the space
$\mk{h}^g$ with the Poincar\'e metric
\begin{align*}
ds^2=\sum_{i=1}^g \frac{dx_i^2+dy_i^2}{y_i^2},
\end{align*}
where $z=(z_i)_{i=1}^g$ and $z_i=x_i+\ii y_i$. The Riemannian volume form of $(\mk{h}^g,ds^2)$ is given in local 
coordinates by
\begin{align}\label{dof}
dV:=\bigwedge_{j=1}^{g} \frac{\ii}{2}\frac{dz_j\wedge d\ov{z}_j}{\Imm(z_j)^2}=\bigwedge_{j=1}^g \frac{dx_j\wedge dy_j}{y_j^2}.
\end{align}
It is also sometimes called the Gau\ss-Bonnet form. 

The group of {\it real analytic isometries} of $\mk{h}^g$, which we denote by $\Isom^{\infty}(\mk{h}^g)$, 
fits into the following short exact sequence:
\begin{align*}
1\rightarrow PSL_2(\RR)^g\rightarrow\Isom^{\infty}(\mk{h}^g)\rightarrow S_g\times S_2^{g}\rightarrow 1.
\end{align*}
Here, the group $S_2^g$ acts on $\mk{h}^g$ 
in the following way: if we let $\sigma_i:=(0,\ldots,1,\ldots,0)\in S_2^g$, where $1$ is placed in the $i$-th position,
then $(\sigma_i(z))_j=z_j$, if $j\neq i$ and $(\sigma_i(z))_i=-\ov{z}_i$. Note that the isometry $\sigma_i\in S_2^g$ changes
(a choice) the orientation of $\mk{h}^g$.

Recall that a discrete subgroup $\ov{\Gamma}\leq PSL_2(\RR)^g$ is said to
be {\it irreducible} if the image of each projection $\pi_j:\ov{\Gamma}\rightarrow PSL_2(\RR)^{g-1}$ (where $\pi_j$ forgets about
the $j$-th coordinate) is dense. We have the following fundamental result of Selberg (see \cite{Sel2}):
\begin{Th}\label{sel}(Selberg)
Let $\ov{\Gamma}\leq PSL_2(\RR)^g$ be an discrete subgroup with $g\geq 2$. Assume that $\ov{\Gamma}$ is
\begin{enumerate}[(i)]
 \item irreducible,
 \item has a fundamental domain of finite volume,
 \item has at least one cusp.
\end{enumerate}
Then, a conjugate of $\ov{\Gamma}$ (inside $PSL_2(\RR)^g$) is commensurable to $PSL_2(\ca{O}_K)$, for
some totally real field $K$ of dimension $g$ over $\QQ$. Here $\ca{O}_K$ denotes the ring of integers
of $K$.
\end{Th}
\begin{Rem}
When $g=1$, it is possible to show that there are (uncountably) many examples of discrete subgroups $\ov{\Gamma}\leq PSL_2(\RR)$ which 
satisfy (ii) and (iii) but which are not commensurable to a conjugate of $PSL_2(\ZZ)$. For one such explicit example, one can
take $\ov{\Gamma}$ to be a {\it Hecke triangle group} where the underlying hyperbolic triangle has one vertex placed 
at infinity and one angle of magnitude $\frac{\pi}{5}$. In this case, $\ov{\Gamma}$ cannot be commensurable to a 
conjugate subgroup of $PSL_2(\ZZ)$. Indeed, since the commensurator of $PSL_2(\ZZ)$ is equal to $PSL_2(\QQ)$,
we must have $g\ov{\Gamma}g^{-1}\leq PSL_2(\QQ)$ for some $g\in PSL_2(\RR)$. Finally, one uses the well-known 
fact that if $\gamma\in PSL_2(\QQ)$ is an element of finite order $n$ then necessarily $n\in\{1,2,3\}$.
\end{Rem}

Now let $\Gamma\leq (GL_2^+(\RR))^g$ be a discrete subgroup. 
In light of Theorem \ref{sel}, we don't loose much by assuming that $\Gamma\leq GL_2^+(K)$ is commensurable to $GL_2(\ca{O}_K)$.
We have the natural projection map 
\begin{align*}
\pi: GL_2^+(K)\rightarrow PSL_2(\RR)
\end{align*}
given by $\gamma\mapsto\frac{\gamma}{\sqrt{\det(\gamma)}}$. We let $\ov{\Gamma}$ be the image of $\Gamma$ under the natural projection map
$\pi:GL_2^+(K)\rightarrow \leq PSL_2(\RR)$, then it follows
that $\ov{\Gamma}$ acts properly discontinuously on $\mk{h}^g$. In particular, for every point $z\in\mk{h}^g$, we have that
$\Stab_{\ov{\Gamma}}(z)$ is a finite group, and, therefore, the quotient space 
\begin{align*}
Y_\Gamma:=\biglslant{\Gamma}{\mk{h}^g}=\biglslant{\ov{\Gamma}}{\mk{h}^g},
\end{align*}
has the structure of a complex analytic orbifold. Since $\ov{\Gamma}\subseteq\Isom^{\infty}(\mk{h}^g)$, the metric $ds^2$
on $\mk{h}^g$ descends to a metric on $Y_{\Gamma}$. Therefore, $Y_{\Gamma}$ is a {\it Riemannian orbifold}.

A parabolic element $\gamma\in\Gamma$, when viewed as acting on $\CC^g$ via its corresponding 
M\"obius transformation, has a unique fixed point $z_0$, which lies necessarily in $\PP^1(K)\subseteq \PP^1(\RR)^g$. 
As it is well-known, the coset $\Gamma z_0\in \biglslant{\Gamma}{\PP^1(K)}$ corresponds geometrically to a {\it cusp},
i.e., an infinite shrinking end inside the space $Y_\Gamma$. It is also easy to see, that for each $\mk{c}\in\PP^1(K)$,
there exists at least one parabolic element $\gamma\in\Gamma$ such that $\mk{c}$ is fixed. 
For the reasons mentioned above, we call the elements in $\PP^1(K)$ simply {\it cusps}. We note that the group $GL_2^+(K)$ acts transitively on $\PP^1(K)$ 
via M\"obius transformation. A right coset in $\biglslant{\Gamma}{\PP^1(K)}$ will be called {\it a relative $\Gamma$-cusp}. 
Usually, a right coset in $\biglslant{\Gamma}{\PP^1(K)}$ containing a cusp 
$\mk{c}$ will be denoted by $[\mk{c}]_\Gamma$. We also define the {\it set of cusps of $Y_{\Gamma}$} as the set of relative $\Gamma$-cusps,
i.e., as the set of cosets $\{[\mk{c}]_{\Gamma}:\mk{c}\in\PP^1(K)\}$.

For $z\in\mk{h}^g$ and $\gamma=\M{a}{b}{c}{d}\in GL_2^+(K)$, we have the formula
\begin{align}\label{tet}
\Imm(\gamma z)=\det(\gamma)\cdot\frac{\Imm(z)}{|cz+d|^2}\in\mk{h}^g. 
\end{align}
Let $\mk{c}\in\PP^1(K)$ be a cusp and choose a writing of $\mk{c}=\frac{a}{b}$ with $a,b\in K$ (we do not require here to choose $a$ and $b$ in $\ca{O}_K$).
For $z\in\mk{h}^g$, following van der Geer (see p. 7 of \cite{Geer}), we define
\begin{align}\label{niz}
\mu(\mk{c},z):=\frac{(\Norm((a,b)))^2\Norm(y)}{|\Norm(-b z+a)|^2}.
\end{align}
Here, $(a,b)=a\ca{O}_K+b\ca{O}_K\subseteq K$ corresponds to the fractional $\ca{O}_K$-ideal generated by $a$ and $b$,
and $\Norm((a,b)):=[\ca{O}_K:(a,b)]\in\QQ_{>0}$; the rational index $[\ca{O}_K:(a,b)]$ being defined as in Section \ref{lat_nu}.
One readily sees that the right-hand side of \eqref{niz} is independent of the choice of the writing of $\mk{c}$.
Moreover, if we choose (arbitrarily) an element $\gamma=\M{a}{*}{b}{*}\in GL_2^+(K)$, such that $\gamma(\infty)=\frac{a}{b}=\mk{c}$,
where $\infty=\frac{1}{0}$, a direct calculation shows that 
\begin{align}\label{tet2}
\mu(\mk{c},z)=\Norm(\det(\gamma))^{-1}\Norm((a,b))^2\Norm(\Imm(\gamma^{-1} z)).
\end{align}
From \eqref{tet2}, one gets that, for all $\eta\in GL_2^+(\ca{O}_K)$, $\mk{c}\in\PP^1(K)$ and $z\in\mk{h}^g$
\begin{align}\label{tet3}
\mu(\eta \mk{c},\eta z)=\mu(\mk{c},z).
\end{align}
Therefore the function $\mu$ is {\it jointly invariant} under the isometries induced by the matrices in $GL_2^+(\ca{O}_K)$.
\begin{Rem}
In general, the function $\mu$ is not necessarily jointly invariant under the isometry of $\mk{h}^g$ induced by a matrix $\eta\in GL_2^+(K)$.  
\end{Rem}

The next proposition describes the main properties of the function $\mu$.
\begin{Prop}\label{keyp}
The function $\mu$ satisfies the following properties:
\begin{enumerate}
 \item Let $\mk{c}_1,\mk{c}_2\in\PP^1(K)$. There exists a positive real number $r$, depending only on $K$, such that
\begin{align*}
\mbox{[$\mu(\mk{c}_1,z)>r$ and $\mu(\mk{c}_2,z)>r$ for all $z\in\mk{h}^g$]} \Longrightarrow \mk{c}_1=\mk{c}_2.
\end{align*}
 \item There exists a positive real number $s$, depending only on $K$, such that for all $z\in\mk{h}^g$
 there exists a cusp $\mk{c}\in\PP^1(K)$ such that $\mu(\mk{c},z)>s$.
 \item For all finite cusp $\mk{c}=\frac{a}{c}\in K\subseteq\PP^1(K)$ (so that $c\neq 0$) and all $z\in\mk{h}^g$, 
 one has that $\mu(\mk{c},z)\leq \frac{1}{|\Norm(y)|}$.
\end{enumerate}
\end{Prop}

{\bf Proof} (1) and (2) are proved on p. 7 and 8 of \cite{Geer}. It remains to prove (3).
Without loss of generality, let us assume that $a,b\in\ca{O}_K$. We have 
\begin{align*}
\mu(\mk{c},z)=\frac{(\Norm((a,b)))^2\Norm(y)}{|\Norm(-b z+a)|^2}\leq \frac{1}{\Norm(y)},
\end{align*}
where the inequality above follows directly from the two obvious inequalities $\Norm((a,b))\leq |\Norm(b)|$ and $|-b^{(j)}z_j+a^{(j)}|
\geq |-b^{(j)}y_j|$ ($1\leq j\leq g$).
The result follows.  \fin
\begin{Rem}
One should view the quantity $\frac{1}{\mu(\mk{c},z)^{1/2}}$ as providing a distance between the cusp $\mk{c}$
and the point $z\in\mk{h}^g$ which is \lq\lq compatible\rq\rq\; with the hyperbolic metric in the sense that $\mu$
is jointly invariant under the discrete group of isometries induced by $GL_2^+(\ca{O}_K)$.
\end{Rem}

\subsubsection{Sphere of influence and neighborhood of a cusp}\label{app_5a}

We use the same notation as in the section above. Recall that $\Gamma\leq GL_2^+(K)$ is a discrete subgroup commensurable 
to $GL_2^+(\ca{O}_K)$ and that $\ov{\Gamma}$ denotes its image under the natural projection $\pi:GL_2^+(K)\rightarrow PSL_2(\RR)$.
Recall also that $Y_\Gamma=\biglslant{\Gamma}{\mk{h}^g}$ is a Riemannian orbifold. Using the function $\mu$ introduced in the previous
section, we first define the \lq\lq sphere of influence\rq\rq\; at a cusp $\mk{c}$ relative to the group $GL_2^+(\ca{O}_K)$ 
(which is a certain subset of $\mk{h}^g$). Secondly, by considering certain translates of these spheres of influence, and 
projecting to $Y_\Gamma$, we construct a certain neighborhood of the relative cusp  $[\mk{c}]_{\Gamma}\subseteq Y_\Gamma$.  

As is explained on p. 8 and 9 of \cite{Geer}, for each cusp $\mk{c}\in\PP^1(K)$, one may define
a {\it sphere of influence} $F_{\mk{c}}\subseteq\mk{h}^g$, where
\begin{align*}
F_{\mk{c}}:=\{z\in\mk{h}^g:\mu(\mk{c},z)\geq \mu(\mk{c}',z) \s\mbox{for all $\mk{c}'\in\PP^1(K)$}\}.
\end{align*}
One may check that for all $\gamma\in GL_2^+(\ca{O}_K)$, $F_{\gamma \mk{c}}=\gamma F_{\mk{c}}$. In particular, 
the set of $GL_2^+(\ca{O}_K)$-translates of $F_{\mk{c}}$ provides a tessellation of $\mk{h}^g$. We call the 
set $F_{\mk{c}}$ the sphere of influence at the cusp $\mk{c}$ relative to the group $GL_2^+(\ca{O}_K)$.

Given a relative cusp $[\mk{c}]_{\Gamma}$, we would like to define a certain neighborhood $B_{\mk{c}}\leq Y_\Gamma$
of $[\mk{c}]_{\Gamma}$. At first, let {\it us assume} that $\Gamma\leq GL_2^+(\ca{O}_K)$ (rather than just in $GL_2^+(K)$), and, as before, we denote its image in $PSL_2(\RR)^g$ by 
$\ov{\Gamma}$. Let $\{\mk{c}_1,\ldots,\mk{c}_h\}\subseteq\PP^1(K)$ be a complete set of representatives of the relative $\Gamma$-cusps.
For each $\mk{c}\in\{\mk{c}_1,\ldots,\mk{c}_h\}$, we let $\Gamma_\mk{c}=\Stab_\Gamma(\mk{c})$. As it is explained on 
p. 8 and 9 of \cite{Geer}, the quotient space $Y_\Gamma=\biglslant{\Gamma}{\mk{h}^g}$ may be decomposed as
\begin{align}\label{tet6}
 Y_{\Gamma}=\bigcup_{j=1}^h B_{[\mk{c}_j]_\Gamma},
\end{align}
where  $B_{[\mk{c}_j]_\Gamma}$ is the image of the natural map $\iota_{\mk{c}_j}:\biglslant{\Gamma_{\mk{c}_j}}{F_{\mk{c}_j}}\rightarrow Y_{\Gamma}$.

Let us now treat the general case. By assumption, $\Gamma$ is commensurable to $GL_2(\ca{O}_K)$.
Therefore, there exists a subgroup $\Gamma'\leq \Gamma$ such that $\Gamma'$ has finite index in $GL_2(\ca{O}_K)$.
In particular, from the previous paragraph, for each relative cusp $[\mk{c}]_{\Gamma'}$, we can associate
a neighborhood $B_{[\mk{c}]_{\Gamma'}}\leq Y_{\Gamma'}$ of the relative cusp $[\mk{c}]_{\Gamma'}$. Since $[\Gamma:\Gamma']<\infty$, we have a natural finite
covering map $p:Y_{\Gamma'}\rightarrow Y_{\Gamma}$. Finally, given a cusp $\mk{c}\in\PP^1(K)$, we define 
\begin{align*}
B_{[\mk{c}]_\Gamma}:=\bigcup_{p([\mk{d}]_{\Gamma'})=[\mk{c}]_\Gamma} p(B_{[\mk{d}]_{\Gamma'}}).
\end{align*}

\subsection{A proof of Proposition \ref{poids} using the point-pair invariant kernel method}\label{app_6}

{\bf Proof} Our proof of Proposition \ref{fat0} is inspired from Kubota's approach presented on p.12 and 13 of \cite{Kubo1},
where the inequality \eqref{chag2} is proved for discrete subgroups of $SL_2(\RR)$ which admit a 
finite covolume and at least one cusp.
Let $\Gamma:=SL_2(\ca{O}_K)$ be the Hilbert modular group of $K$. The group $\Gamma$ acts 
naturally on the left of $\PP^1(K)$ by M\"obius transformations in the following way:
\begin{align*}
\M{a}{b}{c}{d}[x,y]=[ax+by,cx+dy].
\end{align*}
Let 
\begin{align*}
\ca{C}:=\{\mk{c}_1:=[c_1;d_1],\ldots,\mk{c}_h:=[c_h;d_h]\}\in\PP^1(K),
\end{align*} 
be a complete set of representatives of the relative $\Gamma$-cusps, where $c_i,d_i\in\ca{O}_K$ for all $i\in\{1,\ldots,h\}$. 
(see Appendix \ref{app_5} for some basic notions on cusps).
Without loss of generality, we may assume that $\mk{c}_1=[1,0]$, so that $c_1=1$ and $d_1=0$. We may also assume that
each pair $(c_i,d_i)$ is {\it $\ca{O}_K$-reduced} (even though we don't need it in 
the present proof), in the sense that $c_i$ and $d_i$ are not simultaneously divisible by a non-unit element
of $\ca{O}_K$. In general, it is not always possible to choose the elements $c_i$ and $d_i$ to be coprime, since $\ca{O}_K$ may fail to be
a unique factorization domain. 

Let $\Gamma_{\mk{c}_i}=\Stab_{\Gamma}(\mk{c}_i)$ be the stabilizer of the cusp $\mk{c}_i$ and let
$\eta\in\Gamma_{\mk{c}_i}$. Since entries of $\eta$ are algebraic integers, and $\det(\eta)=1$,
it follows that
\begin{align}\label{cari}
\eta\V{c_i}{d_i}=\epsilon_{\eta}\V{c_i}{d_i}, 
\end{align}
for some $\epsilon_{\eta}\in\ca{O}_K^{\times}$. Since $\ca{C}$ is a complete set of representatives 
of the relative $\Gamma$-cusps, for each pair $(c,d)\in\mathcal{R}''$, there exist
$\lambda\in K^{\times}$, $\gamma\in \Gamma$, and an index $i\in\{1,\ldots,h\}$, such that
\begin{align}\label{eyy}
\lambda\cdot\gamma\V{c_i}{d_i}=\V{c}{d}.
\end{align}
Looking more closely at \eqref{eyy}, we see that, in fact, 
$\lambda\in\mk{a}_i^{-1}$ where $\mk{a}_i=c_i\ca{O}_K+d_i\ca{O}_K$. Also, since for 
any $\lambda\in\mk{a}_i^{-1}$, $\V{\lambda c_i}{\lambda d_i}\in\V{\ca{O}_K}{\ca{O}_K}$, it follows from \eqref{cari} and \eqref{eyy} that
\begin{align}\label{gas3}
\bigcup_{i=1}^h\bigcup_{\lambda\in\mk{a}_i^{-1}} 
\left(\bigrslant{\Gamma}{\Gamma_{\mk{c}_i}}\right)\V{\lambda c_i}{\lambda d_i}=\V{\ca{O}_K}{\ca{O}_K}\pmod{\ca{O}_K^{\times}}.
\end{align}

For each cusp $\mk{c}_i=[c_i,d_i]\in\ca{C}$, let us choose, arbitrarily, a matrix 
\begin{align}\label{fus}
\delta_i^{-1}=\M{c_i}{*}{d_i}{*}\in SL_2(K).
\end{align}
In particular, $\delta_i^{-1}\infty=\mk{c}_i$. It follows from \eqref{gas3} that
\begin{align}\label{fas5}
\bigcup_{i=1}^h \bigcup_{\lambda\in\mk{a}_i^{-1}\bs\{0\}}
\left(\bigrslant{\Gamma}{\Gamma_{\mk{c}_i}}\right)\lambda\cdot \delta_i^{-1}=\left\{\M{a}{*_1}{b}{*_2}:a,b\in\ca{O}_K, (a,b)\neq (0,0)\right\}\pmod{\ca{O}_K^{\times}}, 
\end{align}
where the symbols $*_1$ and $*_2$ go over suitable subsets of $K$ which are determined by the above equality. 
Consider now the involution 
\begin{align*}
\mk{s}:M_2(K)&\rightarrow M_2(K)\\
       M=\M{a}{b}{c}{d}&\mapsto M^{\mk{s}}=\M{d}{-b}{-c}{a}
\end{align*} 
The set involution $\mk{s}$ satisfies the following easily verified properties:
\begin{enumerate}[(i)]
 \item For all $M\in M_2(K)$, $M^{\mk{s}}=(M^*)^t$.
 \item For all $M\in M_2(K)$ and $\lambda\in K$, $(\lambda M)^{\mk{s}}=\lambda \cdot M^{\mk{s}}$.
 \item For all $M\in SL_2(K)$, $M^{\mk{s}}=M^{-1}$.
 \item For all $M_1,M_2\in M_2(K)$, $(M_1+M_2)^{\mk{s}}=M_1^{\mk{s}}+M_2^{\mk{s}}$.
 \item For all $M_1,M_2\in M_2(K)$, $(M_1M_2)^{\mk{s}}=M_2^{\mk{s}}M_1^{\mk{s}}$.
\end{enumerate}
In particular, if follows from (ii), (iv) and (v) that $\mk{s}$ is an $K$-linear anti-automorphism of the matrix ring $M_2(K)$. Applying $\mk{s}$ to
each side of \eqref{fas5}, we finally obtain the following equality of sets:
\begin{align}\label{fas6}
\bigcup_{i=1}^h\bigcup_{\lambda\in\mk{b}_i\bs\{0\}} 
\lambda\cdot\delta_i\left(\biglslant{\Gamma_{\mk{c}_i}}{\Gamma}\right)=\left\{\M{*_1}{*_2}{a}{b}:a,b\in\ca{O}_K, (a,b)\neq (0,0)\right\}\pmod{\ca{O}_K^{\times}},
\end{align}
where $\mk{b}_i:=\mk{a}_i^{-1}$. Note that, for each index $i$, if we replace $\delta_i$ by $\delta_i\eta$ in \eqref{fas6}, 
for some $\eta\in\Gamma_{\mk{c}_i}$, then \eqref{fas6} still holds true. 

We would like now to take advantage of the set equality \eqref{fas6} in order to rewrite the \lq\lq almost\rq\rq\; Eisenstein series 
$\ca{E}(z,s)$, which appears in \eqref{chag}, as a finite sum of classical real analytic Poincar\'e-Eisenstein series
of weight $0$. For $1\leq i\leq h$, we define
\begin{align*}
E_i(z,s):=\sum_{\gamma\in \ca{S}_i} \Imm(\delta_i\gamma z)^s,
\end{align*}
where for $1\leq i\leq h$, $\ca{S}_i\subseteq\Gamma$  is a complete set of representatives of $\biglslant{\Gamma_{\mk{c}_i}}{\Gamma}$.
Here $\Gamma_{\mk{c}_i}=\Stab_{\Gamma}(\mk{c}_i)$. One may check that the series $E_i(z,s)$ are well-defined
and that each of them, when viewed as a function in $z$, is a real analytic modular form
of weight zero for the congruence group $\Gamma$.

For a lattice $\mk{a}\leq K$, recall that (see Definition \ref{nor_zeta})
\begin{align}\label{footg}
\wt{Z}_{\mk{a}}(0,0,\omega_{\ov{\bzt}},2s)=\sum_{\{0\neq \lambda\in\mk{a}\}/\ca{O}_K^{\times}(\infty)} \frac{1}{|\Norm(\lambda)|^{2s}},
\end{align}
where $2s\in\Pi_1$. If we let $e=[\ca{O}_K^{\times}:\ca{O}_K^{\times}(\infty)]$, then it follows from \eqref{footg} that
\begin{align}\label{footg2}
\wt{Z}_{\mk{a}}(0,0,\omega_{\ov{\bzt}},2s)=e\cdot \sum_{\{0\neq \lambda\in\mk{a}\}/\ca{O}_K^{\times}} \frac{1}{|\Norm(\lambda)|^{2s}}.
\end{align}
From the definitions of $E_i(z,s)$ and $\wt{Z}_{\mk{a}}(0,0,\omega_{\ov{\bzt}},2s)$, the set equality \eqref{fas6} implies that
\begin{align}\label{poum1}
&\ca{E}(z,s)=\sum_{(c,d)\in\mathcal{R}''}\frac{|y|^{\bu\cdot s}}{|\Norm(cz+d)|^{2s}}\\ \notag
&\leq\frac{1}{e}\left(\wt{Z}_{\ca{O}_K}(0,0,\omega_{\ov{\bzt}},2s)\cdot\left(E_1(z,s)-\Norm(y)^s\right)+
\sum_{i=2}^h \wt{Z}_{\mk{b}_i}(0,0,\omega_{\ov{\bzt}},2s)\cdot E_i(z,s)\right).
\end{align}
Recall here that $\ca{R}''$ is a complete set of representatives of $(\ca{O}_K\bs\{0\})\times\ca{O}_K$ modulo the diagonal action
of $\ca{O}_K^{\times}$. In order to uniformize the notation, we let $\wt{E}_1(z,s):=\left(E_1(z,s)-\Norm(y)^s\right)$
and, for $i\geq 2$, $\wt{E}_i(z,s)=E_i(z,s)$. 

Looking at the defining series of $\wt{Z}_{\mk{b}_i}(0,0,\omega_{\ov{\bzt}},2s)$, 
one readily sees that there exists a constant $D>0$, such that for all $s\in\Pi_1$ and all $i\in\{1,2,\ldots,g\}$,
\begin{align*}
|\wt{Z}_{\mk{b}_i}(0,0,\omega_{\ov{\bzt}},2s)|\leq D.
\end{align*}
Therefore, in order to show the inequality \eqref{chag2} of Proposition \ref{poids},
it is enough to show that, for each $s>1$ fixed, and for all $i\in\{1,\ldots,h\}$, there exists a constant $C_{i,s}>0$ such that 
for all $z\in \mk{h}^g$,
\begin{align}\label{toop}
\wt{E}_i(z,s)\leq C_{i,s}\cdot|\Norm(y)|^{1-s}.
\end{align}
The proof will be a consequence of the following two facts:
\begin{enumerate}[(a)]
 \item The function $z\mapsto\wt{E}_i(z,s)$ is invariant under a co-compact subgroup of the group $\gamma \ca{U}(\RR)\gamma^{-1}$,
 for some $\gamma\in\ca{G}_1(\RR)$, where $\ca{U}(\RR)$ corresponds to the unipotent subgroup of upper triangular matrices.
 \item In the indexing set of the defining sum of $\wt{E}_i(z,s)$, there is no matrix of the form $\M{*}{*}{0}{*}\in\ca{S}_i$ which occurs. 
\end{enumerate}
Note that (a) and (b) are valid for each index $i\in\{1,\ldots,h\}$.
Therefore, in order to prove \eqref{toop}, it is enough to prove it in the special case where $i=1$. Moreover,
replacing $z\in K_{\CC}^{\pm}$ by $\eta z$ in $\wt{E}_i(z,s)$, if necessary, where $\eta$ is a suitable matrix in $\ca{G}(\RR)=GL_2(\RR)^g$ (which does not affect (a) and (b)),
we may assume, without loss of generality, that $z\in\mk{h}^g$.

Now, we need to introduce the key notion of a {\it point-pair invariant kernel}. 
For $z,z'\in\mk{h}^g$, we let $d(z,z')$ be the {\it hyperbolic distance} between $z$ and $z'$.
Given $\epsilon>0$, we consider the point-pair invariant $k_{\epsilon}(z,z')$ given by
\begin{align*}
k_{\epsilon}(z,z')
:=\left\{
\begin{array}{cc}
1 & \mbox{if $d(z,z')<\epsilon$} \\
0 & \mbox{if $d(z,z')\geq\epsilon$}
\end{array}
\right.
\end{align*}
The point-pair invariance means that for
all $\gamma\in \ca{G}_1(\RR)=SL_2(\RR)^g$ and all $z,z'\in\mk{h}^g$, one has that $k_{\epsilon}(z,z')=k_{\epsilon}(\gamma z,\gamma z')$.
Such point-pair invariant kernel gives rise to translation invariant integral operators, see \S 1.3 of \cite{Kubo1}) for more details.

For any $\gamma=\M{*}{*}{c}{d}\in \ca{G}_1(\RR)$, the function
\begin{align}\label{xes}
z\mapsto \Imm(\gamma z)^{\bu\cdot s}=\frac{y^{\bu\cdot s}}{|\Norm(cz+d)|^{2s}},
 \end{align}
is a simultaneous eigenvector, with eigenvalues $s(1-s)$, of the $g$ graded Laplacians $\{\Delta_{j,0}:j=1,\ldots,g\}$. For the definition
of $\Delta_{j,0}$ and a proof that $\Imm(\gamma z)^{\bu\cdot s}$ is a $\Delta_{j,0}$-eigenvector, see Section \ref{car0}. Using a argument
similar to the proof of Theorem 1.3.2 of \cite{Kubo1}, it follows that there exists a positive constant $\Lambda_{\epsilon}>0$, 
{\it independent of the matrix $\gamma\in\ca{G}_1(\RR)$}, such that
\begin{align}\label{eup}
\int_{\mk{h}^g} k(z,z')\cdot \Imm(\gamma z')^{\bu\cdot s} dV(z')=\Lambda_{\epsilon}\cdot\Imm(\gamma z)^{\bu\cdot s}.
\end{align}
Here $dV(z)$ stands for the volume form of $\mk{h}^g$ with respect to the Poincar\'e metric, see Appendix \ref{app_5}. Using the point-pair invariance
of $k_{\epsilon}(z,z')$, a change of coordinates shows that
\begin{align}\label{cire}
\int_{\mk{h}^g} k(z,z')\cdot\Imm(\gamma z')^{\bu\cdot s} dV(z')= \int_{\mk{h}^g} k(\gamma z,z')\Imm(z')^{\bu\cdot s} dV(z').
\end{align}
Note that the integral on the right-hand side of \eqref{cire} corresponds to the integral of the function
$z'\mapsto\Imm(z')^{\bu\cdot s}=\Norm(y')^s$ in a hyperbolic $g$-dimensional disc of radius $\epsilon$ with center $\gamma z$. 

From \eqref{eup}, we have
\begin{align}\label{foot0}
\wt{E}_1(z,s)=\Lambda_{\epsilon}\int_{\mk{h}^g}\sum_{\eta\in \ca{S}_1 } k_{\epsilon}(\eta z,z')\Imm(z')^{\bu\cdot s} dV(z').
\end{align}
Given an $s>1$, it remains to show that the integral on the right-hand side of \eqref{foot0}
is bounded by $C_s\cdot\Norm(y)^{1-s}$, for some
positive constant $C_s$ which independent of $z\in\mk{h}^g$. 

The group $\Gamma_{\infty}\leq \ca{U}(\RR)$ acts naturally on the space $\mk{h}^g$. 
Following Lemma 2.10 of \cite{Frei1}, there are closed parallelotopes $P\subseteq\{x'\in\RR^g\}$ 
(of dimension $g$) and $Q\subseteq\{w'\in\RR^g:\Tr(w')=0\}$ (of dimension $g-1$), such that 
\begin{align*}
\left\{x'+t'e^{w'}\ii\in\mk{h^g}:t'\in\RR_{>0},x'\in P, w'\in Q\right\},
\end{align*}
is a fundamental domain for the action of $\Gamma_{\infty}$ on $\mk{h}^g$. Here 
\begin{align*}
w'=(w_1',\ldots,w_g'):=\frac{1}{\Norm(y')^{1/g}}\left(\log y_1',\ldots,\log y_g'\right).
\end{align*}
By definition of $\mathcal{S}_1$, for each $\eta=\M{a}{b}{c}{d}\in\mathcal{S}_1$, one has $c\neq 0$ with $a,c\in\ca{O}_K$. 
It follows from a similar calculation as the one done in the proof of (iii) of Proposition \ref{keyp} that
\begin{align}\label{footh}
\Imm(\eta z)\leq \frac{1}{\Norm(y)}.
\end{align}
Note that the above inequality hold true precisely because $c\neq 0$.

From \eqref{footh} we find that
\begin{align}\label{ito}
\int_{\mk{h}^g}\sum_{\eta\in \mathcal{S}_1 } k_{\epsilon}(\eta z,z')\Imm(z')^{\bu\cdot s} dV(z')\leq
\int\limits_{ \substack{ z'\in \mathcal{D}}} \Norm(y')^s dV(z'),
\end{align}
where
\begin{align*}
\mathcal{D}:=\left\{x'+t'e^{w'}\ii\in\mk{h^g}:t'\in\left[0,\frac{1}{\Norm(y)}\right],x'\in P, w'\in Q\right\}.
\end{align*}
In order to show that the integral appearing in the right-hand side of \eqref{ito} is bounded
from above, we make a change of variables.
We pose $t:=\Norm(y)$ and $u_i:=\frac{y_i}{t^{1/g}}=e^{w_i}$ for $i\in\{1,\ldots, g\}$. In particular, we have
$\Norm(u)=1$. A direct computation shows that
\begin{align*}
 dy_1\wedge \ldots \wedge dy_{g-1}\wedge dy_g=\frac{u_g}{g}du_1\wedge \ldots \wedge du_{g-1}\wedge dt.
\end{align*}
It thus follows that
\begin{align*}
\int\limits_{ \substack{ z'\in \mathcal{D}}} \Norm(y')^s dV(z')=
&\Big|\int\limits_{ \substack{ z'\in \mathcal{D}}} (t')^{s-2} \frac{u_g}{g}\cdot dt'\wedge dx_1
\wedge\ldots\wedge dx_g\wedge du_1\wedge \ldots \wedge du_{g-1}\Big|\\
&\leq \frac{C}{s-1}\cdot\Norm(y)^{1-s},
\end{align*}
where the constant $C$ {\it depends only} on the parallelotopes $P$ and $Q$ and, therefore, not on the variables 
$s$ and $z\in\mk{h}^g$. This concludes the proof of the inequality \eqref{chag2}. \fin

\bibliographystyle{plain}
\bibliography{\string~/references/biblio2}

\end{document}